\documentclass[12pt]{article}
\usepackage{amssymb,amsmath,amsthm,latexsym}
\usepackage[cp1251]{inputenc}
\usepackage{dsfont}
\usepackage{amssymb}
\usepackage{indentfirst}
\usepackage{graphicx}
\usepackage{caption}
\usepackage{multicol}
\usepackage{enumitem}
\usepackage{hyperref}
\hypersetup{citecolor=red,}
\usepackage{amsfonts, array, hhline}
\usepackage{color}
\usepackage{lmodern}
\usepackage[left=1in,right=1in,top=1in,bottom=1.5in]{geometry}

\DeclareMathOperator{\arcsinh}{arsinh}
\DeclareMathOperator{\arctanh}{artanh}
\DeclareMathOperator{\arccot}{arcot}
\DeclareMathOperator{\sgn}{sgn}
\begin{document}

\newcommand\blfootnote[1]{%
  \begingroup
  \renewcommand\thefootnote{}\footnote{#1}%
  \addtocounter{footnote}{-1}%
  \endgroup
}

\newcommand{\e}{\varepsilon}
\newcommand{\R}{\mathbb{R}}
\renewcommand{\H}{\mathbb{H}}
\renewcommand{\L}{\mathbb{L}}
\renewcommand{\S}{\mathbb{S}}
\newcommand{\E}{\mathbb{E}}
\newcommand{\T}{\mathbb{T}}
\newcommand{\len}{{\rm len}}
\newcommand{\vol}{{\rm vol}}
\newcommand{\diam}{{\rm diam}}
\newcommand{\area}{{\rm area}}
\newcommand{\sarea}{{\rm sarea }}
\newcommand{\inte}{{\rm int}}
\newcommand{\ang}{{\rm ang}}
\newcommand{\inter}{{\rm int}}
\newcommand{\dis}{{\rm dis}}
\newcommand{\dist}{{\rm dist}}
\renewcommand{\sl}{{\rm slope}}
\newcommand{\lev}{{\rm lev}}
\newcommand{\dS}{d\mathbb{S}}
\newcommand{\ds}{d\mathbb{S}}
\newcommand{\PdS}{\mathbb{P}d\mathbb{S}}
\renewcommand{\hat}{\widehat}
\renewcommand{\tilde}{\widetilde}

\makeatletter
\theoremstyle{definition}
\newtheorem*{rep@theorem}{\rep@title}
\newcommand{\newreptheorem}[2]{%
\newenvironment{rep#1}[1]{%
 \def\rep@title{#2 \ref{##1}}%
 \begin{rep@theorem}}%
 {\end{rep@theorem}}}
\makeatother

\numberwithin{equation}{section}

\theoremstyle{definition}
\newtheorem{dfn}{Definition}[subsection]
\newtheorem{exmp}[dfn]{Example}
\newtheorem{thm}[dfn]{Theorem}
\newtheorem{lm}[dfn]{Lemma}
\newtheorem{crl}[dfn]{Corollary}
\newtheorem{cnj}[dfn]{Conjecture}
\newtheorem{prob}[dfn]{Problem}
\newtheorem{thml}{Theorem}
\newreptheorem{lm}{Lemma}
\renewcommand*{\thethml}{\Alph{thml}}
\renewcommand{\proofname}{Proof}
\theoremstyle{remark}
\newtheorem{rmk}[dfn]{Remark}
\newtheorem{prop}{Claim}[dfn]
\renewcommand{\theprop}{\arabic{prop}}

\title{Rigidity of compact Fuchsian manifolds with convex boundary}
\author{Roman Prosanov \thanks{Supported by SNF grant $200021_-169391$ ``Discrete curvature and rigidity''.}}
\date{}
\AtEndDocument{\bigskip{\footnotesize
\par
  \textsc{Technische Universit\"at Wien, Wiedner Hauptstra{\ss}e 8-10/104 A-1040 Wien, Austria} \par
  \textit{E-mail}: \texttt{rprosanov@mail.ru}
}}
\maketitle

\begin{abstract}
A compact Fuchsian manifold with boundary is a hyperbolic 3-manifold homeomorphic to $S_g \times [0; 1]$ such that the boundary component $S_g \times \{ 0\}$ is geodesic. We prove that a compact Fuchsian manifold with convex boundary is uniquely determined by the induced path metric on $S_g \times \{1\}$. We do not put further restrictions on the boundary except convexity. 
\end{abstract}

\setcounter{tocdepth}{2}
\tableofcontents

\section{Introduction}

\subsection{Rigidity of convex bodies in $\E^3$}
\label{Weyl}

The first source of our motivation is the following  result of Pogorelov~\cite{Pog1, Pog2}:

\begin{thm}
\label{pog}
Let $G^1$ and $G^2$ be convex bodies in the Euclidean space $\E^3$ and $f: \partial G^1 \rightarrow \partial G^2$ be an isometry. Then $f$ extends to an isometry of $\E^3$. 
\end{thm}

By a \emph{convex body} in $\E^3$ we mean a compact convex set with non-empty interior. We endow its boundary with the \emph{induced path metric}: the distance between two points is equal to the length of a shortest curve connecting them on the boundary. It is well-known that for convex bodies, even without the smoothness assumption on the boundary, a shortest curve of finite length always exists (although it might be not unique).

Theorem~\ref{pog} was preceded by proofs under additional assumptions that the boundary is smooth (Cohn-Vossen~\cite{CoV} and Zhitomirsky~\cite{Zhi} in the analytic class; Herglotz~\cite{Her} in the $C^3$-class; Sacksteder~\cite{Sac} in the $C^2$-class) or that the boundary is polyhedral (Cauchy if the combinatorial structure is fixed; Alexandrov~\cite{Ale1, Ale2} in the general case). Pogorelov managed to obtain a common generalization of smooth and polyhedral cases.

Theorem~\ref{pog} is frequently stated as \emph{closed convex surfaces in $\E^3$ are globally rigid}. Here by a closed convex surface we mean the boundary of a convex body. Rigidity problems often go along with \emph{isometric realization problems}. One of the most famous instances is \emph{the Weyl problem} asking if every Riemannian metric on the 2-sphere with positive Gaussian curvature can be (smoothly) isometrically realized in $\E^3$. Such a realization should be a closed convex surface. Thus, the smooth instance of Theorem~\ref{pog} says that the realization is unique. The Weyl problem was resolved positively in collective efforts of Weyl~\cite{Wey} (a proof outline); Lewy~\cite{Lew} (the analytic class); Nirenberg~\cite{Nir} (the $C^4$-class); Heinz~\cite{Hei} (the $C^3$-class). It was extended to the case of non-negative Gaussian curvature by  Guan--Li~\cite{GuLi} and Hong--Zuily~\cite{HoZu}.

Meanwhile Alexandrov started to investigate the induced boundary metrics on general convex bodies. He managed~\cite{Ale1, Ale3} to find necessary and sufficient conditions on a metric on the 2-sphere to be realized as the boundary of a convex body. These are so-called \emph{CBB(0) metrics}. The reader can find a similar definition of a CBB($-1$) metric in Subsection~\ref{cbb-1}. First, Alexandrov resolved the problem in the polyhedral case, then he obtained a realization result in the general case by means of polyhedral approximation. About ten years later Pogorelov proved his general rigidity result. Besides, Pogorelov also proved~\cite{Pog2} that a $C^k$-smooth Riemannian metric with positive Gaussian curvature admits a $C^{k-1+\alpha}$-smooth realization for any $\alpha \in (0; 1)$. This allows to obtain a different resolution of the Weyl problem from Alexandrov--Pogorelov works. It is interesting to remark that in the smooth case the rigidity problem is commonly considered easier than the realization problem, but in the general case it seems to be on the contrary.

Pogorelov himself proposed two approaches to Theorem~\ref{pog}, both are quite lengthy and intricate. 
In both of them Pogorelov uses significantly the linear structure of $\E^3$. In the first approach Pogorelov supposes that there exist two convex bodies $G^1$, $G^2$ and an isometry $f: \partial G^1 \rightarrow \partial G^2$ that does not extend to an ambient isometry. Then he considers a supplementary surface defined with the help of the radius functions of $\partial G^1, \partial G^2$. He studies planar sections of the new surface and, with their help, constructs on the boundaries of initial bodies a pair of $f$-isometric simple closed curves bounding regions that need to be intrinsically isometric, but have non-equal generalized Gaussian curvature. This can not happen because Alexandrov proved that the generalized Gaussian curvature is intrinsic (this is a generalization of the Gauss's Theorema Egregium to non-smooth surfaces). In the second approach, from a violation of the global rigidity Pogorelov deduces a violation of the local rigidity and, finally, a violation of the infinitesimal rigidity. He shows that this is not possible. We refer also to the survey of Sen'kin~\cite{Sen} for a discussion of Pogorelov's proof and some further developments.

A natural problem is to quantify the rigidity in Theorem~\ref{pog}. In other words, one can ask: if the boundaries of $G^1$, $G^2$ are sufficiently close in the intrinsic sense, how close are $G^1$, $G^2$ themselves?

\begin{dfn}
Let $f: (M^1, d^1) \rightarrow (M^2, d^2)$ be a homeomorphism between metric spaces. It is called an $\e$-isometry if for any $p, q \in (M^1, d^1)$
$$|d^1(p,q)-d^2(f(p), f(q))|\leq \e.$$
\end{dfn}

\begin{prob}[The Cohn-Vossen problem~\cite{CoV2, Vol}]
\label{c-vprob}
Let $G^1$, $G^2$ be convex bodies in $\E^3$, $\e \geq 0$ and $f: \partial G^1 \rightarrow \partial G^2$ be an $\e$-isometry. Do there exist a continuous function $s: \R_{\geq 0} \rightarrow \R_{\geq 0}$ such that $s(0)=0$ and a constant $C$ that depends only on global geometry of $\partial G^1$ (e.g., its diameter) such that $f$ extends to a $Cs(\e)$-isometry between $G^1$ and $G^2$?
\end{prob}

Volkov claimed a solution to this problem in~\cite{Vol} (for $s(\e)=\e^{1/24}$). An English translation was published as an appendix to the book~\cite{Ale2}. Volkov's proof is completely unrelated to both Pogorelov's approaches, so, in particular, it gives a third way to prove Theorem~\ref{pog}. The paper~\cite{Vol} is highly innovative without any doubts. However, we have some concerns on it. There are some gaps that can be bridged after some non-trivial work. But more important, some essential arguments are written in a very brief and cryptic way making them hard to interpret. We describe our remarks on Volkov's paper in Appendix~\ref{appendix}.

The present paper can be considered as a revival of paper~\cite{Vol} with an application to a rigidity problem for a family of hyperbolic 3-manifolds. We use many important ideas from~\cite{Vol} and its influence on this paper can not be discounted. But in some steps of the proof we have to choose a different road due to our inability to understand some parts of~\cite{Vol}. We think that a proof of Theorem~\ref{pog} can be obtained by the means of the present manuscript, but they are insufficient to resolve Problem~\ref{c-vprob}.

\subsection{Rigidity of hyperbolic 3-manifolds}

The second source of our motivation comes from the importance of hyperbolic 3-manifolds to 3-dimensional topology highlighted in breakthrough works of Thurston. We quickly recall that Thurston proposed \emph{the geometrization program} to understand the topology of 3-manifolds: each 3-manifold can be canonically decomposed into pieces and each piece can be endowed with one of eight canonical geometries (see, e.g. \cite{Thu1, Thu2}). The geometrization program was deeply developed by plenty of researchers with its notable culmination in the works of Perelman. Among all 3-manifolds admitting canonical geometries, hyperbolic manifolds constitute the largest and the most mysterious class. Hence, much of the efforts of researchers were directed toward developing a deeper understanding of hyperbolic 3-manifolds.

One of the foundational properties of hyperbolic manifolds (closed or of finite volume) in dimensions starting from 3 is their rigidity. It was first proved by Mostow~\cite{Mos} for closed manifolds and was extended to the case of finite volume by Prasad~\cite{Pra}. 

\begin{thm}
\label{mos}
Let $f: M^1 \rightarrow M^2$ be a homotopy equivalence between two hyperbolic $n$-dimensional manifolds of finite volume (without boundary) and $n \geq 3$. Then $f$ is homotopic to an isometry.
\end{thm}

One can loosely rephrase it by saying that the geometry of a hyperbolic manifold of finite volume without boundary in dimension at least 3 is completely determined by its topology. It is natural to investigate hyperbolic manifolds with boundary in a similar vein. From the Pogorelov rigidity and the Mostow rigidity one can expect that in the case of convex boundary, the geometry of a hyperbolic 3-manifold should be determined by the topology and the induced boundary metric. In the topologically simplest case when the manifold is the 3-ball a proof was outlined by Pogorelov~\cite{Pog2} and finalized by Milka~\cite{Mil} (see the introduction of~\cite{Mil} for the discussion of proof attempts).

\begin{thm}
\label{pogh}
Let $G^1$ and $G^2$ be convex bodies in $\H^3$ and $f: \partial G^1 \rightarrow \partial G^2$ be an isometry. Then $f$ extends to an isometry of $\H^3$. 
\end{thm}

The proof of Theorem~\ref{pogh} appears to be quite different from the proof of Theorem~\ref{pog} given by Pogorelov, probably, because the latter seems to be quite specific to $\E^3$. For the hyperbolic case Pogorelov proposed a map $\H^3\times \H^3 \rightarrow \E^3 \times \E^3$ that produces a correspondence between pairs of convex surfaces in both spaces and allows to relate the rigidity problem in $\H^3$ to the corresponding problem in $\E^3$.  Pogorelov developed all necessary tools and completed the proof of rigidity in the spherical space, but the details of the proof in the hyperbolic space were furnished in further works, see~\cite{Mil}. 

Next, Schlenker~\cite{Sch4} proved the rigidity of hyperbolic 3-manifolds with smooth strictly convex boundary. Here by strict convexity we mean that the shape operator is positive-definite for the outward choice of the unit normal field. By the Gauss equation this implies that the Gaussian curvature of the boundary is strictly greater than $-1$. We rephrase his result (the uniqueness part of Theorem 0.1 in~\cite{Sch}) as follows:

\begin{thm}
\label{sch}
Let $M^1$ and $M^2$ be compact hyperbolic 3-manifolds with smooth strictly convex boundaries and $f:  M^1 \rightarrow  M^2$ be a homeomorphism such that its restriction to the boundary is isotopic to an isometry. Then $f$ is isotopic to an isometry.
\end{thm}


Previously Labourie~\cite{Lab2} proved the corresponding realization result: if the interior of a compact 3-manifold $M$ with non-empty boundary admits a cocompact hyperbolic metric, then each smooth metric with Gaussian curvature greater than $-1$ on $\partial M$ is induced by a hyperbolic metric on $M$. Schlenker provided another proof of this in~\cite{Sch4}. 

Due to Pogorelov's global rigidity for general bodies, it is natural to expect that the smoothness and strict convexity assumptions are superfluous. One can conjecture

\begin{cnj}
\label{cnj1}
Let $M^1$ and $M^2$ be compact hyperbolic 3-manifolds with convex boundaries and $f:  M^1 \rightarrow  M^2$ be a homeomorphism such that its restriction to the boundary is isotopic to an isometry. Then $f$ is isotopic to an isometry.
\end{cnj}

Here by convex boundary we mean that the boundary is locally modelled on convex subsets of $\H^3$.

The aim of the present paper is to prove this conjecture for a particular family of hyperbolic 3-manifolds with boundary.

\begin{dfn}
\label{fuchsian}
A \emph{compact Fuchsian manifold with boundary} is a hyperbolic 3-manifold homeomorphic to $S_g \times [0; 1]$, where $S_g$ is a closed oriented surface of genus $g>1$, such that $S_g \times \{ 0\}$ is geodesic.
\end{dfn}

In what follows we will omit the word ``compact'' assuming that this is always the case. We will refer to the boundary component $S_g \times \{0\}$ as to the \emph{lower boundary} and to the component $S_g \times \{1\}$ as to the \emph{upper boundary}, and denote them by $\partial_{\downarrow} F$ and $\partial^{\uparrow} F$ respectively, where $F$ is a Fuchsian manifold with boundary.

Fuchsian manifolds with boundary are considered as toy cases in the study of hyperbolic 3-manifolds with boundary, and the interest to them can be traced back to the works of Pogorelov~\cite[Section VI.12]{Pog2} and Gromov~\cite[Section 3.2.4]{Gro}. In particular, Gromov proved the smooth realization result: every smooth metric on $S_g$ of Gaussian curvature $>-1$ is induced on the upper boundary of a Fuchsian manifold. He conjectured that this realization is unique. 

We note that some authors use slightly different definition of a Fuchsian manifold with boundary: they require it to be homeomorphic to $S_g \times [0;1]$ and to contain an embedded geodesic surface isotopic to $S_g \times \{p\}$. If we cut such a manifold along this surface, we obtain two Fuchsian manifolds in the sense of Definition~\ref{fuchsian}. The rigidity of Fuchsian manifolds in the sense of the other definition is equivalent to the rigidity of Fuchsian manifolds from Definition~\ref{fuchsian} with respect to the induced metric on the upper boundary. We find it more convenient to work in the setting of Definition~\ref{fuchsian}.


In the case of polyhedral Fuchsian manifolds Fillastre~\cite{Fil1} proved

\begin{thm}
\label{fil}
Let $F^1$ and $F^2$ be two Fuchsian manifolds with convex polyhedral boundaries and $f: \partial^{\uparrow} F^1 \rightarrow \partial^{\uparrow} F^2$ be an isometry between the upper boundaries. Then $f$ extends to an isometry between $F^1$ and $F^2$.
\end{thm} 

In the case of smooth boundary the same result follows from Theorem~\ref{sch} by the doubling construction. 

Our main result is 

\begin{thml}
\label{main}
Let $F^1$ and $F^2$ be two Fuchsian manifolds with convex boundaries and $f: \partial^{\uparrow} F^1 \rightarrow \partial^{\uparrow} F^2$ be an isometry between the upper boundaries. Then $f$ extends to an isometry between $F^1$ and $F^2$.
\end{thml}

Theorem~\ref{main} implies the following particular case of Conjecture~\ref{cnj1}:

\begin{crl}
Let $M^1$ and $M^2$ be compact hyperbolic 3-manifolds with convex boundaries such that $M^1, M^2$ are homeomorphic to $S_g \times [0; 1]$ and contain geodesic surfaces isotopic to $S_g \times \{p\}$. Let $f:  M^1 \rightarrow  M^2$ be a homeomorphism such that its restriction to the boundary is isotopic to an isometry. Then $f$ is isotopic to an isometry.
\end{crl}

This seems to be the first known rigidity result for general boundary metrics and manifolds more complicated than the ball.

\subsection{Related work}

As we already mentioned, our proof of Theorem~\ref{main} follows the ideas of Volkov. It proceeds by means of polyhedral approximation. However, in order to prove the rigidity by approximation, the rigidity in the polyhedral case is insufficient and one needs to obtain a kind of stability. Assume, e.g., that we have two convex bodies in $\E^3$ that are not ambiently isometric, but have isometric boundaries. We can approximate them by polyhedra. In this way we can obtain two polyhedra that are not close modulo ambient isometries, but have arbitrarily close boundaries. One needs to prove that this can not happen. Volkov's main insight was that this can be achieved by deforming polyhedra through polyhedra with cone-singularities in the interior and by controlling the process with the help of discrete curvature. 

It seems that Volkov's PhD thesis from 1955, recently published as~\cite{VolThe}, was the first place where cone-3-manifolds and their discrete curvature was considered. He provided there a variational proof of Alexandrov's realization and rigidity theorem on Euclidean polyhedra. Volkov also produced a similar proof in the case of so-called \emph{convex polyhedral caps}~\cite{Vol2}. One can find an English translation in the appendix to the book~\cite{Ale2}.

Now we describe some other works that were invaluably impacted by Volkov's ideas.



In~\cite{BoIz} Bobenko and Izmestiev elaborated the ideas of Volkov further and gave another variational proof of the Alexandrov theorem. Based on their work Sechelmann developed~\cite{Sec} a computer program that constructs the polyhedron realizing a given Euclidean cone-metric. This resolved a long-standing problem (original Alexandrov's proof was non-constructive). Izmestiev~\cite{Izm1} also gave a new proof of the realization and rigidity of convex polyhedral caps, which provides a simple introduction to the approach. An investigation of convex polyhedral caps in the Minkowski space was done by Milka~\cite{Mil3} in a way very similar to Volkov. In~\cite{FiIz1},~\cite{FiIz2} Fillastre and Izmestiev used these techniques for hyperbolic and spherical cone-metrics on the torus (realizing them in some hyperbolic and de Sitter manifolds respectively). 

Alexandrov also proved a realization and rigidity result for convex polyhedra in $\H^3$. One might also consider (non-compact) polyhedra in $\H^3$ with \emph{ideal} vertices on the boundary at infinity of $\H^3$ and with \emph{hyper-ideal} vertices outside of it (in the projective model). 
Fillastre~\cite{Fil2} proved an Alexandrov-type result for generalized Fuchsian polyhedra with vertices of these types. In~\cite{Pro1} we investigated ideal Fuchsian polyhedra with the help of cone-manifolds and the discrete curvature. This is of particular interest due to its connection to discrete conformality. See more details, e.g., in~\cite{BPS, GGLSW, Spr, Pro1}. 

The realization counterpart to Theorem~\ref{main} comes from the work of Slutskiy~\cite{Slu}. He considered \emph{compact quasi-Fuchsian manifolds with convex boundary}. They are also homeomorphic to $S_g \times [0;1]$, but the condition on the geodesic boundary is not imposed. The induced metric on a convex surface in a hyperbolic 3-manifold is CBB($-1$). Slutskiy proved that for any pair of CBB($-1$) metrics on $S_g$ there exists a compact quasi-Fuchsian manifold with convex boundary realizing them. He used a smooth approximation and the smooth realization result of Labourie--Schlenker mentioned above. If one replaces the latter by the above mentioned smooth Fuchsian realization of Gromov, then one can obtain the general convex Fuchsian realization result without any difficulties. A similar realization result for CBB($-1$) metrics on the torus was obtained in~\cite{FIV}.

A related and very active research area is the convex realization of metrics and their rigidity in Lorenzian space-forms. We refer to~\cite{Sch2, LaSc, Fil3, FiSl, Bru, Lab, Tam} for works in this direction.

In 70s cone-3-manifolds were rediscovered by Thurston who proposed to use them in order to describe deformations of hyperbolic structures. He demonstrated their usability in his proof of \emph{the hyperbolic Dehn filling theorem}~\cite[Chapter 4]{Thu1}. Hyperbolic cone-3-manifolds were used heavily by Thurston's school of 3-dimensional topology. We refer to~\cite{BLP, HoKe2, Bro3} as just to few examples. 

Mostow's rigidity was also generalized to open hyperbolic 3-manifolds. By an open manifold we mean a non-compact connected 3-manifold without boundary. A subset $F$ of a hyperbolic manifold $M$ is called \emph{totally convex} if $F$ contains every geodesic segment between any two points of $F$. The \emph{convex core} of a hyperbolic 3-manifold $M$ is the intersection of all closed totally convex subsets of $M$. It is non-empty except for a specific family of examples. If $M$ is closed or $M$ is open, but has finite volume, then its convex core is $M$ itself. 

The convex core is an important tool to investigate open hyperbolic 3-manifolds. An open hyperbolic 3-manifold is called \emph{cocompact} if its convex core is compact. It is called \emph{geometrically finite} if the convex core has finite volume. A geometrically finite hyperbolic-3-manifold can be naturally compactified and its boundary at infinity can be endowed with a conformal structure. Marden's rigidity theorem~\cite{Mard} (based on the previous works of Ahlfors, Bers, Kra and Maskit) states that the geometry of a geometrically finite hyperbolic 3-manifold is completely determined by its topology and the conformal structure at infinity. To cover a more general case of all open hyperbolic 3-manifolds with finitely generated fundamental group one needs additional invariants called \emph{ending laminations}. This was established in~\cite{Min, BCM}.

\subsection{Further work directions}

A good direction of further research is Conjecture~\ref{cnj1}. In particular, after Fuchsian manifolds one may attempt to prove it for \emph{compact quasi-Fuchsian manifolds with convex boundary} that were defined above. In this case we have to prescribe the metric on both boundary components. Fuchsian manifolds exhibit a warped product geometry, which quasi-Fuchsian manifolds lack. While remaining topologically still simple, quasi-Fuchsian manifolds already contain all geometric difficulties that appear when we try to proceed further from the Fuchsian case. It might happen that there are not much difficulties to extend the solution from the quasi-Fuchsian case to a more general topological case. 

A nice application of Conjecture~\ref{cnj1} should be another way to determine the geometry of open hyperbolic 3-manifolds. A long standing conjecture is

\begin{cnj}
\label{coco}
The geometry of a geometrically finite hyperbolic 3-manifold is completely determined by the topology and the induced path metric on the boundary of its convex core.
\end{cnj}

The induced metric on the boundary of the convex core is  totally hyperbolic, i.e., locally isometric to the hyperbolic plane. However, it is embedded in a non-smooth way and is bent along a geodesic lamination. This lamination together with an additional data (a transverse measure describing how much the surface is bent) is called \emph{a pleating lamination}. The \emph{pleating lamination conjecture} states that this data provides another way to determine the geometry of a geometrically finite hyperbolic 3-manifold. In some way it is dual to Conjecture~\ref{coco}.

Another important direction is to recover a proof of the Cohn-Vossen problem (Problem~\ref{c-vprob}). As we mentioned in Section~\ref{Weyl}, the Cohn-Vossen problem was presumably solved by Volkov in~\cite{Vol}, but we can not reconstruct some important steps of his proof. We are also unaware of any researchers that claim to understand this proof. One can also formulate the Cohn-Vossen problem, e.g., for Fuchsian manifolds:

\begin{prob}
Let $F^1$, $F^2$ be Fuchsian manifolds with convex boundaries, $\e \geq 0$ and $f: \partial^{\uparrow} F^1 \rightarrow \partial^{\uparrow} F^2$ be an $\e$-isometry. Does $f$ extend to a $Cs(\e)$-isometry between $F^1$ and $F^2$, where $C$ depends on global geometry of $\partial^{\uparrow} F^1$ and $s$ is a continuous function satisfying $s(0)=0$?
\end{prob}


\vskip+0.2cm

\textbf{Acknowledgments.} This work is a part of author's doctoral thesis completed under the supervision of I.~Izmestiev and the author is invaluably grateful to him for his constant attention and advice. The author is also very grateful to F.~Fillastre and M.~Ghomi for their helpful remarks.

\vskip+0.2cm

Unfortunately, our paper is incredibly long so we tried to explain the links between sections as clear as possible. 
In Section~\ref{prelim} we go through preliminaries from metric geometry that we will use. Subsections~\ref{cbb-1}--\ref{polmetr} contain mostly standard material, although we make an emphasis on some specific details that will be important for us. 
In Subsection~\ref{tp} we define \emph{trapezoids} and \emph{prisms} that are our elementary building blocks. With their help in Subsection~\ref{definfcm} we define \emph{Fuchsian cone-manifolds} that are the main heroes of our proof. In Subsection~\ref{disccurv} we formulate some facts about them that are necessary for the proof and define the \emph{discrete curvature functional}. 
In Section~\ref{generalrig} we prove Theorem~\ref{main} modulo several tough lemmas concerning the stability of Fuchsian manifolds with polyhedral boundary. 
Next, in Section~\ref{conechap} we investigate thoroughly properties of Fuchsian cone-manifolds and of the discrete curvature that we use. Section~\ref{stability} is devoted to the proofs of the main lemmas from  Section~\ref{generalrig}. This section is the core of the paper.


\section{Preliminaries}
\label{prelim} 

\subsection{CBB($-1$) metrics}
\label{cbb-1}

First, we briefly sketch the facts that we need about metrics with curvature bounded from below by $-1$ in the sense of Alexandrov (\mbox{CBB($-1$)} for short). For a detailed exposition of CBB($k$) metrics we refer to~\cite{BBI, AKP}. Many properties of CBB($k$) metrics on surfaces are similar to those of CBB(0) metrics, which are treated in details in~\cite{Ale3} (with Chapter XII discussing also CBB($k$) case).

Let $S$ be a connected orientable surface and $d$ be a complete intrinsic metric on $S$. By intrinsic we mean that the distance between each pair of points is equal to the infimum of lengths of all rectifiable paths connecting them. Then this infimum is achieved: there exists a shortest path between any two points. This is a corollary of the Arzela--Ascoli theorem, see~\cite[Theorem 2.5.23]{BBI}.

\begin{dfn}
A \emph{geodesic} is a rectifiable curve (possibly closed) that is locally distance minimizing.
\end{dfn}

Let $\psi$, $\chi$ be two shortest paths in $(S, d)$ sharing an endpoint $p$. Let $q \in \psi$ be the point at distance $x$ from $p$ and $r \in \chi$ be the point at distance $y$ from $p$. Consider the hyperbolic triangle with side lengths $x$, $y$ and $ d(q,r)$ and let $\lambda(x,y)$ be the angle opposite to the side of length $ d(q,r)$. 

\begin{dfn}
\label{cbbdef}
We say that $d$ is a \emph{CBB($-1$) metric} on $S$ if $d$ is complete, intrinsic and for each $p \in S$ there exists a neighbourhood $U \ni p$ such that the function $\lambda(x,y)$ is a nonincreasing function of $x$ and $y$ for every $\psi$, $\chi$ emanating from $p$, in the range $x \in [0; x_0]$, $y \in [0; y_0]$ where the respective points $q, r$ belong to $U$.
\end{dfn}

For other definitions of CBB($k$)-spaces we refer to~\cite[Chapter 4]{BBI} and~\cite[Chapter VII]{Ale3}. 

Definition~\ref{cbbdef} implies that the angle between $\psi$ and $\chi$, which we define as $\lim\limits_{x,y \rightarrow 0} \lambda(x,y)$, exists. Denote it by $\ang(\psi, \chi,d)$.
In this way the angle could be defined between any two rectifiable curves sharing an endpoint, which are possibly not geodesics, but the angle might not exist. 

It is important to note that the locality in Definition~\ref{cbbdef} can be dropped. Namely, for any three points $p, q, r \in (S,d)$ the angle between shortest paths from $p$ to $q$ and to $r$ is at least the respective angle in the \emph{comparison triangle} for $p$, $q$, $r$, i.e., in the hyperbolic triangle with the side lengths $d(p,q)$, $d(q,r)$ and $d(p,q)$. This is called \emph{the Toponogov globalization theorem}, which was proved in the general case by Perelman. We refer to~\cite[Theorem 10.3.1]{BBI}. 

For geodesics in CBB($-1$) spaces the \emph{non-overlapping property} holds: if two geodesics have a segment in common, they can be covered by a larger geodesic. Sometimes we need a slightly stronger version: if two shortest paths have two points in common, then either these are their endpoints or they have a segment in common.


The shortest paths $\psi$ and $\chi$ emanating from a point $p$ divide a sufficiently small neighbourhood of $p$ into two \emph{sectors} $U$ and $U'$. Let $\psi_1, \ldots, \psi_k$ be shortest paths emanating from $p$ belonging to $U$ enumerated in the order from $\psi$ to $\chi$. The angle $\ang(U, d)$ is defined as the supremum of the sums $\ang(\psi, \psi_1,d)+\ldots+\ang(\psi_k, \chi,d)$ over all finite collections of shortest paths from $p$ in $U$. The smallest of $\ang(U, d)$, $\ang(U', d)$ is equal to $\ang(\psi, \chi,d)$. 

\begin{dfn}
The \emph{total angle} $\lambda_p( d)$ of $p$ is equal to $\ang(U,d)+\ang(U',d)$.
\end{dfn}

It does not depend on the choice of initial shortest paths $\psi$, $\chi$.

\begin{dfn}
\label{convexpolyg}
A \emph{geodesic polygon} is a submanifold of $(S, d)$ with piecewise geodesic boundary. It is called \emph{convex} if there is a shortest path between any two of its points that belongs to the polygon. 
\end{dfn}

It might be worth to note~\cite[Chapter II.6]{Ale3}:

\begin{lm}
Assume that $(S, d)$ is a compact CBB($-1$) surface. Then it admits a triangulation consisting of finitely many arbitrarily small convex geodesic triangles.
\end{lm}

The area of a Borel set in $(S,  d)$ is defined intrinsically as its Hausdorff measure. See~\cite[Chapters 1.7, 2.6]{BBI}. 


We sketch the concept of the \emph{intrinsic curvature} $\nu_I$ of a Borel set in $(S, d)$. 

For a point $p \in (S,  d)$ we define $\nu_{I}(p,d):=2\pi-\lambda_p( d)$. 

For a relatively open geodesic segment $\psi$ we define $\nu_I(\psi):=0$. 

For an open geodesic triangle $T$ we define its curvature $\nu_I(T,  d):=\alpha+\beta+\gamma-\pi$, where $\alpha$, $\beta$, $\gamma$ are the angles of $T$.

These three types of sets are called \emph{primitive sets}. An \emph{elementary set} is a set that can be represented as a finite disjoint union of primitive sets. Then its curvature is the sum of the curvatures of these primitive sets. It does not depend on a representation. Then for a closed set in $(S,  d)$ its curvature is defined as the infimum of the curvatures of its elementary supersets. For an open set its curvature is defined as the supremum of the curvatures of its closed subsets. This defines a Borel measure on $(S, d)$: for the details we refer to~\cite[Chapters V, XII.1.6]{Ale3}, \cite[Chapter V]{AlZa} or~\cite{Mac}.

However, we will mostly use the \emph{extrinsic curvature} of a Borel set $B$, namely, $$\nu(B,d):=\nu_I(B,d)+\area(B, d).$$ For CBB($-1$) metrics it is non-negative: see~\cite{Mac}. In what follows we will omit the word extrinsic. 

\subsection{Hyperbolic convex bodies and duality}
\label{duality}

In this subsection we introduce basic facts from convex geometry in $\H^3$ that we will use. Let $G$ be a closed convex set in $\H^3$ with non-empty interior distinct from $\H^3$. Then its boundary $\partial G$ is homeomorphic to an open subset of the sphere $S^2$. First, we recall a fundamental result from convex geometry~\cite{Mil2}:

\begin{lm}[The hyperbolic Busemann-Feller lemma]
Let $p, q \in  \H^3$ and $p', q'$ be the nearest points to them in $G$. Then the hyperbolic distance between $p'$ and $q'$ is at most the hyperbolic distance between $p$ and $q$.
\end{lm}

\begin{crl}
\label{BusFel}
Let $\psi \subset \H^3 \backslash G$ be a rectifiable curve and $\psi'$ be its nearest point projection to the boundary of $G$. Then $\psi'$ is rectifiable and its length is at most the length of $\psi$.
\end{crl}

Another well-known result is~\cite[Chapter XII]{Ale3}:

\begin{lm}
The boundary $\partial G$ equipped with the induced path metric is CBB($-1$).
\end{lm}

In particular, the curvature measure $\nu$ is defined for $\partial G$ as in Section~\ref{cbb-1}.

We will also use the hyperbolic -- de Sitter duality. The \emph{de Sitter space} $\dS^3$ is the space of oriented planes of $\H^3$. It is a Lorenzian manifold of constant curvature 1. More precisely, consider the hyperboloid model. Let $\R^{1,3}$ be the 4-dimensional Minkowski space, i.e., equipped with the scalar product $$\langle   x,   y \rangle := -x_0y_0+x_1y_1+x_2y_2+x_3y_3.$$ We identify $$\H^3 = \{  x \in \R^{1,3}: \langle   x,  x\rangle = -1,~x_0 > 0 \}$$
and define the three-dimensional de Sitter space as
$$\ds^3=  \{  x \in \R^{1,3}: \langle   x,  x\rangle = 1 \}.$$

An oriented plane in $\H^3$ is obtained as an intersection of $\H^3$ with an oriented hyperplane in $\R^{1,3}$ passing through the origin. Naturally, the unit normal to this hyperplane belongs to $\ds^3$. There is a well-known duality between convex sets is $\H^3$ and in $\dS^3$. For a convex set $G \subset \H^3$ we define the dual convex set $G^* \subset \dS^3$ as the set of all planes that do not intersect $\inter(G)$ and are oriented outwards $G$. We refer to~\cite{BeCa, FiSe} for more details.

For each Borel set $U \subset \partial G$ define its dual $U^* \subset \partial G^*$ as the set of all planes tangent to $G$ and passing through points of $U$. It is folklore that
$$\nu(U)=\area(U^*).$$
However, we are unaware of any sources that prove this for non-smooth convex bodies, hence, we are providing a proof by ourselves.

We consider the case when $G$ is a convex body in $\H^3$, i.e. compact, convex and, as before, with non-empty interior. We follow the framework from~\cite{BeCa}. Define $$\dS^3_+:=\{ x\in \dS^3: x_0>0\}.$$
Assume without loss of generality that the point $o=(1,0,0,0)$ is in the interior of $G$. Consider the cone $\mathcal C(G)$ in $\R^{1,3}$ from the origin over $G$. Define its dual cone 
$$\mathcal C(G)^*:=\{x \in \R^{1,3}: \forall y \in \mathcal C(G),~~ \langle x, y \rangle \leq 0 \},$$
which is a convex cone containing the future cone of $\R^{1,3}$ in its interior. Then  $$G^*=\mathcal C(G)^* \cap \dS^3_+$$ (see the details in \cite{BeCa}). 

Define the Gauss map $\mathcal G: \partial G \substack{\rightarrow\\[-1em]  \rightarrow} \partial G^*$ as the multivalued map sending a point on $\partial G$ to the set of its outward unit normals. The set $\partial G^*$ is space-like and comes with a well-defined area measure $\sigma_{\partial G^*}$. This is proven in~\cite[Lemma 2.1]{BeCa}.

Consider $S^2$ as the unit sphere in $T_o\H^3$, let $\mathcal P: S^2 \rightarrow \partial G$ be the radial projection, i.e., a map that sends a point from the sphere to the endpoint on $\partial G$ of the geodesic in the respective direction, and let $r: S^2 \rightarrow \R_{>0}$ be the function measuring the length of the this geodesic. We pull back $\sigma_{\partial G^*}$ to $S^2$ via $\mathcal G \circ \mathcal P$ and denote the obtained measure by $\mu$. We also pull back the area measure $\sigma_{\partial G}$ to $S^2$ via $\mathcal P$ and denote it by $\sigma$ (note that it is not the same $\sigma$ as in \cite{BeCa}).

We pull back our curvature measure $\nu$ to $S^2$ via $\mathcal P$ and, abusing the notation, continue to denote it by $\nu$.

\begin{lm}[\cite{BeCa}, Proposition A.3]
There exists a sequence $G_k$ of convex bodies such that

(1) $G_k$ are smooth and strictly convex;

(2) $r_k$ converge to $r$ uniformly;

(3) $G_k$ converge to $G$ in the Hausdorff sense;

(4) the measures $\mu_k$ converge weakly to $\mu$.
\end{lm}

Now we are ready to prove

\begin{lm}
\label{duallemma}
The measures $\nu$ and $\mu$ coincide. Thus, for each Borel set $U \subset \partial G$
$$\nu(U)=\area(U^*).$$
\end{lm}

\begin{proof}
First, we claim that $\nu_k=\mu_k$. Indeed, for $p \in S^2$ let $K_I(p)$ be the Gaussian curvature of $\mathcal P(p) \in \partial G$ and $K(p)$ be the extrinsic curvature of $\mathcal P(p)$, i.e., the determinant of the shape operator. The Gauss equation says that $K_I(p)+1=K(p)$. Proposition 2.2.1 from \cite{BeCa} states that $\mu_k = K_k d\sigma_k$. We claim that also $\nu_k=K_k d \sigma_k$. Due to the Gauss equation, it is enough to show that $\nu_{I, k}=K_{I,k} d \sigma_k$. Take an open geodesic triangle $T$. The Gauss-Bonnet theorem shows that $\nu_{I, k}(T)=\int_T K_{I, k} d\sigma_k $. The measures of other elementary sets (geodesic arcs and points) are zero for smooth metrics. The proof of $\nu_k=\mu_k$ for any Borel set follows from the definition of $\nu_I$ and elementary properties of the Hausdorff measure $d\sigma$. As $\nu_k=\nu_{I, k}+\sigma_k$, we get $\nu_k=K_k d \sigma_k=\mu_k$.

Theorem 7 from \cite[Chapter 7]{AlZa} says that $\nu_{I, k}$ converge weakly to the intrinsic curvature measure $\nu_I$ of $\partial G$. Theorem 9 from~\cite[Chapter 8]{AlZa} says that the area measures $\sigma_k$ converge weakly to $\sigma$. Thus, $\nu_k$ converge weakly to $\nu$ and, as $\mu_k$ converge weakly to $\mu$, we get $\nu=\mu$.
\end{proof}

\subsection{Cone-metrics and triangulations}
\label{polmetr}

Let $S_g$ be a closed oriented connected surface of genus $g$. In this subsection we deal with hyperbolic cone-metrics on $S_g$ and their geodesic triangulations.

\begin{dfn}
A \emph{topological triangulation} of $S_g$ is a finite vertex set $V$ and a collection of simple disjoint paths with endpoints in $V$ that cut $S_g$ into triangles with vertices in $V$. Two topological triangulations with the same vertex set $V$ are \emph{equivalent} if they are isotopic with respect to $V$ (so the edges are not allowed to pass through points of $V$ during the isotopy). A \emph{triangulation} $\mathcal T$ of $S_g$ is an equivalence class of topological triangulations with the same vertex set $V$.
\end{dfn}

Note that this definition allows loops and multiple edges between two vertices. By $V(\mathcal T)$ we denote the set of vertices of $\mathcal T$ and by $E(\mathcal T)$ we denote the set of edges of $\mathcal T$ considered as isotopy classes of the respective paths.

\begin{dfn}
A \emph{hyperbolic cone-metric} $d$ on $S_g$ is locally isometric to the metric of hyperbolic plane except finitely many points called \emph{conical points}. At a conical point $v$ the metric $d$ is locally isometric to the metric of a hyperbolic cone with angle $\lambda_v(d)\neq 2\pi$. The number $\nu_v(d):=2\pi-\lambda_v(d)$ is called the \emph{curvature} of $v$. We denote the set of conical points of $d$ by $V(d)$. A hyperbolic cone-metric is called \emph{convex} if for every $v \in V(d)$ we have $\lambda_v(d) < 2\pi$. 
Two metrics $d_1$ and $d_2$ are considered equivalent if $V(d_1)=V(d_2)=V$ and there is an isometry fixing $V$ and isotopic to identity with respect to $V$.
\end{dfn}


\begin{lm}[\cite{Ale3}, Chapter XII]
A convex  hyperbolic cone-metric is CBB($-1$).
\end{lm}

From now on we restrict ourselves almost only to the hyperbolic case and omit the word ``hyperbolic'' saying just cone-metrics (except in some special cases).

We note that on a cone-metric there can be multiple shortest paths between two points.

\begin{dfn}
A \emph{geodesic triangulation} is a topological triangulation such that all edges are geodesics. 
\end{dfn}

\begin{dfn}
Let $\mathcal T$ be a triangulation of $S_g$ and $d$ be an intrinsic metric. We say that $\mathcal T$ is \emph{realized} by $d$ if there is a geodesic triangulation of $(S_g, d)$ in the class~$\mathcal T$.
\end{dfn}

Note that if $d$ is a cone-metric, then we do not require in this definition neither $V(\mathcal T) \subseteq V(d)$ nor the converse. We highlight that degenerated triangles are not allowed because the edges are defined up to isotopy with respect to $V(\mathcal T)$. Sometimes when we write $e \in E(\mathcal T)$ and the realization of $\mathcal T$ is evident, then we mean the respective realization of $e$. Similarly, when we consider a triangle $T$ of $\mathcal T$ we frequently mean its realization in a metric.

We recall that a geodesic on a convex cone-metric can not pass through conical points. We will frequently need the following result, which follows from~\cite[Proposition~4]{Izm2}:

\begin{lm}
\label{iztriang}
If $d$ is a cone-metric and ${V\supseteq V(d)}$, then there exists a triangulation $\mathcal T$ realized by $d$ with $V(\mathcal T)=V$. Moreover, any set of disjoint geodesic paths with vertices in $V$ can be extended to a geodesic triangulation.
\end{lm}

In this case the metric $d$ is uniquely determined by $\mathcal T$ and the edge lengths.

\begin{dfn}
Fix $V \subset S_g$. By $\mathfrak D(V)$ we denote \emph{the space of cone-metrics} $d$ on $S_g$ with $V(d) \subseteq V$  up to isometry isotopic to identity with respect to $V$. By $\mathfrak D_c(V)\subset \mathfrak D(V)$ we denote \emph{the subspace of convex metrics}. By $\mathfrak D_{sc}(V)\subset \mathfrak D_c(V)$ we denote \emph{the subspace of strictly convex metrics} with respect to $V$, i.e., $d \in \mathfrak D_{sc}(V)$ if and only if $V(d)=V$ and $d$ is convex.
\end{dfn}

\begin{dfn}
Let $\mathcal T$ be a triangulation of $S_g$ with $V(\mathcal T)=V$. By $\mathfrak D(V, \mathcal T) \subset \mathfrak D(V)$ we denote \emph{the set of cone-metrics realizing $\mathcal T$}. Similarly we denote its subsets $\mathfrak D_c(V, \mathcal T)$ and $\mathfrak D_{sc}(V, \mathcal T)$. 
\end{dfn}

Recall that if $n:=|V|$, then any triangulation of $S_g$ with vertices at $V$ has $2(n+2g-2)$ triangles and $N=3(n+2g-2)$ edges. The edge lengths map $l: {\mathfrak D(V, \mathcal T) \rightarrow \R^{N}}$ is injective. Abusing the notation, we frequently identify $\mathfrak D(V, \mathcal T)$ with its image under this map. Let us study its basic properties.

The set $\mathfrak D(V, \mathcal T)\subset \R^{N}$ is an open polyhedron defined by the strict triangles inequalities for all triangles of $\mathcal T$. For every $v \in V$ the total angle of $v$ defines an analytic function $\lambda_v: \mathfrak D(V, \mathcal T)\rightarrow (0; \infty)$. Then $\mathfrak D_c(V, \mathcal T)$ is the subset of $\mathfrak D(V, \mathcal T)$ satisfying inequalities $\lambda_v(d) \leq 2\pi$. It is a semi-analytic set.

If $d \in \mathfrak D(V)$ realizes two triangulations, then the transition maps are smooth. This endows $\mathfrak D(V)$ with a smooth manifold structure. The set $\mathfrak D_{sc}(V)$ is its open subset. 

\begin{lm}
\label{decreaseangles}
Let $d \in \mathfrak D(V, \mathcal T)$. Define a 1-parameter family of cone-metrics $d_t\in \mathfrak D(V, \mathcal T)$ by $l_e(d_t)=t\cdot l_e(d)$ for each $e \in E(\mathcal T)$, where $t \in [1; +\infty)$. Then $\lambda_v(d_t)$ is strictly monotonously decreasing for every $v\in V$. 
\end{lm}

\begin{proof}
It is clear that all strict triangle inequalities are satisfied after multiplying by $t$, therefore equations $l_e(d_t)=t\cdot l_e(d)$ indeed define a metric $d_t \in \mathfrak D(V, \mathcal T)$. It remains to prove that if all edge lengths of a hyperbolic triangle $ABC$ are multiplied with the same factor $t>1$, then its angles become strictly smaller.

Let $A'B'C'$ be a triangle with the increased side lengths. Consider $B''\in A'B'$ and $C''\in A'C'$ such that $A'B''=AB$ and $A'C''=AC$. It suffices to show that $B''C'' < BC$.

Now let $A'_0B'_0C'_0$ be a Euclidean comparison triangle for $A'B'C'$ and $B''_0 \in A'_0B'_0$, $C''_0 \in A'_0C'_0$ be points such that $A'_0B''_0=AB$ and $A'_0C''_0=AC$. Because of scaling we get $B''_0C''_0=BC$. But $B''_0C''_0 > B''C''$ because of elementary properties of comparison geometry (as the hyperbolic plane has curvature $-1$). 
\end{proof}

In what follows we will need the following corollary. Let $d \in \mathfrak D_c(V, \mathcal T)$, $\sigma>0$ and $\mathfrak B(d, \sigma)$ be the open ball centered at $d$ of radius $\sigma$ in $(\R^{N}, l_{\infty})$. Define $\mathfrak B_c(d, \sigma) = \mathfrak D_c(V, \mathcal T) \cap \mathfrak B(d, \sigma)$ and $\mathfrak B_{sc}(d, \sigma)= \mathfrak D_{sc}(V, \mathcal T) \cap \mathfrak B(d, \sigma)$.

\begin{crl}
\label{connect}
For sufficiently small $\sigma$ the set $\mathfrak B_{sc}(d, \sigma)$ is connected.
\end{crl}

\begin{proof}
If $d \in \mathfrak D_{sc}(V, \mathcal T)$, then the statement is clear as $\mathfrak D_{sc}(V, \mathcal T)$ is an open set.

Assume that $d \in \mathfrak D_c(V, \mathcal T)$, but not in $\mathfrak D_{sc}(V, \mathcal T)$. As $\mathfrak D_c(V, \mathcal T)$ is semi-analytic, it is locally connected, i.e., for sufficiently small $\sigma$ the set $\mathfrak B_c(d, \sigma)$ is connected. Now consider two points in $\mathfrak B_{sc}(d, \sigma)$. They can be joined by a path $d_t$ in $\mathfrak B_c(d, \sigma)$. As $\mathfrak B(d, \sigma)$ is open, then for $t_0>1$  sufficiently close to 1, if we multiply all side lengths of $d_t$ by $t_0$, then the new path $d'_t$ still belongs to $\mathfrak B(d, \sigma)$. On the other hand, $d'_t$ belongs to $\mathfrak D_{sc}(V, \mathcal T)$  by Lemma~\ref{decreaseangles}. It remains to connect the endpoints of the old path with the endpoints of the new one.
\end{proof}

\subsection{Trapezoids and prisms}
\label{tp}

Trapezoids and prisms are the basic building blocks that we use to construct Fuchsian cone-manifolds in the next subsection.

\begin{dfn}
A \emph {trapezoid} is the convex hull of a segment $A_1A_2 \subset \H^2$ and its orthogonal projection to a line such that the segment $A_1A_2$ does not intersect this line. It is called \emph{ultraparallel} if the line $A_1A_2$ is ultraparallel to the second line. 
\end{dfn}

\begin{figure}
\begin{center}
\includegraphics[scale=0.7]{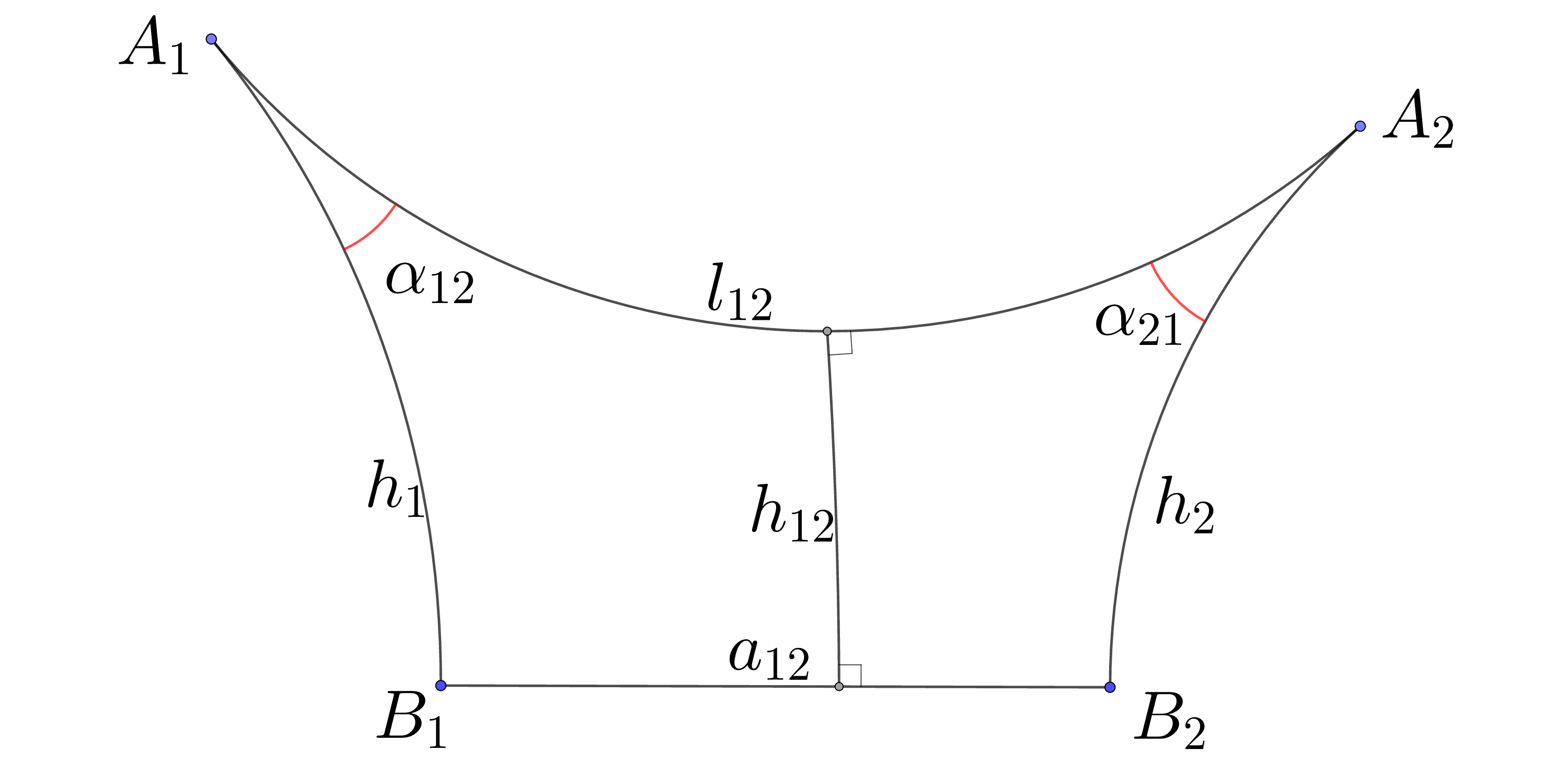}
\caption{A trapezoid.}
\label{Pic3}
\end{center}
\end{figure}

By $B_i$ denote the image of $A_i$ under the projection, $i=1,2$. We refer to $A_1A_2$ as to the \emph{upper edge}, to $B_1B_2$ as to the \emph{lower edge} and to $A_iB_i$ as to the \emph{lateral edges}. The vertices $A_i$ are also called \emph{upper} and $B_i$ are called \emph{lower}. We denote by $l_{12}$ the length of $A_1A_2$, by $a_{12}$ the length of $B_1B_2$, by $h_i$ the length of $A_iB_i$, by $\alpha_{12}$ and $\alpha_{21}$ the angles at vertices $A_1$ and $A_2$ and by $h_{12}$ the distance from the line $A_1A_2$ to the line $B_1B_2$ in the case of an ultraparallel trapezoid. See Figure~\ref{Pic3}.

\begin{dfn}
A \emph {prism} is the convex hull of a triangle $A_1A_2A_3\subset \H^3$ and its orthogonal projection to a plane such that the triangle $A_1A_2A_3$ does not intersect this plane. It is called \emph{ultraparallel} if the plane $A_1A_2A_3$ is ultraparallel to the second plane. 
\end{dfn}

\begin{figure}
\begin{center}
\includegraphics[scale=0.25]{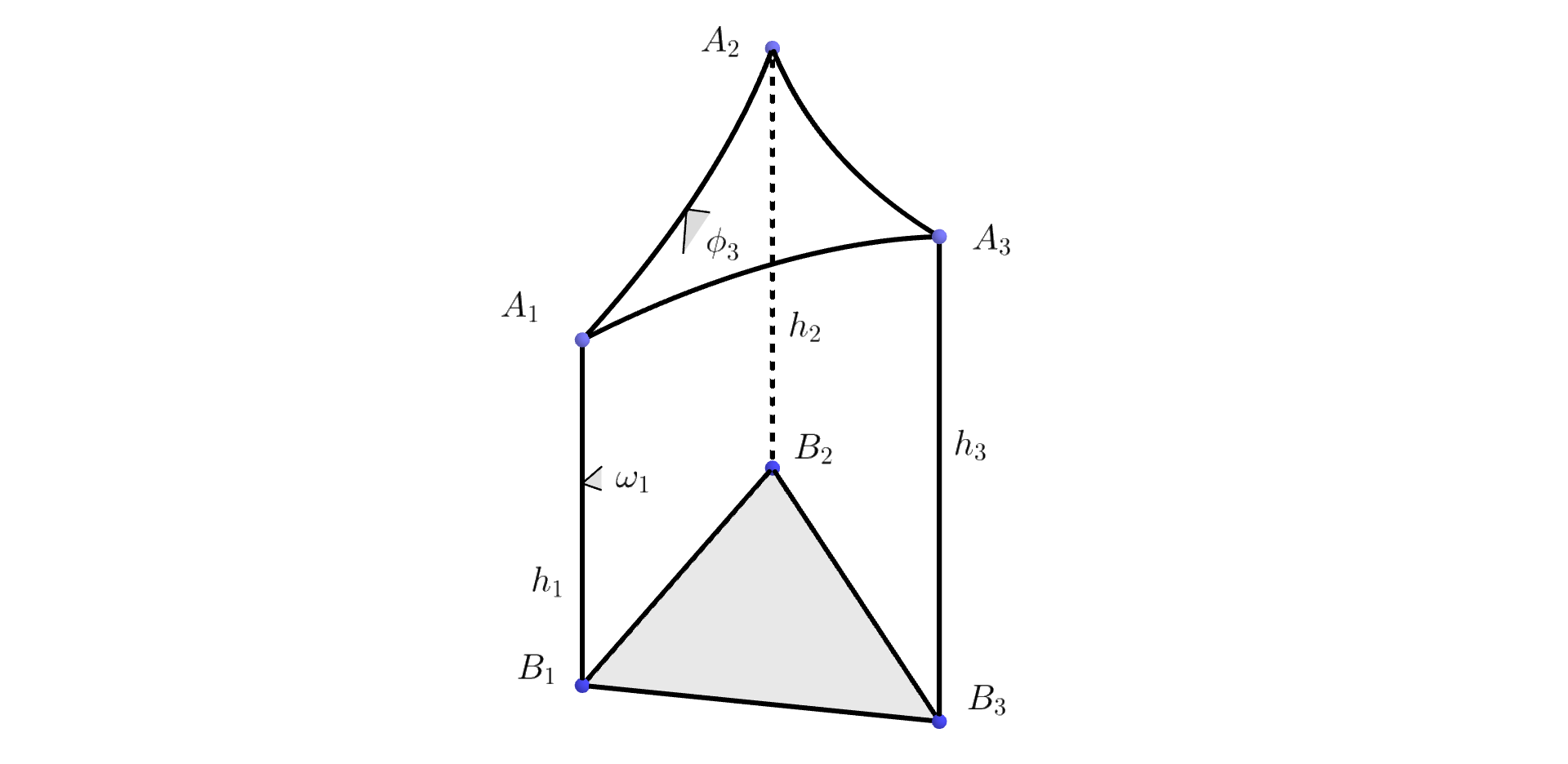}
\caption{A prism.}
\label{Pic2}
\end{center}
\end{figure}

Similarly to trapezoids, by $B_i$ we denote the image of $A_i$ under the projection, $i=1,2,3$, and we distinguish edges and faces of a prism into \emph{upper, lower} and \emph{lateral}. The lateral faces of a prism are trapezoids. The dihedral angles of edges $B_iB_j$ are equal $\pi/2$. The dihedral angles of edges $A_1A_2$, $A_2A_3$ and $A_3A_1$ are denoted by $\phi_3$, $\phi_1$ and $\phi_2$ respectively. The dihedral angle of an edge $A_iB_i$ is denoted by $\omega_i$. Other notation is inherited from trapezoids.  See Figure~\ref{Pic2}.

We do not allow degenerations of the upper triangle, but we consider degenerate prisms with collinear $B_1$, $B_2$ and $B_3$, when the upper plane is orthogonal to the lower one. However, soon we will restrict ourselves only to ultraparallel ones, which do not degenerate.

Once can see a trapezoid (or a prism) as a triangle (respectively, a tetrahedron) with one hyperideal vertex dual to the lower edge (lower face).

The proof of next two lemmas is straightforward: see Thurston's book~\cite[Chapter 2.6]{Thu1} for a general approach allowing to obtain these formulas.

\begin{lm}[The cosine laws]
\label{coslaw}
For a trapezoid we have
$$\cos \alpha_{12}=\frac{\cosh l_{12}\sinh h_1 - \sinh h_2}{\sinh l_{12} \cosh h_1},$$
$$\cosh a_{12} = \frac{\sinh h_1 \sinh h_2 +\cosh l_{12}}{\cosh h_1 \cosh h_2}.$$
\end{lm}



\begin{lm}[The sine law]
\label{sinlaw}
For a trapezoid we have
$$\frac{\sinh a_{12}}{\sinh l_{12}}=\frac{\sin \alpha_{12}}{\cosh h_2}=\frac{\sin \alpha_{21}}{\cosh h_1}.$$
\end{lm}


In particular, for right-angled trapezoids we will use the following identities:

\begin{crl}
\label{rightangl}
Assume that in a trapezoid $\alpha_{21}=\pi/2$. Then
$$\sinh h_1=\cosh l_{12} \sinh h_2,~~~\tanh h_1=\cosh a_{12} \tanh h_2,$$
$$\sin \alpha_{12}=\frac{\cosh h_2}{\cosh h_1},~~~\cot \alpha_{12} = \sinh l_{12} \tanh h_2.$$
\end{crl}

A trapezoid or a prism is clearly determined uniquely by the lengths of the lower and lateral edges. From the cosine laws we also obtain that

\begin{crl}
\label{prismuniq}
A trapezoid or a prism is determined up to isometry by the lengths of the upper and lateral edges. 
\end{crl}

We end the section with a purely computational corollary of cosine/sine laws that will be used further.

\begin{crl}
\label{cotan}
For a trapezoid we have
$$\frac{\cot \alpha_{12}}{\cosh h_1}=\frac{\cosh a_{12} \tanh h_1 - \tanh h_2}{\sinh a_{12}}.$$
\end{crl}

\subsection{Fuchsian cone-manifolds}
\label{definfcm}

Fuchsian cone-manifolds, which we introduce in this subsection, are the main tools of our proof and most of this paper is devoted to understanding their properties.


\begin{dfn}
A \emph{representable triple} is a triple $(d, \mathcal T, h)$, where $\mathcal T$ is a triangulation of $S_g$ with $V(\mathcal T) = V$, $d \in \mathfrak D_c(V, \mathcal T)$  and $h: V \rightarrow \R_{>0}$ is a function on $V$ such that for every triangle $T$ of $\mathcal T$ there exists a prism with the lengths of the upper edges defined by the side lengths of $T$ in $d$ and the lengths of the lateral edges determined by $h$. We write $h(v)$ as $h_v$.
\end{dfn}

\begin{dfn}
Let $(d, \mathcal T, h)$ be a representable triple. Take all prisms determined by $h$, $\mathcal T$ and $d$ and glue them isometrically according to $\mathcal T$. The resulting intrinsic metric space $P$ together with the canonical isometry from $(S_g, d)$ to the upper boundary $\partial^{\uparrow} P$ provided by our construction, is called a \emph{marked Fuchsian cone-manifold with polyhedral boundary}. In what follows we will mostly omit the words ``marked'' and ``with polyhedral boundary'' saying only a \emph{Fuchsian cone-manifold} for short.
\end{dfn}

We say that a representable triple $(d, \mathcal T, h)$ is a \emph{representation} of $P$ and write $P=P(d, \mathcal T, h)$. The upper and lower boundaries of $P$ are defined naturally. By definition, the upper boundary $\partial^{\uparrow} P$ is identified with $(S_g, d)$. We also say that the triangulation $\mathcal T$ is \emph{compatible with} $P$. The function $h$ is called the \emph{height function} of~$P$.

We say that an isometry $f: P^1 \rightarrow P^2$ between two Fuchsian cone-manifolds, where $P^1=P(d^1, \mathcal T^1, h^1)$ and $P^2=P(d^2, \mathcal T^2, h^2)$ with $V(\mathcal T^1)=V(\mathcal T^2)=V$, is a \emph{marked isometry} if its composition with the canonical isometries of $\partial^{\uparrow} P^1$, $\partial^{\uparrow} P^2$ is an isometry between $(S_g, d^1)$ and $(S_g, d^2)$ isotopic to identity with respect to $V$. We will consider Fuchsian cone-manifolds up to marked isometry.

For a Fuchsian cone-manifold $P=P(d, \mathcal T, h)$ and $v \in V(\mathcal T)$ we denote by $\omega_v(P)$ the sum of dihedral angles of the respective lateral edges in all prisms incident to $v$ in $P$ and define the \emph{particle curvature} of $v$ in $P$ as $\kappa_v(P):=2\pi-\omega_v(P)$. If all $\kappa_v(P)=0$, then $P$ is a Fuchsian manifold with polyhedral boundary. 
We will also call it a \emph{polyhedral Fuchsian manifold} for short. 
Let $\mathcal T$ be a triangulation compatible with $P$. For $e \in E(\mathcal T)$ denote by $l_e(d)$ its length in the metric $d$, by $\phi_e(P)$ its dihedral angle in $P$ and by $\theta_e(P):=\pi-\phi_e(P)$ the \emph{curvature} of $e$. A Fuchsian cone-manifold $P$ is called \emph{convex} if the dihedral angles $\phi_e(P)$ of all edges of $P$ are not greater than~$\pi$. Convex Fuchsian cone-manifolds are our main objects and we consider non-convex ones only sometimes in intermediate steps of proofs. 

The upper boundary $\partial^{\uparrow} P$ admits a canonical stratification $\partial^{\uparrow} P= X_0 \sqcup X_1 \sqcup X_2$ as follows. The spherical link of a point $p$ at the boundary of a prism $\Pi \subset \H^3$ is the portion of the unit sphere in $T_p\H^3$ corresponding to the directions to $\Pi$. Now for a point $p \in \partial^{\uparrow} P$ we define naturally its \emph{spherical link} as the gluing of the spherical links of $p$ in all prisms containing it. If this link is a hemisphere, then $p \in X_2$. If the link is a spherical lune, then $p \in X_1$. Otherwise $p \in X_0$. Note that if $p \in X_0$, then $p \in V$ for any representation $P=P(d, V, h)$, so we rather denote it by $v$ in this case. If $v\in X_0$ and $\kappa_v(P)\neq 0$, then the spherical link of $v$ is a spherical polygon with a conical singularity in the interior.

A connected component of $X_2$ is called a \emph{face} of $P$. It would be natural to say that the connected components of $X_1$ (note that they are geodesic segments in $P$) are edges of $P$ and that the points of $X_0$ are  vertices of $P$. However, we want to make another convention. Every compact cone-manifold $P$ under our consideration carries a marked point set $V \subset \partial^{\uparrow} P$ containing $X_0$. We will refer to points of $V$ as to \emph{vertices} of $P$. If we want to emphasize that $v \in X_0$, then we say that $v$ is a \emph{strict vertex} of $P$, and we say that $v$ is a \emph{flat vertex} if $v \in V \backslash X_0$. Next, we call \emph{edge} any geodesic segment in $\partial^{\uparrow} P$ between two points of $V$ (possibly coinciding) that is geodesic in $P$. If it is in $X_1$, then we say that it is a \emph{strict edge}. Otherwise, we call it a \emph{flat edge}.

We denote the set of the faces $P$ by $\mathcal R(P)$ and call it the \emph{face decomposition} of $P$. 
A triangulation $\mathcal T$ with a vertex set $V$ is compatible with $P$ if and only if all edges of $\mathcal T$ are edges of $P$. Note that this means that all strict edges of $P$ are edges of $\mathcal T$.

It might happen that an actual vertex of $P$ is isolated in the sense that there are no strict edges emanating from it. Then its spherical link is a hemisphere with a conical singularity in its center. In case of polyhedral cusps with particles an example is given in~\cite[p. 468]{FiIz1}. It is easy to adapt it to higher genus and we do not do it here. It might also happen that there are homotopically non-trivial faces.

We rephrase the following result of Fillastre~\cite{Fil1}:

\begin{thm}
\label{Fillastre}
For each convex cone-metric $d$ on $S_g$ there exists a unique up to isometry convex polyhedral Fuchsian manifold $P$ with $\partial^{\uparrow} P$ isometric to $(S_g, d)$.
\end{thm}

\begin{rmk}
\label{markeduniq}
We note that the proof of the uniqueness in Theorems~\ref{Fillastre} can be strengthened to the uniqueness up to marked isometry. This is equivalent to the following: if the upper boundary metric $d$ of a convex Fuchsian manifold $P$ admits a self-isometry, then it extends to a self-isometry of $P$. We briefly sketch the argument to convince the reader.

Let $\mathfrak P(V)$ be the set of marked convex polyhedral Fuchsian manifolds with set of vertices $V$ such that all vertices are strict (considered up to isometry isotopic to identity with respect to $V$). For $P \in \mathfrak P(V)$ the induced metric of $\partial^{\uparrow} P$ is in $\mathfrak D_{sc}(V)$. This defines a map $\mathfrak I: \mathfrak P(V) \rightarrow \mathfrak D_{sc}(V)$. After setting a natural topology on $\mathfrak P(V)$, it is proven in~\cite{Fil1} that $\mathfrak I$ is a homeomorphism. Noting that $\mathfrak I$ is equivariant with respect to the action of the mapping class group of $(S_g, V)$ we obtain Theorem~\ref{Fillastre}.
\end{rmk}

One of the most important tools in our study of cone-manifolds is the following:

\begin{dfn}
Let  $P(d, \mathcal T, h)$ be a Fuchsian cone-manifold. The function $\widetilde h: \partial^{\uparrow} P \rightarrow \R_{>0}$ assigning to a point $p \in \partial^{\uparrow} P$ its distance to $\partial_{\downarrow} P$ is called the \emph{extended height function} of $P$.
\end{dfn}

It coincides with $h$ at the  vertices of $\mathcal T$.

\begin{lm}
\label{htriangleinequality}
Take $p,q \in \partial^{\uparrow} P(d, \mathcal T, h)$. Then $\widetilde h(p) \leq \widetilde h(q)+d(p,q).$
\end{lm}

The proof is straightforward.

\subsection{Admissible heights and the discrete curvature}
\label{disccurv}

In this subsection we formulate without proofs a few key facts used in Section~\ref{generalrig} to prove Theorem~\ref{main}. For proofs we refer to Section~\ref{conechap}.

First of all, a convex Fuchsian cone-manifold is completely determined by its intrinsic upper boundary metric and its heights (i.e., we can restore the face decomposition from this data):

\begin{lm}
\label{heightdefine}
Let $\mathcal T^1$ and $\mathcal T^2$ be two triangulations with $$V(\mathcal T^1)=V(\mathcal T^2)=V,$$ $P^1=P(d, \mathcal T^1, h)$ and $P^2=P(d, \mathcal T^2, h)$ be two convex Fuchsian cone-manifolds. Then $P^1$ is marked isometric to $P^2$ with respect to $V$.
\end{lm}

A proof is postponed to Subsection~\ref{heightdefsec}.

This means that while we restrict ourselves to \emph{convex} Fuchsian cone-manifolds, we can drop the triangulation from their definition. However, we need the base set $V$, i.e., we represent a cone-manifold simply as $P=P(d, V, h)$, where $V(d) \subseteq V$.

\begin{dfn}
We call a function $h$ on $V$ to be \emph{admissible} for $(d, V)$ if there exists a convex Fuchsian cone-manifold $P(d, V, h)$. By $H(d,V) \subset \R^{n}$ we denote the set of all admissible $h$ for the pair $(d, V)$, where $n:=|V|$ and functions on $V$ are associated with points in $\R^n$.
\end{dfn}

Thus, $H(d, V)$ can be viewed as the space of convex Fuchsian cone-manifolds for a fixed upper boundary metric $d$ and a marked set $V$. While working with $H(d, V)$ we should take into account the following difficulties: it is non-compact and non-convex. Its closure in $\R^n$ can be obtained by adding the origin (however, it remains non-compact also after this). 

Let $V(d) \subseteq V \subset V'$, $P(d, V, h)$ be a convex Fuchsian cone-manifold and $\widetilde h: \partial^{\uparrow} P \rightarrow \R_{>0}$ be the extension of $h$. The extension gives us a height function $h' \in H(d, V')$. This defines the canonical embedding $H(d,V) \hookrightarrow H(d, V')$. In what follows, when we need such an extension, we will abuse the notation and write just $P=P(d, V', h)$ in the sense that we extend $h$ to $V'$ using $\widetilde h$.

In the opposite direction, assume that $P(d, V', h)$ is a convex Fuchsian cone-manifold, $V \subset V'$ and for any $v \in V' \backslash V$ we have $\nu_v(d)=0$ and $\kappa_v(d)=0$. Let $h|_V$ be the restriction to the set $V$ of $h$. Then $P$ can be represented as $P(d, V, h|_V)$. In such case we will also simply write $P(d, V, h)$ instead. 

Let $P=P(d, \mathcal T, h)$ be a Fuchsian cone-manifold. Define its \emph{discrete curvature} as
$$S(P)=S(d, \mathcal T, h) = -2\vol(P)+ \sum_{v \in V(\mathcal T)}\kappa_v(P) h_v +\sum_{e \in E(\mathcal T)} \theta_e(P) l_e(d).$$

Note that $S$ is independent from the choice of $\mathcal T$ compatible with $P$. Another name of $S$ is \emph{the discrete Hilbert--Einstein functional} as it is a discretization of the Hilbert--Einstein functional from smooth differential geometry: the integral of the scalar curvature over the interior of the manifold plus the integral of the mean curvature over the boundary. We refer to~\cite{Izm3} for a survey.

The discrete curvature $S$ can be considered as a continuous functional over $H(d, V)$. In Subsections~\ref{varapprsec}--\ref{varappr2sec} we will show that it is $C^2$. Also, in Subsection~\ref{varappr2sec} we will prove that $S$ is concave over $H(d, V)$. Moreover, in Subsection~\ref{varapprsec} we will also consider variations of $S$ with respect to the upper boundary.


A foundational result is the following maximization principle: a convex polyhedral Fuchsian manifold (so, without cone-singularities) maximizes the discrete curvature among all convex Fuchsian cone-manifolds with the same upper boundary metric. If we restrict ourselves to $V=V(d)$, then the inequality is strict.

\begin{lm}
\label{maxprinciple}
Let $d \in \mathfrak D_{sc}(V)$, $P^1=P(d, V, h^1)$ be the convex polyhedral Fuchsian manifold realizing $d$ and $P^2=P(d, V, h^2)$ be a convex Fuchsian cone-manifold distinct from $P^1$. Then $S(P^1) > S(P^2)$.
\end{lm}

A proof is given is Subsection~\ref{maxsec}.

\section{Proof of the main theorem}
\label{generalrig}

In this section we prove Theorem~\ref{main} modulo the main lemmas formulated in Subsection~\ref{statements}.

\subsection{Height functions}
\label{cbbheight}

Let $F$ be a Fuchsian manifold with convex boundary. We say that the function $h: \partial^{\uparrow} F \rightarrow \R_{>0}$ assigning to a point of the upper boundary its distance to the lower boundary is the \emph{height function} of $F$. We remark that if $F$ is a polyhedral Fuchsian manifold, then according to the notation from Subsection~\ref{definfcm} we should denote this function by $\widetilde h$ and by $h$ we denote only its restriction to the vertices. However, in this section we consider the height functions always defined over the whole upper boundary and denote them by $h$ slightly abusing the notation.

Consider the universal cover $G$ of $F$ developed to the Klein model of $\H^3$. We consider the Klein model endowed simultaneously with the hyperbolic metric and with the metric of the Euclidean unit ball. We note that convexity is the same in both metrics. Also, if a line is orthogonal to a plane passing through the origin, then the orthogonality also holds in both metrics simultaneously.

Let $Ox_1x_2x_3$ be the Euclidean coordinates. We assume that the lower boundary $\partial_{\downarrow} G$ of $G$, which is a geodesic plane, is developed to the (open) unit disk in the $Ox_1x_2$ plane. By $\rho: \partial^{\uparrow} G \rightarrow \partial_{\downarrow} G$ we denote the orthogonal projection map, which is a homeomorphism due to convexity, and define $$h_{\downarrow}:=h\circ \rho^{-1}: \partial_{\downarrow} G \rightarrow \R_{>0},$$ where $h$ is the height function extended to $\partial^{\uparrow} G$. By $h_{\E}: \partial_{\downarrow} G \rightarrow \R_{>0}$ we denote the Euclidean distance between $x \in \partial_{\downarrow} G$ and $\rho^{-1}(x)$. The convex surface $\partial^{\uparrow} G$ is the graph of $h_{\E}$. Thereby, $h_{\E}$ is a Lipschitz function over every compact subset of $\partial_{\downarrow} G$ and is differentiable almost everywhere. 

Comparing the Euclidean and hyperbolic metrics in the Klein model we compute

\begin{equation}
\label{heightEH}
h_{\E}(x)=h_{\E}(x_1,x_2)=\sqrt{1-x_1^2-x_2^2}\tanh h_{\downarrow}.
\end{equation}

This implies that $h_{\downarrow}$ is also differentiable almost everywhere and is Lipschitz over every compact subset of $\partial_{\downarrow} G$. In particular, if $\psi_{\downarrow}\subset \partial_{\downarrow} G$ is a rectifiable curve, then $\psi=\rho^{-1}(\psi_{\downarrow})$ is also rectifiable. Note that this is similar to the treatment of \emph{horoconvex functions} done in~\cite[Subsection 2.2]{FIV},~\cite{Lab}.




By Corollary~\ref{BusFel} if $\psi \subset \partial^{\uparrow} G$ is rectifiable, then its orthogonal projection to $\partial_{\downarrow} G$ is also rectifiable.

Denote $x_h = \arctanh(x_3/\sqrt{1-x_1^2-x_2^2})$. Let $g_{\downarrow}$ be the metric tensor of $\partial_{\downarrow} G$ in the coordinates $x_1, x_2$ and $g_h$ be the metric tensor of $\H^3$ in the coordinates $x_1, x_2, x_h$. Using the expression of the metric tensor in the Klein model we get

\begin{equation}
\label{heightmetric}
g_h=\cosh^2 x_h g_{\downarrow} + d x_h^2.
\end{equation}

Now we are ready to prove the main result of this subsection:

\begin{lm}
\label{heightrig}
Let $F^1$ and $F^2$ be two Fuchsian manifolds with convex boundaries, $f: \partial^{\uparrow} F^1 \rightarrow \partial^{\uparrow} F^2$ be an isometry and $h^1,  h^2$ be the height functions. Assume that $h^1= h^2 \circ f$. Then $f$ extends to an isometry of $ F^1$ and $ F^2$.
\end{lm}

\begin{proof}
We develop both universal covers $G^1$, $G^2$ to the Klein model such that both lower boundaries coincide with the plane $M=Ox_1x_2$. The isometry $f$ extends to an isometry between $\partial^{\uparrow}G^1$ and $\partial^{\uparrow} G^2$, which we continue to denote by $f$. Let $\rho^1$, $\rho^2$ be the projection maps from the upper boundaries $\partial^{\uparrow} G^1$, $\partial^{\uparrow} G^2$ to $M$ and $h^1_{\downarrow}=h^1\circ (\rho^1)^{-1}$, $h^2_{\downarrow}=h^2\circ(\rho^2)^{-1}$. Define $$f_{\downarrow}:= \rho^2\circ f \circ (\rho^1)^{-1}: M \rightarrow M.$$ It is enough to prove that $f_{\downarrow}$ is an isometry. Indeed, if $f_{\downarrow}$ is a an isometry and $h^1_{\downarrow}=h^2_{\downarrow}\circ f_{\downarrow}$, then the natural extension of $f_{\downarrow}$ to $G^1$ with the help of the orthogonal projection is an equivariant isometry of $G^1$ to $G^2$.

Let $\psi^1_{\downarrow}: [0; \tau] \rightarrow M$ be a rectifiable curve parametrized by length and $\psi^2_{\downarrow}:=f_{\downarrow} \circ \psi^1_{\downarrow}$. As the hyperbolic metric on $M$ is intrinsic, it suffices to prove that the length $l(\psi^1_{\downarrow})=l(\psi^2_{\downarrow})$. Suppose that this is not true. Then there exists a subset $I\subset [0; \tau]$ of strictly positive Lebesgue measure that either $g_{\downarrow}(\dot\psi^2_{\downarrow},\dot\psi^2_{\downarrow})>1$ or $g_{\downarrow}(\dot\psi^2_{\downarrow},\dot\psi^2_{\downarrow})<1$ almost everywhere on $I$. 

Assume that the first case hold, the second is done similarly. Let $\psi^1:=(\rho^1)^{-1}\circ \psi^1_{\downarrow}$ and $\psi^2:=(\rho^2)^{-1}\circ \psi^2_{\downarrow}$. As $f$ is an isometry between $\partial^{\uparrow} G^1$ and $\partial^{\uparrow} G^2$, the length measure $l(\psi^1|_I)$ of $\psi^1$ restricted to $I$ is equal to the length measure $l(\psi^2|_I)$ of $\psi^2$ restricted to $I$. Due to (\ref{heightmetric}) one can compute them as Lebesgue integrals

$$l(\psi^1|_I)=\int_I\bigg(\cosh^2 h_{\downarrow}^1(\psi_{\downarrow}^1) g_{\downarrow}(\dot\psi^1_{\downarrow}, \dot\psi^1_{\downarrow})+(\dot h_{\downarrow}^1(\psi_{\downarrow}^1))^2 \bigg)^{1/2} dt,$$
$$l(\psi^2|_I)=\int_I\bigg(\cosh^2 h_{\downarrow}^2(\psi_{\downarrow}^2) g_{\downarrow}(\dot\psi^2_{\downarrow}, \dot\psi^2_{\downarrow})+(\dot h_{\downarrow}^2(\psi_{\downarrow}^2))^2 \bigg)^{1/2} dt.$$

For all $t \in I$ (almost all $t$ for the third inequality) we have
$$h_{\downarrow}^1(\psi_{\downarrow}^1(t))=h_{\downarrow}^2(\psi_{\downarrow}^2(t)),$$ $$g_{\downarrow}(\dot\psi^1_{\downarrow}(t), \dot\psi^1_{\downarrow}(t))=1,$$ $$g_{\downarrow}(\dot\psi^2_{\downarrow}(t),\dot\psi^2_{\downarrow}(t))>1.$$ Thus, $l(\psi^1|_I)\neq l(\psi^2|_I)$ and we get a contradiction. 
\end{proof}

Choose an arbitrary homeomorphism $f^1: S_g \rightarrow \partial^{\uparrow} F^1$ and let $d$ be the pull-back by $f^1$ of the intrinsic metric of $\partial^{\uparrow} F^1$. Then $(S_g,  d)$ is a CBB($-1$) metric space and $f^1$ is an isometry. The pair $(F^1, f^1)$, where $f^1$ is defined up to isotopy, is called a \emph{marked Fuchsian manifold with boundary}. Define $f^2:=f\circ f^1$, which is also an isometry. The height functions $ h^1$ and $ h^2$ can be pulled back to the functions on $S_g$, which we continue to denote by $ h^1$ and $ h^2$. Lemma~\ref{heightrig} shows that if $ h^1$ and $ h^2$ coincide as functions on $S_g$, then Theorem~\ref{main} follows. We denote the marked manifolds $ (F^1, f^1)$ and $ (F^2, f^2)$ as $F(d,  h^1)$ and $F(d,  h^2)$ respectively.

\subsection{Flat points}
\label{flatpoints}

Recall that the curvature measure $\nu$ is non-negative for CBB($-1$) metrics. We call a point $p \in (S_g,  d)$ \emph{flat} if there exists an open neighborhood $U \ni p$ such that $\nu(U)=0$. One can show (e.g., from the arguments below) that in this case $U$ is isometric to an open subset of the hyperbolic plane. In the other case we say that $p$ is \emph{non-flat}. In particular, if $\nu_p( d) \neq 0$, then $p$ is non-flat. The main result of this subsection is

\begin{lm}
\label{flatlm}
Let $F(d, h^1)$ and $F(d, h^2)$ be two marked Fuchsian manifolds with convex boundaries. Assume that for each non-flat $p \in (S_g,  d)$ we have $h^1(p)= h^2(p)$. Then this is true for any $p$. 
\end{lm}

It is reminiscent of the fact from Euclidean convex geometry that the support of $m$-th curvature measure is the closure of the set of $m$-extreme points: see~\cite[Theorem 4.5.1]{Sch}. 

First we need to establish some auxiliary facts (note that they are not used outside of the proof of Lemma~\ref{flatlm}).

\begin{lm}
\label{nonflat}
Let $G \subset \H^3$ be a convex body and $M$ be a supporting plane such that $R=M\cap G$ has non-empty interior relative to $M$. Define $\overline R = \partial G \backslash R$. Then $\nu(\overline R)>0$.
\end{lm}

\begin{proof}
We prove it with the help of duality from Section~\ref{duality} and Lemma~\ref{duallemma}. In particular, we use the notation from there.

Define the equator $$(S^2)^*:=\{x \in \dS^3: x_0=0\}$$ and the projection map $\mathcal Q: \partial G^* \rightarrow (S^2)^*$ sending a point $p \in \partial G^*$ to the endpoint of a geodesic from $p$ orthogonal to $(S^2)^*$. It is easy to see that it is a homeomorphism: for the details see~\cite[Section 1]{BeCa}. Due to Lemma 2.1 from~\cite{BeCa} and the duality $\nu=\mu$ from Lemma~\ref{duallemma}, to prove $\nu(\overline R)>0$ it is enough to show that $\mathcal Q((\overline R)^*)$ contains an open set. 

Recall that the point $o$ with coordinates $(1,0,0,0)$ in $\R^{1,3}$ is an interior point of $G$. Let $C$ be the convex cone in $\H^3$ with the apex $o$ over $\partial R$. Consider a plane $K$ through $o$ that is tangent to $C$, but not to $\partial R$, oriented outwards. The cone $C$ is strictly convex, therefore the set of dual points $K^*$ for all such planes $K$ form an open subset of $(S^2)^*$. Consider a geodesic segment in $\dS^3$ connecting $K^*$ and $\mathcal Q^{-1}(K^*) \in \partial G^*$. It corresponds to a path of mutually ultraparallel planes all orthogonal to the same ray from $o$ directed towards $\overline R$. Thus, $K^* \in \mathcal Q((\overline R)^*)$. This finishes the proof.
\end{proof}

We need the last Lemma for a proof of the following fact:

\begin{lm}
\label{flatconv}
Let $ F=F(d,  h)$ be a Fuchsian manifold with convex boundary. Then all flat points of $(S_g, d)$ belong to the convex hull in $F$ of non-flat points.
\end{lm}

\begin{proof}
Consider the universal cover $G$ of $F$ in $\H^3$. We say that a segment in $\partial^{\uparrow} G$ is \emph{extrinsically geodesic} if it is a geodesic segment in $\H^3$. Note that $\partial^{\uparrow} G$ does not contain an extrinsically geodesic ray. Indeed, if such a ray intersects with the boundary at infinity of $\partial_{\downarrow} G$, then $h$ goes to zero along this ray. Otherwise, $h$ goes to infinity. Both conclusions contradict to the compactness of $F$.

Recall that a point $p \in \partial^{\uparrow} G$ is \emph{extreme} if it does not belong to the relative interior of an extrinsically geodesic segment in $\partial^{\uparrow} G$. A point $p$ is called \emph{exposed} if there exists a plane $M$ in $\H^3$ such that $M \cap G = \{p\}$. As $\partial^{\uparrow} G$ does not contain extrinsically geodesic rays, it is straightforward that $\partial^{\uparrow} G$ is contained in the convex hull of its extreme points. Clearly, being extreme or exposed depends only on the preimage of $p$ from $\partial^{\uparrow} F$. 

Consider $\H^3$ in the Klein model. Note that the notions of extreme and exposed points are purely affine, hence, if we consider $G$ as a Euclidean convex set, then extreme and exposed points  remain the same. For Euclidean closed convex sets the Straszewicz theorem~\cite{Str} says that the the extreme points belong to the closure of the set of exposed points. One can see, e.g., \cite[Theorem 1.4.7]{Sch} for a modern proof under additional assumption of compactness. To make $G$ closed we only need to add the boundary at infinity of $\partial_{\downarrow} G$, which does not change the result because $\partial^{\uparrow} G$ does not contain extrinsically geodesic rays.

Our plan is to prove that an exposed point is non-flat. The set of non-flat points is closed by definition, hence, extreme points are also non-flat and this finishes the proof due to the discussion above. We use some ideas from Olovyanishnikov~\cite{Olo}.



Suppose that $p \in \partial^{\uparrow} G$ is a flat exposed point. Let $U$ be a neighbourhood of $p$ such that $\nu(U)=0$ and $M$ be a supporting plane to $G$ such that $M\cap G=\{p\}$. We push $M$ slightly inside by a hyperbolic isometry with the axis passing through $p$ orthogonally to $M$ and obtain the plane $M'$ with $M' \cap G = \chi$ where $\chi$ is a closed curve bounding a compact set $U' \subset U$. Lemma~\ref{nonflat} implies that $\nu(U')>0$. This is a contradiction with $\nu(U)=0$.
\end{proof}

Now we need the notion of a $\mathcal F(-1)$-concave function due to Alexander and Bishop~\cite{AlBi}. 

\begin{figure}
\begin{center}
\includegraphics[scale=1]{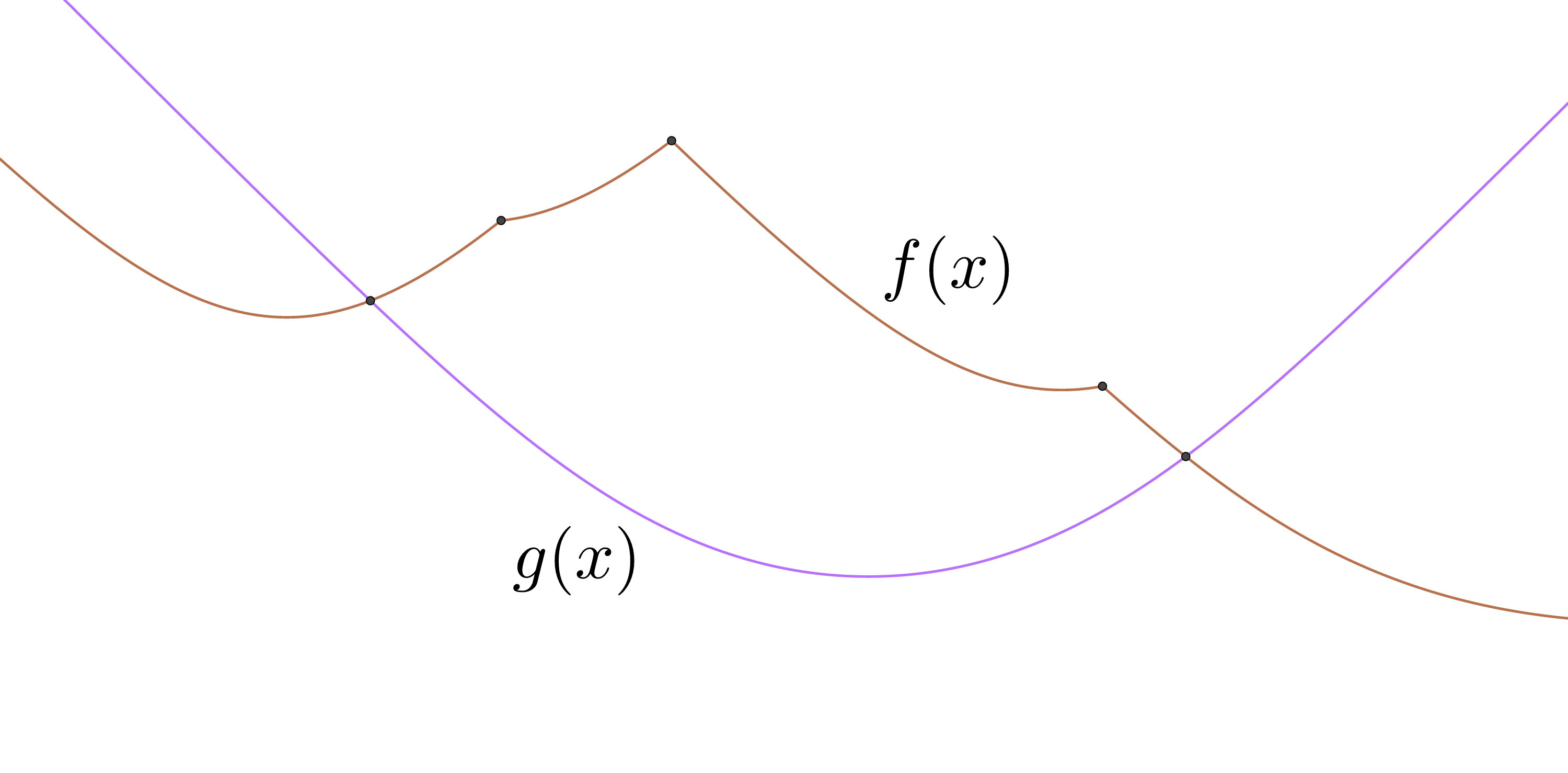}
\caption{An $\mathcal F(-1)$-concave function.}
\label{Pic5}
\end{center}
\end{figure}

\begin{dfn}
\label{f-1concave}
By $\mathcal F(-1)$ we denote the set of twice continuously differentiable functions $g: \R \rightarrow \R$ satisfying $g''=g$. A continuous function $f: I \rightarrow \R$ is \emph{$\mathcal F(-1)$-concave} if for each $[x_1; x_2] \subseteq I$ we have $f \geq g$, where $g \in \mathcal F(-1)$, $g(x_1)=f(x_1)$ and $g(x_2)=f(x_2)$ (see Figure~\ref{Pic5}). 

A continuous function $f: (S_g, d) \rightarrow \R$ is \emph{$\mathcal F(-1)$-concave} if it becomes $\mathcal F(-1)$-concave when restricted to every unit speed geodesic.
\end{dfn}

An element of $\mathcal F(-1)$ is a linear combination of $\sinh x$ and $\cosh x$ and is defined uniquely by values at two points. 

The last tool we need is:

\begin{lm}
\label{heightbend}
Let $ F=F(d,  h)$ be a Fuchsian manifold with convex boundary. Then the function $\sinh h$ is $\mathcal F(-1)$-concave.
\end{lm}

\begin{proof}
Let $f: I \rightarrow \R$ be a function. We skip proofs of two elementary facts below, which can be proven in the same way as in the theory of concave functions. 

\begin{prop}
\label{loctoglob}
Assume that for every $x \in I$ there exists a neighbourhood of $x$ in $I$ for which $f$ is $\mathcal F(-1)$-concave. Then $f$ is $\mathcal F(-1)$-concave over $I$.
\end{prop}

\begin{prop}
\label{tangent}
Assume that $f$ has the left and the right derivatives at every point. Then $f$ is $\mathcal F(-1)$-concave if and only if for every sufficiently close $t$ and $t'$ with $t<t'$ (resp. $t>t'$) we have $f(t') \leq g(t')$, where $g$ is a unique $\mathcal F(-1)$ function such that $f(t)=g(t)$ and $\dot g(t)$ is equal to the right (resp. the left) derivative of $f$ at $t$. 
\end{prop}

The proof of the next Claim is also straightforward.

\begin{prop}\cite[Lemma 2.1]{AlBi}
\label{f(-1)}
Let $\psi$ be a line in $\H^2$ and $g$ be the hyperbolic sine of the distance to $\psi$. Then the restriction of $g$ to every unit speed geodesic is in $\mathcal F(-1)$.
\end{prop}

We adapt Liberman's method~\cite{Lib},~\cite[Lemma 1 in Chapter IV.6]{Ale3}. Let $\tilde\psi: [0; \tau] \rightarrow (S_g, d)$ be a unit speed geodesic and $\tilde L$ be the union of all segments connecting points of the image of $\tilde\psi$ with their orthogonal projections to $\partial_{\downarrow} F$. By abuse of notation, we denote $h\circ\tilde \psi(t)$ by $h(t)$. The surface $\tilde L$ can be developed isometrically to a subset $L$ of $\H^2$ by straightening its lower boundary curve. We denote its boundary components by $\psi$ and $\psi_{\downarrow}$.

The curve $\psi_{\downarrow}$ is a geodesic segment. The set $L$ is convex. Indeed, otherwise there are two arbitrarily close points $p, q \in \psi$ such that the geodesic segment connecting them lies above $\psi$. Its length is smaller than the length of $\psi$ between $p$ and $q$. On the other hand, this segment corresponds to a curve that lies above $\partial^{\uparrow}F$ in the ambient Fuchsian manifold $\overline F$ and connects the preimages $\tilde p$ and $\tilde q$ of $p$ and $q$ under the developing map. Due to Corollary~\ref{BusFel}, the length of this curve is at least the length of its orthogonal projection to $\partial^{\uparrow}F$. Hence, it is at least $d(\tilde p, \tilde q)$. We can choose $p$ and $q$ so that the segment of $\tilde \psi$ between $\tilde p$ and $\tilde q$ is a shortest path. This is a contradiction.

The number $h(t)$ is equal to the distance from $\psi(t)$ to $\psi_{\downarrow}$. Let $p=\psi(0)$ and $p_{\downarrow}$ be the base of the perpendicular from $p$ to $\psi_{\downarrow}$. Due to convexity of $L$, there exists the right half-tangent $\chi$ to $\psi$ at $p$. This is well-known in the Euclidean case and extends to the hyperbolic case with the help of the Klein model. Parametrize $\chi$ with the unit speed and denote the distance from $\chi(t)$ to $\psi_{\downarrow}$ by $h^+(t)$. By $\alpha$ denote the angle between $\chi$ and the segment $pp_{\downarrow}$.

\begin{prop}
The right derivatives of $h$ and $h^+$ exist at zero and both are equal to $\cos\alpha$.
\end{prop}

Indeed, define $l(t):=d_{\H^2}(\psi(t), \psi(0))$. The convexity of $L$ implies 
$$\lim\limits_{t \rightarrow 0+} \frac{l(t)}{t}=1.$$
In the Euclidean case this is Lemma 2 from~\cite[Chapter IV.6]{Ale3}. If we consider the Klein model and put $\psi(0)$ in the origin, then we see that this is true in the hyperbolic case also. By $\alpha(t)$ we denote the angle between segments $p\psi(t)$ and $pp_{\downarrow}$. With the help of Lemma~\ref{coslaw} we get
$$\lim_{t \rightarrow 0+} \frac{h(t)-h(0)}{t}=\lim_{t \rightarrow 0+} \frac{\sinh h(t)-\sinh h(0)}{t \cosh h(0)} =$$
$$=\lim_{t \rightarrow 0+} \frac{\sinh h(t)-\sinh h(0)}{ l(t) \cosh h(0)}=
\lim_{t \rightarrow 0+} \frac{\cosh l(t) \sinh h(t)-\sinh h(0)}{ \sinh l(t) \cosh h(0)}=$$
$$= \lim_{t \rightarrow 0+} \cos \alpha(t) = \cos\alpha.$$
For $h^+(t)$ the same computations hold.
\vskip+0.2cm

Thereby, the functions $h$ and $h^+$ have the same right derivatives at 0, $h(0)=h^+(0)$ and $\sinh h^+\in \mathcal F(-1)$ due to Claim~\ref{f(-1)}. Let us now prove that for sufficiently small $t$ we have $h(t) \leq h^+(t)$. One can prove similar statement for any sufficiently close points $p$ and $q$ on $\psi$. Then it follows from Claim~\ref{loctoglob} and Claim~\ref{tangent} that $\sinh h$ is $\mathcal F(-1)$-concave.

Let $q = \psi(t)$ be a point such that $\psi$ is a shortest path between $p$ and $q$. By $q_{\downarrow}$ denote the base of perpendicular from $q$. Suppose that $\alpha \geq \pi/2$. Define $t'=d_{\H^2}(  p,   q)\leq t$. We have $h^+(t') \geq h(t)=d_{\H^2}(  q,   q_{\downarrow})$. Indeed, to get $\chi(t')$ we need to rotate the segment $pq$ away from the line $\psi_{\downarrow}$. As $t'\leq t$ and $\alpha \geq \pi/2$, we get $h^+(t)\geq h^+(t') \geq h(t)$.

Now we suppose that $\alpha < \pi/2$. Extend the segment $  q   q_{\downarrow}$ to the intersection with $  \chi$ and denote the intersection point by $  q'=  \chi(t')$. We have $h^+(t') \geq h(t)$. Define $t''=t'+d_{\H^2}(q', q)$. Due to convexity, $t'' \geq t$. Consider the point $  q''=  \chi(t'')$. As the triangle $  q''  q'   q$ is isosceles, we get $\angle   q_{\downarrow}  q   q'' > \pi/2$. Therefore, $h^+(t'') \geq h(t).$

We need to get $h^+(t) \geq h^+(t'')$. As $\alpha < \pi/2$, if $\chi$ is not ultraparallel to $  \psi_{\downarrow}$, then this is true for any $t'' \geq t$. Otherwise, we assume that $t$ is sufficiently small, so $  q''$ lies on the same side with $  q'$ from the closest point of $  \chi$ to $  \psi_{\downarrow}$. Thus, $h^+(t) \geq h^+(t'') \geq h(t)$ and the proof is finished.
\end{proof}

Now we have all the ingredients that we need.

\begin{proof}[Proof of Lemma~\ref{flatlm}.]
Let $p \in (S_g, d)$ be a flat point. Due to Lemma~\ref{flatconv}, $p$ belongs to the convex hull of non-flat points. The Caratheodory theorem (normally it is stated for Euclidean convex bodies, but this does not matter because of the Klein model) implies that it belongs to the convex hull of at most four non-flat points. If this number can not be reduced to at most three, then $p$ must be an interior point of $F$. Thereby, $p$ belongs either to the relative interior of a segment $\psi$ with non-flat endpoints $q_1$ and $q_2$ or to a triangle $T$ with non-flat vertices $q_1$, $q_2$ and $q_3$  that is extrinsically geodesic in $F^1$. Consider the first case. Parametrize $\psi$ by length. Note that $\sinh h^1$ restricted to an extrinsically geodesic segment in $F^1$ is $\mathcal F(-1)$: this is Claim~\ref{f(-1)} from the previous proof. As $h^1$ coincides with $h^2$ at the endpoints of $\psi$, by Lemma~\ref{heightbend} we have $h^2|_{\psi} \geq h^1|_{\psi}$. In particular $h^2(p)\geq h^1(p)$. In the second case, let $\psi$ be the edge of $T$ opposite to $q_1$. Similarly, we have $h^2|_{\psi} \geq h^1|_{\psi}$. Connect $q_1$ with a point $q\in \psi$ by a geodesic $\chi$ passing through $p$. This is possible as $T$ is extrinsically geodesic in $F^1$, so it is isometric to a hyperbolic triangle. Applying the same arguments to $\chi$ we still get that $h^2(p)\geq h^1(p)$. In the same way one can show that $h^1(p) \geq h^2(p)$. Thus, $h^1(p)=h^2(p)$.
\end{proof}

We showed that to prove Theorem~\ref{main} it is enough to show that $h^1$ coincides with $h^2$ at all non-flat points.

\subsection{Polyhedral approximation}
\label{triang}

Let $F=F(d, h)$ be a marked Fuchsian manifold with convex boundary, so $(S_g, d)$ is a CBB($-1$) metric space. In this subsection we will understand how it can be approximated by convex polyhedral Fuchsian manifolds.

We will consider various convex surfaces in $F$. We mark all of them with the help of the vertical projection map (along the perpendiculars to $\partial_{\downarrow} F$). Their intrinsic metrics and height functions are pulled back to $S_g$ with the help of marking maps. In particular, this induces a marking on $\partial_{\downarrow} F$.



The embedding of the universal cover of $F$ as $G \subset \H^3$ determines a representation of $\pi_1(S_g)$ as a Fuchsian subgroup $\Gamma$ of isometries of $\H^3$: it leaves invariant the plane $\partial_{\downarrow} G$. The quotient by this action of the closed invariant half-space of $\H^3$ containing $G$ is called \emph{the extended Fuchsian manifold} and is denoted by $\overline F$. It is homeomorphic to $S_g \times [0; \infty)$.

We need also to define several classes of triangulations:

\begin{dfn}
\label{geodtriang}
A geodesic triangulation $\mathcal T$ of $(S_g, d)$ is called \emph{short} if all triangles are convex, all edges are shortest paths and all angles are strictly smaller than $\pi$. It is called \emph{strictly short} if additionally all edges are unique shortest paths between their endpoints.
\end{dfn}

\begin{dfn}
\label{fine}
A geodesic triangulation $\mathcal T$ of $(S_g, d)$ is called \emph{$\delta$-fine} if for every triangle $T$ we have $\diam(T, d)< \delta$ and $\nu(T, d) < \delta$.
\end{dfn}

\begin{lm}
\label{goodtriang}
For every $\delta >0$ there exists a strictly short $\delta$-fine triangulation $\mathcal T$ of $(S_g, d)$. Moreover, each triangulation $\mathcal T'$ can be refined to a strictly short $\delta$-fine triangulation (by adding new vertices).
\end{lm}

\begin{proof}
The proof is the same as in CBB(0) case~\cite[Lemma in Chapter III.1, p. 134]{Pog2}. However, Pogorelov uses the intrinsic curvature $\nu_I$ instead of $\nu$. This can be easily fixed as the area of sufficiently small triangles uniformly diminishes. 
\end{proof}




Consider a sequence of positive numbers $\mu_m \rightarrow 0$ and a sequence of finite $\mu_m$-dense sets $V_m\subset (S_g, d)$. Take the convex hull of $V_m$ in $F$ and obtain a convex Fuchsian manifold $F_m$ with polyhedral boundary. By $d_m$ we denote its boundary metric. Here the convex hull is considered in the totally convex sense: the inclusion-minimal closed totally convex set containing $V_m$.

\begin{lm}
\label{convuni}
The metrics $d_m$ converge uniformly to $d$.
\end{lm}

\begin{proof}
We follow the ideas from~\cite[Lemma 3.7]{FIV} and~\cite[Lemma 10.2.7]{BBI}. Unfortunately, some details are less straightforward in our setting.


Let $h$, $h_m$ be the height functions of $\partial^{\uparrow} F$, $\partial^{\uparrow} F_m$. The sequence $\partial^{\uparrow} F_m$ converges to $\partial^{\uparrow} F$ in the Hausdorff sense. This implies that $h_m$ converges to $h$ uniformly, i.e., $$\zeta_m :=\sup_{p \in S_g} (h(p)-h_m(p))\rightarrow 0$$ (note that $h_m \leq h$).

\begin{prop}
\label{majorate}
Let $\partial^{\uparrow} F^1$, $\partial^{\uparrow} F^2$ be two convex surfaces in an extended Fuchsian manifold $\overline F$, $h^1$ and $h^2$ be their height functions and $d^1$, $d^2$ be their intrinsic metrics. Assume that $h^1 \leq h^2$ and $$\zeta:=\sup_{p \in S_g} (h^2(p)-h^1(p)).$$ Then $d^1\leq d^2 + 2\zeta$.
\end{prop}

The proof easily follows from Corollary~\ref{BusFel} and is the same as in~\cite[Lemma 2.12]{FIV}, we do not provide it here. Claim~\ref{majorate} and $h_m \leq h$ imply $d_m \leq d + 2\zeta_m$ with $\zeta_m \rightarrow 0$. 

Next we claim that for each $t$ the function $h'_t$ defined by the equation $$\tanh h'_t = e^{-t} \tanh h$$ is the height function of another convex surface in $F$. Indeed, as $\partial^{\uparrow} F$ is convex, its Euclidean height function $h_{\E}$ is concave. Due to (\ref{heightEH}), the Euclidean height function of the surface defined by $h'_t$ is different from $h_{\E}$ just by multiplication by $e^{-t}$. Thus, it is also concave and its graph is a convex surface. 

The transformation $h \rightarrow h'_t$ is an analogue of homothety in our setting.

We can choose a sequence $t_m \rightarrow 0$ such that $h'_{t_m} \leq h_m$ and $h'_{t_m}$ converges uniformly to $h$. Denote $h'_{t_m}$ by $h'_m$, the respective surface by $\partial^{\uparrow} F'_m$ and let $d'_m$ be the intrinsic metric of $\partial^{\uparrow}  F'_m$. From Claim~\ref{majorate} we obtain $d'_{m} \leq d_m + 2 \zeta'_m$ with $ \zeta'_m \rightarrow 0$. 

It remains to relate $d$ with $d'_m$. First, we check how the distances in $F$ change. Fix $m$, let $p, q \in \partial^{\uparrow} F$ and $p', q' \in \partial^{\uparrow} F'_m$ be the corresponding points. Define $l:=d_F(p,q)$, $l':=d_F(p', q')$. We want to show that there exist numbers $\xi_m \rightarrow 0$ independent of $p$, $q$ such that $l\leq l'(1+ \xi_m).$ This follows from

\begin{prop}
\label{rat}
Let $A_1A_2B_1B_2$ and $A'_1 A'_2 B_1 B_2$ be two ultraparallel trapezoids with $$\tanh h'_1 = e^{-t}\tanh h_1,$$ $$\tanh h'_2 = e^{-t}\tanh h_2$$ (here we use the notation from Subsection~\ref{tp} and mark the parameters of the trapezoid $A'_1 A'_2 B_1 B_2$ with prime symbol). By $h_0$ we denote the distance between lines $A_1A_2$ and $B_1B_2$. For $0<m<M$ by ${\rm Trap}(t, m, M)$ we denote the space of such pairs of trapezoids with $h_0 \geq m$ and $h_1, h_2 \leq M$. Define the function ${\rm rat}: {\rm Trap}(t, m, M) \rightarrow \R_{>0}$ as $l_{12}/l'_{12}.$ Fix $m, M$ and let $t \rightarrow 0$. Then $$\sup_{{\rm Trap}(t, m, M)} {\rm rat}\rightarrow 1.$$
\end{prop}

First we finish the proof of the Lemma and then give a proof of Claim~\ref{rat}. As $F$ is compact, the function $h$ is bounded from below and from above by positive constants. Let $\tilde p,\tilde q \in \partial^{\uparrow} G \subset \H^3$ be the lifts of $p,q$ to the universal cover such that $d_{\H^3}(\tilde p,\tilde q)=d_F(p,q)$. The line passing through $\tilde p$, $\tilde q$ is ultraparallel to the lower plane $\partial_{\downarrow} G$. Let $r$ be the closest point from this line to $\partial_{\downarrow} G$. If $r$ lies outside the segment $\tilde p\tilde q$, then it lies above $\partial^{\uparrow} G$ and $d_{\H^3}(r, \partial_{\downarrow} G)$ is at least the infimum of $h$. If $r$ lies between $\tilde p$ and $\tilde q$, then due to Corollary~\ref{rightangl} we have $$\sinh h(\tilde p) = \sinh d_{\H^3}(r, \partial_{\downarrow} G) \cosh d_{\H^3} (\tilde p, r).$$ As $d_{\H^3} (\tilde p, r)\leq d_F(p,q)\leq \diam(F)$, we get that there exists $m=m(F)>0$ such that $d_{\H^3}(r, \partial_{\downarrow} G) \geq m$. Thus, we can apply Claim~\ref{rat} and get $l \leq l'(1+ \xi_m)$ with $\xi_m \rightarrow 0$.

Now connect $p'$ and $q'$ by a shortest path $\psi' \subset \partial^{\uparrow} F'_m$. Let $\psi \subset \partial^{\uparrow} F$ be its image under the vertical projection. It is rectifiable. Consider a polygonal approximation of $\psi'$ with sufficiently small segments. It corresponds to a polygonal approximation of $\psi$ of total length multiplied by at most $1+\xi_m$. As the lengths of $\psi$ and $\psi'$ are the suprema of the lengths of their polygonal approximations, we get $d \leq d'_m(1+\xi_m).$

In total, we obtain $|d -d_m|\rightarrow 0$ uniformly.

\end{proof}

\begin{proof}[Proof of Claim~\ref{rat}.]
The space ${\rm Trap}$ of all ultraparallel trapezoids up to isometry can be parametrized by $h_0$, $h_1$ and $a_{12}$ with $h_0 >0$, $h_1 \geq h_0$ and $a_{12} \in \R\backslash \{0\}$ as follows. First we choose two ultraparallel lines at distance $h_0$. Let $A$ be the closest point on the upper line to the lower line. We choose $A_1$ to the right from $A$ such that its distance to the lower line is $h_1$. Then we choose $A_2$ to the right from $A_1$ if $a_{12}$ is positive and to the left if $a_{12}$ is negative. One can see that two trapezoids are isometric if and only if the parameters are the same.

The trapezoid $A'_1 A'_2 B_1 B_2$ is determined by $A_1A_2B_1B_2$ and $t$, hence for each $t$ we can consider ${\rm Trap}(t, m, M)$ as a subset of ${\rm Trap}$. The closure $\overline {\rm Trap}(t, m, M)$ is compact (note that $m \leq h_0 \leq h_1 \leq M$) and is obtained by adding the degenerated trapezoids with $a_{12}=0$. We want to show that the function ${\rm rat}$ extends there continuously. Indeed, by Lemma~\ref{sinlaw} we get
$$\lim_{a_{12}\rightarrow 0} {\rm rat} (h_0, h_1, a_{12}) =\lim_{a_{12}\rightarrow 0} \frac{l_{12}}{l'_{12}}=\lim_{a_{12}\rightarrow 0} \frac{\sinh l_{12}}{\sinh l'_{12}}=\lim_{a_{12}\rightarrow 0} \frac{\cosh h_1 \sin \alpha'_{12}}{\cosh h'_1 \sin  \alpha_{12}}.$$ 

Applying Lemma~\ref{sinlaw} the second time we eliminate the angles:
$$\lim_{a_{12}\rightarrow 0} {\rm rat} (h_0, h_1, a_{12}) = \frac{\cosh^2 h_1 \cosh h'_0}{\cosh^2 h'_1 \cosh h_0}>0.$$

Now for arbitrary $\tau>0$ we consider the set 
$${\rm TRAP}(\tau, m, M)=\{ \cup {\rm Trap}(t, m, M): 0\leq t \leq \tau\}$$
and its closure $\overline{\rm TRAP}(\tau, m, M)$. The function ${\rm rat}$ is a continuous function over this compact and is equal to 1 as $t=0$. The claim follows.
\end{proof}

\begin{crl}
\label{diamareaglob}
$\diam(S_g, d_m) \rightarrow \diam(S_g, d),~~~\area(S_g, d_m) \rightarrow \area(S_g, d).$
\end{crl}

\begin{proof}
The first statement is trivial. For the second we use Theorem 9 of~\cite[Chapter 8]{AlZa}.
\end{proof}

Let $T$ be a geodesic triangle in $(S_g,  d)$. By $T( d) = (l_1, l_2, l_3, \lambda_1, \lambda_2, \lambda_3)( d)$ denote the 6-tuple of its side lengths and angles in the metric $ d$. We will always consider these 6-tuples as points of $\R^6$ endowed with $l_{\infty}$ metric.

\begin{lm}
\label{refin}
Let $\mathcal T$ be a strictly short triangulation of $(S_g, d)$. Consider a sequence of finite $\mu_m$-dense sets $V_m\subset (S_g, d)$ with $\mu_m \rightarrow 0$ such that all of them contain $V(\mathcal T)$ and none of them contains a point in the interior of an edge of $\mathcal T$. By $d_m$ denote the pull back of the upper boundary metric of the convex hull of $V_m$ in $F$ defined as above.

Then $\mathcal T$ is realized by infinitely many $d_m$. Moreover, the realizations can be chosen so that they are short and for every triangle $T$ we have (after taking a subsequence)

(1) $T(d_m) \rightarrow T(d)$;

(2) $\diam(T, d_m) \rightarrow \diam(T, d)$;

(3) $\area(T, d_m) \rightarrow \area(T, d).$
\end{lm}

\begin{proof}
Consider $u, v \in V(\mathcal T)$ connected by an edge $e\in E(\mathcal T)$. Connect $u$ and $v$ by a shortest path $\psi_m(e)$ in $\partial^{\uparrow} F_m$. By $\psi(e)$ we denote the realization of $e$ in $\partial^{\uparrow} F$. By~\cite[Chapter II.1, Theorems 4 and 5]{Ale3} applied to the intrinsic metric space $F$, for every $e\in E(\mathcal T)$ the sequence $\{\psi_m(e)\}$ converges uniformly (up to taking a subsequence) in $F$ as parametrised curves to a rectifiable path $\psi'(e) \subset \partial^{\uparrow} F$ of length at most $l_e(d)$. As $\psi(e)$ is the unique shortest path in $\partial^{\uparrow} F$ between its endpoints, we have $\psi(e)=\psi'(e)$.

Choose $\xi_1>0$ such that $\xi_1$-neighborhoods of vertices of $\mathcal T$ in $F$ do not intersect. Then we choose $\xi_2>0$ such that for every pair of edges $e'$ and $e''$, $\xi_2$-neighborhoods of $\psi(e'), \psi(e'')$ intersect only if $e'$ and $e''$ share an endpoint $v$ and if they do, then the intersection lies in the $\xi_1$-neighborhood of $v$. For sufficiently large $m$, $\psi_{m}(e)$ belongs to $\xi_2$-neighborhood of $\psi(e)$. This means that $\psi_{m}(e')$ and $\psi_{m}(e'')$ can intersect only if $e'$ and $e''$ share an endpoint. But in this case they can not intersect except at this endpoint by the non-overlapping property: if two shortest paths have two points in common, then either these are their endpoints or they have a segment in common, see~\cite[Chapter II.3, Theorem 1]{Ale3}.

It is clear now that the union of $\psi_m(e)$ gives a realization of $\mathcal T$ in $d_m$. We will continue to work with this realization. The convergence of side lengths of each triangle is already shown. Now we proceed to diameters and areas.

Let $\rho: F \rightarrow \partial^{\uparrow} F$ be the vertical projection map. We claim that $\rho(\psi_m(e))$ converge uniformly to $\psi(e)$ in the metric $d$. To see this, first project everything to $\partial_{\downarrow} F$. Denote the images by $\psi_{m\downarrow}(e)$ and $\psi_{\downarrow}(e)$. The curves $\psi_m(e)$ converge uniformly to $\psi(e)$ in $F$ and the projection to $\partial_{\downarrow} F$ contract $F$-distances. Therefore, $\psi_{m\downarrow}(e)$ converge uniformly to $\psi_{\downarrow}(e)$ in $\partial_{\downarrow} F$. It remains to recall that $\partial^{\uparrow} F$ is the graph of a Lipschitz function $h$, hence $d \leq Cd_{\downarrow}$ for some constant $C$, where  $d_{\downarrow}$ is the intrinsic metric of $\partial_{\downarrow}F$.

The convergence of the diameters follows easily from this and the uniform convergence of the metrics. For the convergence of areas we apply Theorem 8 of~\cite[Chapter~8]{AlZa}. 


Let us prove now the convergence of angles. Consider $v \in V(\mathcal T)$ and let $U_1, \ldots, U_k$ be the sectors between the edges of $\mathcal T$ emanating from $v$ in the cyclic order. We claim that a proof of the bound $$\ang(U_i, d) \leq \liminf\limits_{m \rightarrow \infty} \ang(U_i, d_m)$$ follows line by line the proof in the Euclidean case, see~\cite[Chapter IV.4, Theorem 2]{Ale3}.

Assume that for $\xi>0$ and some $i$ we have $$\ang(U_i, d) < \liminf\limits_{m \rightarrow \infty}  \ang(U_i, d_m)-\xi.$$ Then for some $m$ we have $\lambda_v(d) < \lambda_v(d_m)-\xi$. This contradicts to 

\begin{prop}
\label{coneins}
For each $v \in V(\mathcal T)$ and each $m$ we have $\lambda_v(d_m) \leq \lambda_v(d).$
\end{prop}

It remains to prove Claim~\ref{coneins} and the proof of Lemma~\ref{refin} is over.

\end{proof}

\begin{proof}[Proof of Claim~\ref{coneins}.]
Consider the universal covers $\partial^{\uparrow} G$, $\partial^{\uparrow} G_m \subset \H^3$ of $\partial^{\uparrow} F$ and $\partial^{\uparrow} F_m$. We use the Klein model of $\H^3$ embedded into $\E^3$ as the interior of the unit Euclidean ball $B(1)$. By $d_{\E}$ and $d_{\H}$ denote the Euclidean and the hyperbolic metrics on $B(1)$ respectively. We assume that $v$ is developed to the origin. Both $\partial^{\uparrow} G$ and $\partial^{\uparrow} G_m$ are also convex surfaces in the Euclidean metric. By $d'$ we denote the intrinsic metric on $\partial^{\uparrow} G$ induced by $d_{\E}$. We want to show that $\lambda_v(d)=\lambda_v(d').$

Let $B(r)$ be the Euclidean ball with the radius $r<1$ and center $v$. Comparing the Euclidean and the hyperbolic metric one can see that the identity map between $(B(r), d_{\E})$ and $(B(r), d_{\H})$ is bi-Lipschitz with constant $C(r)$ such that $C(r) \rightarrow 1$ as $r \rightarrow 0$, where $d_{\E}$ and $d_{\H}$ are the Euclidean and the hyperbolic metrics respectively on $B(r)$. 

Let $\psi$ and $\chi$ be two geodesics on $\partial^{\uparrow} G$ in the metric $d$ emanating from $v$. Due to the bi-Lipschitz equivalence of small neighborhoods of $v$ it is easy to see that the angle $\ang(\psi, \chi, d')$ is also defined and is equal to $\ang(\psi, \chi, d)$. Let $p \in \psi$, $q \in \chi$ be two points, connect them by shortest paths $\psi'$, $\chi'$ with $v$ in the metric $d$. By definition, $\ang(\psi', \chi', d') \rightarrow \ang(\psi, \chi, d')$. Using the definition of the total angle as the supremum of sums of angles between consecutive geodesics, we get $\lambda_v(d) \leq \lambda_v(d')$. The converse inequality could be obtained in the same way.

Consider the tangent cone to $\partial^{\uparrow} G$ at $v$ in the Euclidean metric (see~\cite[Chapter IV.5]{Ale3}). Its total angle is equal to $\lambda_v(d')=\lambda_v(d)$~\cite[Chapter IV.6, Theorem 3]{Ale3}. Similarly there exists the tangent cone to $\partial^{\uparrow} G_m$ at $v$ with the total angle $\lambda_v(d_m).$ The latter tangent cone is inscribed in the former. We note that for the Euclidean convex cones $2\pi$ minus the total angle is equal to the area of the Gaussian image. This shows the desired inequality.
\end{proof}

\subsection{Statements of the stability lemmas}
\label{statements}

Our proof of Theorem~\ref{main} is based on several lemmas, but even their statements are slightly cumbersome. Hence, we are going to discuss each one before we formulate it.\vskip+0.2cm

From Lemma~\ref{maxprinciple} we know that if $P^1=P(d, V, h^1)$ is a convex polyhedral Fuchsian manifold (so it does not have cone-singularities) and $P^2=P(d, V, h^2)$ is another convex (polyhedral) Fuchsian cone-manifold with isometric upper boundary, then $S(P^1)>S(P^2)$. We need to find a better quantitative lower bound on $S(P^1)-S(P^2)$. Naturally, such bound should depend on $\max_{v \in V} |h^2_v-h^1_v|$ and on global geometry of $P^1$. To our purposes it is enough to treat only the case when there exists $v \in V$ such that $h^2_v > h^1_v$. 

A difficulty arises on our way: the bound we obtain depends on the curvature $\nu_v$ of a point $v$ with $h^2_v \neq h^1_v$. This implies additional complications in our proof of Theorem~\ref{main}. Imagine that for two Fuchsian manifolds with isometric convex boundaries we have $h^2_v \neq h^1_v$ in a point with zero curvature. Then when we approximate these manifolds by polyhedral ones, the curvature of $v$ in the approximation goes to zero and we are unable to use our quantitative result. Luckily, from results of Subsection~\ref{flatpoints} we can choose $v$ to be non-flat. However, this does not mean that the curvature of $v$ stays bounded from zero in the approximations: imagine the case of smooth boundaries. But it guarantees that the curvature of a small neighbourhood of this point actually stays bounded from zero. This means that we should be able to handle the case when we have a lower bound on the curvature of a set of points where $h^1$ is distinct from $h^2$, but not at any precise point. This leads us to the following technical statement: 

\begin{lm}[Main Lemma I]
Let $d$ be a convex cone-metric on $S_g$, $V=V(d)$, $P^1=P(d, V, h^1)$ be a convex polyhedral Fuchsian manifold and $P^2=P(d, V, h^2)$ be a convex (polyhedral) Fuchsian cone-manifold. Let $W \subseteq V$ be a subset of vertices. Define $$\nu:=\sum\limits_{w \in W} \nu_w(d),~~m:=\min\limits_{w \in W} \tanh h^1_w,~~\tau:=\min\limits_{w \in W} \ln \left(\frac{\sinh h^2_w}{\sinh h^1_w}\right),$$ $$M:=\max\limits_{w \in W} \cosh^2\Big(\arcsinh\big(e^{\tau}\sinh h^1_w\big)\Big).$$ Assume that for all $w \in W$ we have $h^2_w > h^1_w$, thereby, $\tau > 0$. Then
$$S(P^1)-S(P^2) \geq \nu m \left( M e^{\tau/M}-\tau\right).$$
\end{lm}

Note that if $W=\{w\}$, then $M=\cosh^2 h^2_w$ and the statement is easier to perceive. At the first reading we advice the reader to consider this case.

We conjecture that the presence of $\nu$ in our Main Lemma I is excessive and it should be possible to give a bound without it. This would have simplified our exposition and also is helpful for a resolution of the Cohn-Vossen problem for Fuchsian manifolds.

\begin{rmk}
\label{monot}
The function $\nu m \left( M e^{\tau/M}-\tau\right)$ is positive and is increasing in $\nu$, $m$, $M$ and $\tau$ in the range $\nu>0$, $m>0$ and $M>\tau>0$. The definition of $M$ implies that $M>\tau$.
\end{rmk}

Now we discuss Main Lemma II. A natural way to approximate a metric $d$ by cone-metrics (not relying on a convex isometric embedding) is to take a sufficiently fine geodesic triangulation $\mathcal T$ and to replace each triangle by a triangle from a model space (e.g., $\H^2$) with the same side lengths. However, it turns out to be insufficient for our purposes. We propose a new way to approximate CBB metrics: we  replace each triangle of $\mathcal T$ by a triangle with the same side lengths, the same angles and one conical point in the interior. This is the subject of Main Lemma II. A big advantage is that the curvature of each vertex and each triangle of $\mathcal T$ coincides in both approximating and approximated metrics. We hope that Main Lemma II with its proof are of independent interest.

Recall from Subsection~\ref{triang} that $T(d)$ means the 6-tuple of side lengths and angles of a triangle $T$ in metric $d$. We consider these 6-tuples as elements of $\R^6$ endowed with $l_{\infty}$-norm. We also recall definitions of strictly short (Definition~\ref{geodtriang}) and $\delta$-fine (Definition~\ref{fine}) triangulations.

\begin{dfn}
Let $d$ be a cone-metric on $S_g$ and $\mathcal T$ be a triangulation of $S_g$. We call $d$ \emph{swept with respect to} $\mathcal T$ if $d$ realizes $\mathcal T$ and each triangle $T$ of $\mathcal T$ in $(S_g, d)$ has at most one conical point in the interior.
\end{dfn}

\begin{dfn}
Let $(S_g, d)$ be a CBB($-1$) metric space and $\mathcal T$ be its geodesic triangulation. A convex cone-metric $\hat d$ on $S_g$ is called a \emph{sweep-in} of $d$ with respect to $\mathcal T$ if 

(1) $\hat  d$ is swept with respect to $\mathcal T$;

(2) for each triangle $T$ of $\mathcal T$ we have $T(d)=T(\hat d)$.
\end{dfn}

\begin{lm}[Main Lemma II]
There exists an absolute constant $\delta>0$ such that if $d$ is a convex cone-metric on $S_g$ and $\mathcal T$ is its $\delta$-fine triangulation into convex triangles, or $d$ is a CBB($-1$) metric and $\mathcal T$, in addition, is strictly short, then there exists a unique sweep-in $\hat d$ of $d$ with respect to $\mathcal T$. 
\end{lm}

Next we turn to Main Lemma III, which appears to be decomposed into three parts. It is the core of our proof strategy. Assume that we have two cone-metrics $d^1$ and $d^2$ that are very close in some way and there is a convex polyhedral Fuchsian manifold $P^1=P(d^1, V(d^1), h^1)$. We claim that there is a convex Fuchsian cone-manifold $P^{21}=P(d^2, V(d^2), h^{21})$ such that $\widetilde h^{21}(p) \geq \widetilde h^1(p)$ at some ``essential'' points $p$ and $S(P^{21})$ is not significantly smaller than $S(P^1)$. The idea behind this is that we can transform $d^1$ to $d^2$ with the help of some elementary operations and modify $P^1$ along this transformation so that (1) we control the change of heights of some points; (2) we control the decrement of $S$.

Our transformation of $d^1$ to $d^2$ goes roughly: (1) do a sweep-in of $d^1$; (2) transform the metric in ``the class'' of sweep-ins; (3) restore $d^2$ from its sweep-in. The transformation of cone-manifolds corresponding to the first part is Main Lemma IIIA, to the second part is Main Lemma IIIC, to the third part is Main Lemma IIIB.

\begin{lm}[Main Lemma IIIA]
There exists an absolute constant $\delta>0$ with the following properties.
Let $d$ be a convex cone-metric on $S_g$, $\mathcal T$ be a short $\delta$-fine triangulation of $(S_g, d)$, $\hat d$ be a sweep-in of $d$ with respect to $\mathcal T$ and $P=P(d, V(d), h)$ be a convex Fuchsian cone-manifold. There exists $\hat h \in H(\hat d, V(\hat d))$ such that

(1) $\hat h_v \geq h_v$ for each $v \in (V(\mathcal T)\cap V(d))$;

(2) $S(\hat P)=S(\hat d, V(\hat d), \hat h) \geq S(P).$
\end{lm}

\begin{lm}[Main Lemma IIIB]
For any numbers $\e, A, D>0$ there exists $\delta=\delta(\e, A, D)$ with the following properties. Let $d$ be a convex cone-metric on $S_g$ with  $$\diam(S_g, d)<D,~~~\area(S_g, d)<A;$$ 
$\mathcal T$ be a $\delta$-fine triangulation of $(S_g, d)$,
$\hat d$ be a sweep-in of $d$ with respect to $\mathcal T$ and  $\hat P=P(\hat d, V(\hat d), \hat h)$ be a convex Fuchsian cone-manifold. Then there exists $h \in H(d, V(d))$ such that

(1) $h_v \geq \hat h_v$ for each $v \in (V(\mathcal T)\cap V(d))$;

(2) $S(P)=S(d, V(d), h) \geq S(\hat P) - \e.$
\end{lm}

\begin{lm}[Main Lemma IIIC]
Let $\hat d$ be a convex cone-metric on $S_g$ swept with respect to a triangulation $\mathcal T$. For each $\e>0$ there exists $\delta=\delta(\e, \hat d, \mathcal T)>0$ with the following properties. Let $d^1$ and $d^2$ be convex cone-metrics on $S_g$ swept with respect to $\mathcal T$
and for each triangle $T$ of $\mathcal T$ we have 
$$||T(d^1)-T(\hat d)||_{\infty}<\delta,~~~||T(d^2)-T(\hat d)||_{\infty}<\delta.$$
Let $P^1=P(d^1, V(d^1), h^1)$ be a convex Fuchsian cone-manifold. Then there exists $h^{21} \in H(d^2, V(d^2))$ such that

(1) $h^{21}_v \geq h^1_v$ for each $v \in V(\hat d)$;

(2) $S(P^{21})=S(d^2, V(d^2), h^{21}) \geq S(P^1) - \e.$
\end{lm}

\subsection{Proof of Theorem~\ref{main}}
\label{proofmain}

Suppose the contrary. We identify the upper boundaries of $F^1$ and $F^2$ with the help of $f$ and further identify them with $(S_g, d)$ for a CBB($-1$) metric $d$ (see the end of Subsection~\ref{cbbheight}). Let $h^1$ and $h^2$ be the corresponding height functions. Because the heights of non-flat points together with the upper boundary metric determine uniquely a Fuchsian manifold with convex boundary (Lemma~\ref{heightrig} and Lemma~\ref{flatlm}), there exists a non-flat point $u \in (S_g, d)$ such that $h^1(u)\neq h^2(u)$. Without loss of generality assume that ${h^2(u)-h^1(u)>0}$. Define $$H:=\frac{\min\{h^1(u), h^2(u)-h^1(u)\}}{3}.$$

Let $U$ be an (open) geodesic triangle with ${\diam(U, d)<H}$ containing $u$ in the interior. Define 
$$D:=2\diam(S_g, d),~~~A:=2\area(S_g, d),~~~\nu:=\nu(U,d)/2>0$$
$$m:=\tanh (h^1(u)-H)>0,~~\tau:=\ln \left(\frac{\sinh (h^2(u)-H)}{\sinh (h^1(u)+H)}\right)>0,$$ $$M:= \cosh^2\Big(\arcsinh\big(e^{\tau}\sinh(h^1(u)-H)\big)\Big)>0,$$
$$\e=\frac{1}{5}\nu m\left(Me^{\tau/M}-\tau\right)>0.$$

Take $\delta_1$ as the minimum of $\delta$ from Main Lemma II, IIIA and from IIIB for these $\e, A, D$. Construct a strictly short $\delta_1$-fine geodesic triangulation $\mathcal T$ of $(S_g, d)$ with the help of Lemma~\ref{goodtriang}. 
We can assume that

(1) $u \in V(\mathcal T)$,

(2) $\mathcal T$ subdivides $U$ in the sense that there is a subset $\mathcal W$ of triangles of $\mathcal T$ such that $U^c=\cup_{T \in \mathcal W} T^c$, where $U^c$, $T^c$ are the closures of $U$, $T$.

Indeed, if this is not the case, we can refine $\mathcal T$ to a triangulation $\mathcal T'$ for which this is true. It might happen that the new triangulation is not strictly short. But due to  Lemma~\ref{goodtriang}, we can refine it further, finally obtaining a triangulation fulfilling all our demands, which we continue to denote by $\mathcal T$. 

By $W'$ denote the set of those vertices of triangles from $\mathcal W$ that belong to $U$ (recall that $U$ is open). So we have $\sum_{T\in \mathcal W} \nu (T, d)+\sum_{w \in W'}\nu(w, d)=\nu(U, d)=2\nu$. 

Due to Main Lemma II, there exists a sweep-in $\hat d$ of $d$ with respect to $\mathcal T$. Let $\delta_2$ be $\delta$ from Main Lemma IIIC for $\e$, $\hat d$ and $\mathcal T$.

We have two convex Fuchsian manifolds with boundary $F^1=F(d, h^1)$ and $F^2=F(d, h^2)$. Choose a sufficiently dense set $V \subset (S_g, d)$ containing $V(\mathcal T)$. Define $P^1=P(d^1, V(d^1), h^1)$ and $P^2=P(d^2, V(d^2), h^2)$ to be convex polyhedral Fuchsian manifolds obtained by taking the convex hull of $V$ in $F^1$ and $F^2$ respectively. Here $d^1$ and $d^2$ are induced metrics on the boundaries pulled back to $S_g$ with the help of the vertical projections, and, abusing the notation, we denote restrictions of $h^1$, $h^2$ to $V(d^1)$, $V(d^2)$ still by $h^1$ and $h^2$. 

Because $u$ is non-flat in $(S_g, d)$, we have $u \in (V(d^1)\cap V(d^2))$.

Due to the convergence of inscribed polyhedral manifolds (Lemma~\ref{refin} and Corollary~\ref{diamareaglob}) we can choose $V$ such that 

(1) both $d^1$ and $d^2$ realize $\mathcal T$;

(2) $\diam(S_g, d^1)<D,~~~\diam(S_g, d^2)<D,~~~\area(S_g, d^1)<A,~~~~\area(S_g, d^2)<A;$

(3) $\mathcal T$ is $\delta_1$-fine on $d^1$, $d^2$;

and for each triangle $T$ of $\mathcal T$ we have 

(4) $||T(d^1)-T(d)||_{\infty}<\delta_2,$ $||T(d^2)-T(d)||_{\infty}<\delta_2.$

Because $T(d)=T(\hat d)$ we can rewrite (4) as

(4) $||T(d^1)-T(\hat d)||_{\infty}<\delta_2,$ $||T(d^2)-T(\hat d)||_{\infty}<\delta_2.$


Denote $W:=V\cap U$. As all the angles and the areas of the triangles of $\mathcal T$ in $d^1$ can be chosen to be arbitrarily close to those in $d$, we can assume also that
$$\nu':=\sum_{T\in \mathcal W} \nu (T, d^1)+\sum_{w \in W'}\nu(w, d^1)=\nu(U, d^1)=\sum_{w \in W} \nu_w(d^1)> \nu.$$

Let $\hat d^1$, $\hat d^2$ be sweep-ins of $d^1$, $d^2$ with respect to $\mathcal T$. They exist due to Main Lemma II as $d^1$, $d^2$ are $\delta_1$-fine. 

First we apply Main Lemma IIIA to $P^2$ and $\hat d^2$. We get $\hat h^2 \in H(\hat d^2, V(\hat d^2))$ such that

(1) $\hat h^2_u \geq h^2_u$;

(2) $S(\hat P^2)=S(\hat d^2, V(\hat d^2), \hat h^2) \geq S(P^2).$

Now we apply Main Lemma IIIC to $\hat P^2$ and $\hat d^1$. We get $\hat h^{12} \in H(\hat d^1, V(\hat d^1))$ such that

(1) $\hat h^{12}_u \geq \hat h^2_u$;

(2) $S(\hat P^{12})=S(\hat d^1, V(\hat d^1), \hat h^{12}) \geq S(\hat P^2) - \e.$

Finally, we take $\hat P^{12}$, $d^1$ and apply Main Lemma IIIB. We get $h^{12} \in H( d^1, V(d^1))$ such that

(1) $h^{12}_u \geq \hat h^{12}_u$;

(2) $S(P^{12})=S(d^1, V(d^1), h^{12}) \geq S(\hat P^{12}) - \e.$

Summarizing all last three steps together we see:

(1) $h^{12}_u \geq h^{2}_u$;

(2) $S(P^{12})=S(d^1, V(d^1), h^{12}) \geq S(P^2)-2\e$.

We get $h^{12}_u \geq  h^2_u \geq h^1_u+ 3H$.

Now we do the same but starting from $P^1$ and transforming $d^1$ to $\hat d^1$ to $\hat d^2$ to $d^2$. We obtain $h^{21} \in H(d^2, V(d^2))$ such that (here we need only properties (2))

$$S(P^{21})=S(d^2, V, h^{21}) \geq S(P^1)-2\e.$$

In total, we summarize our main transformations in the following diagram:
$$F^1=F(d, h^1) \longrightarrow P^1=P(d^1, V(d^1), h^1) \longrightarrow P^{21}=P(d^2, V(d^2), h^{21})$$
$$F^2=F(d, h^2) \longrightarrow P^2=P(d^2, V(d^2), h^2) \longrightarrow P^{12}=P(d^1, V(d^1), h^{12}).$$

Now let $$m':=\min\limits_{w \in W} \tanh (h^1_w),~~\tau':=\min\limits_{w \in W} \ln \left(\frac{\sinh h^{12}_w}{\sinh h^1_w}\right),$$ $$M':=\max\limits_{w \in W} \cosh^2\Big(\arcsinh\big(e^{\tau}\sinh(h^1_w)\big)\Big)$$
and $\nu'$ is defined above.
Then Main Lemma I implies that
$$S(P^1)-S(P^{12}) \geq \nu' m'\left(M'e^{\tau'/M'}-\tau'\right).$$

Recall that we chose $U$ such that $\diam(U, d)<H$. Thus, for each $w \in W$ we have $d(u, w)<H$. We know $d^1 \leq d$, $d^2 \leq d$ due to Corollary~\ref{BusFel}. Thus, for each $w \in W$ we have $d^1(u, w)<H$, $d^2(u, w)<H$. From this we conclude the bounds on heights with the help of Lemma~\ref{htriangleinequality}: for every $w \in W$ we obtain $$h^{12}_w\geq h^{12}_{u}-H\geq h^2_u-H,$$ $$h^1_u-H\leq h^1_w \leq  h^1_{u}+H.$$ 
From this we deduce that that $m' \geq m$, $\tau' \geq \tau$ and $M' \geq M$. Recall also that $\nu' \geq \nu$. Thus, due to Remark~\ref{monot} we get
$$\nu' m'\left(M'e^{\tau'/M'}-\tau'\right)\geq \nu m\left(Me^{\tau/M}-\tau\right)=5\e.$$

By Lemma~\ref{maxprinciple} we have
$$S(P^2)-S(P^{21})\geq 0.$$

Summing this up we have
$$4\e = 2\e+2\e \geq \Big(S(P^{2}) - S(P^{12})\Big)+ \Big(S(P^1) - S(P^{21})\Big)=$$
$$=\Big(S(P^2)-S(P^{21})\Big) + \Big(S(P^1)-S(P^{12})\Big) \geq 0+5\e=5\e.$$

Thus, $\e \leq 0$ and we obtain a contradiction.

\section{Properties of Fuchsian cone-manifolds}
\label{conechap}


In this section we examine some important properties of Fuchsian cone-manifolds.

\subsection{Ultraparallelism} 

\begin{lm}
\label{ultrapar2}
Let $P=P(d, \mathcal T, h)$ be a convex Fuchsian cone-manifold. Then each prism of $P$ is ultraparallel. 
\end{lm}

\begin{proof}
Suppose the contrary. Then we construct an infinite sequence of pairwise distinct prisms of $P$.

Suppose that there is a prism $\Pi$ of $P$ developed to $\H^3$ as $A_1A_2A_3B_1B_2B_3$ such that its upper boundary plane $M^{\uparrow}$ intersects the lower boundary plane $M_{\downarrow}$ in line $L$ under dihedral angle $\phi$.

We call the distance function $w$ to $L$ over the triangle $A_1A_2A_3$ the \emph{weight function} and its minimum the \emph{weight} $W(\Pi)$. It is easy to see that the point $p$, where $W(\Pi)$ is attained, is unique and belongs to the boundary of $A_1A_2A_3$. We have
$$\sinh \widetilde h(p)=\sinh w(p) \sin \phi.$$

Suppose that $p$ is an interior point of an edge (say, $A_2A_3$). Then the line $A_2A_3$ is ultraparallel to $L$. We observe that that the dihedral angle of $\Pi$ at $A_2A_3$ is greater than $\pi/2$. Take next the prism $\Pi'$ adjacent to $\Pi$ at the edge corresponding to $A_2A_3$. Embed it to $\H^3$ adjacent to $A_1A_2A_3B_1B_2B_3$. Then its upper boundary plane also intersects the lower boundary plane $M_{\downarrow}$ in line $L'$. Let $\phi'$ be the respective dihedral angle. Due to the convexity condition, the angle of $\Pi'$ at $A_2A_3$ is smaller than $\pi/2$ and $\phi \leq \phi' \leq \pi/2$. To see this one considers the orthogonal plane section to lines $L$ and $A_2A_3$, which exists as they are ultraparallel. We get $w'(p) \leq w(p)$, where $w'$ is the weight function of $\Pi'$ defined as the distance to $L'$, and $w'(p)=w(p)$ if and only if $\phi=\phi'$.

If the weight of $\Pi'$ is attained at $p$, then its dihedral angle at the edge containing $p$ is also greater than $\pi/2$. This can not happen as the total angle is at most $\pi$. Therefore, it can not be attained at $p$ and $W(\Pi')<W(\Pi)$, where $W(\Pi')$ is the minimum of $w'$.

Now suppose that $W(\Pi)$ is attained at a vertex of $\Pi$ (say, $A_3$). We claim that among the edges incident to $A_3$ there exists at least one with the dihedral angle greater than $\pi/2$. Indeed, either the line containing one edge does not intersect $L$ and $A_1A_2A_3$ lies from the different side of this line than $L$, then this edge is the desired. Or there is an edge such that $A_1A_2A_3$ lies in the non-obtuse angle formed by $L$ and the line containing this edge. Let this edge be $A_2A_3$.

Now take the prism $\Pi'$ adjacent to the edge corresponding to $A_2A_3$ and embed $\Pi'$ to $\H^3$ adjacent to $A_1A_2A_3B_1B_2B_3$. Similarly to the case before, convexity implies that the upper plane of $\Pi'$ intersects $M_{\downarrow}$ in line $L'$ with the dihedral angle $\phi'$ satisfying $\phi \leq \phi' \leq \pi/2$. If $\phi'>\phi$, then $w'(A_3)<w(A_3)$ and, therefore, $W(\Pi')<W(\Pi)$. Suppose that $W(\Pi')=W(\Pi)$, then $W(\Pi')$ is also attained at $A_3$ and $\phi'=\phi$, so the upper plane of $\Pi'$ coincides with $M^{\uparrow}$.

By $L_0$ denote the line in $M^{\uparrow}$ through $A_3$ orthogonal to $L$ oriented outwards to $L$. Assume that $A_3A_2$ lies to the left of $L_0$ and $\gamma$ be the angle between $L_0$ and $A_3A_2$. The condition that $W(\Pi)$ is attained at $A_3$ implies that $\gamma \geq \pi/2$. 

Let $A_3A_4$ be the edge of $\Pi'$ incident to $A_3$. As $W(\Pi)=W(\Pi')$, then $A_3A_4$ also lies to the left of $L_0$ and for its angle $\gamma'$ with $L_0$ we have $\gamma > \gamma' \geq \pi/2$.

In this way we construct an infinite sequence of prisms of $P$. Two prisms in this sequence are distinguished either by their weight, which decreases monotonously, or, if the weights are equal, by the angle $\gamma$, which strictly decreases as long as the weight remains constant. Clearly, both $W(\Pi)$ and $\gamma$ are independent of the embedding. This shows a contradiction. 

Now it remains to consider the case when the initial prism $\Pi$ is asymptotically parallel, i.e., the upper boundary plane $M^{\uparrow}$ intersects the lower boundary plane $M_{\downarrow}$ only at infinity. Let $L$ be a horocycle in $M^{\uparrow}$ through the intersection point. We can similarly define the weight function as the distance to $L$. All arguments above remain valid. Moreover, either we will find in our sequence a prism such that its upper boundary plane intersects the lower boundary plane in a line, and then we are reduced to the case above, or all prisms in the sequence after development will have $M^{\uparrow}$ as their boundary plane.
\end{proof}

Let  $P(d, \mathcal T, h)$ be convex and $\widetilde h$ be the extended height function. Take $T \in \mathcal T$ and embed the prism $\Pi$ containing $T$ to $\H^3$. Lemma~\ref{ultrapar2} say that this prism is ultraparallel. Let $A$ and $B$ be the closest points on the upper and lower boundary planes respectively. Due to Corollary~\ref{rightangl}, $\widetilde h$ satisfies for every $p \in T$
\begin{equation}
\label{distfunc}
\sinh \widetilde h(p) = \sinh AB \cosh({\rm dist}(p,A)).
\end{equation}
In particular, its restriction to a geodesic segment in $P$ has the form
\begin{equation}
\label{pd}
\arcsinh(b\cosh(x-a)).
\end{equation}
We highlight that $A$ and $B$ may lie outside of the image of $\Pi$ in $\H^3$. However, neither ${\rm dist}(x, A)$ nor the distance $AB$ depend on the choice of representation or embedding. 

\subsection{Heights define a convex cone-manifold}
\label{heightdefsec}

In this subsection we prove Lemma~\ref{heightdefine}. We recall Definition~\ref{f-1concave} of an $\mathcal F(-1)$-concave function.

Consider a convex Fuchsian cone-manifold $P(d, \mathcal T, h)$ and a unit speed geodesic $\psi: [0; \tau] \rightarrow (S_g, d)$. Let $x_1, \ldots, x_k \in (0;\tau)$ be its intersection points with strict edges of $P$. We set $x_0:=0$ and $x_{k+1}:=\tau$. Due to~(\ref{pd}), on each segment $[x_i; x_{i+1}]$ the restriction of $\widetilde h$ has the form
$$\widetilde h(x)=\arcsinh(b_i\cosh(x-a_i))$$
for some positive $b_i$. Moreover, the convexity implies that for each kink point $x_i$ we have the left derivative greater than the right derivative (see $f(x)$ on Figure~\ref{Pic5}).

\begin{lm}
\label{distf-1conc}
The function $\sinh \widetilde h$ is $\mathcal F(-1)$-concave.
\end{lm}

\begin{proof}
Let $g: [x^0; x^1] \rightarrow \R$ be in $\mathcal F(-1)$ and $f: [x^0; x^1] \rightarrow \R$ be a function such that

(1) $g(x^0)=f(x^0)$, $g(x^1)=f(x^1)$;

(2) there is a subdivision $$[x^0; x^1] = [x_0; x_1]\cup\ldots\cup[x_k; x_{k+1}]$$ ($x^0=x_0$ and $x_{k+1}=x^1$) such that the restriction of $f$ to each $[x_i; x_{i+1}]$ is $\mathcal F(-1)$;

(3) at every point $x_i$ the left derivative of $\tilde{h}(x)$ is greater than the right derivative.

We claim that $f(x) \geq g(x)$ for all $x \in [x^0; x^1]$. Due to the discussion above, this will imply that $\sinh\widetilde h$ is $\mathcal F(-1)$-concave. 

\begin{prop}
Two distinct $\mathcal F(-1)$-functions on $\R$ can not coincide at more than one point.
\end{prop}

This follows from the definition of the $\mathcal F(-1)$-class.

\begin{prop}
Let $f_i(x): \R \rightarrow \R$ be the extension to $\R$ of the restriction $f|_{[x_i; x_{i+1}]}$. For all $x \not\in [x_i; x_{i+1}]$ we have $f_i(x) > f(x)$.
\end{prop}

Indeed, consider $j>i$ and prove by induction over $j-i$ that for $x \in [x_j; x_{j+1}]$ we have $f_i(x) > f_j(x)$. The base case $j=i+1$ follows from Claim 1 and the fact that in the kink point $x_{i+1}$ the derivative of $f_i$ is greater or equal than the derivative of $f_{i+1}$. The inductive step is obvious. The case $j<i$ is the same.

\begin{prop}
For all $x \in (x_0; x_1]$ we have $f(x) \geq g(x)$.
\end{prop}

Indeed, Claim 1 and the assumption $f(x^0)=g(x^0)$ imply that the sign of the difference $f(x) - g(x)$ is constant on $(x_0; x_1]$. Assume that it is negative, i.e., for every $x > x^0$ we have $f_0(x) < g(x)$. Also assume $k \neq 0$. Substituting $x=x^1$ and using Claim 2 we obtain that $f(x^1)=f_k(x^1) < f_0(x^1) < g(x^1)$, which contradicts the statement.
\vskip+0.2cm

Now we may assume that the sign of the difference $f(x) - g(x)$ is positive on $(x_0; x_1]$: if it is zero, then we just cut off $[x_0; x_1]$ and proceed in the same way. Under this assumption, suppose that for some $i$ and $x' \in (x_i; x_{i+1}]$, $x' \neq x^1$, we have $f(x') = g(x')$. Due to Claim 1 and the assumption, we have $f_i(x) < g(x)$ for all $x > x'$. Thus, due to Claim 2 we have $f(x) < g(x)$  for all $x > x'$. We obtain a contradiction. 

It implies that $f(x) \geq g(x)$ over the interval $[x^0; x^1]$.
\end{proof}

We note a simple property of $\mathcal F(-1)$-concave functions similar to a property of ordinary concave functions. Let $f: I \rightarrow \R$ be a $\mathcal F(-1)$-concave function and $g \in \mathcal F(-1)$ coincides with $f$ at points $x_1, x_2 \in I$. Then for all $x \not\in [x_1; x_2]$ we get $f(x) \leq g(x)$. If for some $x \in (x_1; x_2)$ we have $f(x) > g(x)$, then $f(x)>g(x)$ for all $x \in (x_1; x_2)$ and for all $x \not\in [x_1; x_2]$ we get $f(x) < g(x)$.

Now we are going to prove that a convex Fuchsian cone-manifold is uniquely determined by $d$ and $h$. If we know $d$ and $h$, then it remains only to restore its face decomposition. 

\begin{replm}{heightdefine}
Let $\mathcal T^1$ and $\mathcal T^2$ be two triangulations with $$V(\mathcal T^1)=V(\mathcal T^2)=V,$$ $P^1=P(d, \mathcal T^1, h)$ and $P^2=P(d, \mathcal T^2, h)$ be two convex Fuchsian cone-manifolds. Then $P^1$ is marked isometric to $P^2$ with respect to $V$.
\end{replm}

\begin{proof}
Take the intersection point $v'$ of two edges $e^1 \in E(\mathcal T^1)$ and $e^2 \in E(\mathcal T^2)$. The restriction of $\sinh \widetilde h^1$ to $e^1$ is $\mathcal F(-1)$. Due to Lemma~\ref{distf-1conc}, the restriction of $\sinh \widetilde h^2$ to $e^1$ is $\mathcal F(-1)$-concave. As they coincide at the endpoints of $e^1$, we get $\sinh \widetilde h^2(v') \geq \sinh \widetilde h^1(v')$. Applying the same reasoning to $e^2$, we get $\sinh \widetilde h^2(v') \leq \sinh \widetilde h^1(v')$. Hence, $ \widetilde h^2(v') =  \widetilde h^1(v')$.

Let $V'$ be $V$ together with all the intersection points between the edges of $\mathcal T^1$ and $\mathcal T^2$. We obtained $\widetilde h^1|_{V'}=\widetilde h^2|_{V'}$. The union of geodesic triangulations $\mathcal T_1$ and $\mathcal T_2$ decomposes the metric space $(S_g, d)$ into geodesic polygons with vertices in $V'$. Let $\mathcal T'$ be a geodesic triangulation refining this decomposition. The convex Fuchsian cone-manifold $P(d, \mathcal T', \widetilde h^1|_{V'})$ is marked isometric to both $P^1$, $P^2$ due to Corollary~\ref{prismuniq}. This finishes the proof. 
\end{proof}

\subsection{The space of admissible heights}
\label{compactlemmas}

Recall from Subsection~\ref{disccurv} that $H(d,V) \subset \R^n$ is the set of all admissible heights for the pair $(d, V)$, i.e., it is a set of convex Fuchsian cone-manifolds for fixed $d$ and $V$. In this subsection we collect key facts about $H(d, V)$ that will be used in the proofs of the main lemmas.

\begin{lm}
\label{movingup}
Let $d \in \mathfrak D_c(V)$ and $h \in H(d, V)$. Define $h_t$ by $\sinh h_{t, v}=e^t \sinh h_v$ for all $v \in V$. Then $h_t \in H(d, V)$ for all $t \in \R$, and all Fuchsian cone-manifolds $P_t=P(d, V, h_t)$ have the same face decomposition. 
\end{lm}

\begin{proof}
Consider two adjacent prisms such that the dihedral angle of the adjacent upper edge is $\pi$. Develop them to $\H^3$ in the half-space model, let $M_{\downarrow}$ be the plane containing the lower boundaries and $M^{\uparrow}$ --- upper. By $A \in M^{\uparrow}$ and $B \in M_{\downarrow}$ denote the closest points. Assume that the line $AB$ intersects the ideal boundary at point $C$ and $B$ is between $C$ and $A$. The homothety $\rho$ with center $C$ is a hyperbolic isometry. 
Consider the images of all upper vertices of the prisms and their projections to $M_{\downarrow}$. We obtain two new prisms with upper boundaries isometric to those of former prisms and the dihedral angle of the common upper edge remains equal to $\pi$. 

Assume that $\cosh(\rho(A)B)=e^{t}\cosh AB.$ Let $h_{\rho}$ be the distance function of the new vertices. Due to formula~(\ref{distfunc}), $\sinh h_{\rho, v} = e^{t}\sinh h_v$ for each vertex $v$ of the prisms.

Let $\mathcal T$ be a triangulation compatible with $P$, so $P=P(d, \mathcal T, h)$. We saw that all flat edges remain flat in all $P(d, \mathcal T, h_t)$. But for sufficiently small $t$ all strict edges remain to be strict. Thus if $t$ is sufficiently small, then $P(d, \mathcal T, h_t)$ is convex and $h_t \in H(d, V)$. However, if a strict edge becomes flat, then it is flat in all $P_t$. Thus, $h_t \in H(d, V)$ for all $t \in \R$ and all $P_t$ have the same face decomposition.
\end{proof}

\begin{lm}
\label{finite}
For each $V \subset S_g$ and $d \in \mathfrak D_c(V)$ there are finitely many $\mathcal T$ with $V(\mathcal T)=V$ for which there is $h \in  H(d, V)$ such that $\mathcal T$ is compatible with $P(d, V, h)$.
\end{lm}

\begin{proof}
Due to Lemma~\ref{movingup}, it is enough to prove that there are finitely many $\mathcal T$ compatible with convex cone-manifolds $P=P(d, V, h)$ with $h_v \geq m$ for some $m>0$ and all $v \in V$. 

Let $e$ be an edge of $P$. Develop the trapezoid containing $e$ as $A_1A_2B_1B_2 \subset \H^2$. Let $A$ and $B$ be the closest points on the lines containing its upper and lower boundaries. Assume that $A_1$ lies between $A_2$ and $A$. Denote the distance $AB$ by $h$, $AA_1$ by $l_1$ and $AA_2$ by $l_2$. Then from Corollary~\ref{rightangl} we get
$$\sinh h_1=\cosh l_1 \sinh h,~~~\sinh h_2=\cosh l_2 \sinh h.$$
Then we have 
$$\frac{\sinh h_2}{\sinh h_1} = \cosh(h_2-h_1)+\sinh(h_2-h_1)\coth h_1 = \frac{\cosh l_2}{\cosh l_1}\geq \cosh l_{12}.$$

Due to Lemma~\ref{htriangleinequality}, we have $h_2-h_1 \leq \diam(S_g, d)$. Then we get
$$\cosh l_{12} \leq \cosh\big(\diam(S_g, d)\big)(1+\coth m).$$

If $A$ lies between $A_1$ and $A_2$, then one of $l_1$, $l_2$ is at least $l_{12}/2$. Assume that this is $l_1$. Note that $h_1-h \leq \diam(S_g, d)$. Considering the trapezoid $AA_1BB_1$, we get similarly
$$\cosh \frac{l_{12}}{2} \leq \cosh\big(\diam(S_g, d)\big)(1+\coth m).$$

Hence there exists a constant $M=M(d, m)$ such that the lengths of edges of $P$ are at most $M$. But the lengths of geodesics between points of $V$ in $(S_g, d)$ form a discrete set by an argument from~\cite[Proposition 1]{ILTC}. Thus there are finitely many of them that can appear as edges of some $P$. Therefore, there are finitely many realizable triangulations.
\end{proof}

Define
$$H(d, V, m, M) = \{h \in H(d,V): \min_{v \in V} h_v \geq m, \max_{v \in V} h_v \leq M\}.$$ 

\begin{lm}
\label{H_d comp}
For every $d \in \mathfrak D_c(V)$ and $0<m \leq M$ the set $H(d, V, m, M)$ is compact.
\end{lm}

\begin{proof}
From Lemma~\ref{finite} there are finitely many triangulations compatible with $P(d, V, h)$ for $h \in H(d, V)$. For every triangulation $\mathcal T$ the set of heights $h$ bounded between $m$ and $M$, such that $(d, \mathcal T, h)$ is a representable triple, is compact. Its subset of admissible heights $h$, such that $\mathcal T$ is compatible with the convex cone-manifold $P(d, V, h)$, is defined by considering also inequalities $\phi_e(P) \leq \pi$ for all edges of $\mathcal T$. Thus, it is a closed subset, hence, also compact. Finally, $H(d, V, m, M)$ is compact as a finite union of compact sets.
\end{proof}


Now we are going to establish two fundamental properties used in several parts of our proof: even if $h$ is on the boundary of $H(d, V)$, the height of a point with non-positive particle curvature can always be decreased (Corollary~\ref{decrease}) and the heights of all points with positive particle curvatures can always be simultaneously increased (Lemma~\ref{increase}). Together with this we show how non-positive/positive particle curvature impacts on the upper boundary geometry.

\begin{lm}
\label{negcurv}
Let $P(d, V, h)$ be a convex Fuchsian cone-manifold, $v \in V$ and either $\kappa_v(P) < 0$ or $\kappa_v(P)=0$ and $\nu_v(d)>0$. Then all angles at $v$ of faces incident to $v$ are strictly less than $\pi$.
\end{lm}

\begin{proof}
After passing to the spherical link of $v$ this becomes a lemma from spherical geometry, which is essentially due to Volkov~\cite{Vol}. A proof can also be found in~\cite[Lemma 5.3]{Izm1}.
\end{proof}

\begin{lm}
\label{dihed}
Let $A_1A_2A_3B_1B_2B_3$ and $A_2A_3A_4B_2B_3B_4$ be two adjacent compact prisms in $\H^3$ such that the total dihedral angle $\phi$ of the common edge $A_2A_3$ is $\pi$. By $\lambda$ denote the total angle of $A_2$ in the upper faces. Then
$$\sgn \frac{\partial \phi}{\partial h_3} = \sgn(\lambda-\pi).$$
\end{lm}

\begin{proof}
We highlight that $\lambda$ is the face angle opposite to $A_3$! For the angle $\alpha_{23}$ at $A_2$ in the trapezoid $A_2A_3B_2B_3$ from Lemma~\ref{coslaw} we have $$\frac{\partial \alpha_{23}}{\partial h_3}>0.$$ Consider the spherical link of $A_2$ in the union of both prisms. It consists of two spherical triangles glued altogether. When $h_3$ changes, only the common edge $\alpha_{23}$ changes and the derivative is positive. Let $\phi^+$, $\phi^-$ be the dihedral angles of $A_2A_3$ in two prisms and $\lambda^+$, $\lambda^-$ be the angles at $A_2$ of their upper boundaries. Then
$$\frac{\partial \phi^+}{\partial \alpha_{23}}=-\frac{\cot\lambda^+}{\sin\alpha_{23}},~~~ \frac{\partial \phi^-}{\partial \alpha_{23}}=-\frac{\cot\lambda^-}{\sin\alpha_{23}}.$$
So we get
$$\frac{\partial \phi}{\partial h_3} = -\frac{\cot\lambda^++\cot\lambda^-}{\sin\alpha_{23}}\cdot\frac{\partial \alpha_{23}}{\partial h_3}.$$
This implies the desired.
\end{proof}

\begin{lm}
\label{decrease0}
Let $P=P(d, V, h)$ be a convex Fuchsian cone-manifold and $v \in V$ be a vertex that has angle less than $\pi$ in all faces. For $\xi>0$ define $h'$ by $h'_v:=h_v-\xi$ and $h'_u:=h_u$ for all $u \in V$, $u \neq v$. Then $h' \in H(d, V)$ provided $\xi$ is sufficiently small.
\end{lm}

\begin{proof}
The proof ideas are essentially due to Volkov~\cite{Vol}. It was also reproduced for cusps with particles in~\cite{FiIz1}. We sketch it here because the ideas from the proof are used further.

\begin{figure}
\begin{center}
\includegraphics[scale=0.4]{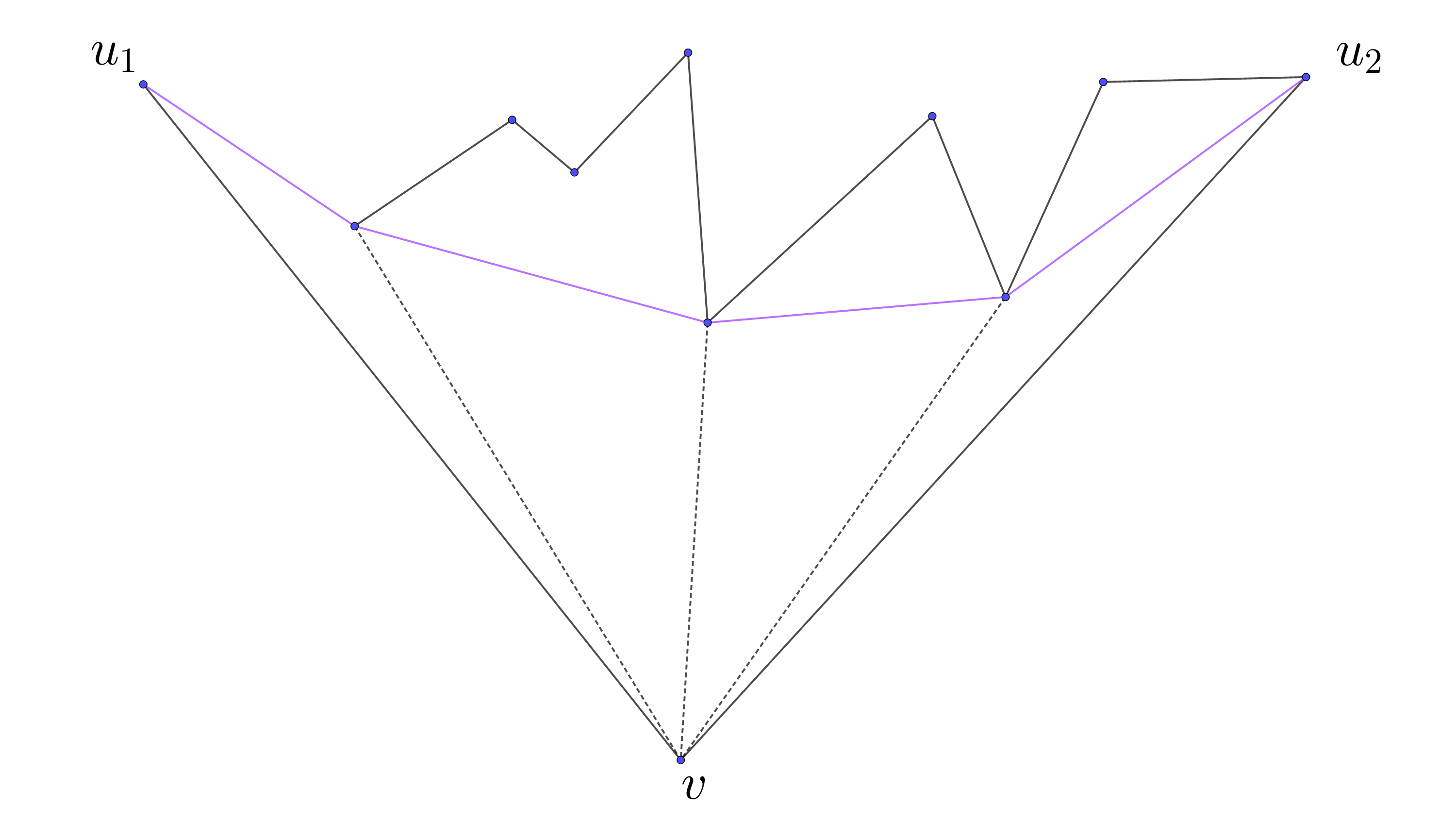}
\caption{To the proof of Lemma~\ref{decrease0}.}
\label{Pic6}
\end{center}
\end{figure}

Let $v$ be a vertex such that it has angle less than $\pi$ in all faces of $P$. We claim that each face can be triangulated such that all new flat edges become convex after we decrease the height of $v$. Let $R$ be a face incident to $v$ and $u_1$, $u_2$ be two vertices incident to $v$. Consider the shortest path in $R$ connecting $u_1$ and $u_2$ and homotopic to the polygonal curve $u_1vu_2$ (see Figure~\ref{Pic6}). This path is a polygonal curve that cuts off $R$ a polygon $R'$ such that it has angles greater than $\pi$ at all vertices except $v$, $u_1$, $u_2$. Triangulate it by diagonals from $v$. Lemma~\ref{dihed} shows that this triangulation is a desired one.
\end{proof}

Combining Lemma~\ref{negcurv} and Lemma~\ref{decrease0} we get

\begin{crl}
\label{decrease}
Let $P=P(d, V, h)$ be a convex Fuchsian cone-manifold and $v \in V$ be a vertex with either $\kappa_v(P) < 0$ or $\kappa_v(P)=0$ and $\nu_v(d)>0$. For $\xi>0$ define $h'$ by $h'_v:=h_v-\xi$ and $h'_u:=h_u$ for all $u \in V$, $u \neq v$. Then $h' \in H(d, V)$ provided $\xi$ is sufficiently small.
\end{crl}

Recall that a vertex of $P$ is called isolated if there are no strict edges emanating from it.

\begin{lm}
\label{isolvertex}
If $v$ is an isolated vertex of a Fuchsian cone-manifold $P=P(d, V, h)$, then $\kappa_v(P)=\nu_v(d)\geq 0.$
\end{lm}

\begin{proof}
Consider a triangulation $\mathcal T$ compatible with $P$ and start developing all prisms containing $v$ to $\H^3$ one by one preserving the incidence, each prism developed only once. We get a chain of prisms with all lower boundaries belonging to the same plane $M_{\downarrow}$ and all upper boundaries belonging to the same plane $M^{\uparrow}$ as all upper edges are flat. It is easy to see that as $\nu_v(d)\geq 0$, then $\kappa_v(P) \geq 0$ with the equality in one case implying the equality in the second case. Thereby, our chain of prisms is not self-intersecting. Assume that $\nu_v(d)>0$. The first and the last upper edge of the development are images of the same edge of $P$, therefore the distance function restricted to these edges is the same. Due to condition that the sum of the dihedral angles at this edge is $\pi$, we see that $v$ is developed to the closest point from $M^{\uparrow}$ to $M_{\downarrow}$. This gives the desired equality. 
\end{proof}

\begin{lm}
\label{increase}
Let $P=P(d, V, h)$ be a convex Fuchsian cone-manifold and $V=V(d)$. For $\xi>0$ define $h'$ by $\sinh h'_v:=e^{\xi}\sinh h_v$ if $\kappa_v(P)>0$ and $h'_u:=h_u$ otherwise. Then $h' \in H(d, V)$ provided $\xi$ is sufficiently small.
\end{lm}

\begin{proof}

Lemma~\ref{finite} shows that there are finitely many triangulations compatible with $P$. Consider $\xi$ sufficiently small such that for any $\mathcal T$ compatible with $P$ all strict edges of $P$ remain strict in $P(d, \mathcal T, h')$. We start from an arbitrary triangulation $\mathcal T$ and use the \emph{flip algorithm} to find a triangulation $\mathcal T'$ such that $P(d, \mathcal T', h')$ is convex. 
The flip algorithm will also be used in some further parts of the proof.

We simply take a concave edge $e$ of $P(d, \mathcal T, h')$ and flip it, i.e., replace by the other diagonal in the quadrangle formed by two triangles incident to $e$. There are two questions: why the algorithm can not run infinitely and why a concave edge can be flipped.

By induction it is easy to see that each triangulation that appears is compatible with $P$. Indeed, if this is true for the current triangulation, then by the choice of $\xi$ all strict edges of $P$ are strict in $P(d, \mathcal T, h')$. Therefore, they can not be flipped and the new triangulation is also compatible with $P$. The extended height function $\widetilde h$ is pointwise non-decreasing after each flip. Moreover, it strictly increases at all the interior points of the quadrangle defining the flip. This means that no triangulation appears twice during the algorithm. As there are finitely many of them, the algorithm can not run infinitely.

There are two cases, when an edge can not be flipped. First, it can be an edge of a triangle with two sides glued altogether, so there is no quadrangle at all. Second, it can be a diagonal in a concave quadrangle. We illustrate this by Figure~\ref{Pic7}.

\begin{figure}
\begin{center}
\includegraphics[scale=0.6]{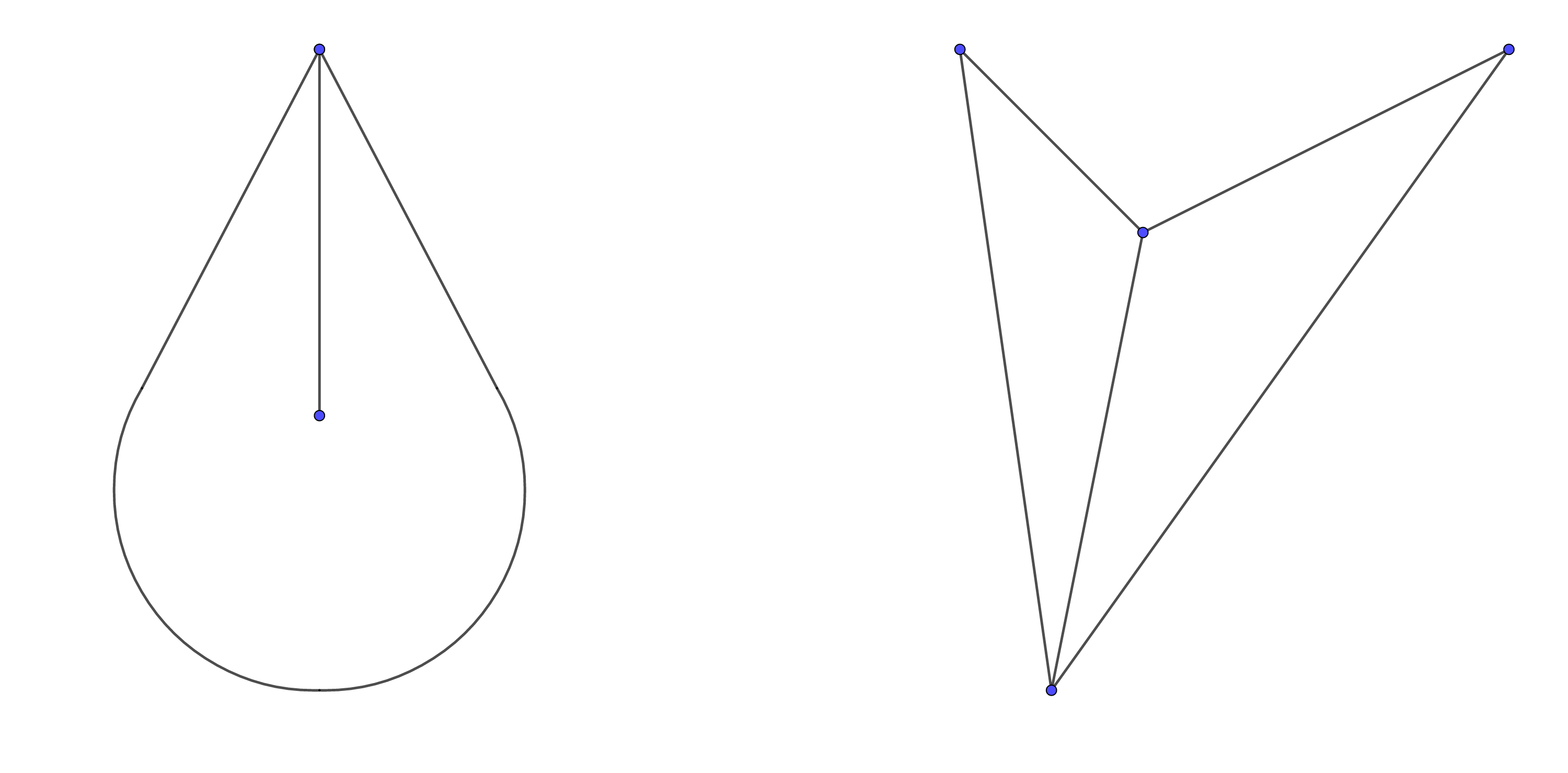}
\caption{Two cases when an edge of a triangulation can not be flipped.}
\label{Pic7}
\end{center}
\end{figure}

We prove that under our conditions this can not happen. Suppose that for some $P(d, \mathcal T, h')$ an edge $e$ is incident only to one triangle $T$, $v$ is the vertex incident only to $e$ and $u$ is the other vertex of $e$. Then $\nu_v(d)>\pi$. The edge $e$ is flat in $P$, thus, $v$ is isolated in $P$. Lemma~\ref{isolvertex} implies that $\kappa_v(P)>\pi$, so $\sinh h'_v=e^{\xi}\sinh h_v$. If also $\sinh h'_u=e^{\xi}\sinh h_u$, then $e$ remains flat in $P(d, \mathcal T, h')$ (see the proof of Lemma~\ref{movingup}). Consider the case $h'_u=h_u$. Note that as $T$ is isosceles, then both angles of $T$ at $u$ are smaller than $\pi/2$. Therefore, $e$ becomes convex in $P(d, \mathcal T, h')$ due to Lemma~\ref{dihed}.

Suppose that for some $P(d, \mathcal T, h')$ an edge $e$ is a diagonal of a concave quadrilateral and $v$ is the vertex of this quadrilateral with total angle at least $\pi$. The edge $e$ is flat in $P$. Lemma~\ref{negcurv} shows that $\kappa_v(P)>0$ (as $V=V(d)$, we have $\nu_v(d)>0$) and, therefore, $\sinh h'_v=e^{\xi}\sinh h_v$. If this is true for all vertices of the quadrilateral, then $e$ remains flat in $P(d, \mathcal T, h')$. We prove that in all other cases $e$ becomes convex.

Indeed, denote three other vertices by $u_1$, $u_2$, $u_3$ such that $e$ connects $v$ and $u_3$. If only $v$ increased its height, then the claim follows from Lemma~\ref{dihed}. Otherwise, we first increase the heights of all $u_1$, $u_2$, $u_3$ by $\sinh h'_u=e^{\xi}\sinh h_u$. The edge $e$ stays flat. If we decrease now the height of $u_3$, then $e$ becomes convex due to Lemma~\ref{dihed}. If we decrease the heights of $u_1$ or $u_2$, then it is clear that the dihedral angle of $e$ decreases. This also applies if we need to decrease heights of several of them. Thus, all cases are considered.
\end{proof}

\begin{lm}
\label{polygons}
Let $P=P(d, V, h)$ be a convex Fuchsian cone-manifold such that all faces are strictly convex hyperbolic  polygons. Then $h$ is in the interior of $H(d, V)$.
\end{lm}

\begin{proof}
The proof is similar to~\cite[Proposition 3.15(2)]{FiIz1}.
Lemma~\ref{finite} shows that there are finitely many triangulations compatible with $P$. We can choose a small enough neighbourhood $N$ of $h$ in $\R^n$ such that for any such triangulation $\mathcal T$ and any $h'\in N$ all strict edges of $P$ remain strict in $P(d, \mathcal T, h')$. Then the flip algorithm from the proof of Lemma~\ref{increase} performs inside each face of $P$ separately. Therefore, it finds in finitely many steps a triangulation $\mathcal T'$ such that $P(d, \mathcal T', h')$ is convex.
\end{proof}

\begin{lm}
\label{homotopnontriv}
If a face $R$ contains a curve, which is not homotopic to a point in $\partial^{\uparrow} P$, then it contains a face angle greater than $\pi$.
\end{lm}

\begin{proof}
We cut $R$ along geodesic segments until it becomes simply connected (the existence of such set of cuts follows, e.g., from Lemma~\ref{iztriang}). Before the last cut we have an annulus $R'$. Let $\psi_1$ and $\psi_2$ be two its boundary components, $\psi$ be the geodesic segment between them that defines the last cut and $R''$ be the resulting polygon. Each cut only decreases the values of the angles, hence it is enough to prove that there is an angle at least $\pi$ before the last cut.

If there are no angles greater than $\pi$ in $R''$, then $R''$ is convex. Develop the prisms having upper boundaries in $R''$ to $\H^3$. Note that as $R''$ is convex, the development is not self-intersecting. Let $\chi_1$, $\chi_2$ be geodesic segments, which are images of $\psi$, $M_{\downarrow}$ be the plane containing the lower boundaries and $M^{\uparrow}$ --- upper. Denote by $A \in M^{\uparrow}$ and $B \in M_{\downarrow}$ their closest points. The distance function from $\chi_1$ and $\chi_2$ to $M_{\downarrow}$ coincides at the respective points. Therefore, the trapezoids containing $\chi_1$ and $\chi_2$ are isometric. There exists a unique orientation-preserving isometry of $\H^3$ mapping one trapezoid to the other. Because $M_{\downarrow}$ is orthogonal to both trapezoids planes, it is preserved by this isometry. Because $\psi$ corresponds to a flat edge of $P$, $M^{\uparrow}$ also must be preserved by this isometry. This means that this is a rotation along the line $AB$ at the angle smaller than $2\pi$, which maps $\chi_1$ to $\chi_2$. Because $R''$ is a convex polygon, when we consider the restriction of this rotation to $M^{\uparrow}$, we get an embedding of $R'$ to a convex hyperbolic cone. Assume that $\psi_2$ is embedded closer to the apex than $\psi_1$. Write the formula for the area of the cone-polygon bounded by $\psi_2$ in terms of the angles. We see that $\psi_2$ contains in $R'$ an angle greater than $\pi$.
\end{proof}

From Lemmas~\ref{negcurv},~\ref{isolvertex},~\ref{polygons} and~\ref{homotopnontriv} we conclude

\begin{crl}
\label{inter}
Let $P=P(d, V, h)$ be a convex Fuchsian cone-manifold with $V=V(d)$ such that for all $v \in V$ we have $\kappa_v(P)\leq 0$. Then the graph of the strict edges of $P$ is connected and all faces are strictly convex hyperbolic polygons. In particular, $h \in \inter(H(d,V))$.
\end{crl}

\begin{rmk}
A crucial difficulty for some of our arguments is given in the fact that $H(d, V)$ is not convex for a cone-metric $d$. However, it becomes convex if we do the coordinate change $h_v \rightarrow \sinh h_v$. But in these coordinates the discrete curvature functional does not have the magical properties that we need. 
\end{rmk}

\subsection{First variation formulas}
\label{varapprsec}

Let $P=P(d, \mathcal T, h)$ be a convex Fuchsian cone-manifold. Recall from Subsection~\ref{disccurv} that its discrete curvature is
$$S(P)=S(d, \mathcal T, h) = -2\vol(P)+ \sum_{v \in V(\mathcal T)}\kappa_v(P) h_v +\sum_{e \in E(\mathcal T)} \theta_e(P) l_e(d).$$
(This makes sense not only for convex Fuchsian cone-manifolds, but we are interested only in convex ones.) Clearly, $S$ is independent from the choice of $\mathcal T$ compatible with $P$ and is a continuous functional on $H(d, V) \subset \R^n$ in the coordinates $h_v$. 

Fix $V \subset S_g$. For a subset $\mathfrak U \subseteq \mathfrak D_c(V)$ define 
$$\mathfrak H(\mathfrak U) = \{ (d, h): d \in \mathfrak U, h \in H(d,V)\}.$$ 

This set is endowed with the topology induced from the product topology. Any point of $\mathfrak H(\mathfrak U)$ can be considered as a convex Fuchsian cone-manifold with the vertex set~$V$. We can define $S$ over $\mathfrak H(\mathfrak U)$. In order to understand its behaviour we need to introduce charts on $\mathfrak H(\mathfrak U)$ related to geodesic triangulations of $(S_g, d, V)$.

First, for a triangulation $\mathcal T$ with the vertex set $V(\mathcal T)=V$ denote by $H(d, V, \mathcal T) \subset  H(d, V)$ the set of all admissible $h \in \R^n$ such that $\mathcal T$ is compatible with $P(d, V, h)$. This defines a subdivision of $H(d, V)$ into subsets corresponding to different triangulations. It is evident that the boundary of $H(d, V, \mathcal T)$ is piecewise analytic.

Next, define the set $\mathfrak H\big(\mathfrak D_c(V), \mathcal T\big)$ as the set of pairs $(d, h)$, where $d \in \mathfrak D_c(V, \mathcal T)$, $h \in H(d, V, \mathcal T)$. The boundary of the set $\mathfrak H\big(\mathfrak D_c(V), \mathcal T\big)$ is piecewise analytic. In the chart $\mathfrak H(\mathfrak D_c(V), \mathcal T)$ on $\mathfrak H(\mathfrak D_c(V))$ the functional $S$ is a continuous function of heights $h_v$ and lengths of the upper edges $l_e(d)$.

\begin{lm}
\label{S1der}
$S$ is continuously differentiable over $\mathfrak H(\mathfrak D_c(V))$ and 
$$\frac{\partial S}{\partial h_v}=\kappa_v,~~~\frac{\partial S}{\partial l_e}=\theta_e.$$
\end{lm}

\begin{rmk}
\label{boundcontdiff}
The set $\mathfrak H(\mathfrak D_c(V))$ has boundary. By continuously differentiable at the boundary points we mean that all partial derivatives admit a continuous extension there.
\end{rmk}

\begin{rmk}
The partial derivatives $\frac{\partial S}{\partial l_e}$ are not defined for flat edges of cone-manifolds corresponding to boundary points of $\mathfrak H(\mathfrak D_c(V), \mathcal T)$. For these points they make sense only as directional derivatives along paths for  which these edges become strict.
\end{rmk}

\begin{proof}

The set $\mathfrak H\big(\mathfrak D_c(V)\big)$ is locally modelled in $\R^{n+N}$, where $n=|V|$ and $N=3(n+2g-2)$ is the number of edges of any triangulation of $S_g$ with $n$ vertices.

Assume first that $(d, h)$ belongs to the interior of some $\mathfrak H\big(\mathfrak D_c(V), \mathcal T\big)$. Then the combinatorics of Fuchsian cone-manifolds in a neighbourhood of $(d, h)$ does not change and all the dihedral angles can be written as sums of dihedral angles in the same prisms. For a single prism $\Pi=A_1A_2A_3B_1B_2B_3$ each dihedral angle is a smooth function of its lengths and by the Schl\"affli formula we get
$$-2d\vol(\Pi)=h_1d\omega_1+h_2d\omega_2+h_3d\omega_3+l_{12}d\phi_{12}+l_{13}d\phi_{13}+l_{23}d\phi_{23}.$$

Summing this for all prisms we obtain
$$-2d\vol(P)=-\sum_{v\in V} h_vd\kappa_v-\sum_{e \in E(\mathcal T} l_e d\theta_e.$$
Therefore,
$$dS(P)=\sum_{v\in V}\kappa_v d h_v + \sum_{e\in E(\mathcal T)} \theta_e d l_e.$$

This shows the first derivative formula. As the angles are smooth functions of the lengths, we see that $S$ is smooth in this case.

Now assume that $(d, h)$ is not in the interior of any $\mathfrak H\big(\mathfrak D_c(V), \mathcal T\big)$. This means that either some vertices of $V$ have the total angle $2\pi$ or there are flat edges in $P$. In any case, all boundary conditions are piecewise analytic. Let $E_s(d, h)$ be the set of strict edges of $P(d, V, h)$. By $\mathcal T_0$ we denote any triangulation with vertices at $V$ realized by $d$ containing all edges from $E_s(d, h)$, so it is compatible with $P(d,V, h)$. (By Lemma~\ref{iztriang} such a triangulation exists.) It provides a chart over $(d, h)$ in $\R^{n+N}$. Now let $\bar x \in \R^{n+N}$ be the point corresponding to $(d, h)$ and $\bar \xi$ be a vector such that $\bar x + t\bar \xi \in \inter\Big(\mathfrak H\big(\mathfrak D_c(V)\big)\Big)$ for all sufficiently small $t>0$. Then, as all boundary conditions are piecewise analytic and the number of triangulations compatible with $P(d,V,h)$ is finite, there exists $\mathcal T$ such that $\bar x + t\bar \xi \in \mathfrak H(\mathfrak D_c(V), \mathcal T)$ for all sufficiently small $t$. The triangulation $T$ might be different from $T_0$, but it contains all edges from $E_s(d, h)$. 

From the previous argument, we get
$$\left.\frac{\partial S}{\partial \bar \xi}\right|_{\bar x}=\sum_{v\in V}\kappa_v \xi_v + \sum_{e\in E_s(d, h)} \theta_e \xi_e,$$
where $\xi_v$, $\xi_e$ are the corresponding coordinates of $\xi$. The second sum is only over strict edges as flat edges give zero directional derivatives at $\bar x$. Thus, all partial derivatives exist and are continuous in a neighbourhood of $(d, h)$. We see that $S$ is continu ously differentiable and its first derivatives are the respective curvatures.
\end{proof}

\subsection{Second variation formulas}
\label{varappr2sec}

Now we are going to investigate the second derivatives of $S$ over $H(d, V)$ (we do not consider the second derivatives with respect to the upper boundary). To this purpose we need a new notation. 

By $\vec E(\mathcal T)$ we denote the set of \emph{oriented} edges of $\mathcal T$ in the sense that each edge $e \in E(\mathcal T)$ gives rise to two oriented edges with respect to two different possible orientations. By $\vec E_v (\mathcal T)$ we denote the set of oriented edges starting at $v$ and by $\vec E_{vu}(\mathcal T)$ we denote the set of oriented edges starting at $v$ and ending in $u$. In particular, $\vec E_{vv}(\mathcal T)$ is the set of loops from $v$ to itself and each loop is counted twice. For $P$ compatible with $\mathcal T$ and $\vec e \in \vec E(\mathcal T)$ by $\phi^+_{\vec e}$ and $\phi^-_{\vec e}$ we denote the dihedral angles of $\vec e$ in the right and the left prisms incident to $\vec e$ respectively. By $\alpha_{\vec e}$ we denote the angle of the trapezoid containing $\vec e$ at the vertex, where $\vec e$ starts. By $l_{\vec e}$ and $a_{\vec e}$ we continue to denote the lengths of $\vec e$ in the upper and the lower boundaries of $P$ respectively.

\begin{lm}
\label{S2der}
Let $d \in \mathfrak D_c(V)$. Then $S$ is twice continuously differentiable over $H(d, V)$ and
$$v \neq u:~~~\frac{\partial^2 S}{\partial h_v\partial h_u}=\frac{\partial \kappa_u}{\partial h_v}=\frac{\partial \kappa_v}{\partial h_u}=\sum_{\vec e \in \vec E_{vu}(\mathcal T)}\frac{\cot \phi^+_{\vec e} + \cot \phi^-_{\vec e}}{\sin \alpha_{\vec e}}\cdot \frac{1}{\cosh h_v \sinh a_{\vec e}}\geq 0;$$
$$v=u:~~~\frac{\partial ^2 S}{\partial h^2_v}=\frac {\partial \kappa_v}{\partial h_v} = -\sum_{\vec e \in \vec E_v(\mathcal T)} \frac{\cot \phi^+_{\vec e} + \cot \phi^-_{\vec e}}{\sin \alpha_{\vec e}}\cdot \frac{\coth a_{\vec e}}{\cosh h_v }+$$
$$+\sum_{\vec e \in \vec E_{vv}(\mathcal T)}\frac{\cot \phi^+_{\vec e} + \cot \phi^-_{\vec e}}{\sin \alpha_{\vec e}}\cdot \frac{1}{\cosh h_v \sinh a_{\vec e}} \leq$$
$$\leq -\sum_{u \neq v}\frac{\partial \kappa_u}{\partial h_v}-\sum_{\vec e \in \vec E_{vv}(\mathcal T)}\frac{\cot \phi^+_{\vec e} + \cot \phi^-_{\vec e}}{\sin \alpha_{\vec e}}\cdot \frac{\cosh a_{\vec e}-1}{\cosh h_v \sinh a_{\vec e}} \leq 0.$$
\end{lm}



\begin{proof}
Similarly to the proof of Lemma~\ref{S1der}, it is enough to prove the claim for $h$ in the interior of a set $H(d, V, \mathcal T)$ of heights such that $\mathcal T$ is compatible with the respective cone-manifold. This is reduced to computations for a single prism $\Pi=A_1A_2A_3B_1B_2B_3$.

Consider a spherical link at the vertex $A$. From the cosine rule we get
$$\cos \omega_1 = \frac{\cos \lambda_1 - \cos \alpha_{12}\cos\alpha_{13}}{\sin\alpha_{12}\sin\alpha_{13}}.$$ 

Differentiating the spherical cosine rule, we obtain
$$\frac{\partial \omega_1}{\partial \alpha_{12}}=-\frac{\cot \phi_{12}}{\sin \alpha_{12}},~~~\frac{\partial \omega_1}{\partial \alpha_{13}}=-\frac{\cot \phi_{13}}{\sin \alpha_{13}}.$$

Differentiating the cosine rule for trapezoids from Lemma~\ref{coslaw}, we see
$$\frac{\partial \alpha_{12}}{\partial h_1}=-\frac{\coth a_{12}}{\cosh h_1},~~~\frac{\partial \alpha_{13}}{\partial h_1}=-\frac{\coth a_{13}}{\cosh h_1},~~~\frac{\partial \alpha_{12}}{\partial h_2}=\frac{1}{\cosh h_1 \sinh a_{12}}.$$

So we have
\begin{equation}
\label{der11}
\frac{\partial \omega_1}{\partial h_1}=\frac{\partial \omega_1}{\partial \alpha_{12}}\cdot\frac{\partial \alpha_{12}}{\partial h_1}+\frac{\partial \omega_1}{\partial \alpha_{13}}\cdot\frac{\partial \alpha_{13}}{\partial h_1}=\frac{\cot \phi_{12}\coth a_{12}}{\sin \alpha_{12}\cosh h_1}+\frac{\cot \phi_{13}\coth a_{13}}{\sin \alpha_{13}\cosh h_1},
\end{equation}
\begin{equation}
\label{der12}
\frac{\partial \omega_1}{\partial h_2}=\frac{\partial \omega_1}{\partial \alpha_{12}}\cdot\frac{\partial \alpha_{12}}{\partial h_2}=-\frac{\cot \phi_{12}}{\sin \alpha_{12}\cosh h_1 \sinh a_{12}}.
\end{equation}

Consider $v \neq u$ and sum up a formula of type (\ref{der12}) for each prism incident to both $v$ and $u$. Then we get the desired formulas for $\frac{\partial \kappa_v}{\partial h_u}=-\frac{\partial \omega_v}{\partial h_u}$. For $v=u$ we sum up the formulas of type (\ref{der11}) in all prisms incident to $v$ plus a formula of type (\ref{der12}) for each prism incident to an oriented loop and get $\frac{\partial \kappa_v}{\partial h_v}$.
This finishes the proof.
\end{proof}


\begin{lm}
\label{concave} 
Fix $d \in \mathfrak D_c(V)$ and identify $H(d, V)$ with a subset of $\R^n$. The Hessian of $S$ over $H(d, V)$ is non-positive definite. Moreover, its kernel is spanned by vectors $\bar x^u=(x^u_v)_{v \in V}$ defined by $x^u_v=1$ if $v = u$ and 0 otherwise, where $u$ is an isolated vertex of $P$.

\end{lm}

\begin{proof}
We have
$$\bar x^T{\rm Hess}(S)\bar x=\sum_{u,v \in V}\frac{\partial \kappa_u}{\partial h_v} x_ux_v=-\sum_{\substack{u,v \in V\\ u\neq v}}\frac{\partial \kappa_u}{\partial h_v}(x_u-x_v)^2+\sum_{v\in V}x^2_v\sum_{u \in V}\frac{\partial \kappa_u}{\partial h_v}$$

From Lemma~\ref{S2der} we see 
$$\sum_{u \in V}\frac{\partial \kappa_u}{\partial h_v}= \frac{\partial \kappa_v}{\partial h_v}+\sum_{\substack{u \in V\\ u\neq v}}\frac{\partial \kappa_u}{\partial h_v}=-\sum_{\vec e \in \vec E_{v}(\mathcal T)}\frac{\cot \phi^+_{\vec e}+ \cot \phi^-_{\vec e}}{\sin \alpha_{\vec e}}\cdot \frac{\cosh a_{\vec e}-1}{\cosh h_v \sinh a_{\vec e}}\leq 0.$$

Hence, the Hessian is non-positive definite. Also $\frac{\partial \kappa_u}{\partial h_v} = 0$  if and only if there are no strict edges in $P$ between $u$ and $v$. Thus, it is easy to see that $\bar x^T{\rm Hess}(S)\bar x=0$ if and only if $\bar x$ has non-zero coordinates only on isolated vertices of $P$.
\end{proof}

Thus, $S$ is concave over $H(d, V)$.

\begin{rmk}
Although the Hessian of $S$ can be degenerate and, moreover, critical points might be non-isolated, one can show that it can not be degenerate along any segment in $H(d, V)$. Thus, $S$ is strictly concave over $H(d, V)$. In what follows we will not need it, so we do not prove this. This would lead to a proof of local rigidity of convex Fuchsian cone-manifolds.
\end{rmk}

Together with Lemma~\ref{isolvertex} we get

\begin{crl}
\label{hessnond}
Let $P=P(d, V, h)$ be a convex Fuchsian cone-manifold with $V=V(d)$ such that $\kappa_v(P)\leq 0$ for all $v \in V$. Then ${\rm Hess}(S)$ is non-degenerate at $P$.
\end{crl}

\subsection{Behaviour of $S$ at infinity}
\label{maxsec}

The set $H(d, V)$ is non-compact and we need to understand the behaviour of $S$ outside a compact set. In particular, we will prove the maximization principle Lemma~\ref{maxprinciple}.

\begin{lm}
\label{compactabove}
Let $d \in \mathfrak D_c(V)$. For each $K \in \R$ there exists $M=M(K, d, V)>0$ such that if $h \in H(d, V)$ and for some vertex $v \in V$ one has $h_v \geq M$, then we get $S(d, V, h) \leq K$.
\end{lm}

\begin{proof}

\begin{prop}
For each $\e >0$ there exists $\delta > 0$ such that if the side lengths of a hyperbolic triangle are at most $\delta$, then the sum of its angles is at least $\pi-\e$.
\end{prop}

This is because the area of a hyperbolic triangle is equal to the sum of its angles minus $\pi$. By computing the Euler characteristic we get

\begin{prop}
For each $\e>0$ there exists $\delta >0$ such that if for a Fuchsian cone-manifold $P=P(d, \mathcal T, h)$ for each $e \in E(\mathcal T)$ the length $a_e$ of the projection of $e$ to $\partial_{\downarrow} P$ is at most $\delta$, then $$\sum_{v \in V} \kappa_v(P) \leq 2\pi(2-2g) + \e.$$
\end{prop}

The number of triangulations compatible with some $P(d, V, h)$ is finite for fixed $d$, $V$ due to Lemma~\ref{finite}. Due to Lemma~\ref{coslaw}, if a trapezoid has fixed upper edge and growing heights, then its lower edge becomes smaller. Thus, we conclude

\begin{prop}
For each $\delta>0$ there exists $M_0=M_0(\delta, d, V)>0$ such that if $h\in H(d, V)$ and for each $v\in V$ we have $h_v \geq M_0$, then for each edge $e$ of $P=P(d,V,h)$ we have $a_e \leq \delta$.
\end{prop}

From the finiteness of the set of compatible triangulations we also get

\begin{prop}
There exists $C=C(d, V)>0$ such that for any convex Fuchsian cone-manifold $P=P(d, V, h)$ compatible with a triangulation $\mathcal T$ we have $$\Big|\sum_{e \in E(\mathcal T)} \theta_e(P) l_e(d)\Big|\leq C.$$
\end{prop}

Finally, the triangle inequality for trapezoids (Lemma~\ref{htriangleinequality}) gives 

\begin{prop}
\label{allincrease}
For every $M_0>0$ there exists $M=M(M_0, d)>0$ such that if $h\in H(d, V)$ and for some $v \in V$ we have $h_v \geq M$, then for each $u \in V$ we get $h_u \geq M_0$.
\end{prop}

Combining all these claims together we see that for every $\e>0$ and $M_0>0$ we can choose $M>0$ such that if $h \in H(d, V)$ and $h_v \geq M$ for some $v \in V$, then we get
$$S(d, V, h) \leq -2\vol(P) + M_0(2\pi(2-2g)+\e) + C.$$
As $2-2g<0$ and $M_0$ is arbitrary, this finishes the proof.
\end{proof}

Define 
$$H_S(d, V, m, K)=\{h \in H(d, V): \min_{v \in V}h_v \geq m, S(d, V, h) \geq K\}.$$
Lemma~\ref{H_d comp} and Lemma~\ref{compactabove} together give

\begin{crl}
\label{H_S comp}
For every $d \in \mathfrak D_c(V)$, $m \in \R_{>0}$ and $K \in \R$ the set $H_S(d, V, m, K)$ is compact.
\end{crl}

Now we turn to cone-manifolds with small heights.

\begin{lm}
\label{alldecrease}
\label{nottozero}
If $h^n \in H(d, V)$ is a sequence such that for some $w \in V$ we have $h^n_w \rightarrow 0$, then for all $v \in V$ we get $h^n_v \rightarrow 0$ and $\kappa^n_v \rightarrow \nu_v(d) \geq 0$, where $\kappa^n_v$ are particle curvatures in $P^n=P(d, V, h^n)$. Also for every edge $e$ of all $P^n$ we have $\phi^n_e \rightarrow \pi$, where $\phi^n$ is the dihedral angle of $e$ in $P^n$.
\end{lm}

\begin{proof}
Lemma~\ref{finite} implies that there are finitely many triangulations $\mathcal T$ compatible with some $P^n$. Hence, up to taking a subsequence we may suppose that the same triangulation $\mathcal T$ is compatible with all $P^n$.

Take $v$ connected with $w$ by an edge $e$ of $\mathcal T$. If for some $m>0$ we have $h^n_v \geq m$, then the trapezoid corresponding to $e$ becomes not ultraparallel for large enough $n$. This contradicts Lemma~\ref{ultrapar2}. Hence, $h^n_v \rightarrow 0$. By the connectivity of the edge graph of $\mathcal T$, this extends to all vertices.

Take an oriented edge $\vec e \in \vec E_v(\mathcal T)$ emanating from $v \in V$. Develop the trapezoid containing $\vec e$ from $P^n$ to $\H^2$ as $A_1A_2B_1B_2$. Let $A$ and $B$ be the closest points on the lines containing the upper and the lower boundary respectively. We have $A_1B_1=h^n_v$. Denote $AB$ by $h^n_{\vec e}$. Clearly, $h^n_{\vec e} \leq h^n_v$, hence $h^n_{\vec e} \rightarrow 0$. From Corollary~\ref{rightangl} we get
$$\sin \alpha_{\vec e}=\frac{\cosh h^n_{\vec e}}{\cosh h^n_v} \rightarrow 1.$$

Hence, $\alpha_{\vec e} \rightarrow \pi/2$. Now we take a triangle $T$ of $\mathcal T$ and the prism from $P^n$ containing $T$. Let $T$ be incident to $v \in V$, $\lambda^n_{v, T}$ be the angle of $T$ at $v$, $\omega^n_{v, T}$ be the dihedral angle of the respective lateral edge and $\phi^n_{e, T}$ be the dihedral angle of $e$ in the prism. Consider the spherical link of $v$. As all $\alpha_{\vec e} \rightarrow \pi/2$, we obtain $\omega^n_{v, T} \rightarrow \lambda^n_{v, T}$ and $\phi^n_{e, T} \rightarrow \pi/2$. Thereby, $\omega_v(P^n) \rightarrow \lambda_v(d)$, $\kappa_v(P^n) \rightarrow \nu_v(d) \geq 0$ and $\phi_e(P^n) \rightarrow \pi$.
\end{proof}

Thus, when some heights tend to zero, the upper boundary falls down to the lower boundary as we could expect. Lemma~\ref{alldecrease} together with Lemma~\ref{finite} imply

\begin{crl}
\label{zerocont}
Let $o \in \R^n$ be the origin. The functional $S$ with its partial derivatives can be continuously extended to $o$ over $H(d, V)$ by putting $S(o)=0$, $\frac{\partial S}{\partial h_v}=\nu_v$.
\end{crl}

Now we are ready to prove the maximization principle from Subsection~\ref{disccurv}:

\begin{replm}{maxprinciple}
Let $d \in \mathfrak D_{sc}(V)$, $P^1=P(d, V, h^1)$ be the convex polyhedral Fuchsian manifold realizing $d$ and $P^2=P(d, V, h^2)$ be a convex Fuchsian cone-manifold distinct from $P^1$. Then $S(P^1) > S(P^2)$.
\end{replm}

\begin{proof}
We want to show that there exists at least one global maximum of $S$ over $H(d, V)$ and there are $0<m<M$ such that all global maxima are contained in $H^o(d, V, m, M)$, where 
$$H^o(d, V, m, M):=\{h \in H(d, V): m<h_v<M {\rm~for~all~}v\in V\}.$$
Indeed, an upper bound follows from Lemma~\ref{compactabove}. If there is no lower bound, then we have a sequence of global maxima converging to $o$. Lemma~\ref{alldecrease} implies that all $\kappa_v$ become positive at some point $P$. Lemma~\ref{movingup} shows that we can increase all heights of $P$, which increases $S$, so $P$ was not a maximum point. This means that we can find $0<m<M$ such that $H^o(d, V, m, M)$ contains all global maxima and at least one global maximum exists.

Let $P=P(d, V, h)$ be a convex Fuchsian cone-manifold that is a local maximum (not necessarily strict) of $S$ over $H^o(d, V, m, M)$. We claim that it is actually a convex polyhedral Fuchsian manifold. Indeed, if there is $v \in V$ such that $\kappa_v(P) < 0$, then due to Corollary~\ref{decrease} we can decrease $h_v$ while staying in $H^o(d, V, m, M)$. Due to Lemma~\ref{S1der} this increases $S$. If there is a vertex $v \in V$ such that $\kappa_v(P)>0$, then Lemma~\ref{increase} shows that we can increase simultaneously the heights $h_v$ of all such vertices and stay in $H^o(d, V, m, M)$. Again, this increases $S$.

Hence, $P$ is a convex polyhedral Fuchsian manifold. Theorem~\ref{Fillastre} and Remark~\ref{markeduniq} show that $P=P^1$ is unique.  
\end{proof}

Here we appealed to Theorem~\ref{Fillastre} because the set $H^o(d, V, m, M)$ may be non-convex. However, one can see that it is contractible. With the help of elementary Morse theory one can obtain an alternative proof of uniqueness and, thus, reprove Theorem~\ref{Fillastre}.

\section{Stability lemmas}
\label{stability}

\subsection{Proof of Main Lemma I}

\subsubsection{Dual area}

Let $P=P(d, V, h)$ be a convex Fuchsian cone-manifold. In Section~\ref{definfcm} we defined the spherical link of a vertex $v \in V$. Intrinsically it is a disk with a spherical cone-metric that has (at most) one conical point in the interior, piecewise geodesic boundary and the angles of all kink points of the boundary are less than $\pi$. Geometry of spherical links is highly connected with the behaviour of the discrete curvature. In order to get the quantitative bound of Main Lemma I, first we prove an inequality concerning spherical links, which will be used further in estimating the derivatives of the discrete curvature.

\begin{lm}
\label{spherarea}
Let $R$ be a disk with a spherical cone-metric that has a single conical point $Z$ in the interior of total angle $\omega$ and of curvature $\kappa:=2\pi-\omega$, and has a piecewise geodesic boundary with kink points $A_1, \ldots, A_m$ denoted in the cyclic order. By $\alpha_i$ denote the length $ZA_i$, by $\phi^+_i$ and $\phi^-_i$ denote the angles $ZA_iA_{i+1}$ and $ZA_iA_{i-1}$ respectively. Assume that $\kappa \leq 0$, $m\geq 1$, for each $i=1\ldots m$ we have $\phi^-_i+\phi^+_i < \pi$ and the perimeter $\lambda$ of $R$ is (strictly) less than $ 2\pi$. Define $\nu:=2\pi-\lambda$. Then 
\begin{equation}
\label{orthoineq}
\sum_{i=1}^m\frac{\cot \phi^+_i+\cot \phi^-_i}{2\sin\alpha_i}\cdot \cot \alpha_i \geq \nu  - \kappa>0.
\end{equation}
\end{lm}

\begin{proof}
We start with the case $\kappa=0$, so $\omega=2\pi$. Then there is no conical point, but we assume that $Z$ is an arbitrary marked point in the interior of $R$. 

In this case we develop the polygon $R$ to the unit sphere so that $Z$ coincides with the south pole. Then $R$ determines a solid convex cone with the apex at the origin. Consider its polar cone. It defines the polar spherical polygon $\tilde R$, which belongs to the open hemisphere centered at the north pole $\tilde Z$. Project it from the center of sphere to the tangent plane at $\tilde Z$. We get a Euclidean polygon $\bar R$. If all $\alpha_i < \pi$, then $\tilde Z$ is in the interior of both $\tilde R$, $\overline R$. If there are some $\alpha_i > \pi$, then $\tilde Z$ is outside from both of them.

Regardless the location of $\tilde Z$, one can compute that the left-hand side of (\ref{orthoineq}) is exactly the area of $\bar R$ and the right-hand side is the area of $\tilde R$ (in particular this will follow from the general discussion below). Hence, (\ref{orthoineq}) is true because the central projection to a tangent plane increases areas.

In the general case $\kappa \neq 0$ one can still consider similar polar polygons (which contain cone-singularities) as long as all $\alpha_i <\pi$. Otherwise the polar interpretation breaks down and $\bar R$, $\tilde R$ can be defined only in a virtual sense described below.

We call a spherical triangle $O=ZAB$ a \emph{(spherical) orthoscheme} if $\angle B=\pi/2$, the sides $ZA$ and $ZB$ are either both at most $\pi/2$ or both at least $\pi/2$ and $O$ is equipped with plus or minus sign. If $ZA$ is equal to $\pi/2$, then the orthoscheme is called \emph{singular}. If $ZB=\pi/2$, then also automatically $ZA=\pi/2$ and we call $O$ \emph{double-singular}. From simple trigonometry we conclude that if $O$ is non-singular, then  
$$\sgn(\pi/2-\angle A)=\sgn(\pi/2-ZA)=\sgn(\pi/2-ZB)\neq 0,$$  
$$AB< \pi/2,~~~\angle Z < \pi/2.$$

Define its \emph{dual} orthoscheme $\tilde O=\tilde Z\tilde A\tilde B$ as a spherical triangle with $\angle \tilde A = \pi/2$, $\angle \tilde Z = \angle Z$ and $\tilde Z \tilde A = |\pi/2 - ZA|$. One can conclude $$\tilde Z\tilde B = |\pi/2-ZB|,~~~\tilde A \tilde B = |\pi/2-\angle A|,~~~\angle \tilde B = \pi/2-AB.$$
See Figure~\ref{orth}.
For a non-singular orthoscheme $O$ the sign of $\tilde O$ is defined as the sign of $O$ if $ZA, ZB < \pi/2$ or as the converse sign otherwise. For a singular orthoscheme, the dual orthoscheme is degenerate and we assign zero as its sign. 

\begin{figure}
\begin{center}
\includegraphics[scale=0.20]{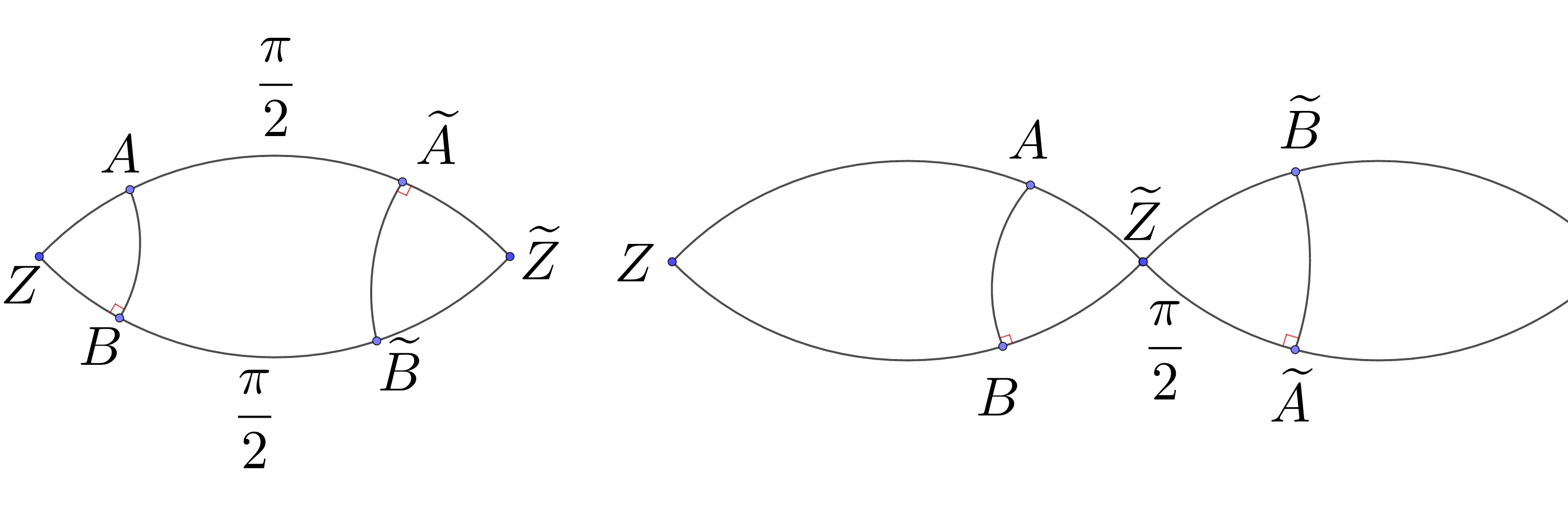}
\caption{Orthoscemes and their duals in cases $ZA<\pi/2$ and $ZA>\pi/2$.}
\label{orth}
\end{center}
\end{figure}

This duality is not involutive. By definition, only orthoschemes with both $ZA, ZB < \pi/2$ can be obtained as duals and each such non-degenerate orthoscheme is obtained as the dual of exactly two orthoschemes.

Let $O$ be an orthoscheme. Consider the unit sphere and place the dual orthoscheme $\tilde O$ such that $\tilde Z$ becomes the north pole. Then $\tilde O$ belongs to the open northern hemisphere. Project it from the center of the sphere to the tangent plane at $\tilde Z$. The image is called \emph{the Euclidean dual orthoscheme} to $O$ and is denoted by $\bar O=\bar Z\bar A \bar B$ (where $\tilde Z=\bar Z$). One can compute the signed areas, i.e., the areas taken with the sign of the dual orthoshemes:
$$\sarea~\tilde O=\angle Z - AB,~~~\sarea~\bar O=\frac{\cot \angle A}{2\sin ZA}\cot ZA.$$

Let $ZA_1A_2$ be a spherical triangle with both $ZA_1$, $ZA_2$ at least $\pi/2$ or both at most $\pi/2$, but at least one (say, $ZA_1$) not equal to $\pi/2$. Put it on the unit sphere and consider the geodesic great circle containing segment $A_1A_2$. Then there are two perpendiculars from $Z$ to it. Let $ZB$ be the one such that $\sgn(\pi/2-ZB)=\sgn(\pi/2-ZA_1)$. The characteristic function of triangle $ZA_1A_2$ is the signed sum of the characteristic function of triangles $ZBA_1$, $ZBA_2$. Consider the latter triangles as orthoschemes endowed with signs coming from this sum. We say that $ZA_1A_2$ is \emph{canonically decomposed} into the orthoschemes $ZBA_1$, $ZBA_2$. In the case of $ZA_1=ZA_2=\pi/2$ the canonical decomposition of $ZA_1A_2$ into orthoschemes is $ZA_1A_2$ itself taken with the plus sign (then it is a double-singular orthoscheme).

Now we turn to the spherical cone-polygon $R$. Since now we call the kink points \emph{vertices} of $R$. Cut $R$ into triangles $ZA_iA_{i+1}$. If there are $i$ such that $ZA_i<\pi/2$ and $ZA_{i+1}>\pi/2$ (or the converse), then there exists a unique $A \in A_i A_{i+1}$ such that $ZA=\pi/2$ and we add $A$ to the set of vertices. Hence, we may assume that for all $i$ either both $ZA_i$, $ZA_{i+1}$ are at least $\pi/2$ or both are at most $\pi/2$. We decompose each triangle canonically into orthoschemes and call it \emph{the canonical decomposition} of $R$. By $\mathcal O$ denote the set of orthoschemes of the canonical decomposition. The total angle $\omega$ is the signed sum of angles at $Z$ in all orthoschemes, similarly is the perimeter $\lambda$. Thereby,
$$\sum_{O\in \mathcal O} \sarea ~\tilde O= \omega-\lambda=\nu-\kappa$$
and the inequality~(\ref{orthoineq}) is equivalent to
$$\sum_{O \in \mathcal O}\sarea ~\bar O \geq \sum_{O\in \mathcal O} \sarea~\tilde O.$$
Define $\phi_i := \phi_i^+ + \phi_i^-$. We have all $\phi_i < \pi$, except possibly some $i$ with $\alpha_i = \pi/2$ (which correspond to recently added vertices). 

Take $A_i$ such that $\alpha_i \neq \pi/2$ and consider two orthoschemes incident to $ZA_i$. Denote them by $O^+_i$ and $O^-_i$ respectively. Note that due to convexity at least one of them is positive. Let $B_i$ be the base of perpendicular from $Z$ to the geodesic $A_iA_{i+1}$ in the canonical decomposition. Develop the dual orthoschemes $\tilde O^+_i$ and $\tilde O^-_i$ to the sphere so that they have the segment $\tilde Z \tilde A_i$ in common, with opposite orientations if they have the same sign and with the same orientation otherwise (see Figure~\ref{orthtr}). We say that the triangle $\tilde Z \tilde B_{i-1} \tilde B_i$ obtained this way is the sum $\tilde O^+_i + \tilde O^-_i$ of the dual orthoschemes. We equip it with the sign determined from the sum of the characteristic functions of $\tilde O^+_i$ and $\tilde O^-_i$ taken with the respective signs.

\begin{figure}
\begin{center}
\includegraphics[scale=0.3]{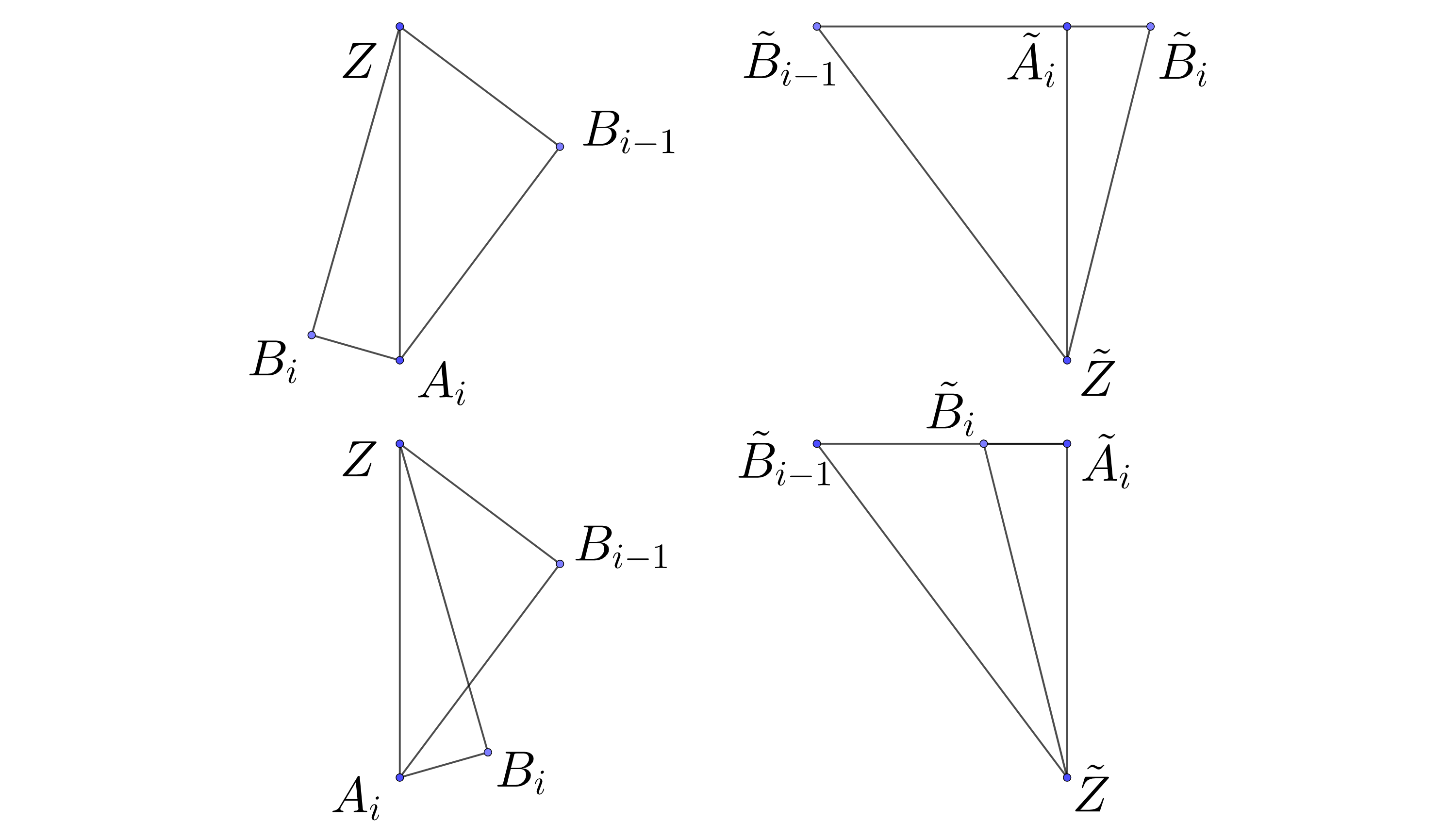}
\caption{The transformation of a pair of adjacent orthoschemes to a triangle.}
\label{orthtr}
\end{center}
\end{figure}

One can see that if $\alpha_i < \pi/2$, then the sum of the dual orthoschemes $\tilde O^+_i + \tilde O^-_i$ is a positive spherical triangle $\tilde Z \tilde B_{i-1} \tilde B_i$ with $\tilde B_{i-1}\tilde B_i=\pi-\phi_i$. We call it a positive pair. We obtain that if $\alpha_i < \pi/2$, then
$$\sarea~\bar O^+_i+\sarea~\bar O^-_i=\sarea\left(\bar O^+_i+\bar O^-_i\right) \geq$$
$$\geq \sarea\left(\tilde O^+_i+\tilde O^-_i\right) = \sarea~\tilde O^+_i+\sarea~\tilde O^-_i$$
because the central projection increases areas.

Similarly, if $\alpha_i=\pi/2$, then 
$$\sarea~\bar O^+_i+\sarea~\bar O^-_i=\sarea~\tilde O^+_i+\sarea~\tilde O^-_i=0.$$

Therefore, if all $\alpha_i \leq \pi/2$, then the proof is finished.

It remains to consider the case, when some $\alpha_i > \pi/2$, which is much more subtle. 

First, we investigate the boundary structure. Let $S^1$ be an oriented circle of length $\omega$. It naturally parametrizes the set of geodesics emanating from $Z$. We say that $0 \in S^1$ corresponds to the segment $ZA_1$ and positive direction corresponds to the order $A_1, A_2, A_3, \ldots$. If $ZA$ is the segment corresponding to $t \in S^1$, where $A$ belongs to the boundary of $R$, then $t$ is the angle of the sector between $ZA_1$, $ZA$ in the positive direction. 
Let $\alpha(t)$ be the length of $ZA$. It defines a continuous function $\alpha$ on $S^1$. We now consider the case when for some $t \in S^1$ we have $\alpha(t) > \pi/2$. 

Lemma~\ref{negcurv} implies that no geodesic segment in the boundary of $R$ has length at least $\pi$. This clearly means that no geodesic segment in $R$ has length at least $\pi$. Otherwise we can extend this segment until it intersects the boundary and cut off the part of the polygon that does not contain $Z$. Then we obtain the cone-polygon with non-positive curvature, small perimeter and geodesic segment in the boundary of length at least $\pi$, which can not happen. 

In particular this means that for each $t \in S^1$ we get $\alpha(t)+\alpha(t+\pi) < \pi$. Then there exists $t \in S^1$ such that $\alpha(t)<\pi/2$. We want to show

\begin{prop}
\label{pbstructure}
There exist exactly two $t_1, t_2 \in S^1$ such that $\alpha(t_1)=\alpha(t_2)=\pi/2$. 
\end{prop} 

With the preliminary discussion this means that $S^1$ is subdivided by $t_1$, $t_2$ into two complementary open intervals $I, J$ such that $\alpha|_I > \pi/2$ and $\alpha|_J < \pi/2$.

Indeed, consider $t$ such that $\alpha(t) > \pi/2$. Let $I$ be the maximal interval containing $t$ for which $\alpha|_I > \pi/2$. We saw that $I \neq S^1$. Let $t_1, t_2$ be the left and the right endpoints of $I$, so $\alpha(t_1)=\alpha(t_2)=\pi/2$. Let $A_{i_2}$ be the vertex of $R$ corresponding to $t_2$. Then $\phi^-_{i_2}>\pi/2$. Note that the strict inequality here is because $\alpha(t)>\pi/2=\alpha(t_2)$ for $t$ close to $t_2$ from the left. Consider the interval $I^+=(t_2; t_2+\pi]$. Develop the part of $R$ corresponding to $I^+$ to the unit sphere. From convexity, $\phi^-_{i_2}>\pi/2$ and $\alpha(t_2)=\pi/2$ one can see that $\alpha|_{I^+}<\pi/2$.

Assume that there are more points such that $\alpha(t) = \pi/2$. Take the closest of them to $t_2$ in positive direction and denote it by $t_1'$. We showed that $t'_1-t_2 > \pi$. Together with $\alpha(t_2)=\alpha(t'_1)=\pi/2$ this implies that the part of perimeter of $R$ between the segments from $Z$ corresponding to $t_2$ and $t'_1$ is greater than $\pi$.  Similarly, we consider $t_2'$ that is the closest to $t_1$ in the negative direction such that $\alpha(t_2')=\pi/2$. We get that the segments corresponding to $t_2'$ and $t_1$ also cut off the part of perimeter of length greater than $\pi$. The intervals $[t_2; t_1')$ and $(t_2'; t_1]$ do not intersect. Hence, the perimeter of $R$ is greater than $2\pi$, which is a contradiction. Claim~\ref{pbstructure} is proven.
\vskip+0.2cm

By $A_{i_1}, A_{i_2}$ we denote the points at the boundary of $R$ corresponding to $t_1$, $t_2$. Due to our previous agreement, they are vertices even if their angles are $\pi$. By $\mathcal A_I, \mathcal A_J$ we denote the sets of vertices corresponding to the segments $I, J \subset S^1$.

We are going to discuss now what happens with dual orthoschemes. We focus on Euclidean duals, but we note that exactly the same happens with spherical duals. We think that our discussion is easier to imagine on the Euclidean plane.

Take $A_i \in \mathcal A_I$, we get $\alpha_i> \pi/2$. For two orthoschemes incident to $ZA_i$ the sum of their Euclidean duals $\bar O_i^++\bar O_i^-$ is a negative triangle $\bar Z \bar B_{i-1} \bar B_i$ (recall Figure~\ref{orthtr}). We call it a negative pair. 

Consider all $A_i \in \mathcal A_I$ and start developing all dual triangles $\bar Z \bar B_{i-1} \bar B_i$ one by one to the Euclidean plane. The sides $\bar B_{i-1} \bar B_i$ constitute a polygonal curve $\bar B_{i_1}\bar B_{i_1+1}\ldots\bar B_{i_2-1}$, which we denote by $\bar B_{\smile}$. See Figure~\ref{Pic8-1}. Define 
$$\omega_{\smile}:=\sum_{i_1}^{i_2-2} \angle \bar B_i\bar Z \bar B_{i+1}.$$
The angle $\bar B_{i-1}\bar B_i \bar B_{i+1}$ as seen from $\bar Z$ is
$$\angle \bar B_{i-1} \bar B_i\bar Z + \angle  \bar Z\bar B_i \bar B_{i+1} = \pi + \lambda_i>\pi,$$
where $\lambda_i = A_iA_{i+1} < \pi$. This means that the polygonal curve $\bar B_{\smile}$ is convex as seen from $\bar Z$. From this one can show that $\omega_{\smile} < \pi$ and $\bar B_{\smile}$ is not self-intersecting. 

\begin{figure}
\begin{center}
\includegraphics[scale=0.85]{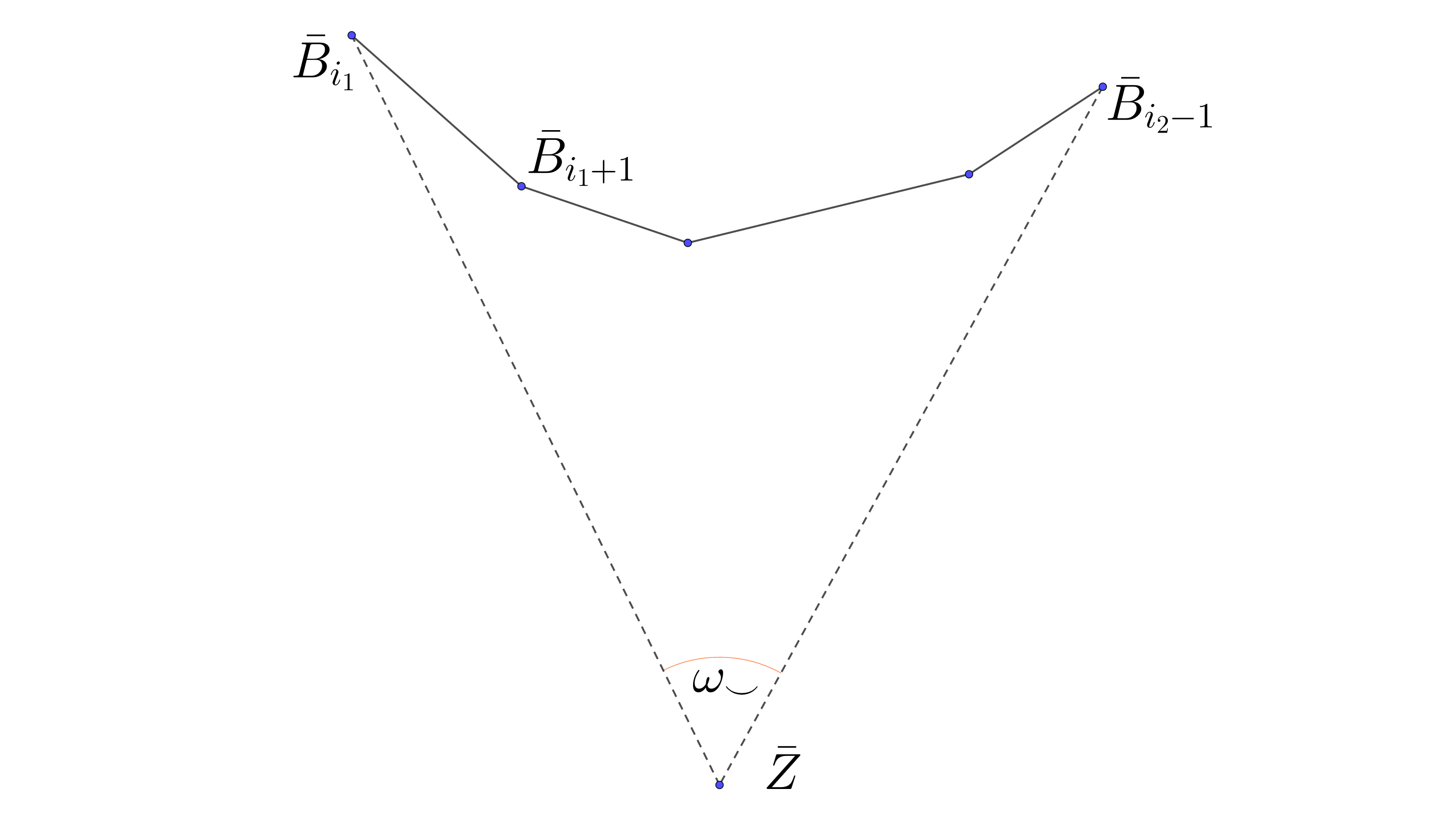}
\caption{Polygonal curve $\bar B_{\smile}$.}
\label{Pic8-1}
\end{center}
\end{figure}

We remark that $\omega_{\smile} \neq \angle A_{i_1} Z A_{i_2}$ or $\angle A_{i_1+1}Z A_{i_2-1}$. The angle $\omega_{\smile}$ can be seen as follows. Take all orthoschemes for triangles $ZA_{i_1} A_{i_1+1}, \ldots ZA_{i_2-1}A_{i_2}$ except the first and the last ones, which are singular. The (signed) sum of their angles at $Z$ is negative and is equal to $-\omega_{\smile}$.

Our plan is to cover the negative polygon $\bar Z \bar B_{i_1}\ldots \bar B_{i_2-1}$ by positive triangles constructed from remaining orthoschemes. This will finish the proof: we obtain (\ref{orthoineq}) for the sum of all negative pairs with the positive pairs, which we used to cover. All remained orthoschemes are split into positive pairs, for which the inequality also holds. If this is not possible, then we will prove $\lambda > 2\pi$, which is a contradiction. 

If we do the same construction as above, but with all $A_i \in \mathcal A_J$, then we obtain a concave polygonal curve $\bar B_{\frown} = \bar B_{i_2}\ldots \bar B_{i_1-1}$ seen from $\bar Z$ under total angle $\omega_{\frown}$. Then $\omega_{\frown} \geq \omega_{\smile}$. Indeed, this follows from
$$\pi-\omega_{\smile}+\pi+\omega_{\frown} = \omega \geq 2\pi$$
that we obtain from the decomposition of $R$ into orthoschemes. Note that $\omega_\frown$ may be greater than $2\pi$ and $\bar B_\frown$ may be self intersecting when developed to the plane.

Now we return to $\bar B_{\smile}$ and start to develop positive triangles $\bar Z \bar B_{i_1-1} \bar B_{i_1-2}, \ldots$ to the right from the ray $\bar Z \bar B_{i_1}$. We note that if $\phi_{i_1}=\pi$, then $\bar B_{i_1-1}$ coincides with $\bar B_{i_1}$, otherwise, $\bar B_{i_1}$ lies in the interior of the segment $\bar Z \bar B_{i_1-1}$. The new polygonal curve is concave (it is a part of polygonal curve $\bar B_{\frown}$). As $\omega_{\frown}\geq \omega_{\smile}$, there are two possibilities: either the new polygonal curve is entirely above $\bar B_{\smile}$ (except possibly $\bar B_{i_2-1}$) or it intersects $\bar B_{\smile}$ transversely. 

\begin{figure}
\begin{center}
\includegraphics[scale=0.6]{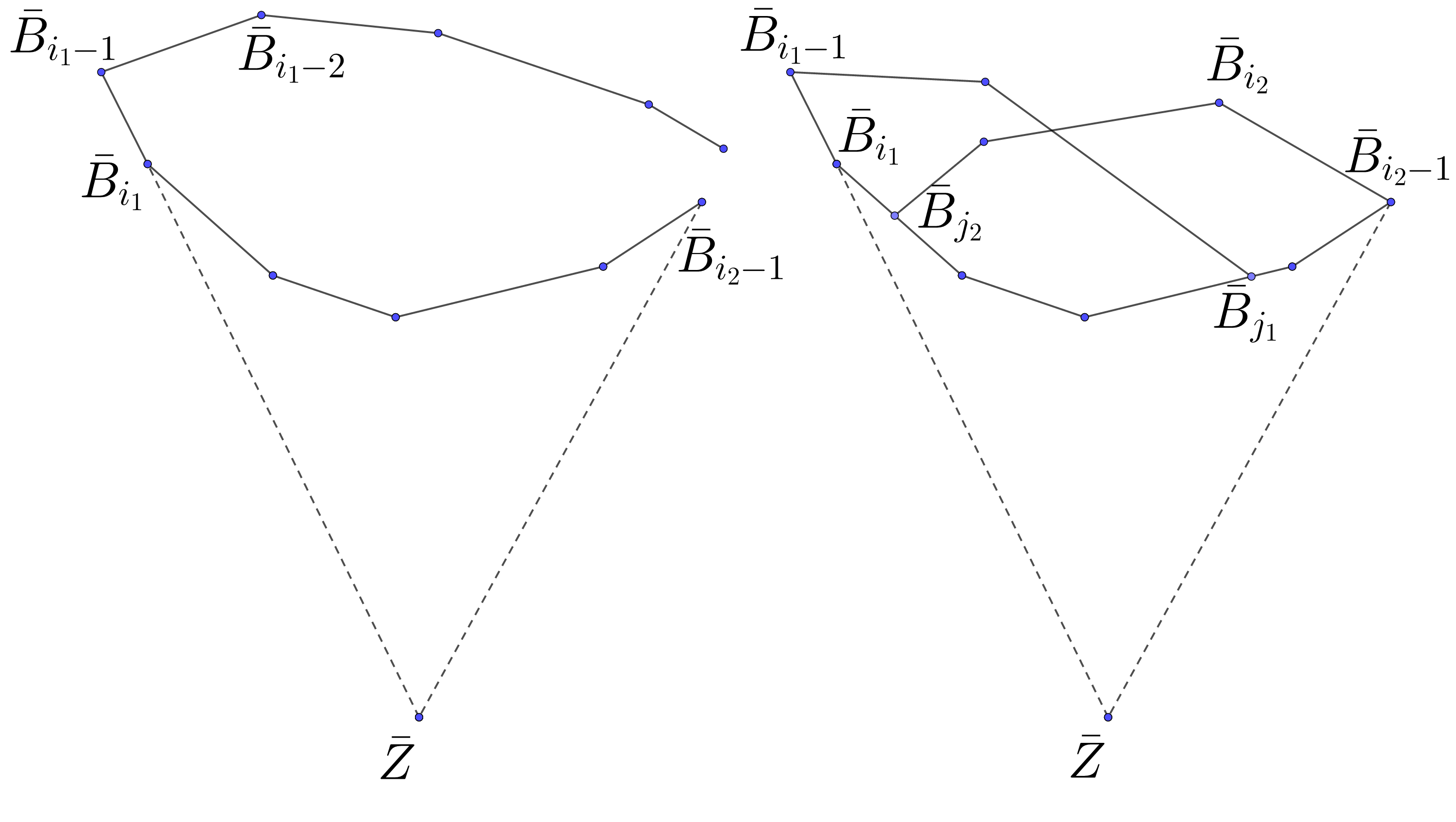}
\caption{Left: curve $\bar B_{\frown 1}$ is above $\bar B_{\smile}$. \\ Right: Curves $\bar B_{\frown 1}$ and $\bar B_{\frown 2}$ intersect.}
\label{Pic8-2}
\end{center}
\end{figure}

In the first case the proof is finished. In the second case if the intersection point is not a vertex of the upper polygonal curve, then we declare it a vertex. Denote it by $B_{j_1}$. Then we start to develop the triangles $\bar Z \bar B_{i_2} \bar B_{i_2+1}, \ldots$ to the left from the ray $\bar Z \bar B_{i_2-1}$. Similarly, either the new polygonal curve is above $\bar B_{\smile}$ or it intersects it. Then declare the intersection point to be a vertex and denote it by $\bar B_{j_2}$. 

We need to consider the case when both upper polygonal curves intersect $\bar B_{\smile}$. If $\bar B_{j_2}$ is to the left from $\bar B_{j_1}$ (or they coincide), then the upper polygonal curves intersect (possibly non-transversely). Let $\bar B$ be an intersection point. 
The polygon $\bar Z \bar B_{i_1} \ldots \bar B_{i_2-1}$ is covered by the union of polygons $\bar Z \bar B_{i_1-1} \ldots \bar B$ and $\bar Z \bar B \ldots \bar B_{i_2}$. They are made from distinct dual orthoschemes because $\angle \bar B_{i_1-1} \bar Z \bar B + \angle \bar B \bar Z \bar B_{i_2} = \omega_{\smile} \leq \omega_{\frown}.$

The last remaining case is when $\bar B_{j_1}$ lies to the left from $\bar B_{j_2}$. We denote the polygonal curve $\bar B_{i_1-1} \ldots \bar B_{j_1}$ by $\bar B_{\frown 1}$ and the polygonal curve $\bar B_{j_2}\ldots \bar B_{i_2}$ by $B_{\frown 2}$. We will prove that $\lambda > 2\pi$. Our plan is to find points $A_j$, $A_k$ on the boundary of $R$ such that

(1) the angles of both sectors between $ZA_j$ and $ZA_k$ are at least $\pi$;

(2) $ZA_j+ZA_k> \pi$.

Indeed, this means that both parts of the boundary between $A_j$ and $A_k$ have lengths greater than $\pi$. This will finish the proof. 

Assume that $\bar B_{j_1}$ lies on the segment $\bar B_{l_1-1} \bar B_{l_1}$ of $\bar B_{\smile}$ (possibly $\bar B_{l_1}$ coincides with $\bar B_{j-1}$ as a point, however we distinguish then $\bar B_{j-1}$ as a vertex of $\bar B_{\frown 1}$ and $\bar B_{l_1}$ as a vertex of $\bar B_{\smile}$). We claim that the angle subtended by $A_{j_1+1}A_{j_1+2}\ldots A_{l_1}$ is greater than $\pi$.

Indeed, we consider the polygon formed by $\bar B_{\smile}$, $\bar B_{\frown 1}$ in the tangent plane to the unit sphere at the north pole. Consider the solid cone  from the origin over this polygon and take its polar cone. The polar cone intersects the sphere in a polygon consisting of two parts: one is the development of $ZA_{j_1+1} A_{j_1+2} \ldots A_{l_1}$ from $R$ and the second is a geodesic spherical triangle formed by connecting $ZA_{l_1}$ with $ZA_{j_1+1}$ via a geodesic segment. The length of this segment is smaller than $\pi$, hence it is seen from $Z$ with angle smaller than $\pi$. This shows that the angle subtended by $A_{j_1+1}A_{j_1+2}\ldots A_{l_1}$ is greater than $\pi$.

Similarly, let $\bar B_{j_2}$ belong to segment $\bar B_{l_2-1} \bar B_{l_2}$ on $\bar B_{\smile}$. Then the angle subtended by $A_{l_2}A_{l_2+1}\ldots A_{j_2}$ is at least $\pi$.

Now take an arbitrary point of $\bar B_{\smile}$ between $\bar B_{j_1}$ and $\bar B_{j_2}$. Let it belong to a segment $\bar B_{k-1} \bar B_k$ on $\bar B_{\smile}$. We see that the angles of both sectors between $ZA_{j_1+1}$ and $ZA_k$ are at least $\pi$; similarly for $ZA_{j_2}$ and $ZA_k$. It remains to show that either $ZA_{j_1+1}+ZA_k> \pi$ or $ZA_{j_2}+ZA_k> \pi$. We look at three lines: $L_{j_1+1}$ through $\bar B_{j_1}, \bar B_{j_1+1}$; $L_{j_2}$ through $\bar B_{j_2}, \bar B_{j_2-1}$ and $L_k$ through $\bar B_{k-1}, \bar B_k$. Due to their location with respect to $\bar Z$, either $L_{j_1+1}$ or $L_{j_2}$ intersects the shortest segment from $\bar Z$ to $L_k$. Without loss of generality, assume that this is $L_{j_2}$. Then the distance from $\bar Z$ to $\bar L_{j_2}$ is smaller than the distance from $\bar Z$ to $L_k$. Going to the spherical duals recall that for every $i$ the distance from $\tilde Z$ to the geodesic $\tilde B_i \tilde B_{i-1}$ is $|\pi/2-ZA_i|$. Then our conclusion on the lines is equivalent to
$$ZA_k - \pi/2 > \pi/2 - ZA_{j_2},$$
which means $ZA_{j_2}+ZA_k> \pi$ as desired.
\end{proof}

\subsubsection{Comparing discrete curvature}


Now we investigate how maxima of the discrete curvature behave under some constraints. Let $d \in \mathfrak D_{sc}(V)$ and $P^0=P(d,V, h^0)$ be the convex polyhedral Fuchsian manifold. By Lemma~\ref{maxprinciple} $P^0$ maximizes $S$ over $H(d, V)$. Consider $W\subseteq V$. For each $t \geq 0$ define the set 
$$H( t)=\{h \in H(d,V): \sinh h_w \geq e^{ t} \sinh h^0_w \textrm{~for~every~}w \in W\},$$
which is non-empty due to Lemma~\ref{movingup}.

The core of our proof of Main Lemma I is the following

\begin{lm}
\label{maxconstrain}
The maximum of $S$ over $H(t)$ is attained at a unique point $h(t)$. Moreover, the map $t \rightarrow h(t)$ is $C^1$ and satisfies 

(1) for each $v \in V \backslash W$ we have $\kappa_v(t)=0$ and $0 < \dot h_v(t) <\tanh h_v(t)$;

(2) for each $w \in W$ we have $\sinh h_w(t)=e^t \sinh h^0_w$ and $\dot h_w(t)=\tanh h_w(t)$. 
\end{lm}

Here $\kappa_v(t)$ are curvatures $\kappa(P(t))$ in the cone-manifold $P(t)=P(d, V, h(t))$. The property $\dot h_w(t)=\tanh h_w(t)$ follows directly from $\sinh h_w(t)=e^t \sinh h^0_w$.

\begin{proof}
For reader's convenience we advise to assume during the first read that $W$ is a single point $w$. 

First, there is at least one maximum of $S$ over $H(t)$. The argument is the same as in Lemma~\ref{maxprinciple}. There are constants $0<m<M$ such that the supremum of $S$ over $H(t)$ is equal to the supremum over $H(t) \cap H(d, V, m, M)$, which is compact (see Lemma~\ref{H_d comp}). Indeed, $M$ exists due to Lemma~\ref{compactabove}, $m$ exists due to Lemma~\ref{alldecrease}. 

Let $\check h( t)$ be any of maximum points of $S$ over $H( t)$. By $\check \kappa_v( t)$ we denote the particle curvatures in the cone-manifold $\check P( t):=P(d, V, \check h( t))$. 

\begin{prop}
For every $v \in V\backslash W$ we have $\check \kappa_v( t)\geq 0$.
\end{prop}

Suppose the contrary. Because of Corollary~\ref{decrease} the height of $v$ can be decreased. This deformation stays in $H( t)$ and increases $S$ due to Lemma~\ref{S1der}, which contradicts to the choice of $\check h( t)$.
\vskip+0.2cm

\begin{prop}
For every $v \in V$ we have $\check \kappa_v( t)\leq 0$. 
\end{prop}

The proof is the same as of the previous Claim, but we use Lemma~\ref{increase}.
\vskip+0.2cm

We obtain that for each $\check h(t)$ and each $v \in V\backslash W$, we have $\check \kappa_v( t)=0$. For $w \in W$ we have $\check \kappa_w(t) \leq 0$. Also if $t>0$, then there exists $w \in W$ such that $\check \kappa_w(P) < 0$ (we use the fact that the critical point of $S$ is unique due to Theorem~\ref{Fillastre} and Remark~\ref{markeduniq}). Corollary~\ref{decrease} implies that if for some $w \in W$ we have $\check \kappa_w(P)<0$, then $\sinh \check h_w( t)=e^{ t}\sinh h^0_w$; otherwise we can increase $S$ while staying in $H(t)$.

For all $t\ge0$ Corollary~\ref{inter} implies that $\check h(t) \in \inter(H(d,V))$ and Corollary~\ref{hessnond} implies that the Hessian of $S$ is non-degenerate at $\check h(t)$.

Take $t_0 > 0$ and some $\check h(t_0)$. Let us show that there is a unique piecewise $C^1$-curve $h(t): (-\infty; t_0] \rightarrow H(d, V)$ such that:

(0) $h(t_0) = \check h(t_0)$;

(1) for any $v \in V \backslash W$, we have $\kappa_v(t)=0$ (here we write $\kappa_v(t)$ for particle curvatures of the cone-manifold $P(t)=P(d, V, h(t))$);

(2) for $w \in W$ if for some $t'>t$ we have $\kappa_w(t')=0$, then $\kappa_w(t)=0$; otherwise $$\sinh  h_w(t)=e^{t-t_0}\sinh \check h_w(t_0)=e^{t}\sinh h^0_w.$$ 

Due to this definition, if for some $t$ we have $h(t)=h^0$, then $h(t)$ stabilizes, i.e., for each $t' \leq t$ we have $h(t')=h^0$.

First, we see that if $h(t)$ is defined for some $t \leq t_0$ and $h(t) \neq h^0$, then there exists $\Delta t>0$ such that $h$ is defined uniquely for all $t' \in (t-\Delta t; t]$. Indeed, due to Corollary~\ref{hessnond} and Corollary~\ref{inter}, the Hessian of $S$ is non-degenerate at $h(t)$ and $h(t) \in \inter(H(d,V))$ (note that the $\kappa_w$ can not be positive for $w \in W$ as they are non-negative for $t=t_0$ and once they become 0, they stay equal to 0). Then our claim follows from the Implicit Function Theorem.

Assume that $h(t)$ is defined uniquely by conditions (0)-(2) over the interval $(t'; t_0]$ for some $t'<t_0$, and $h(t) \neq h^0$ for all $t \in (t'; t_0]$. The second assumption implies that there exists $w \in W$ such that $\kappa_w(t) <0$ for all $t \in (t'; t_0]$. Thus, $$\lim_{t \rightarrow t'} \sinh h_w(t)=e^{t'}\sinh h^0_w.$$ By Lemma~\ref{alldecrease} there exists $m>0$ such that for all $t \in (t'; t_0]$ and all $v \in V$ we have $h_v(t) \geq m$. By Claim~\ref{allincrease} from the proof of Lemma~\ref{compactabove} there exists $M>0$ such that for  all $t \in (t'; t_0]$ and all $v \in V$ we get $h_v(t) \leq M$. Then for all $t \in (t'; t_0]$ we have $h(t) \in H(d, V, m, M)$, which is compact by Lemma~\ref{H_d comp}. Then there exists a unique point $h(t') \in H(d, V)$ that is a continuous extension of $h(t)$. Moreover, by continuity for all $v \in V$ we have $\kappa_v(t') \leq 0$. Then the Hessian of $S$ is non-degenerate at $h(t')$ and $h(t') \in \inter(H(d,V))$. Due to the Implicit Function Theorem, $h(t)$ is $C^1$ in some right neighbourhood of $t'$.

Thereby, $h(t)$ is uniquely defined for all $t \leq t_0$. There might be kink points corresponding to those $t$, when $\kappa_w$ becomes 0 for a vertex $w \in W$. 

It follows that either there exists $\check t_0 < t_0$ such that $h(\check t_0) = h^0$ or there is $w \in W$ such that $\kappa_w(t) < 0$ for all $t \in (-\infty; t_0]$. In the second case $h_w(t) \rightarrow 0$ as $t \rightarrow -\infty$. Then Lemma~\ref{nottozero} excludes this case as it shows that $\kappa_w$ becomes positive for $t$ close to $-\infty$.

Thus, there exists $\check t_0<t_0$ such that $h(\check t_0) = h^0$. We assume that $\check t_0$ is the maximal with this property. There is $w \in W$ such that for all $t > \check t_0$ we have $\kappa_w(t) < 0$. This means that 
$$\sinh h_w(\check t_0)=e^{\check t_0}\sinh h^0_w.$$ 
As it is equal to $\sinh h^0_w$, then $\check t_0 = 0$.

Our current aim is

\begin{prop}
\label{noWzero}
For all $w \in W$ and all $0 < t \leq t_0$ we have $\kappa_w(t) < 0$.
\end{prop}

This requires some preparation. Fix $t$ that is not a kink point of $h(t)$ and define
$$W(t) = \{w \in W: \kappa_w(t)<0\}.$$
If $w \in W(t)$, then $\dot h_w(t)=\tanh h_w(t)$.

\begin{prop}
\label{range}
For each $v\in V$ we have $\dot h_v( t)>0$.
\end{prop}

Indeed, we already saw this for $v \in W(t)$. If $v \not\in W(t)$, then $\kappa_v( t)$ is identically 0. Lemma~\ref{S2der} gives 
$$-\dot\kappa_v=\frac{1}{\cosh h_v}\sum_{\vec e=\overrightarrow{vu} \in \vec E_v} \frac{\cot\phi^+_{\vec e}+\cot\phi^-_{\vec e}}{\sin \alpha_{\vec e}}\cdot \frac{\cosh a_{\vec e} \dot h_v-\dot h_u}{\sinh a_{\vec e}}=0.$$

Let $v$ be a vertex, where the minimum of $\dot h_v(t)$ is attained over the vertices not from $W(t)$. Because of Lemma~\ref{isolvertex} there exists at least one strict edge emanating from $v$, so the sum is non-empty. If $\dot h_v<0$, then is easy to see that each summand is negative, which is a contradiction. If $\dot h_v=0$, then looking at the sum we see that for every $u$ connected with $v$ we also have $\dot h_u=0$. Repeating this and noting that the graph of strict edges is connected due to Corollary~\ref{inter}, we obtain a contradiction. This proves Claim~\ref{range}.
\vskip+0.2cm

\begin{prop}
\label{range2}
For every $v \in V\backslash W(t)$, we have $\dot h_v < \tanh h_v$
\end{prop}

Define $x_v(t):=\tanh h_v( t)$. For every $v\in V$ we get
$$-\frac{\partial \kappa_v}{\partial h_v} x_v - \sum_{u \neq v}\frac{\partial \kappa_v}{\partial h_u}x_u=$$
$$= \frac{1}{\cosh h_v}\sum_{\vec e=\overrightarrow{vu} \in \vec E_v} \frac{\cot\phi^+_{\vec e}+\cot\phi^-_{\vec e}}{\sin \alpha_{\vec e}}\cdot \frac{\cosh a_{\vec e} \tanh h_v - \tanh h_u}{\sinh a_{\vec e}} .$$

In this sum, $\vec e = \overrightarrow{vu}$ is an oriented edge from $v$ to $u$ and the sum is over all edges starting at $v$. With the help of Lemma~\ref{cotan} we obtain 
$$-\frac{\partial \kappa_v}{\partial h_v} x_v - \sum_{u \neq v}\frac{\partial \kappa_v}{\partial h_u}x_u=
\frac{1}{\cosh^2 h_v}\sum_{\vec e \in \vec E_v} \frac{\cot\phi^+_{\vec e }+\cot\phi^-_{\vec e }}{\sin \alpha_{\vec e }}\cdot \cot \alpha_{\vec e} .$$

Recall that $V=V(d)$, thus we have $\nu_v(d)>0$. From Lemma~\ref{spherarea} we see
\begin{equation}
\label{1line}
-\frac{\partial \kappa_v}{\partial h_v} x_v - \sum_{u \neq v}\frac{\partial \kappa_v}{\partial h_u}x_u > 0.
\end{equation}
If $v \not\in W(t)$, then we also have
\begin{equation}
\label{2line}
-\frac{\partial \kappa_v}{\partial h_v} \dot h_v - \sum_{u \neq v}\frac{\partial \kappa_v}{\partial h_u}\dot h_u = 0.
\end{equation}
Multiplying (\ref{1line}) with $\frac{\dot h_v}{x_v}>0$ and subtracting (\ref{2line}) from it we write
$$\sum_{u \neq v}\frac{\partial \kappa_v}{\partial h_u}\left(\dot h_u-x_u\frac{\dot h_v}{x_v}\right) > 0.$$
The coefficients $\frac{\partial \kappa_v}{\partial h_u}$ are non-negative. Thus, there exists $u$ connected with $v$ such that $$\dot h_u-x_u\frac{\dot h_v}{x_v}>0,$$
or equivalently $$\frac{\dot h_u}{x_u}>\frac{\dot h_v}{x_v}.$$

Consider the maximum of the expression $\frac{\dot h_v}{x_v}$ over $V.$ We obtain that it is attained at points of $W(t)$, where it is equal 1. Therefore, for every $v \not\in W(t)$ we get $\dot h_v < \tanh h_v$ and Claim~\ref{range2} is proven.
\vskip+0.2cm

Let $t$ is a kink point of $h(t)$, so for some $w \in W$, $\kappa_w(t)$ is 0 and it is not zero for all larger $t$. Then the right derivative of $h_w$ at $t$ is $\tanh h_w(t)$, but the left derivative is smaller. For other $v$ the derivatives of $h_v$ remain continuous.

Suppose that there exist $t\in (0; t_0]$ such that for some $w \in W$ we have $\kappa_w(t)=0$. Take $t$ to be the largest with this property.  We have $\sinh h_w(t) \geq e^{t} \sinh h^0_w$. However, for all $t' \in [0; t)$, we get $\dot h_w < \tanh h_w$. This gives $\sinh h_w(0) > \sinh h^0_w$, which is a contradiction. This proves Claim~\ref{noWzero}.

From all the above discussion we can see that the maximum point $\check h(t)$ over $H(t)$ is unique for all $t \geq 0$ and these points form a $C^1$-curve $h(t) \in H(d, V)$ that is determined by the conditions:

(0) $h(0) = h^0$;

(1) for any $v \in V$, $v \not\in W$, we have $\kappa_v(t)=0$;

(2) for any $w \in W$ we have $$\sinh  h_w(t)=e^{t}\sinh h^0_w.$$ 

Moreover, for any $w \in W$ we obtained $\dot h_w(t)=\tanh h_w(t)$ and for any $v \not\in W$ we got $0<\dot h_v < \tanh h_v$. 

This finishes the proof.
\end{proof}

Main Lemma I follows from Lemma~\ref{maxconstrain} quite easily.

\begin{proof}[Proof of Main Lemma I]
Let $h(t)$, $H(t)$ and $P(t)=P(d, V, h(t))$ be taken from Lemma~\ref{maxconstrain}  (note that $P^1$ from the notation of Main Lemma I is $P(0)=P^0$ in the notation of Lemma~\ref{maxconstrain}). Define $S(t)=S(P(t)).$ We have $h^2 \in H(\tau)$. As $h(t)$ is the maximal point of $S$ over $H(t)$, we see that $S(P^2) \leq S(\tau)$.

From Lemma~\ref{S1der} we obtain 
$$-\dot S( t)=-\sum\limits_{w\in W}\kappa_w\tanh h_w( t).$$
We have $$S(P^1) - S(P^2) \geq S(0) - S(\tau)= - \int_0^{\tau}d S=$$
$$= - \int_0^{\tau} \sum\limits_{w\in W} \kappa_w \tanh h_w d  t\geq  - m\cdot \int_0^{\tau} \sum\limits_{w\in W} \kappa_wd  t.$$

From Lemma~\ref{maxconstrain} for any $w \in W$ we get
$$-\dot\kappa_w \geq \frac{1}{\cosh h_w}\sum_{\vec e=\overrightarrow{wu} \in \vec E_w} \frac{\cot\phi^+_{\vec e}+\cot\phi^-_{\vec e}}{\sin \alpha_{\vec e}}\cdot \frac{\cosh a_{\vec e} \tanh h_w-\tanh h_u}{\sinh a_{\vec e}}$$

Lemma~\ref{cotan} transforms the right hand side as
$$\frac{1}{\cosh h_w}\sum_{\vec e=\overrightarrow{wu} \in \vec E_w} \frac{\cot\phi^+_{\vec e}+\cot\phi^-_{\vec e}}{\sin \alpha_{\vec e}}\cdot \frac{\cosh a_{\vec e} \tanh h_w-\tanh h_u}{\sinh a_{\vec e}}=$$
$$= \frac{1}{\cosh^2 h_w}\sum_{\vec e \in \vec E_w} \frac{\cot\phi^+_{\vec e}+\cot\phi^-_{\vec e}}{\sin \alpha_{\vec e}}\cdot \cot\alpha_{\vec e}  $$

Combining this with Lemma~\ref{spherarea} we get
$$-\dot\kappa_w  \geq\frac{\nu_w-\kappa_w}{M}.$$

Define $\kappa:=\sum\limits_{w\in W} \kappa_w \leq 0$. We obtain $-\dot \kappa\geq \frac{\nu-\kappa}{M}.$
Rewrite it as $$\frac{d (\nu - \kappa( t))}{ \nu - \kappa( t)} \geq \frac{d  t}{M}.$$
Integrating this we see $$-\kappa( t)\geq \nu\left(\exp\left(\frac{ t}{M}\right)-1\right).$$

Finally, we get $$S(P^1) - S(P^2) \geq -  m\cdot \int_0^{\tau} \kappa d  t\geq$$ $$\geq\nu m \cdot \int_0^{\tau}\left(\exp\left(\frac{ t}{ M}\right)-1\right) d  t =$$ $$= \nu m \left( M \exp\left(\frac{\tau}{M}\right)-\tau\right) .$$

This finishes the proof.
\end{proof}

\newpage
\subsection{Proof of Main Lemma II}

Main Lemma II states that if we have a sufficiently fine strictly short triangulation of $(S_g, d)$, then we can replace each triangle by a triangle with a swept cone-metric and the same sides and angles. To this purpose it is enough to treat each triangle separately. We summarize our main definitions restricted to a single triangle:

\begin{dfn}
By a (topological) \emph{triangle} $T$ we mean a topological disc with three marked points at the boundary called vertices. By a metric on $T$ we always mean a metric such that the boundary is geodesic except (possibly) the vertices. 

Cone-metrics or CBB($-1$) metrics on $T$ are defined naturally. We call a cone-metric $d$ on $T$ \emph{convex} if, in addition to the condition $\nu_p(d) \leq 2\pi$ for the interior points of $T$, also the angles of the vertices of $T$ are at most $\pi$. One can see that this requirement holds for convex triangles in CBB($-1$) metrics (Definition~\ref{convexpolyg}). 

A metric on $T$ is called \emph{short} if the sides of $T$ are shortest paths and the angles of vertices are strictly smaller than $\pi$. It is called \emph{strictly short} if, in addition, the sides are unique shortest paths. A cone-metric $d$ on $T$ is called \emph{swept} if $|V(d)|\leq 1$. Here we note that by $V(d)$ we mean the set of interior singularities of $d$, i.e., we do not include the vertices of $T$ in $V(d)$. Metric $\hat d$ on $T$ is a \emph{sweep-in} of metric $d$ if $\hat d$ is swept and $T$ has the same side lengths and angles in $d$ and $\hat d$.
\end{dfn}

Main Lemma II follows from two lemmas:

\begin{lm}
\label{mlIIv}
There exists an absolute constant $\delta>0$ such that if $d$ is a convex cone-metric on a triangle $T$ with $\diam(T, d) < \delta$ and $\nu(T, d) < \delta$, then there exists a unique sweep-in $\hat d$ of $d$ on $T$.
\end{lm}

\begin{lm}
\label{mlII}
There exists an absolute constant $\delta>0$ such that if $d$ is a strictly short CBB($-1$) metric on a triangle $T$ with $\diam(T, d) < \delta$ and $\nu(T, d) < \delta$, then there exists a unique sweep-in $\hat d$ of $d$ on $T$.
\end{lm}

The plan is first to prove in Subsection~\ref{mergesec} the existence in Lemma~\ref{mlIIv} via the \emph{curvature merging process}. We also give some quantitative estimates on $\hat d$ that will be crucial further in the proof of Main Lemmas IIIA-B. Next, in Subsection~\ref{conesec} we prove the uniqueness in both lemmas. To this purpose we only need to prove that a swept convex metric on a triangle is determined by side lengths and angles. The proof relies on elementary geometry. Last, in Subsection~\ref{cbbswept} we prove the existence in Lemma~\ref{mlII} with the help of approximation by cone-metrics. Note that the condition on $d$ to be strictly short in the statement of Lemma~\ref{mlII} seems to be non-essential. It is used to avoid some annoying degenerations that we want to exclude from the present manuscript in order to simplify the exposition. 

\subsubsection{Curvature merging}
\label{mergesec}

The main proof idea of Lemma~\ref{mlIIv} is that we can take two cone-singularities, cut the metric along a geodesic connecting them and glue-in there a piece of a hyperbolic cone so that the number of cone-singularities is reduced. The main difficulty is that we can not do this in all circumstances. In order to perform this operation, the length of the geodesic and the curvatures of the singularities should satisfy some mutual conditions. Hence, we should control carefully what happens with the metric on each step in order to be able to perform the next one.

We start from two simple lemmas on hyperbolic triangles.

\begin{lm}
\label{exist}
Let $\Delta > 0$ and define
$$\eta=\eta(\Delta)=2\arcsin\frac{1}{\sqrt 2 \cosh\frac{\Delta}{2}} \in (0; \pi/2).$$

Then for any numbers $a \in (0; \Delta)$ and $\beta, \gamma >0$ such that $\beta+\gamma < \eta$ there exists a hyperbolic triangle $ABC$ with $BC=a$, $\angle B = \beta$ and $\angle C = \gamma$. Moreover, $\angle A> \pi/2$.
\end{lm}

\begin{proof}
Suppose that the triangle $ABC$ exists. Write the dual cosine law for $\angle A=\alpha$:
\begin{equation}
\label{dualcos}
\cos\alpha=-\cos\beta\cos\gamma+\sin\beta\sin\gamma\cosh a.
\end{equation}

One can prove that $ABC$ exists if and only if the right hand side is in the range of cosine of a non-zero angle, i.e., it belongs to $(-1; 1)$. Due to elementary trigonometry, the lower bound is always true. We are interested in $\angle A>\pi/2$, i.e., when the right hand side of (\ref{dualcos}) is less than 0. 

Let us fix $a$ and fix the sum of $\beta$ and $\gamma$. Then if one considers the right hand side of (\ref{dualcos}) as the function of $\beta$, by taking the derivative it follows that the right hand side is maximized when $\beta=\gamma$. 
Now let us check what is the value of $\beta=\gamma$ making the right hand side zero. It is easy to compute that
$$\beta=\gamma=\arcsin\frac{1}{\sqrt 2 \cosh\frac{a}{2}}=\frac{\eta(a)}{2}.$$

This discussion shows that for any $\beta, \gamma>0$ such that $\beta+\gamma<\eta(a)$ the right hand side is smaller than zero, so the triangle $ABC$ exists and $\angle A>\pi/2$. It remains only to note that $\eta(a)$ decreases as $a$ increases.
\end{proof}

\begin{lm}
\label{curvatbound}
Let $ABC$ be a hyperbolic triangle with sides $a,b,c$ and angles $\alpha, \beta, \gamma$ respectively such that $\alpha>\pi/2$. Then
$$\pi-\alpha<\beta+\gamma\cosh a.$$
\end{lm}

\begin{proof}
Define
$$\alpha' = \pi - \beta-\gamma.$$
We have$$\cos\alpha'=\sin\beta\sin\gamma-\cos\beta\cos\gamma,$$ $$\cos\alpha=\sin\beta\sin\gamma\cosh a-\cos\beta\cos\gamma,$$
$$\cos \alpha - \cos \alpha'=\sin \beta \sin \gamma (\cosh a -1).$$
We have $\pi>\alpha' > \alpha > \pi/2$. From the concavity of cosine over $[\pi/2; \pi]$ we obtain $$\frac{\alpha'-\alpha}{\cos\alpha-\cos\alpha'} < \frac{1}{\sin\alpha'}.$$
From this we get
$$\alpha'-\alpha < \frac{\cos\alpha-\cos\alpha'}{\sin \alpha'}=\frac{\sin \beta \sin \gamma (\cosh a -1)}{\sin(\beta+\gamma)}.$$
Note that $\beta<\beta+\gamma<\pi/2$ and $\sin \gamma<\gamma$. Therefore,
$$\alpha'-\alpha< \gamma(\cosh a -1).$$
It follows that
$$\pi-\alpha<\beta+\gamma+\gamma(\cosh a - 1)=\beta+\gamma\cosh a.$$
\end{proof}

Fix $\Delta>0$ and choose $\eta=\eta(\Delta)$ from Lemma~\ref{exist}. Then choose $\delta_{\nu}>0$ such that

\begin{equation}
\label{choice1}
\delta_{\nu}<\min\{\frac{2\eta}{\cosh\Delta},
 \frac{\Delta}{\sinh\Delta}\}.
\end{equation}

After that we choose $\delta_D>0$ such that 
\begin{equation}
\label{choice2}
\delta_D  < \Delta-\delta_{\nu}\sinh\Delta.
\end{equation}

Now we formulate and prove a quantitative version of the existence part of Lemma~\ref{mlIIv}:

\begin{lm}
\label{merge}
Let $d$ be a convex cone-metric on $T$ with $\diam(T, d)<\delta_D$ and $\nu(T, d)<\delta_{\nu}$. Then there exists a sweep-in $\hat d$ of $d$ on $T$ such that $$\diam(T, d)\leq \diam(T, \hat d) <\Delta,~~~\nu(T, d)\leq \nu(T, \hat d)<2\eta.$$
\end{lm}

\begin{proof}
Put $T_0=T$, $d_0=d$. At first step take two points $v_0, v_1 \in V(d_0)$ and connect them by a shortest path $\chi_1$. 

Consider the hyperbolic triangle $ABC$ such that $$BC = \len(\chi_1, d_0)=d_0(v_0, v_1),~~~\angle B =  \nu_{v_0}(d_0)/2,~~~\angle C=\nu_{v_1}(d_0)/2.$$ 
We have $$BC \leq \diam(T_0, d_0)<\delta_D<\Delta,$$ $$\angle B+\angle C=\frac{\nu_{v_0}(d_0)+\nu_{v_1}(d_0)}{2}\leq \frac{\nu(T, d)}{2}<\frac{\delta_{\nu}}{2}<\eta.$$
Thus, Lemma~\ref{exist} implies that $ABC$ exists and $BC$ is its greatest side. 

We take two copies of $ABC$ and glue isometrically the sides corresponding to $AB$ and $AC$. We call the obtained figure  a \emph{bigon} and denote it by $X_1$. It has the piecewise geodesic boundary consisting of two segments equal to $BC$, two kink points at the boundary, which we continue to denote by $B$ and $C$ abusing the notation, and a conical point in the interior, which we continue to denote by $A$.

Cut $T_0$ along $\chi_1$. Glue isometrically the boundary of $X_1$ with the boundary of the cut: the vertex $B$ is glued with $v_0$, $C$ is glued with $v_1$. Denote the resulting topological space by $T_1$ and the induced path metric $d_1$. Since now we consider $X_1$ as a subset of $T_1$. We have a natural map $$g_1: T_0 \backslash \chi_1 \rightarrow T_1 \backslash X_1,$$ which is a local isometry with respect to metrics $d_0$, $d_1$.

We identify $V(d_0)$ with its $g_1$-image. The point $A$ will be also denoted by $w_1 \in T_1$. Then $$\nu_{v_0}(d_1)=\nu_{v_1}(d_1)=0,$$ $$\nu_{w_1}({d_1}) > \nu_{v_0}(d_0)+\nu_{v_1}(d_1),$$ $$V(d_1)=(V(d_0)\cup\{w_1\})\backslash\{v_0,v_1\}.$$ 

\begin{prop}
\label{distincr}
The map $g_1: (T_0\backslash \chi_1, d_0)\rightarrow (T_1\backslash X_1, d_1)$ does not decrease the distances.
\end{prop}

Indeed, let $p', p'' \in T_1\backslash X_1$. Connect them by a shortest path in $(T_1, d_1)$. If it does not intersect $X_1$, then the distance between them did not decrease. If it intersects $X_1$, then the intersection is a segment. This is because of the strong non-overlapping property (see Section~\ref{cbb-1}). Let the shortest path intersect the boundary of $X_1$ in points $q', q'' \in T_1$ such that $q'$ is closer to $p'$ and $q''$ is closer to $p''$. See Figure~\ref{Pic9}. Choose points $q'_0, q''_0 \in \chi_1\subset T_0$ such that $$d_0(v_1, q'_0)=d_1(v_1, q'), $$ $$d_0(v_1, q''_0)=d_1(v_1, q'').$$

\begin{figure}
\begin{center}
\includegraphics[scale=0.6]{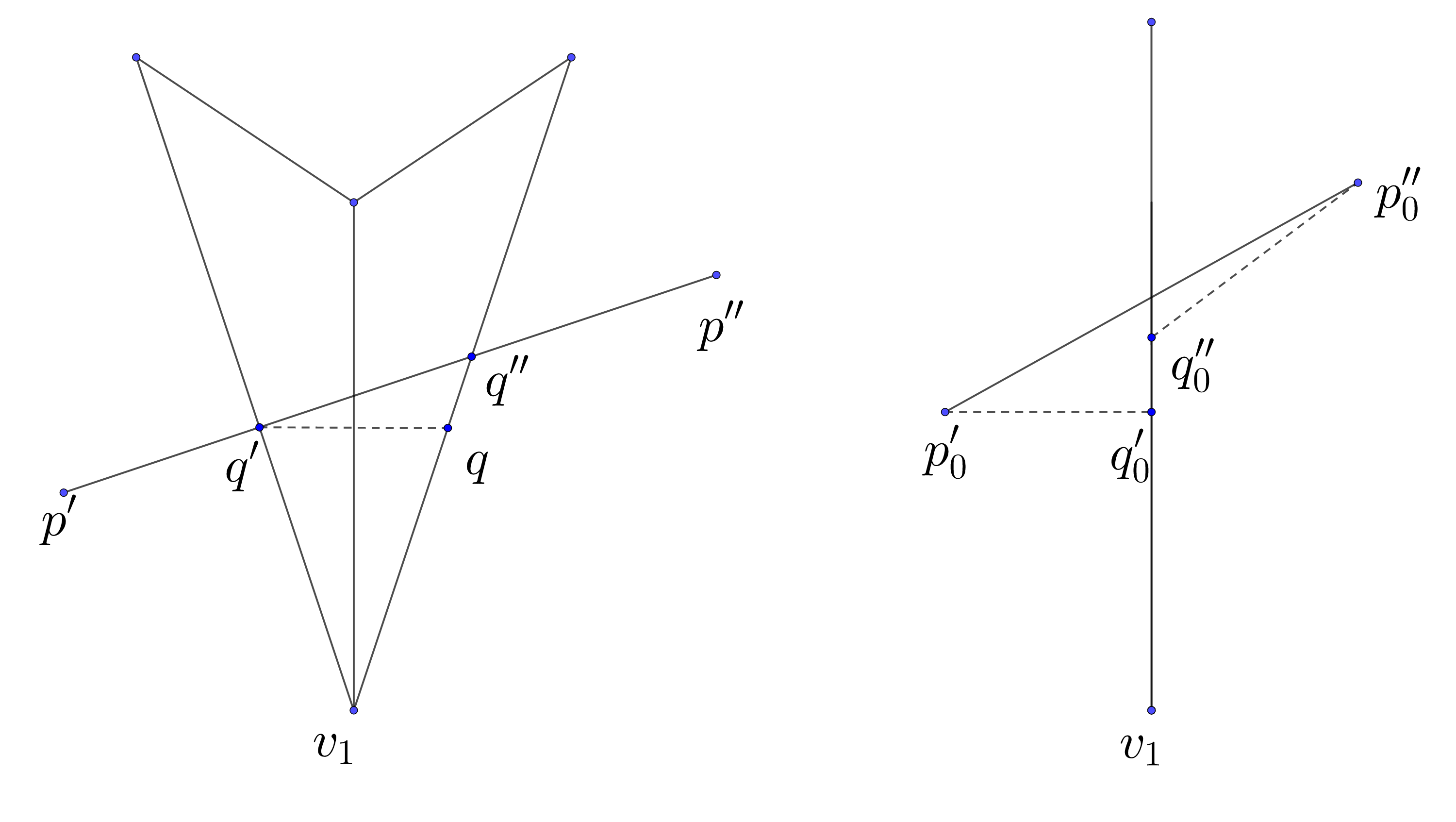}
\caption{To the proof of Claim~\ref{distincr}.}
\label{Pic9}
\end{center}
\end{figure}

Without loss of generality we assume that $q'_0$ lies between $v_1$ and $q''_0$. Then choose a point $q \in T_1$ on the boundary of $X_1$ between $v_1$ and $q''$ such that $d_1(v_1, q)=d_1(v_1, q')$. It is clear that $d_1(q'', q)\leq d_1(q', q'')$.  Define $p'_0=g_1^{-1}(p')$, $p''_0=g_1^{-1}(p'')$. Then
\begin{multline*}
d_0(p'_0, p''_0) \leq d_0(p'_0, q'_0) + d_0(q'_0, q''_0) + d_0(q''_0, p''_0) \leq \\ \leq d_1(p', q') + d_1(q, q'') + d_1(q'', p'')\leq \\ \leq d_1(p', q') + d_1(q', q'') + d_1(q'', p'') =d_1(p', p'').
\end{multline*}
Claim~\ref{distincr} is shown. As a corollary we get
\vskip+0.2cm

\begin{prop}
\label{diambound}
$\diam(T_0, d_0)\leq\diam(T_1, d_1).$
\end{prop}

Now we need to show an upper bound.

\begin{prop}
\label{distbound}
$\diam(T_1, d_1) < \diam(T_0, d_0) + \nu_{v_1}(d_0)\sinh(\diam(T_0, d_0)).$
\end{prop}

Indeed, consider in $\H^2$ two triangles $AB'C$ and $AB''C$ sharing the edge $AC$, each is isometric to the triangle $ABC$. As $\angle B'AC=\angle B''AC>\pi/2$, the point $A$ belongs to the interior of the triangle $B'B''C$.

Let $\psi$ be the (unique) shortest path connecting $w_1$ and $v_0$ in $X_1$. Cut $X_1$ along $\psi$ and consider the developing map $${\rho: X_1 \backslash \psi \rightarrow CB'AB''}.$$ One can see that it does not decrease the distances. Then $$\diam(X_1) \leq \diam(B'B''C) = \max(B'B'', B'C).$$ Note that $B'C \leq \diam(T_0, d_0)$ and $$B'B'' < 2\sinh(B'B''/2) = 2\sin(\nu_{v_1}(d_0)/2)\sinh B'C <$$ $$< \nu_{v_1}(d_0)\sinh(\diam(T_0, d_0)).$$

Consider $p' \in T_1 \backslash X_1$ and $p'' \in X_1$. Connect them by a shortest path. Let it intersect the boundary of $X_1$ in a point with $\rho$-image in $CB'$. In the triangle $B'B''C$ draw the segment through $\rho(p'')$, which cuts off equal segments from the segments $CB'$ and $CB''$. Consider its intersection point with $CB'$ and let $q \in T_1$ be its preimage under $\rho$. Let $q_0 \in \chi_1\subset T_0$ be a point such that $d_0(v_0, q_0) = d_1(v_0, q)$ and $p'_0=g^{-1}(p')\in T_0 \backslash \chi_1.$ See Figure~\ref{Pic10-5}. Then $$d_{1}(p', p'') \leq d_{1}(p', q) +d_{1}(p'', q) \leq d_0(p'_0, q_0) + B'B'' <$$ $$< \diam(T_0, d_0) + \nu_{v_1}(d_0)\sinh(\diam(T_0, d_0)).$$

\begin{figure}
\begin{center}
\includegraphics[scale=0.6]{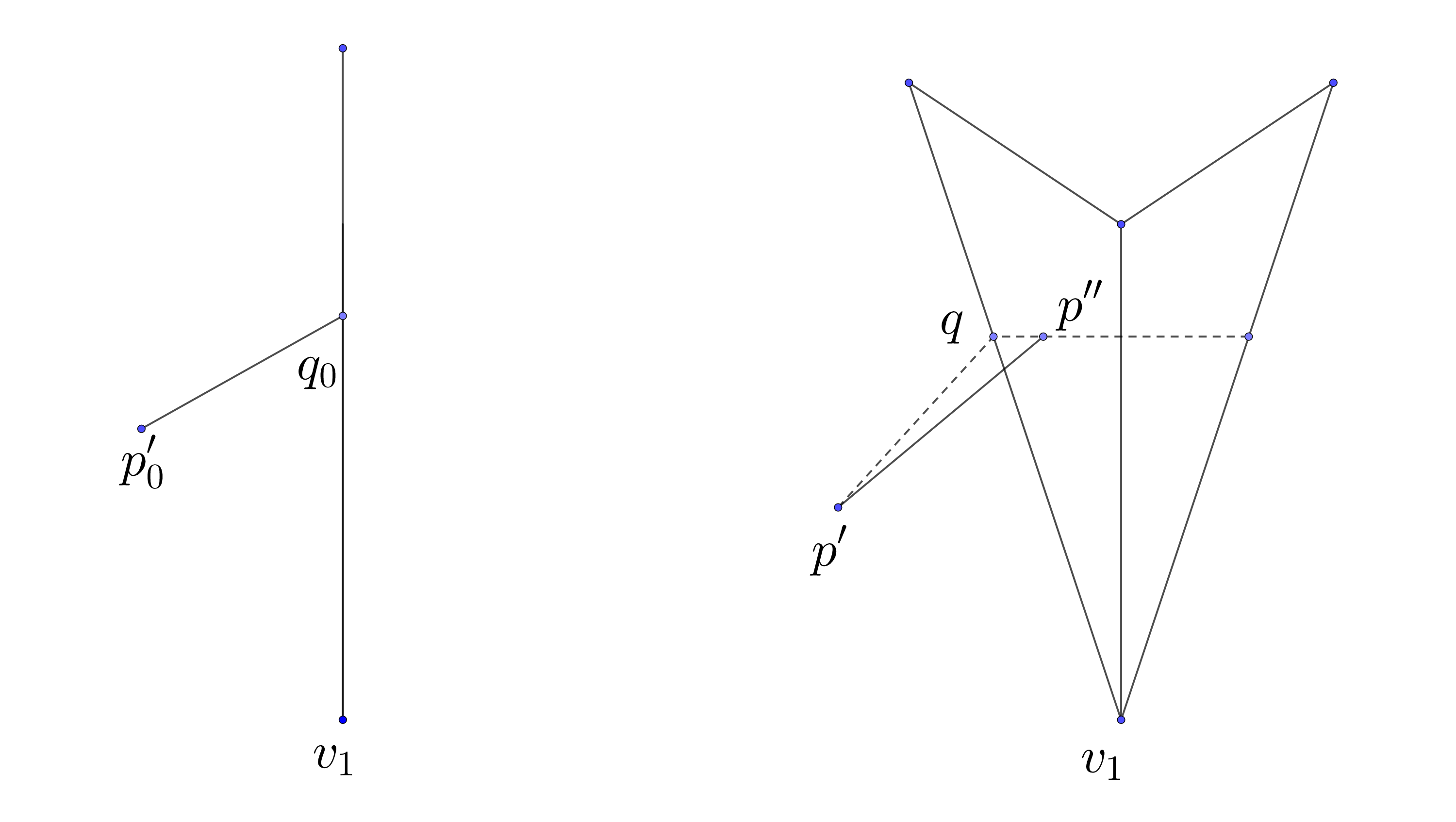}
\caption{To the proof of Claim~\ref{distbound}.}
\label{Pic10-5}
\end{center}
\end{figure}


Now we consider $p', p'' \in T_1 \backslash X_1$. Let $p'_0 = g_1^{-1}(p') \in T_0$ and $p''_0=g_1^{-1}(p'')\in T_0$. If $d_0(p_0', p_0'')\neq d_1(p', p'')$, then all shortest paths connecting $p'_0$ and $p''_0$ in $d_0$ intersect $\chi_0$. Consider an arbitrary shortest path between $p'_0$ and $p''_0$. Let it intersect $\chi_0$ at the point $q$. We consider the points $q', q''\in T_{1}$ at the borders of $X_1$ such that $d_{1}(q', v_1)=d_{1}(q'', v_1)=d_0(q, v_1)$. See Figure~\ref{Pic10}. Then $$d_{1}(p', p'') \leq d_{1}(p', q')+d_{1}(q', q'') + d_{1}(q'', p'') \leq$$ $$\leq d_0(p'_0, p''_0) + B'B'' < \diam(T_0, d_0)+ \nu_{v_1}(d_0)\sinh(\diam(T_0, d_0)).$$

\begin{figure}
\begin{center}
\includegraphics[scale=0.6]{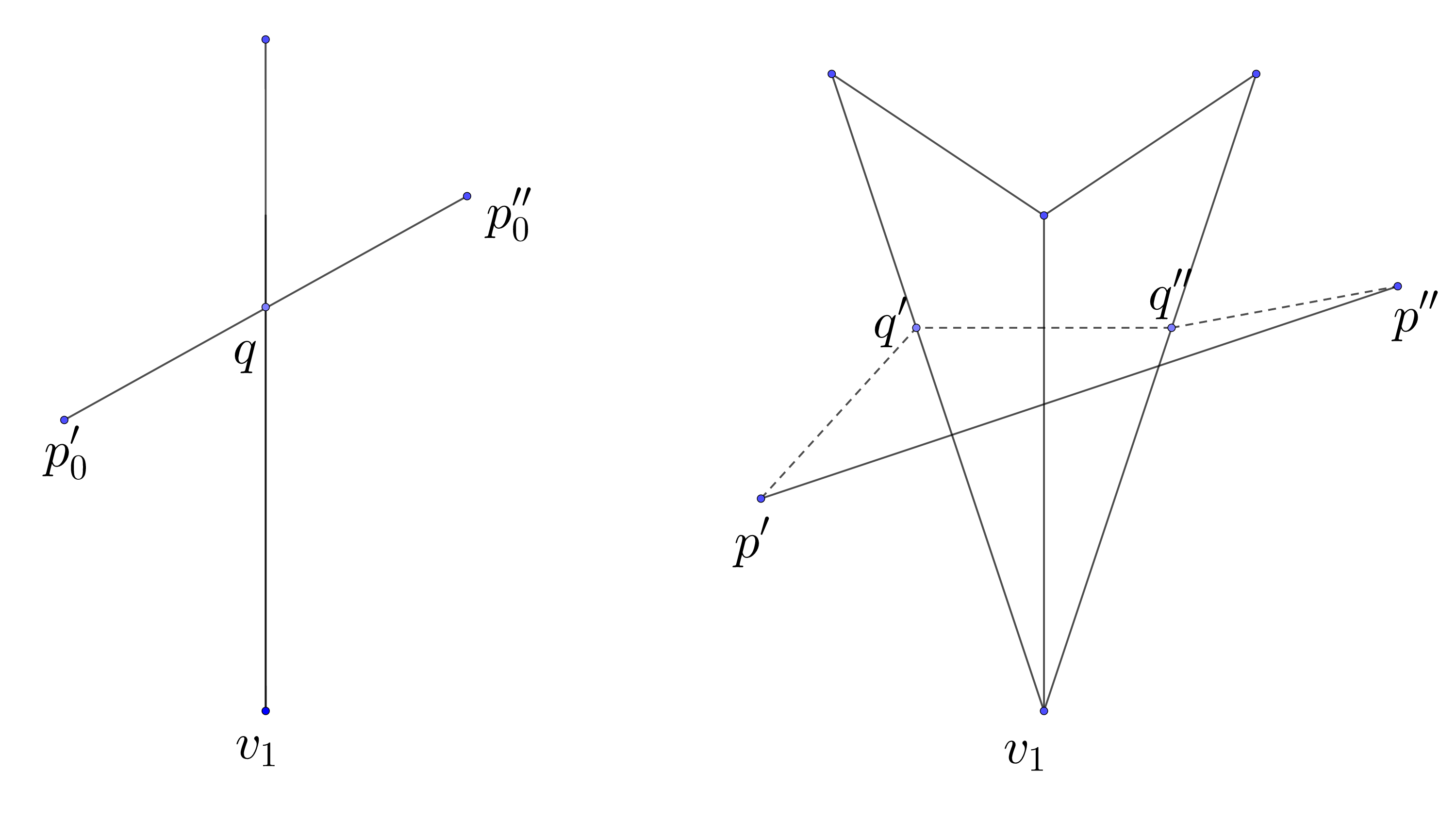}
\caption{To the proof of Claim~\ref{distbound}.}
\label{Pic10}
\end{center}
\end{figure}

This finishes the proof of Claim~\ref{distbound}.
\vskip+0.2cm

\begin{prop}
\label{curvbound}
$\nu_{v_0}(d_0)+\nu_{v_1}(d_0)<\nu_{w_1}(d_1) <\nu_{v_0}(d_0)+\nu_{v_1}(d_0)\cosh(\diam(T_0, d_0)).$
\end{prop} 

Indeed, the lower bound is trivial, the upper is just Lemma~\ref{curvatbound} applied to triangle $ABC$.
\vskip+0.2cm

From Claim~\ref{distbound} and (\ref{choice2}) we see that
$$\diam (T_1, d_1)<\diam(T_0, d_0) + \nu_{v_1}(d_0)\sinh(\diam(T_0, d_0))<\delta_D+\delta_{\nu}\sinh\Delta<\Delta.$$

From Claim~\ref{curvbound} we obtain
$$\nu_{w_1}(d_1)<(\nu_{v_0}(d_0)+\nu_{v_1}(d_0))\cosh\Delta.$$

We continue the process of merging the curvature by gluing new bigons. At every step we take the point that we obtained after the last step, take another conical point and cut the metric along a shortest path connecting both points. Let us check that there always exists a bigon to glue in the cut.

To this purpose we show by induction that after $i$ steps
\begin{equation}
\label{induction}
\diam(T_i, d_i)<\Delta,~~~\nu_{w_i}(d_i)<\left(\sum_{j=1}^i\nu_{v_j}(d_0)\right)\cosh\Delta.
\end{equation}
Assume that these bounds hold. Let $v_{i+1} \in V(d_i)$ be distinct from $w_i$. Connect them by a shortest path $\chi_{i+1}$. Note that $v_{i+1} \in V(d_0)$ and $\nu_{v_{i+1}}(d_i)=\nu_{v_{i+1}}(d_0)$. We have
$$\len(\chi_{i+1}, d_i) = d_i(w_i, v_{i+1})\leq \diam(T_i, d_i) < \Delta,$$
\begin{multline*}
\nu_{w_i}(d_i) < \left(\sum_{j=1}^i\nu_{v_j}(d_0)\right)\cosh\Delta<\\
< \left(\sum_{j=1}^{i+1}\nu_{v_j}(d_0)\right)\cosh\Delta-\nu_{v_{i+1}}(d_0) \leq\\
\leq \nu(T, d_0)\cosh\Delta-\nu_{v_{i+1}}(d_0)<\\
<\delta_{\nu}\cosh\Delta-\nu_{v_{i+1}}(d_0)<2\eta - \nu_{v_{i+1}}(d_0). 
\end{multline*}

Then Lemma~\ref{exist} shows that the bigon $X_{i+1}$ exists to glue with the cut along $\chi_{i+1}$. We glue-in the bigon and obtain a new metric space $(T_{i+1}, d_{i+1})$ and a local isometry $g_{i+1}$. We identify $V(d_i)$ in $T_i$ with its $g_{i+1}$-image in $T_{i+1}$. Let $w_{i+1}$ be the interior vertex of $X_{i+1}$. Then $$V(d_{i+1})=(V(d_i)\cup\{w_{i+1}\} )\backslash \{w_i, v_{i+1}\}.$$
The proofs of Claim~\ref{distincr}--\ref{curvbound} remain to be valid. (To repeat the proof of Claim~\ref{distbound} we use that the boundary segments of $X_{i+1}$ are greater than the distances from boundary vertices to the interior vertex.) From these Claims,~(\ref{choice1}), (\ref{choice2}) and the induction hypothesis we show 
\begin{multline*}
\diam(T_{i+1}, d_{i+1}) <\diam(T_i, d_i)+\nu_{v_{i+1}}(d_i)\sinh(\diam(T_i, d_i))<\\
<\diam(T_{i-1}, d_{i-1})+\nu_{v_i}(d_{i-1})\sinh(\diam(T_{d_{i-1}}, d_{i-1}))+\nu_{v_{i+1}}(d_0)\sinh\Delta<\\
<\diam(T_0, d_0)+\left(\sum_{j=0}^{i+1}\nu_{v_j}(d_0)\right)\sinh\Delta<\\
<\diam(T_0, d_0)+\nu(T,d_0)\sinh\Delta<\Delta,
\end{multline*}
\begin{multline*}
\nu_{w_{i+1}}(d_{i+1})< \nu_{w_i}(d_i)+\nu_{v_{i+1}}(d_i)\cosh(\diam(T_i, d_i))<\\
<\left(\sum_{j=1}^{i}\nu_{v_j}(d_0)\right)\cosh\Delta+\nu_{v_{i+1}}(d_0)\cosh\Delta=\\
=\left(\sum_{j=1}^{i+1}\nu_{v_j}(d_0)\right)\cosh\Delta.
\end{multline*}

Thus, we established the induction step for (\ref{induction}). Finally we obtain a metric space $(T_n, d_n)$ with just one conical point in the interior. As both $T$, $T_n$ are triangles, we may consider $d_n$ as a metric $\hat d$ on $T$. The side lengths and angles do not change during the process. Thus, we obtain a sweep-in $\hat d$ of $d$ on $T$. From (\ref{choice1}) and (\ref{induction}) we get $\diam(T, \hat d)<\Delta$ and 
$$\nu(T, \hat d)<\nu(T, d)\cosh\Delta<\delta_{\nu}\cosh\Delta<2\eta.$$
\end{proof}

Now we make a brief sidestep from metrics on $T$ to metrics on $S_g$. We claim that if we perform the sweeping process described in this section for each triangle of $T$, then the diameter of the resulting metric can be estimated in exactly the same way as for a single triangle. In other words, we state

\begin{lm}
\label{diamglobal}
Let $\Delta$ be a positive number and $\delta_{\nu}, \delta_{D}$ be positive numbers satisfying (\ref{choice1}), (\ref{choice2}). If $(S_g, d)$ is a cone-metric space and $\mathcal T$ is its $\delta$-fine triangulation into convex triangles, where $\delta:=\min\{\delta_{\nu}, \delta_D\}$, then there exists a unique sweep-in $\hat d$ of $d$ with respect to $\mathcal T$. Moreover, $\diam(S_g, \hat d)$ satisfies
$$\diam(S_g, d) \leq \diam(S_g, \hat d)  < \diam(S_g, d) + \nu(S_g, d)\sinh\Delta.$$
\end{lm}

The existence here follows from Lemma~\ref{merge} and the uniqueness will be proven in the next section. Thereby, what we remark now is the diameter estimate. The proof is exactly the same as the proof of similar bounds in Lemma~\ref{merge}. We only need to note that proofs of Claim~\ref{distincr}, Claim~\ref{distbound} from the proof of Lemma~\ref{merge}, which are responsible for these bounds, work without changes if instead of points in $T$ we consider points in $S_g$.

\subsubsection{Swept triangles}
\label{conesec}

\begin{dfn}
Since now we call a metric space $(T, \hat d)$, where $T$ is a triangle and $\hat d$ is a swept metric, a \emph{swept triangle}. It is called \emph{convex} if $\hat d$ is convex. It is called \emph{strict} if $|V(\hat d)|=1$ and the curvature of the conical point is positive. Note that a strict swept triangle might be non-convex: this happens if some boundary angles are greater than $\pi$.
\end{dfn}

Let $(T, \hat d)$ be a swept triangle. We also denote it as $OA_1A_2A_3$, where $A_1, A_2, A_3$ are vertices of $T$ and $O$ is the interior cone-point of $\hat d$. In case $V(\hat d)=0$ we assume that $O$ is a marked point in the interior of $T$, i.e., a swept triangle is always a triangle with a marked point inside. By $l_i$ denote the length of $A_{i-1}A_{i+1}$ (we assume that $A_i$ are enumerated modulo three) and by $\lambda_i$ denote the angle $A_{i-1}A_iA_{i+1}$. We naturally call $l_i$ the \emph{side lengths} of $OA_1A_2A_3$ and $\lambda_i$ its \emph{angles}. There are unique geodesics from $O$ to each vertex. By $x_i$ denote the length of the geodesic $OA_i$ and by $\beta_i$ denote the sector angle $A_{i-1}OA_{i+1}$. One can see that $0<\beta_i < \pi$.

Clearly, a swept triangle is uniquely determined by $x_i$, $\beta_i$ (up to isometry). This gives us a parametrization of the set of strict swept triangles up to isometry by the following subset of $\R^6$:
$$SCT=\{(x_1, x_2, x_3, \beta_1, \beta_2, \beta_3): x_i>0,~0<\beta_i<\pi,~\beta_1+\beta_2+\beta_3 < 2\pi\}.$$

We prove the uniqueness of a sweep-in in Lemmas~\ref{mlIIv}, \ref{mlII}. It is easy to see that if a swept triangle has the side lengths and angles coinciding to the side lengths and angles of a hyperbolic triangle, then it is isometric to this hyperbolic triangle. The rest  follows from

\begin{lm}
\label{conerigid}
A strict swept triangle is uniquely determined by the side lengths and angles (up to isometry)
\end{lm}

\begin{proof}
We show that the distances $x_i$ in a strict swept triangle $OA_1A_2A_3$ are uniquely determined by the side lengths and angles. This also shows that $\beta_i$ are uniquely determined, which implies the claim.

Cut $OA_1A_2A_3$ along the geodesic $OA_1$. Develop it in $\H^2$. Let $A'_1A_2A_3A''_1O$ be the resulting polygon. One can show that it is not self-intersecting even if some boundary angles are greater than $\pi$.

The side lengths $l_1$, $l_2$, $l_3$ and the angles $\lambda_2$ and $\lambda_3$ determine the polygonal line $A'_1A_2A_3A''_1$ in $\H^2$ up to isometry.
The point $O$ belongs to the perpendicular bisector to segment $A'_1A''_1$. We can see that as $O$ moves along this bisector, the sum of angles $OA'_1A_2$ and $OA''_1A_3$ changes monotonously. Then there exists no more than one position where this sum is equal to the angle $\lambda_1$. This finishes the proof.
\end{proof}

Consider the map $\Theta: SCT \rightarrow \R^6$, which associates to a strict swept triangle the 6-tuple of its side lengths and angles, i.e.,
$$\Theta(x_1,x_2,x_3,\beta_1,\beta_2,\beta_3)=(l_1,l_2,l_3,\lambda_1,\lambda_2,\lambda_3).$$
Lemma~\ref{conerigid} states that this map is injective. Clearly it is continuous, thereby by the Brouwer Invariance of Domain Theorem we get

\begin{crl}
\label{homeo}
The image $\Theta(SCT)$ is open and $\Theta$ is a homeomorphism onto its image. 
\end{crl}

The map $\Theta$ will be used prominently in the next subsection.

\subsubsection{Sweeping-in CBB($-1$) metrics}
\label{cbbswept}

As we wrote before, the proof of the existence in Lemma~\ref{mlII} is done via an approximation by cone-metrics. Assume that $(T, d)$ is a triangle with a {CBB($-1$)} metric. Let $(T, d_m)$ be a sequence of convex cone-metrics converging to $(T, d)$. There is a sweep-in $\hat d_m$ for each metric $d_m$. We need to investigate the convergence of swept triangles $(T, \hat d_m)$. It might happen that such a sequence diverges. For instance, it might happen that some $x_i=0$ in the limit. The tough part of the proof is to show that such degeneracies can not happen when $d_m$ converges to a CBB($-1$) metric. 

The approximation can be done as in Subsection~\ref{triang}. To this purpose first we need to realize $(T, d)$ on a convex surface. This can be done with the help of Alexandrov's gluing theorem and Alexandrov's realization theorem:

\begin{thm}[The gluing theorem]
\label{gluingthm}
Let $(M^1, d^1)$ and $(M^2, d^2)$ be two surfaces with boundaries and CBB($-1$) metrics. Let $f: \partial M^1 \rightarrow \partial M^2$ be an isometry (it is possible to consider only some boundary components of $M^1$). Then the gluing of $(M^1, d^1)$ and $(M^2, d^2)$ along $f$ is a CBB($-1$) metric space.
\end{thm}

A proof is identical to the CBB(0) case proven by Alexandrov in~\cite[Section VIII.1]{Ale3}. Another proof in the CBB($k$) case for arbitrary dimension can be found in~\cite{Pet}.

\begin{thm}[The realization theorem]
\label{alex2h}
For every CBB($-1$) metric $d$ on $S^2$ there exists a convex body $G \subset \H^3$ such that $(S^2, d)$ is isometric to the boundary of $G$. 
\end{thm}

The proof is also identical to the CBB(0)-case~\cite[Section VII.1, Theorem 1]{Ale3}. An explanation of differences between CBB(0) and CBB($-1$) cases is given in~\cite[Section XII.2]{Ale3}.

Now consider a convex body $G \subset \H^3$. Let its boundary be identified with a CBB($-1$) metric space $(S^2, d)$. Similarly to Subsection~\ref{triang}, we consider a sequence of $\mu_m$-dense sets $V_m \subset (S^2, d)$ with $\mu_m \rightarrow 0$. Let $G_m$ be the convex hull of $V_m$. Pick a point $o \in \inter (G_m)$ for all $m$. The metrics of $\partial G_m$ can be pulled back to convex cone-metrics $d_m$ on $S^2$ with the help of central projection from $o$. 

We can prove an analogue of Lemma~\ref{convuni} and Lemma~\ref{refin}:

\begin{lm}
\label{refin2}
Let $\mathcal T$ be a strictly short triangulation of $(S^2, d)$. Consider a sequence of finite $\mu_m$-dense sets $V_m\subset (S^2, d)$ with $\mu_m \rightarrow 0$ such that all of them contain $V(\mathcal T)$ and none of them contains a point in the interior of an edge of $\mathcal T$. By $d_m$ denote the metrics on $S^2$ defined as above.

Then $d_m$ converge to $d$ uniformly and the triangulation $\mathcal T$ is realized by infinitely many $d_m$. Moreover, the realizations can be chosen so that they are short and for every triangle $T$ of $\mathcal T$ we have (after taking a subsequence)

(1) $T(d_m) \rightarrow T(d)$;

(2) $\diam(T, d_m) \rightarrow \diam(T, d)$;

(3) $\area(T, d_m) \rightarrow \area(T, d).$
\end{lm}

A proof of Lemma~\ref{refin2} can be done exactly in the same way as our proof of Lemma~\ref{refin}. Instead of distance functions to the lower boundary one should use the distance function to the point $o$. We do not provide here a proof of Lemma~\ref{refin2}. We note that if we want to prove only Main Lemma II, then it is enough to prove the existence of the sweep-in only for metrics coming from strictly short triangles on the upper boundary of a convex compact Fuchsian manifold. In this case we may just use below Lemma~\ref{refin} instead of Lemma~\ref{refin2}. However, in a more general Lemma~\ref{mlII} we do not know if the triangle comes this way and we use Lemma~\ref{refin2} to overcome this. 

Now we turn to a triangle $T$ with a strictly short CBB($-1$) metric $d$. Take two copies of $(T, d)$ and glue along the boundary. We obtain a metric $\tilde d$ on the 2-sphere $S^2$. Theorem~\ref{gluingthm} implies that it is CBB($-1$). Two copies of $(T, d)$ in $(S^2, \tilde d)$ form a strictly short triangulation of the latter. We fix an inclusion $\rho: (T, d) \rightarrow (S^2, \tilde d)$. 

By Theorem~\ref{alex2h}, there exists a convex body $\tilde G \subset \H^3$ with boundary isometric to $(S^2, \tilde d)$. In order to prove Lemma~\ref{mlII} we need to consider triangulations of $(T, d)$. Let $\mathcal T$ be a strictly short triangulation of $(T, d)$ with no additional vertices at the edges of $T$. Transfer it to $(S^2, \tilde d)$ via $\rho$. Together with the second copy of $(T, d)$ in $(S^2, \tilde d)$ it forms a strictly short triangulation $\tilde{\mathcal T}$ of $(S^2, \tilde d)$. Now we can approximate $(S^2, \tilde d)$ by cone-metrics realizing $\tilde{\mathcal T}$ with the help of Lemma~\ref{refin2}. After cutting off from these metrics the redundant triangle of $\tilde{\mathcal T}$, we obtain approximations of $(T, d)$ realizing $\mathcal T$. In other words, we get

\begin{crl}
\label{refin3}
Let $\mathcal T$ be a strictly short triangulation of $(T,  d)$. There exists a sequence of convex cone-metrics $d_m$ on $T$ such that they converge uniformly to $d$, all of them realize $\mathcal T$ by short triangulations and for every triangle $T'$ of $\mathcal T$ we have

(1) $T'(d_m) \rightarrow T'(d)$;

(2) $\diam(T', d_m) \rightarrow \diam(T', d)$;

(3) $\area(T', d_m) \rightarrow \area(T', d).$
\end{crl}


Now we are ready to prove Lemma~\ref{mlII}.

\begin{proof}[Proof of Lemma~\ref{mlII}]
It remains only to prove the existence for metrics that are not cone-metrics. 

Let $\Delta$ be an arbitrary positive number and $\eta=\eta(\Delta)$ be from Lemma~\ref{exist}. Take $\delta_D$, $\delta_{\nu}$ satisfying~(\ref{choice1}) and~(\ref{choice2}). Define $\delta$ as the minimum of $\delta_D$, $\delta_{\nu}$. Let $d$ be a short CBB($-1$) metric on $T$ with $\diam(T, d)<\delta$, $\nu(T,d)<\delta$ as in the statement of Lemma~\ref{mlII}. If $T(d)$ happens to be the 6-tuple of side lengths and angles of a hyperbolic triangle, then we just let $\hat d$ to be the metric of this triangle. Consider the other case.

Corollary~\ref{refin3} shows that there exists a sequence of convex cone-metrics $d_m$ on $T$ such that $\diam(T, d_m)<\delta_D$, $\nu(T, d_m)<\delta_{\nu}$ and $T(d_m) \rightarrow T(d)$. Moreover, we can demand that $T(d_m)$ are not 6-tuples of side lengths and angles of hyperbolic triangles.

For each $m$ by Lemma~\ref{merge} we construct a convex cone-metric $\hat d_m$ on $T$ such that $(T, \hat d_m)$ is a convex strict swept triangle and $T(\hat d_m) = T(d_m)$, so $T(d_m) \in \Theta(SCT)$. (Recall the definitions of $SCT \subset \R^6$ and the map $\Theta$ from Subsection~\ref{conesec}.) It follows that  $T(d) \in \overline{\Theta(SCT)}.$ We need to understand what happens in the limit. Define $Z_m := \Theta^{-1}(T(d_m))$. By $x_{m, i}$ and $\beta_{m, i}$ we denote the respective components of $Z_m$.

We claim that after taking a subsequence $Z_m$ converges in $\R^6$. Indeed, otherwise for some $i$ we have $x_{m,i} \rightarrow \infty$. This can not happen because Lemma~\ref{merge} says that $\diam(T, \hat d_m) < \Delta$. Let $Z_m \rightarrow Z \in \overline{SCT}.$ By $x_i$, $\beta_i$ we denote the respective components of $Z$. If $Z \in SCT$, then $T(d) = \Theta(Z)$, so there exists a convex swept triangle with side lengths and angles equal to those of $(T, d)$ and the proof is finished. It remains to assume that $Z \in \partial SCT$ and consider different types of degeneracies that can happen with $Z$. We divide them in the following cases:

(i) $\beta_i=0$ for some $i$;

(ii) $x_i=0$ for at least two $i$;

(iii) $x_i =0$ for one $i$; 

(iv) $\beta_1 +\beta_2 +\beta_3  = 2\pi$ and $x_i >0$ for all $i$; 

(v) $\beta_i =\pi$ for some $i$, but $\beta_1 +\beta_2 +\beta_3  < 2\pi$ and $x_i >0$ for all~$i$. 

Case (iii) turns out to be the most complicated. So first we rule out the four others.

(i) Let $\beta_1=0$. Consider a strict swept triangle  $OA_1A_2A_3$ like in Subsection~\ref{conesec} with $OA_1, OA_2 < \Delta$. Then for every $\xi>0$ there exists $\mu >0$ such that if $\angle A_1OA_2 \leq \mu$, then either $A_1A_2<\xi$ or $\angle OA_1A_2 > \pi-\xi$ or $\angle OA_2A_1 > \pi-\xi$. Now take $\xi$ smaller than the length of any edge of $(T, d)$ and than the difference of $\pi$ and any angle of $(T, d)$. Choose $\mu$ for this $\xi$. For sufficiently large $m$ we see that all $\beta_{m,1} > \mu$. Thus, $\beta_1 > 0$.

(ii) Let $x_1=x_2=0$. Consider a strict swept triangle  $OA_1A_2A_3$. For every $\xi>0$ there exists $\mu >0$ such that if $OA_1, OA_2 \leq \mu$, then $A_1A_2 \leq \xi$. Take $\xi$ smaller than the length of any edge of $(T, d)$. Choose $\mu$ for this $\xi$. For sufficiently large $m$ we see that either $x_{m, 1}$ or $x_{m, 2}$ is greater than $\mu>0$. Thus, either $x_1$ or $x_2$ is greater than 0.

(iv) In this case $T(d)$ has the side lengths and angles of a hyperbolic triangle.

(v) Let $\beta_1=\pi$. Consider the convex quadrilateral with sides $x_2$, $x_3$, $l_2$, $l_3$ and angle $\beta_2+\beta_3 < \pi$ between the sides of lengths $x_2$, $x_3$. The angle $\lambda_1$ of $(T, d)$ is equal to the angle of this quadrilateral between the sides of lengths $l_2$, $l_3$. Then $\lambda_1$ is smaller than the respective angle in the comparison hyperbolic triangle with lengths $l_1=x_2+x_3$, $l_2$ and $l_3$. This contradicts to the Toponogov globalization theorem (see Section~\ref{cbb-1}) and the condition that sides of $(T, d)$ are shortest paths.

(iii) If $\nu(T,d)=0$, then $T(d)$ has the side lengths and angles of a hyperbolic triangle. This can be seen from approximating by cone-metrics as above and noting that $2\pi-\beta_{m,1}-\beta_{m,2}-\beta_{m,3}$ converges to $\nu(T,d)$. So from now assume that $\nu(T,d)>0$. The idea is the following. If $x_1=0$ and $d'$ is any cone-metric sufficiently close to $d$, then the respective component $x_1'$ of $\Theta^{-1}(T(\hat d'))$ is sufficiently small. However, we are going to construct arbitrarily close convex cone-metrics to $d$, for which $x_1'$ is bounded away from zero. This is done by paying the attention to how we do the approximation. Subdivide $(T, d)$ into several triangles so that one of them contains almost all curvature and is not adjacent to the vertex $v_1$ of $T$. Then we find a cone-metric $d'$ close to $d$ realizing this new triangulation. When we merge the curvature inside $T$ in $d'$, the resulting swept triangle does not depend on the process of merging. Hence, we can first merge all curvature inside ``the biggest part'' of $T$ and then all other curvature. We can show that the resulting conical point after the final merging is not close to $v_1$.

Now we proceed to the details. 
Let $x_1=0$. Denote the vertices of $T$ by $v_1$, $v_2$, $v_3$ and the respective edges by $e_1$, $e_2$, $e_3$. 
Fix $t>0$ and choose points $p_2 \in e_3$ and $p_3 \in e_2$ at distance $t$ from $v_1$. Connect $p_2$ with $v_3$ and $p_3$ with $v_2$ by shortest paths. If there are several shortest paths from $p_2$ to $v_3$, then we choose the one that makes the smallest angle with $e_3$ at $v_3$ (it is easy to see that this determines the shortest path uniquely). Note that except endpoints it is contained in the interior of $T$ as otherwise it has to be contained in the boundary and then the angle of $T$ at $v_1$ is $\pi$. Similarly choose the shortest path from $p_3$ to $v_2$.

Let the chosen paths intersect in the point $v$. By $\psi_2$ and $\psi_3$ we denote the segments from $v$ to $v_2$ and $v_3$ respectively. See Figure~\ref{Pic11-5}. It clear from the definitions that geodesics in a CBB($-1$) metric do not pass through points $p$ with $\nu_p(d)\neq 0$. Therefore, $\nu_v(d)=0$.  Let $U(t)$ be the (open) region bounded by $e_2$, $e_3$, $\psi_2$, $\psi_3$ and $T'(t)$ be the (open) triangle bounded by $e_1$, $\psi_2$ and $\psi_3$. Define $$\nu'(t):=\nu(T'(t), d),~~~\nu''(t):=\nu(U(t), d).$$
We have $\nu'(t)+\nu''(t)=\nu(T, d)$. If $t_1<t_2$, then $U(t_1) \subset U(t_2)$. The intersection of all $U(t)$ is empty. Thus, one can see that $\nu''(t)$ tends to zero as $t$ tends to zero. We choose $t_0 > 0$ such that $$\frac{2\sinh\Delta}{\Delta}\cdot\frac{\tan(\nu''(t_0)\cosh\Delta)}{\nu'(t_0)} < \frac{1}{2}.$$

\begin{figure}
\begin{center}
\includegraphics[scale=0.3]{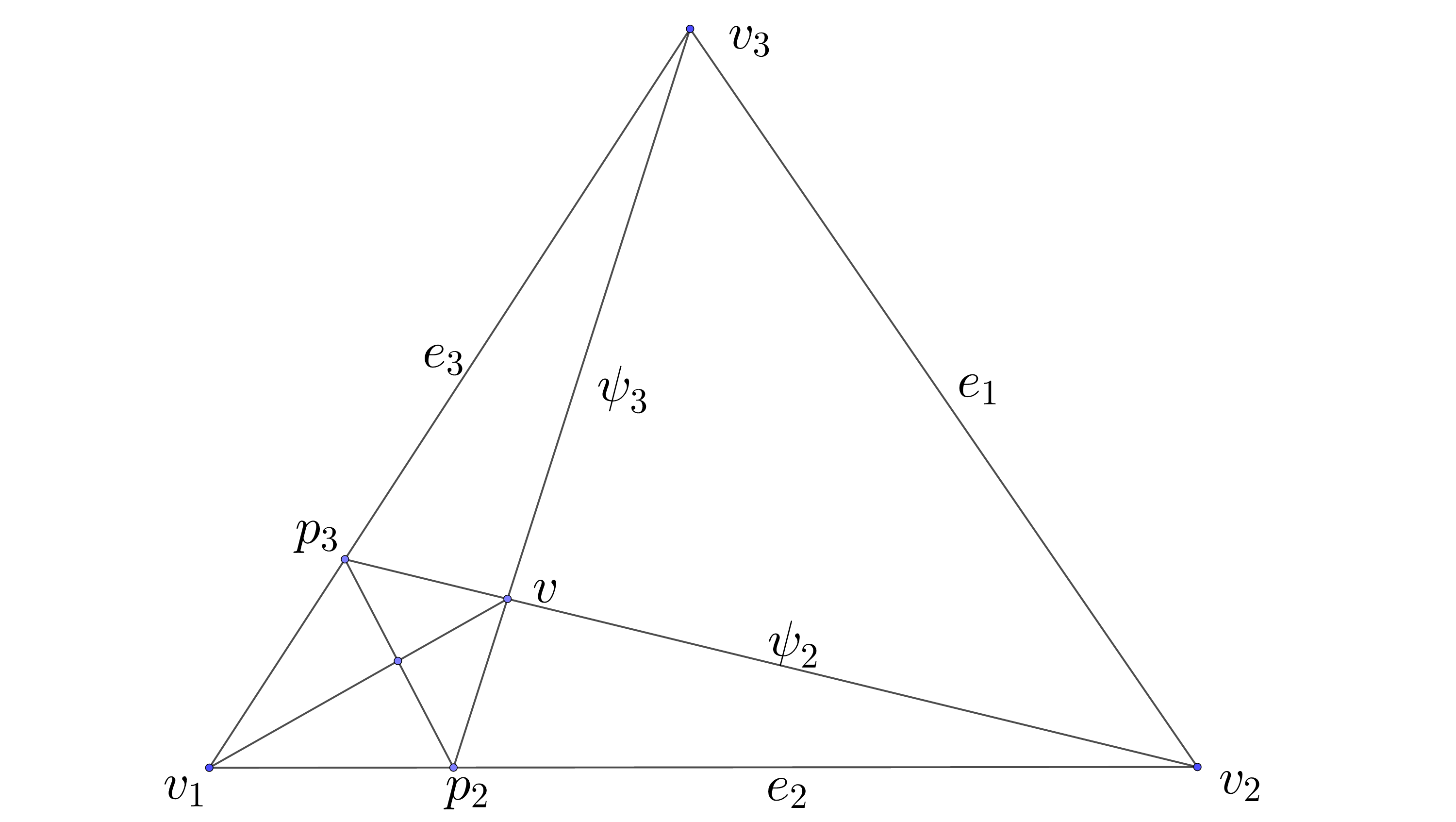}
\caption{Positions of points.}
\label{Pic11-5}
\end{center}
\end{figure}

From now on we denote $T'(t_0)$, $U(t_0)$, $\nu'(t_0)$ and $\nu''(t_0)$ by $T'$, $U$, $\nu'$ and $\nu''$ . Also we define $$ y := \inf_{p \in T'} d(v_1, p).$$

All sides of the triangle $T'$ are unique shortest paths due to the strong non-overlapping property. It is very easy to extend our current decomposition of $T$ to a strictly short triangulation. To this purpose we just draw two diagonals in the quadrangle $v_1p_2vp_3$. From the strong non-overlapping property, all edges of the obtained triangulation are unique shortest arcs. Therefore, we obtain a strictly short triangulation $\mathcal T$ of $T$ containing $T'$ (angles are less than $\pi$ as no two geodesics can form zero angle, the convexity of triangles follows from the convexity of $(T,d)$). Due to Corollary~\ref{refin3} there exists a sequence of convex cone-metrics $\{d_n\}$ on $T$ converging to $d$ uniformly and realizing $\mathcal T$. Define $$\nu'_n:=\nu(T', d_n),~~~\nu''_n:=\nu(U,d_n),~~~y_n := \inf_{p \in T'}  d_n(v_1, p).$$
Here by $U$, $T'$ we mean their realizations in $d_n$. From Corollary~\ref{refin3} we can choose $\{d_n\}$ so that $\nu'_n \rightarrow {\nu'}$, $\nu''_n \rightarrow  \nu''$ and $T(d_n)\rightarrow T(d)$ as $n$ goes to infinity. From the uniform convergence one can get $y_n \rightarrow y$.

Now choose $n$ such that 
$$\diam(T, d_n) < \delta_D,$$ 
$$\nu(T, d_n)<\delta_{\nu},$$ 
$$y_n >  \frac{y}{2},$$ 
$$\frac{2\sinh\Delta}{\Delta}\cdot\frac{\tan(\nu''_n\cosh\Delta)}{\nu'_n} < \frac{1}{2},$$ 
$$\lambda_{1, n}< \frac{\lambda_1+\pi}{2}<\pi,$$
\begin{equation}
\label{choosexn}
x_{n,1} < \min\left\{\frac{y}{4}, \arctanh\left(\cos\frac{\lambda_1+\pi}{2}\cdot\tanh \frac{y}{2} \right)\right\}.
\end{equation}




Here $\lambda_{n,1}$, $\lambda_1$ are angles at the vertex $v_1$ in $(T, d_n)$, $(T, d)$ respectively and $x_{n,1}$ is the distance from $v_1$ to the cone point in the sweep-in of $(T, d_n)$. Since now we denote $d_n$, $y_n$, $\nu'_n$ and $\nu''_n$ by $\tilde d$, $\tilde y$, $\tilde \nu'$ and $\tilde \nu''$.

By Lemma~\ref{conerigid} we see that if we merge all curvature inside $(T, \tilde d)$ in one point (by gluing bigons), then the resulting swept triangle does not depend on the choice of a merging algorithm. Lemma~\ref{merge} says that its diameter is smaller than $\Delta$. Claim~\ref{diambound} of Lemma~\ref{merge} shows that after each merging the diameter of $T$ is not decreasing, therefore $\Delta$ is an upper bound on the diameter of $T$ during the whole process. 

First, we merge all curvature inside $T'$ to point $w'$. By Claim~\ref{distincr} of Lemma~\ref{merge} the distance from $w'$ to $v_1$ is at least $y$. Now we merge all curvature in $U$ to a point $w''$. We note that $U$ is not convex, hence we do not expect that $w''$ stays in $U$. Let $ d_*$ be the resulting cone-metric on $T$. We would like to glue a bigon along the shortest path $\psi$ connecting $w'$ and $w''$. From Claim~\ref{curvbound} of Lemma~\ref{merge} we see that $\tilde\nu'\leq\nu_{w'}( d_*)<\tilde\nu'\cosh\Delta$ and $\tilde\nu''\leq\nu_{w''}( d_*)<\tilde\nu''\cosh\Delta$. Then from inequalities~(\ref{choice1}) we get $\nu_{w'}( d_*)+\nu_{w''}( d_*)<2\eta$. Together with $\len(\psi,  d_*)<\Delta$, it follows from Lemma~\ref{exist} that the last bigon exists.

Let $\hat w$ be the resulting curvature point, $\hat  d_*$ be the resulting metric and $w$ be the base of perpendicular from $\hat w$ to the boundary segment of the glued bigon that is intersected by a shortest path from $\hat w$ to $v_1$. See Figure~\ref{Pic11}. Then $\hat  d_*(v_1, \hat w) \geq  d_*(v_1, w)$. 

\begin{figure}
\begin{center}
\includegraphics[scale=0.6]{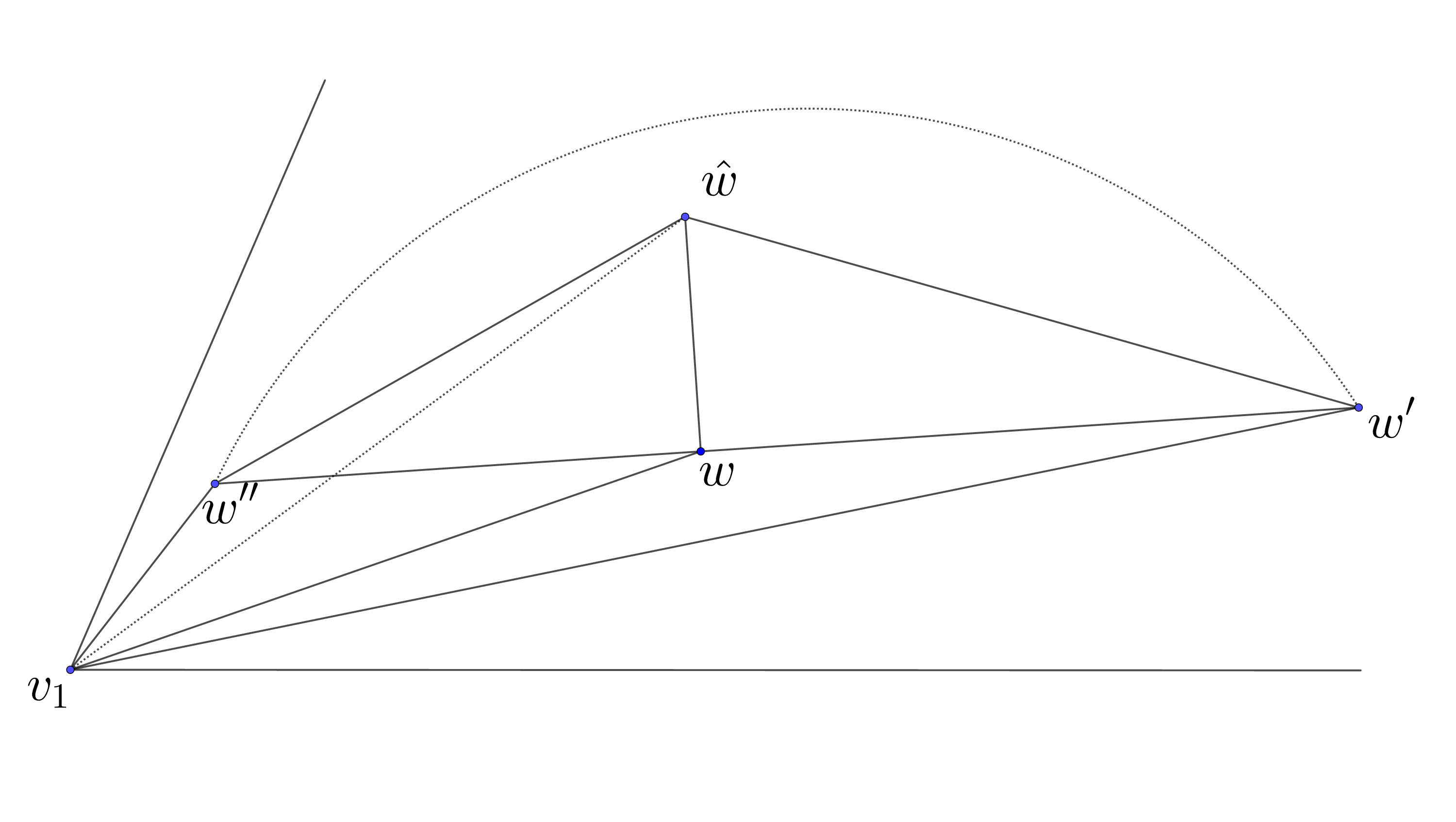}
\caption{Positions of points.}
\label{Pic11}
\end{center}
\end{figure}

Note that $\hat  d_*(v_1, \hat w)=x_{n,1}.$ By the choice of $n$ there is the upper bound~(\ref{choosexn}) on $x_{n,1}$. We are going to obtain a contradiction with this bound. We consider two cases. First, suppose that $ d_*(v_1, w'') \geq y/2$. Then we claim that 
$$\tanh d_*(v_1, w) \geq \cos\frac{\lambda_1+\pi}{2}\cdot\tanh \frac{y}{2} .$$
Indeed, elementary trigonometry shows that for each hyperbolic triangle $ABC$ with $\angle A \leq \frac{\lambda_1+\pi}{2}<\pi$ and $AB$, $AC \geq \frac{y}{2}$ we have the hyperbolic tangent of the distance from $A$ to the segment $BC$ is at least $\cos\frac{\lambda_1+\pi}{2}\cdot\tanh \frac{y}{2} $. We apply this to the triangle $v_1w'w''$. This contradicts to~(\ref{choosexn}).

In the other case, we suppose $ d_*(v_1, w'') \leq y/2\leq d_*(v_1, w')$. Write $$ d_*(v_1, w) \geq  d_*(v_1, w') -  d_*(w', w).$$ Now we bound $ d_*(w', w)$. Recall that $\nu_{w'}( d_*)\geq \tilde \nu'$.

$$d_*(w',w)\leq \sinh( d_*(w',w)) = \sinh( d_*(w, w'')) \frac{\tan(\nu_{w''}( d_*)/2)}{\tan(\nu_{w'}( d_*)/2)} \leq $$ $$\leq
\sinh( d_*(w', w'')) \frac{\tan(\nu_{w''}( d_*)/2)}{\tan(\tilde\nu'/2)}\leq  d_*(w', w'')\cdot\frac{\sinh\Delta}{\Delta} \cdot\frac{\tan(\nu_{w''}( d_*))}{\tilde \nu'}.$$

Recall that $\nu_{w''}( d_*)<\tilde \nu'' \cosh\Delta.$

Also $ d_*(w', w'') \leq  d_*(v_1, w') + d_*(v_1, w'') \leq 2 d_*(v_1, w') $ due to our assumption.

In total we get 
$$
 d_*(v_1, w) \geq   d_*(v_1, w') -  d_*(w', w)\geq$$ 
 $$
 \geq d_*(v_1, w') \left(1-\frac{2\sinh\Delta}{\Delta}\cdot\frac{\tan(\tilde\nu''\cosh\Delta)}{\tilde\nu'} \right)>\frac{y}{4}>x_{n,1}.
$$

This gives the last contradiction with $x_{n,1}=\hat  d_*(v_1, \hat w) \geq  d_*(v_1, w)$. It follows that $x_1>0$ and the proof is over.
\end{proof}
\newpage
\subsection{Proof of Main Lemmas III}

The rough idea of our proofs of Main Lemmas IIIA-C is as follows. Assume that we have a path of metrics $d_t$ parametrized by $t \in [0; \tau]$ and a convex Fuchsian cone-manifold $P_0=P(d_0, V, h_0)$. We want to transform this cone-manifold continuously along the path $d_t$ so that (1) heights of the vertices can move only in one direction (up); (2) we control the change of discrete curvature $S$ (for our purposes it is enough to control it also only in one direction). 
There are two aspects impacting the discrete curvature: the change of heights and the change of the upper boundary metric. A great observation is that we can move up only vertices with positive curvature. This always increases $S$ and it remains to bound the change of $S$ under the metric change. Another idea is that we may consider separately two types of small deformations of $P_t$: (1) transform heights with fixed boundary metric (Section~\ref{height}); (2) transform boundary metric with fixed heights. Our transformation is decomposed virtually into a discrete process. A ``continuous induction'' argument (Lemma~\ref{collect}) allows us to combine this into a continuous one. We note that for technical reasons we can control efficiently only heights of vertices with $\nu_v(d)\neq 0$. Thereby, we pay plenty of attention to the change of the singularity set $V(d_t)$. 

Last, we make remarks on metric paths $d_t$. For proofs of Main Lemmas IIIA-B we turn the operations of gluing/cutting a bigon from the proof of Lemma~\ref{merge} into continuous processes. This is done in Subsection~\ref{gluing}. For our proof of Main Lemma IIIC we introduce a process of dissolving/creating cone-points of small curvature. This is described in Subsection~\ref{flat}. The rest of the proof of Main Lemma IIIC follows from an indirect transformation based on a study of local properties of the space $\mathfrak H(\mathfrak D_c(V))$.

\subsubsection{Changing heights with fixed boundary}
\label{height}

Recall that $\mathfrak D_{sc}(V)$ and $\mathfrak D_{c}(V)$ denote the sets of convex cone-metrics on $S_g$ with $V(d)=V$ and with $V(d)\subseteq V$ respectively.

\begin{lm}
\label{open}
Let $P=P(d, V, h)$ be a convex Fuchsian cone-manifold with $d\in \mathfrak D_{sc}(V)$. Define $\mathfrak N(d, \zeta) \subset \mathfrak D_{c}(V)$ as the set of convex cone-metrics $d'$ such that for every $u, v \in V$, $$|d'(u,v) - d(u, v)|<\zeta.$$ Then there exist $\zeta >0$ and $h' \in H(d, V)$ such that:

(1) $h'_v \geq  h_v$ for every $v \in V$;

(2) $S(d, V, h') \geq S(d, V, h)$;

(3) $h' \in H(d', V)$ for every metric $d' \in \mathfrak N(d, \zeta)$.
\end{lm}

\begin{proof}
Consider the face decomposition $\mathcal R=\mathcal R(P)$. By a triangle of $\mathcal R$ we mean a triangle of a triangulation $\mathcal T$ compatible with $P$ such that $V=V(\mathcal T)$. Recall that due to Lemma~\ref{finite}, there are finitely many triangulations compatible with $P$ and, therefore, finitely many triangles of $\mathcal R$.

First we see that in some cases we can simply take $h'=h$. Consider a triangulation $\mathcal T$ compatible with $P$ and the Fuchsian cone-manifold $P'=P(d', \mathcal T, h)$. It might happen that $P'$ is not convex. Let us use the flip algorithm described in Lemma~\ref{increase}: take an arbitrary concave edge and flip it. Provided $\zeta$ is sufficiently small, strict edges of $P$ remain strict in all intermediate cone-manifolds during the algorithm and, therefore, can not be flipped. We refer to the proof of Lemma~\ref{increase} for details. Thus, all triangulations that appear are compatible with $P$. We also recall that the extension $\widetilde {h}$ of the height function is pointwise non-decreasing during the process and is increased at least in one point. Hence, the algorithm can not run infinitely.

If every concave edge can always be flipped, then the algorithm terminates with a triangulation $\mathcal T'$ such that the Fuchsian cone-manifold $P''=P(d', \mathcal T', h)$ is convex. Recall from the proof of Lemma~\ref{increase} that there are two types of edges that can not be flipped. Motivated by this we define the notion of \emph{obstructive} vertices. 

(a) A vertex $v \in V$ is called \emph{obstructive of type I} if $v$ is isolated and there exists a triangle $T$ of $\mathcal R$ with two edges adjacent to $v$ glued together. We say that $T$ is \emph{associated} to $v$.

(b) A vertex $v \in V$ is called \emph{obstructive of type II} if there exist two triangles $T_1$ and $T_2$ of $\mathcal R$ sharing a flat edge incident to $v$ and the total angle of $v$ at this edge is at least $\pi$ (it might happen that this edge is incident to $v$ twice, then we count only one incidence). We say that $T_1$ and $T_2$ are \emph{associated} to $v$.

Otherwise a vertex is called \emph{non-obstructive}. Suppose that all vertices are non-obstructive in $P$. It is possible to choose $\zeta$ small enough so in any $P'=P(d', \mathcal T, h)$ all vertices remain non-obstructive. The discussion above shows that in this case any concave edge can be flipped, the flip algorithm terminates and we can take $h'=h$.

Now we assume that there are obstructive vertices. Our aim is to find $h' \in H(d, V)$ such that properties (1) and (2) from the statement of Lemma~\ref{open} are satisfied and there are no obstructive vertices in $P'=P(d, V, h')$. 

There might be plenty of triangles associated to an obstructive vertex. We remark that being obstructive of type II is a stronger property than having an angle at least $\pi$ in a face of $\mathcal R$. Indeed, take a hyperbolic polygon $A_1A_2B_2B_1$ with $A_1B_1=A_2B_2$, $\angle A_1=\angle A_2 > \pi/2$ and $\angle B_1 = \angle B_2 < \pi/2$. Glue the sides $A_1B_1$ and $A_2B_2$. The angle of vertex $A$ obtained from $A_1$ and $A_2$ is greater than $\pi$. However, one can see that it is not obstructive of type II, i.e., for any triangulation of the respective surface, any two adjacent triangles have total angle smaller than $\pi$ at $v$.

Because an obstructive vertex $v$ of type I is isolated and $V=V(d)$, Lemma~\ref{isolvertex} implies that $\kappa_v(P)>0$ (one can see even more that $\kappa_v(P)>\pi$ in this case). Lemma~\ref{negcurv} shows that $\kappa_v(P)>0$ also for an obstructive vertex $v$ of type II.

\begin{prop}
\label{plainloop}
Let $v$ be obstructive of type I, $T$ be the associated triangle and $u$ be the other vertex of $T$. Then $h_u > h_v$.
\end{prop}

Indeed, develop the prism containing $T$ as the prism $A_1A_2A_3B_1B_2B_3$ in $\H^3$, where $A_1$ corresponds to $v$, $A_2$ and $A_3$ correspond to $u$. 

The lengths $A_1A_2$ and $A_1A_3$ are equal. The distance to the lower plane is the same at the corresponding points of $A_1A_2$ and $A_1A_3$. It implies that the prism has a plane symmetry with respect to the plane orthogonal to $B_1B_2B_3$, passing through $A_1$ and the midpoint of $A_2A_3$. Then the dihedral angles of $A_1A_2$ and $A_1A_3$ are equal to $\pi/2$. Hence, $A_1$ is the closest point from the upper plane to the lower. From this we see that $h_u=A_2B_2>A_1B_1=h_v$.
\vskip+0.2cm

\begin{prop}
Let $v$ be obstructive of type II and $T_1$, $T_2$ be the respective associated triangles. Then there exists a vertex $u$ of $T_1$ or $T_2$ such that $h_u > h_v$.
\end{prop}

Indeed, develop the prisms containing $T_1$ and $T_2$ as a polyhedron in $\H^3$. Let $A$ be the vertex of this polyhedron corresponding to $v$. It has angle at the upper face at least $\pi$, therefore it lies in the convex hull of other upper boundary vertices. By the convexity of the distance function, one boundary vertex has height strictly greater than $h_v$.
\vskip+0.2cm

Together these claims imply

\begin{prop}
\label{obstr}
There exists a triangle associated to an obstructive vertex and containing a non-obstructive vertex.
\end{prop}

Indeed, consider the set of pairs $(v, u)$, where $v$ is obstructive and $u$ belongs to a triangle associated with $v$. Take a pair $(v, u)$ among them such that $u$ has the maximum height. Then $u$ is non-obstructive.
\vskip+0.2cm

Now we describe a change of heights satisfying the conditions (1) and (2) of the Lemma such that all strict edges remain strict and at least one new strict edge appears. This new edge is a flat edge of $P$. There are finitely many of them, therefore, the process terminates with a Fuchsian cone-manifold without obstructive vertices.

The deformation is the same as described in Lemma~\ref{increase}, but it is applied now to the set of obstructive vertices, which is a subset of vertices with positive curvature. Namely, for $\xi>0$ define $h'$ by $\sinh h'_v:=e^{\xi}\sinh h_v$ if $v$ is obstructive and $h'_v:=h_v$ otherwise. One can see that the argument from the proof of Lemma~\ref{increase} works without modifications to show that for small enough $\xi$ we get $h' \in H(d, V)$. Since $\kappa_v(P)>0$ for every obstructive $v$, this deformation increases $S$.
 
Consider $P'=P(d, V, h')$. It is clear that for sufficiently small $\xi$, strict edges of $P$ remain strict in $P'$ and non-obstructive vertices remain non-obstructive. 
By Claim~\ref{obstr} there is a triangle $T$ associated to an obstructive vertex that contains a non-obstructive vertex. In the proof of Lemma~\ref{increase} it is shown that a flat edge of $T$ necessarily becomes strict.
\end{proof}

The assumption $d \in \mathfrak D_{sc}(V)$ was essential in our proof of Lemma~\ref{open}, but we need also to deal with $d \in \mathfrak D_c(V)$. In this case we can guarantee the following weaker version:

\begin{lm}
\label{openV}
Let $P=P(d, V, h)$ be a convex Fuchsian cone-manifold, $d \in \mathfrak D_c(V)$. Define $\mathfrak N(d, \zeta) \subset\mathfrak D_c(V)$ as the set of convex cone-metrics $d'$ such that  for every $u, v \in V,$ $$|d'(u,v) - d(u, v)|<\zeta.$$ Then for every $\mu>0$ there exist $\zeta >0$ and $h' \in H(d, V)$ such that:

(1) $ h'_v \geq  h_v$ for every $v \in V$;

(2) $S(d, V, h') \geq S(d, V, h)-\mu$;

(3) for every metric $d' \in \mathfrak N(d, \zeta)$ we have $h' \in H(d', V)$.
\end{lm}

\begin{proof}
The proof is the same to Lemma~\ref{open}. There is one thing that goes wrong: there can be obstructive vertices of type II with $\kappa_v(P)=\nu_v(d)=0$. This tiny aspect appears to be important: when we move such a vertex up, $\kappa_v$ becomes negative due to Lemma~\ref{2pi} below. Thus, moving up such vertices may produce an unpredictable negative impact on $S$. Nevertheless, we can control it. At each step of moving obstructive points higher (see the end of the proof of Lemma~\ref{open}) we can choose $\xi$ small enough such that if $S$ decreases, then it decreases by arbitrarily small amount. We have an upper bound for the number of steps (all edges that appear are flat edges of $\mathcal R$ and there are finitely many of them), therefore the decrement of $S$ can be kept under $\mu$.
\end{proof}

Now we consider deformations of the boundary metric.

\subsubsection{Gluing a bigon continuously}
\label{gluing}

We return to the operation described in the proof of Lemma~\ref{merge}: modifying a cone-metric by gluing a bigon. We turn it into a continuous process. After this we will learn how to transform a Fuchsian cone-manifold along with gluing a bigon to its upper boundary.

Let $d_0 \in \mathfrak D_{sc}(V)$, $u, v \in V$,
$\chi$ be a shortest path between $u$ and $v$ and
\begin{align}
\label{conditions}
\begin{aligned}
\len(\chi, d_0)=d_0(u,v) < \Delta, \\
\nu_u(d_0) + \nu_v(d_0) < 2\eta.
\end{aligned}
\end{align}

Here $\Delta>0$ is arbitrary and $\eta=\eta(\Delta)$ is taken from Lemma~\ref{exist}. On $\chi$ choose a point $w$ such that $$\sinh(d_0(u,w))\tan(\nu_u(d_0)/2)=\sinh(d_0(w,v))\tan(\nu_v(d_0)/2).$$ Define $V' = V \cup \{w\}$. Now we consider $d_0$ as an element of $\mathfrak D_c(V')$.

As in the proof of Lemma~\ref{merge}, we consider a hyperbolic triangle $ABC$ such that $$BC= d_0(u,v),~~~\angle B =  \nu_u(d)/2,~~~\angle C=\nu_v(d_0)/2.$$ Lemma~\ref{exist} implies that $ABC$ exists and $BC$ is its greatest side. Let $AH$ be the height of $ABC$ from point $A$, so $H \in BC$. Then 
$$\sinh BH \tan\angle B=\sinh CH \tan \angle C.$$

Define $\tau=AH$. Note that 
\begin{equation}
\label{heightbound}
\tau \leq \sinh\tau \leq \sin(\nu_v(d_0)/2)\sinh(d_0(u,v)) < \nu_v(d_0)\sinh\Delta.
\end{equation}
For every $t \in [0; \tau]$ let $A_t \in AH$ be the point such that $HA_t = t$. 

Cut $(S_g, d_0)$ along $\chi$. Glue there a bigon $X_t$ as in Section~\ref{mergesec} consisting of two copies of $A_tBC$. We obtain a metric $d_t$ on $S_g$. The set $V$ can be naturally defined in all metric spaces $(S_g, d_t)$. We associate the remaining point $w \in V'$ with the vertex $A_t$ of $X_t\subset (S_g, d_t)$. Thus we get a continuous path of convex cone-metrics $d_t \subset \mathfrak D_c(V')$ parametrized by $t\in [0; \tau]$. Moreover, for each $t \in (0; \tau)$\ we have $d_t \in \mathfrak D_{sc}(V')$. Also note that $V(d_{\tau}) = V'\backslash \{u, v\}$ and $\nu_u(d_0)+\nu_v(d_0)<\nu_w(d_{\tau})$.


By $\psi_1, \psi_2 \subset X_t$ we denote the unique shortest paths from $w$ to $u$ and $v$ in each metric $d_t$ respectively. Remark that for any $0\leq t' < t'' \leq \tau$ the metric $d_{t''}$ can be obtained from $d_{t'}$ by cutting $d_{t'}$ along $\psi_1$ and $\psi_2$ and gluing there the conical 4-gon $X_{t',t''}$ obtained from two copies of the 4-gon $A_{t'}BA_{t''}C$. Treating this operation similarly as in the proof of Claim~\ref{distincr} of Lemma~\ref{merge} one can show that the diameter of $d_t$ is non-decreasing with respect to $t$.

Assume that we have a convex Fuchsian cone-manifold $P_0=P(d_0, V, h_0)$. 
We want to transform it into a convex Fuchsian cone-manifold $P_{\tau}=P(d_{\tau}, V', h_{\tau})$ and estimate the change of $S$. 

\subsubsection{Changing boundary with fixed heights under gluing a bigon}
\label{boundchange}

All the notation is taken from the previous subsection. Here we are going to estimate the change of $S$ under a small deformation of the boundary metric with the fixed heights. We do it modulo several quantitative results that we prove separately in Subsection~\ref{slopes}. After this in Subsection~\ref{simult} we combine it with the results of Subsection~\ref{height} to get the full deformation.

For a metric $d$ on $S_g$ define $\alpha(d):=\arccot(2\cosh(\diam(S_g, d))).$ 

\begin{lm}
\label{metricch}
Let $0\leq t' < t'' \leq \tau$. Assume that there is 
$$h \in \bigcap_{t \in [t'; t'']} H(d_t, V')$$
and a triangulation $\mathcal T_h$ with $V(\mathcal T_h)=V'$ such that every convex Fuchsian cone-manifold $P_t=P(d_t, V', h)$ is compatible with $\mathcal T_h$. Then
\begin{align}
\label{tempfinal}
\begin{split}
0 \leq S(d_{t''}, V', h) - S(d_{t'}, V', h)\leq\\
\leq\left(\frac{(\pi+\Delta)\sinh\Delta}{\Delta}+1\right)\frac{2\pi(t''-t')}{\sin^{2}\alpha(d_{\tau})}.
\end{split}
\end{align}
\end{lm}

\begin{proof}
By $S(t)$ denote $S(P_t)$. As all $P_t$ are compatible with $\mathcal T_h$, $S(t)$ is a function of the upper lengths $l_e(d_t)=l_e(t)$ for $e \in E(\mathcal T_h)$. From Lemma~\ref{S1der} $$\frac{\partial S}{\partial l_{e}}    = \theta_e =\pi-\phi_e .$$

We have $$S(t'')-S(t')=\int_{t'}^{t''} \langle {\rm grad~} S,~\dot l \rangle dt = \int_{t'}^{t''} \sum_{e \in E(\mathcal T_h)} (\pi  - \phi_e )\dot l_e  dt.$$

First we suppose that the edges of $E(\mathcal T_h)$ changing their lengths over $[t'; t'']$ are split into three groups:

(1) edges that intersect $\psi_1$ once;

(2) edges that intersect $\psi_2$ once;

(3) edges with $w$ as one (and only one) of the endpoints.

We call this case \emph{simple}. In the general case it might happen that an edge belongs to several types, or intersects $\psi_1$ or $\psi_2$ more than once, or has $w$ as both endpoints. First we give a proof in the simple case when neither of this happens, and then show how to reduce the general case to the simple one.

We denote the set of edges of each type by $E_1$, $E_2$ and $E_3$ respectively. 

Consider $(S_g, d_t)$ for $t \in (0; \tau]$. Triangulate the bigon $X_t \subset (S_g, d_t)$ in the natural way and extend it to a triangulation $\mathcal T_0$ in any way. This is possible because of Lemma~\ref{iztriang}, which allows to extend any partial triangulation to a full triangulation. Note that for all $t \in (0; \tau]$ the metric $d_t$ realizes $\mathcal T_0$. In the metric $d_0$ our triangulation degenerates, but this provides no additional difficulties to the discussion below. The deformation $d_t$ can be described easily in the chart $\mathfrak D(V', \mathcal T_0)$ as only two edges of $\mathcal T_0$ ($\psi_1$ and $\psi_2$) change their lengths in $d_t$.

Take $e \in E_1$. Develop the path of triangles of $\mathcal T_0$ intersecting $e$ in the moment $t \in [t'; t'']$ to $\H^2$. The development results in a (possibly non-convex and self-intersecting) triangulated polygon $R_0$ in $\H^2$. See Figure~\ref{Pic12}.  Denote the image of the endpoints of $e$ in the development by $D_1$ and $D_2$. By construction, $D_1D_2$ is in the interior of $R_0$ and intersects each interior edge of its triangulation. Note that a triangle of $\mathcal T_0$ can be repeated several times in the triangulation of $R_0$.

\begin{figure}
\begin{center}
\includegraphics[scale=0.55]{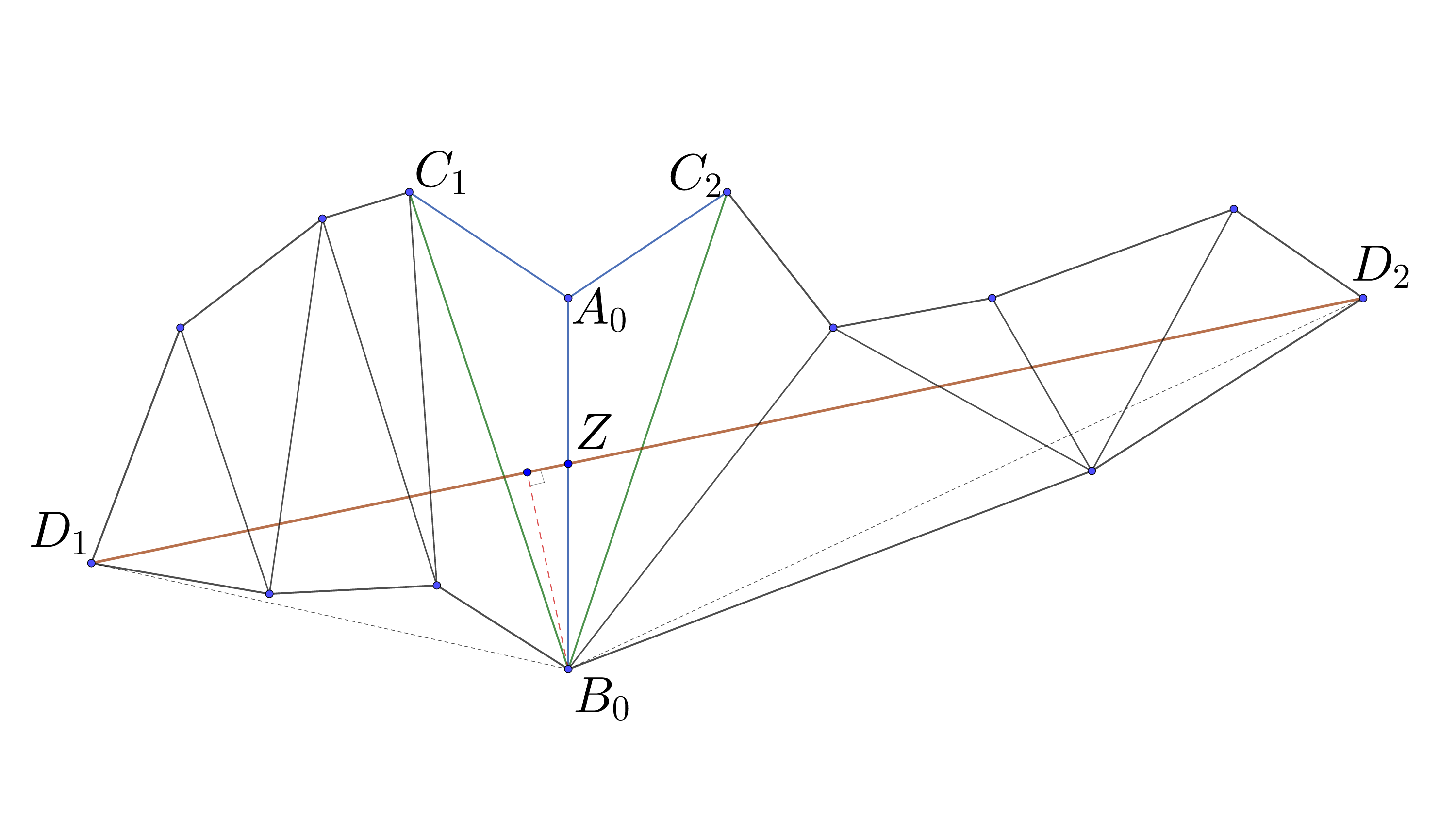}
\caption{Polygon $R_0$ for an edge of type (1).}
\label{Pic12}
\end{center}
\end{figure}

As $e \in E_1$ and the situation is simple, there exists exactly one interior edge of $R_0$ that is the image of $\psi_1$ under the developing map. Denote its endpoints by $A_0$ and $B_0$ so that $A_0$ corresponds to $w$, $B_0$ to $u$, and by $C_1$, $C_2$ denote the points corresponding to $v$. Thereby, the triangles $A_0B_0C_1$ and $A_0B_0C_2$ are isometric to the triangle $A_tBC$ from the bigon $X_t$, the distances $B_0D_1$ and $B_0D_2$ do not change during the deformation and $l_e =D_1D_2$.

Let $y_e $ be the hyperbolic sine of the length of perpendicular from the point $B_0$ to $D_1D_2$, $ \rho $ be the angle $C_1B_0C_2$ and $ \rho_e $ be the angle $D_1B_0D_2$. We claim 
\begin{equation}
\label{templine}
\frac{\partial l_e}{\partial  \rho}  = y_e .
\end{equation}
Indeed, differentiating the cosine rule for the triangle $D_1B_0D_2$ we obtain $$\frac{\partial l_e}{\partial  \rho_e}  = y_e .$$
However, $\frac{\partial  \rho_e}{\partial  \rho} =1,$ therefore we get~(\ref{templine}).
It is easy to see that $\dot \rho \geq 0$, therefore $\dot l_e  \geq 0$.

Let $D_1D_2$ and $A_0B_0$ intersect in a point $Z$. By $x_e $ denote $\sinh B_0Z$. By $\gamma_e $ denote the smallest angle at which $D_1D_2$ intersects $A_0B_0$. We have 
\begin{equation}
\label{templine2}
y_e =x_e\sin\gamma_e.
\end{equation}

In $P_t$ consider the geodesic ray from the point corresponding to $Z$ orthogonal to $\partial_{\downarrow} P_t$. Let this ray make the angles $\beta^-_e$ and $\beta^+_e$ with the segments of $\psi_1$. Define $\beta_e:=\beta^-_e+\beta^+_e$. Corollary~\ref{slopepoint} shows that
$$\alpha(d_t)\leq \beta^-_e, \beta^+_e \leq \pi-\alpha(d_t).$$
Together with Corollary~\ref{slopeedge} and Lemma~\ref{impact} this gives
$$\pi-\phi_e    \leq \frac{\pi-\beta_e}{\sin^2\alpha(d_t)\sin\gamma_e }.$$
Combined this with~(\ref{templine}) and~(\ref{templine2}) we get
\begin{equation}
\label{templine3}
(\pi-\phi_e  )\frac{\partial l_e}{\partial  \rho}   \leq \frac{\pi-\beta_e}{\sin^2\alpha(d_t)\sin\gamma_e }x_e\sin\gamma_e  \leq \frac{\pi-\beta_e}{\sin^2\alpha(d_{\tau})}x_e .
\end{equation}

Consider the union of all segments connecting points of $\psi_1$ with their orthogonal projections to $\partial_{\downarrow} P_t$. This set is isometric to a convex hyperbolic polygon, which is the union of trapezoids with the angles at the upper boundary equal to $\beta_e$ (here $e$ varies in $E_1$). The total length of the upper boundary is equal to $A_tB$. Apply Corollary~\ref{bend} and get 
\begin{equation}
\label{templine4}
\sum_{e \in E_1} (\pi-\beta_e) \leq \pi+A_tB < \pi+\Delta.
\end{equation}
Here we used $A_{t}B<\Delta$ for all $t$.

Define $c(t)=\sinh A_tB$. We have $$x_e  \leq\sinh A_0B_0= \sinh A_tB \leq \sinh A_{t''}B=c(t'').$$ 
From this,~(\ref{templine3}) and~(\ref{templine4}) we obtain
$$\sum_{e \in E_1}(\pi-\phi_e )\frac{\partial l_e}{\partial  \rho} \leq \frac{1}{\sin^2\alpha(d_{\tau})} \sum_{e \in E_1} (\pi-\beta_e)x_e  < \frac{(\pi + \Delta)c(t'')}{\sin^2\alpha(d_{\tau})}.$$

The map $t \rightarrow \rho$ is monotone. Substitute $\rho$ in the integral $$0 \leq \int_{t'}^{t''} \sum_{e \in E_1} (\pi-\phi_e ) \dot l  dt = \int_{\rho(t')}^{\rho(t'')} \sum_{e \in E_1} (\pi-\phi_e ) \frac{\partial l_e}{\partial  \rho} d  \rho <$$ $$< \frac{\pi + \Delta}{\sin^2\alpha(d_{\tau})}c(t'')( \rho(t'')- \rho(t')).$$ 

Now we consider the triangles $A_{t'}BC$, $A_{t''}BC$ embedded in the triangle $ABC$. Note that $\angle A_{t'}BA_{t''}=\frac{1}{2}( \rho(t'')- \rho(t'))<\pi/2$ and $\sinh A_{t''}B=c(t'').$ Therefore, $$c(t'')( \rho(t'')- \rho(t'))= 2\sinh A_{t''}B \cdot \angle A_{t'}BA_{t''}\leq$$ $$\leq \pi \sinh A_{t''}B \cdot \sin\angle A_{t'}BA_{t''}\leq \pi \sinh A_{t'}A_{t''} \leq$$ $$\leq \frac{\pi\sinh\Delta}{\Delta}A_{t'}A_{t''}=\frac{\pi\sinh\Delta}{\Delta}(t''-t').$$
Here we used the concativity of $\sin$ and the convexity of $\sinh$.

In total we get $$0 \leq \int_{t'}^{t''} \sum_{e \in E_1} (\pi-\phi_e ) \dot l  dt \leq \frac{\pi(\pi+\Delta)\sinh\Delta}{\Delta\sin^2\alpha(d_{\tau})}(t''-t').$$

The same bound holds for the edges of the second type: $$0\leq \int_{t'}^{t''} \sum_{e \in E_2} (\pi-\phi_e ) \dot l  dt \leq\frac{\pi(\pi+\Delta)\sinh\Delta}{\Delta\sin^2\alpha(d_{\tau})}(t''-t').$$

Now consider edges of the third type. By Lemma~\ref{vertexangle} and Lemma~\ref{section} we have $$\sum_{e \in E_3} (\pi - \phi_e )\leq \omega_w(P_t) \leq \frac{2\pi}{\sin\alpha(d_t)}\leq \frac{2\pi}{\sin^{2}\alpha(d_{\tau})}.$$
Recall here that $\omega_w(P_t)$ is the total cone angle of $w$ in the cone-manifold $P_t$.

Let $e \in E_3$. Develop the path of triangles of $\mathcal T_0$ intersecting $e$ in the moment $t$ to $\H^2$. Let $A_0$ be the image of $w$ under the developing map, $D$ be the image of the second endpoint of $e$, $B_0C_0$ be the image of the boundary segment of $X_t$ intersected by $e$. 

Consider the moment $t+\Delta t$ and develop this path of triangles so that the image of each triangle coincides with its image under the previous development, except the first one, which is developed to $A_1B_0C_0$. Then $l_e(t)=A_0D$, $l_e(t+\Delta t)=A_1D$. Note that the angle $\angle DA_0A_1$ is not acute, therefore $l_e(t+\Delta t)\geq l_e(t) $, which implies $\dot l_e  \geq 0$.

From the triangle inequality we see $$l_e(t+\Delta t)-l_e(t)  = DA_1-DA_0 \leq A_0A_1= \Delta t.$$ Hence, $\dot l_e  \leq 1$. Then $$0 \leq \int_{t'}^{t''} \sum_{e \in E_3} (\pi - \phi_e ) \dot l_e  dt \leq \int_{t'}^{t''} \sum_{e \in E_3} (\pi - \phi_e ) dt \leq$$ $$\leq \frac{2\pi}{\sin^{2}\alpha(d_{\tau})}(t''-t').$$

Summing up all three estimates we get exactly the estimate~(\ref{tempfinal}). This finishes the proof in the simple case.

Consider the general case, when an edge $e$ may intersect each of $\psi_1$ and $\psi_2$ several times and in the same time may have endpoints at $w$. We subdivide our exposition further. We suppose now that each edge $e \in E(\mathcal T_h)$, when it crosses $X_t$, intersects either $\psi_1$ or $\psi_2$ and then leaves $X_t$. This means that we exclude the case when after intersecting $\psi_1$ our edge $e$ intersects $\psi_2$ or vice versa. We call this case \emph{semi-simple}. 

As before, we develop the path of triangles intersecting $e$ in the moment $t$ to $\H^2$ as a triangulated polygon $R_0$. Orient $e$ arbitrarily and enumerate all the triangles of $R_0$ with respect to the orientation.
By $k$ denote the number of edges in $R_0$ (including the interior edges of the triangulation of $R_0$) and by $z_{1}(t), \ldots, z_{k}(t)$ denote their lengths in the moment $t$. The lengths $z_{1}(t), \ldots, z_{k}(t)$ and the combinatorics determine $R_0$ up to isometry. Thus, the length $l_e $ can be considered as the composite function $l_e(t)=l_e(z_{1}(t), \ldots, z_{k}(t)).$ We consider the derivative $\dot l_e $ as the composite derivative $$\dot l_e =\sum \frac{\partial l_e}{\partial z_{i}}\frac{\partial z_{i}}{\partial t} $$ and decompose this sum into \emph{elementary deformations}, which have types (1), (2) and (3) similar to the edge types in the simple case.

(1) Assume that $T_{i-1}$ and $T_{i}$ are two subsequent triangles of $R_0$ that are the images of $X_t$ under the developing map, and they share an edge that is the image of $\psi_1$ intersected by $e$.
Let $z_{i_1}$ be the length of this edge and $z_{i_2}$, $z_{i_3}$ be the lengths of the images of $\psi_2$ (as $e$ leaves $X_t$ immediately, the edges of $T_{i-1}$, $T_i$ at the boundary of $R_0$ are the images of $\psi_2$). Then we say that the sum of coordinate derivatives 
$$L(e, i_1, i_2, i_3) =\frac{\partial l_e}{\partial z_{i_1}}\frac{\partial z_{i_1}}{\partial t} + \frac{\partial l_e}{\partial z_{i_2}}\frac{\partial z_{i_2}}{\partial t}  + \frac{\partial l_e}{\partial z_{i_3}}\frac{\partial z_{i_3}}{\partial t} $$
is an elementary deformation of the first type.

An elementary deformation of the first type is analyzed in the same way as the deformation of an edge of the first type in the simple case. Let $y(e, i_1, i_2, i_3) $ be the hyperbolic sine of the length of perpendicular from the image of $u$ in $T_{i-1}\cup T_i$ to the image of $e$ and consider $z_{i_s} $ as the function of $ \rho$, which is the angle at vertex $u$ in the bigon $X_t$. Then 
$$L(e, i_1, i_2, i_3) =y(e, i_1, i_2, i_3) \frac{\partial  \rho}{\partial t} .$$

(2) An elementary deformation of the second type is defined similarly, but swapping $\psi_1$ and $\psi_2$.

(3) Assume that $w$ is the starting vertex of $e$ and $z_{i_1}$, $z_{i_2}$ are the lengths of the edges of $T_1$ that are the images of $\psi_1$ and $\psi_2$. Then we say that the sum of coordinate derivatives 
$$L(e, i_1, i_2) =\frac{\partial l_e}{\partial z_{i_1}}\frac{\partial z_{i_1}}{\partial t} + \frac{\partial l_e}{\partial z_{i_2}}\frac{\partial z_{i_2}}{\partial t}$$
is an elementary deformation of the third type. Similarly, if $w$ is the ending vertex of $e$, then the sum of coordinate derivatives of side lengths of the last triangle $T_k$ also form an elementary deformation of the third type.

As in the case of an edge of the third type, we have $$0\leq L(e, i_1, i_2)  \leq 1.$$

We see that every coordinate derivative $$\frac{\partial l_e}{\partial z_{i}}\frac{\partial z_{i}}{\partial t} $$ either belongs to exactly one elementary deformation or is equal to zero. Therefore, the total derivative $\dot l_e $ can be decomposed into the sum of elementary deformations.

We do this for every edge $e$ of $\mathcal T_h$ and group elementary deformations of each type together for all edges. Then we decompose the sum $$\sum_{e \in E(\mathcal T_h)} (\pi  - \phi_e )\dot l_e $$ into three summands corresponding to elementary deformations of each type. The integral of each summand is estimated in the same way as in the simple case. This finishes the proof in the semi-simple case.

\begin{figure}
\begin{center}
\includegraphics[scale=0.55]{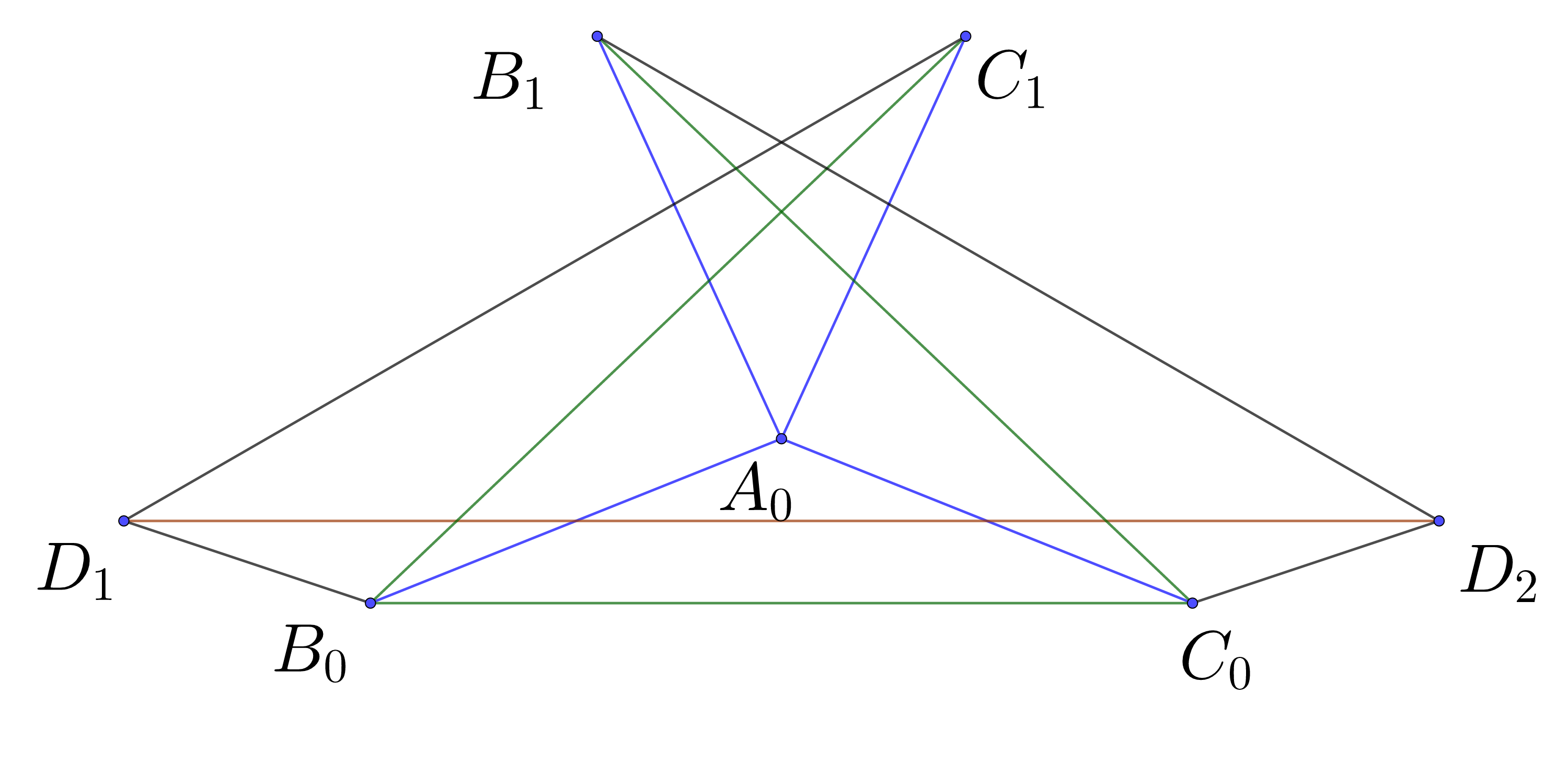}
\caption{Polygon $R_0$ when $e$ crosses $\psi_1$ and then $\psi_2$.}
\label{Pic12-1}
\end{center}
\end{figure}

Now we proceed to the last remaining case when, in addition to the previously described events, $e$ may enter $X_t$ and cross $\psi_1$ and then $\psi_2$ or vice versa. Note that after this $e$ needs to leave $X_t$ as otherwise it self-intersects inside $X_t$. Let the development of $X_t$ in $R_0$ corresponding to this event be the pentagon $A_0C_1B_0C_0B_1$. It is self-intersecting as the total angle $\lambda_w(d_t)$ is greater than $\pi$. See Figure~\ref{Pic12-1}. Here $A_0B_0$ and $A_0B_1$ are images of $\psi_1$, $A_0C_0$ and $A_0C_1$ are images of $\psi_2$. By $z_{i_1}, \ldots, z_{i_4}$ we denote their lengths. By $\rho_1$ we denote the angle of $X_t$ at $u$, the endpoint of $\psi_1$, and by $\rho_2$ we denote the angle of $X_t$ at $v$, the endpoint of $\psi_2$.

In this case we need to decompose the partial derivatives further. Define
$$L(e, i_1, i_2, i_3, i_4, j)(t):=\sum_{s=1}^{4}\frac{\partial l_e}{\partial z_{i_s}}\frac{\partial z_{i_s}}{\partial \rho_j}\frac{\partial \rho_j}{\partial t}(t).$$
We say that the elementary deformation $L(e, i_1, i_2, i_3, i_4, 1)$ belongs to the type (1) and $L(e, i_1, i_2, i_3, i_4, 2)$ belongs to the type (2). 

Let $y(e, i_1, i_2, i_3, i_4, 1)$ be the hyperbolic sine of the length of perpendicular from $B_0$ to the image of $e$, $y(e, i_1, i_2, i_3, i_4, 2)$ be the hyperbolic sine of the length of perpendicular from $C_0$. We have $\rho_1=\angle C_1B_0C_0$, $\rho_2=\angle B_0C_0B_1$. Thus, it is easy to see that
$$L(e, i_1, i_2, i_3, i_4, 1) =y(e, i_1, i_2, i_3, i_4,1) \frac{\partial  \rho_1}{\partial t},$$
$$L(e, i_1, i_2, i_3, i_4, 2) =y(e, i_1, i_2, i_3, i_4,2) \frac{\partial  \rho_2}{\partial t},$$
so indeed $L(e, i_1, i_2, i_3, i_4, 1)$, $L(e, i_1, i_2, i_3, i_4, 2)$ can be treated in the same way as elementary deformations of types (1) and (2) that were defined for the semi-simple case. Hence, the bound for $$\sum_{e \in E(\mathcal T_h)} (\pi  - \phi_e )\dot l_e $$ remains the same and we finish the proof.
\end{proof}

\subsubsection{Slopes and angle estimates}
\label{slopes}

All results of this subsection are auxiliary for the proof of Lemma~\ref{metricch} and a similar Lemma~\ref{flatSbound} below.

Let  $P=P(d, V, h)$ be a convex Fuchsian cone-manifold, $p \in (S_g, d)$ be a point and $T$ be a face triangle containing $p$. Embed the prism containing $T$ to $\H^3$. Lemma~\ref{ultrapar2} says that this prism is ultraparallel. By $A_p$ we denote the image of $p$, by $B_p$ its projection to the lower plane, by $AB$ denote the common perpendicular between the upper and the lower boundary planes. Define
$$\sl(p, P, T):=\angle B_pA_pA.$$

It is clear that it does not depend on the choice of an embedding and if $T$ is contained in a face $R$ of $P$, then $\sl(p, P, T)$ does not depend on $T$, so we can denote it as $\sl(p, P, R).$ We have $0<\sl(p, P, R)\leq \pi/2$.

Recall that $\alpha(d):=\arccot(2\cosh(\diam(S_g, d))).$

\begin{lm}
\label{slope}
For each $p\in (S_g, d)$ and a face $R \in \mathcal R(P)$ with $p \in \overline R$ we have $$\alpha(d) \leq \sl(p, P, R).$$
\end{lm}

\begin{proof}
In the notation above we get (using Corollary~\ref{rightangl})
$$\cot (\sl(p, P, R)) = \sinh AA_p \tanh AB\leq \cosh AA_p \tanh AB =$$
$$=\frac{\sinh A_pB_p}{\cosh AB}=\frac{\sinh (A_pB_p-AB+AB)}{\cosh AB}=$$ $$=\sinh(A_pB_p-AB) + \cosh(A_pB_p-AB)\tanh AB \leq$$
$$\leq 2\cosh (A_pB_p-AB).$$

We claim that there exists a point $q \in S_g$ such that $\widetilde h(q) \leq AB$. Indeed, we can assume to this purpose that $p \in R$, i.e, $p$ is not a conical point and does not belong to a strict edge, because this does not change $AB$. Consider the unit speed geodesic $\psi\subset (S_g, d)$ starting from $p$ in the direction $-{\rm grad_p}\widetilde h$. Suppose that $\psi$ extends to the distance $AA_p$ from $p$ (it may not happen in the only case when $\psi$ ends in a conical point before). Let $q \in \psi$ be the point at the distance $AA_p$. Consider the function $f$ on $\psi$ defined by
$$\sinh f(x)= \cosh(x-AA_p)\sinh AB.$$

We recall the definition of $\mathcal F(-1)$-function and $\mathcal F(-1)$-concave function from Section~\ref{heightdefsec}. By definition, $\sinh f$ is $\mathcal F(-1)$; $\sinh \widetilde h$ is $\mathcal F(-1)$-concave due to Lemma~\ref{distf-1conc} and $f(x)$ coincides with the restriction of $\widetilde h$ to $\psi$ in a neighbourhood of $p$. Then $$\sinh \widetilde h(q) \leq \sinh f(q)=\sinh AB.$$

In the case if $\psi$ ends in a conical point before moving at distance $AA_p$, we can perturb slightly the direction of $\psi$ to avoid this. By the argument above, we obtain $q_{\e}$ such that $\widetilde h(q_{\e}) \leq AB+\e$ for some $\e>0$. We can make $\e$ arbitrarily small and consider a limit point $q$ of the set $\{q_{\e}\}$. By continuity of $\widetilde h$, we have $\widetilde h(q) \leq AB$.

We get $$\cot (\sl(p, P, R)) \leq 2\cosh(\widetilde h(v)-\widetilde h(q)) \leq 2\cosh(\diam(S_g,d))$$ due to Lemma~\ref{htriangleinequality}. This proves the desired inequality.
\end{proof}

The following corollary easily follows after a local development.

\begin{crl}
\label{slopepoint}
Let $p\in (S_g, d)$ and $\alpha$ be the angle between a geodesic ray from $p$ in $\partial^{\uparrow} P$ and the geodesic ray from $p$ orthogonal to $\partial_{\downarrow} P$. Then
$$\alpha(d) \leq \alpha \leq \pi-\alpha(d).$$
\end{crl}

The next corollary is just slightly more involved.

\begin{crl}
\label{slopeedge}
Let $e$ be an edge of $P$, $\phi^-_e$, $\phi^+_e$ be the dihedral angles of $e$ in the prisms containing $e$. Then $$\alpha(d) \leq \phi^-_e, \phi^+_e \leq \pi-\alpha(d).$$
\end{crl}

\begin{proof}
Let $\Pi$ be the prism containing $e$ with dihedral angle $\phi^-_e$ and $T$ be the upper boundary triangle of $\Pi$. Embed it to $\H^3$ and let $A$ be as above. By $A_e$ denote the closest point to $A$ on the line containing $e$ and by $B_e$ denote its projection to the lower boundary plane. Then $\angle AA_eB_e$ is either $\phi^-_e$ or $\pi-\phi^-_e$ and $\angle AA_eB_e\leq \pi/2$. The only issue is that $A_e$ may lie outside of $e$. Then take any $p \in e$ and note that $\sl(p, P, T)\leq \angle AA_eB_e$. Apply Lemma~\ref{slope} and finish the proof.
\end{proof}






\begin{lm}
\label{vertexangle}
Let $P=P(d, V, h)$ be a convex Fuchsian cone-manifold. Then for every $v \in V$ we have $$\omega_v(P) \leq \frac{2\pi}{ \sin \alpha(d)}.$$
\end{lm}

\begin{proof}
Let $M^{\uparrow}$ and $M_{\downarrow}$ be ultraparallel hyperbolic planes in $\H^3$, $A\in M^{\uparrow}$ and $B \in M_{\downarrow}$ be the closest points, $A_1 \in M^{\uparrow}$ be a point, $B_1$ be its projection to $M_{\downarrow}$ and $\alpha=\angle B_1A_1A$.

Consider an angle of value $\lambda$ in $M^{\uparrow}$ with the vertex $A_1$. Let $\omega$ be the value of its orthogonal projection to $M_{\downarrow}$. We claim that 
\begin{equation}
\label{slte1}
\omega\leq \frac{\lambda}{ \sin \alpha}.
\end{equation}

First, assume that $A$ belongs to one of the two rays bounding $\lambda$ and $\lambda \leq \pi/2$. We call it an angle of the first type. By looking at the spherical link we get
\begin{equation}
\tan\omega=\frac{\tan \lambda}{\sin \alpha}.
\end{equation}


Let $f(x)=\arctan\left(\frac{\tan x}{\sin \alpha}\right).$ We have $f(0)=0$ and $$f'(x)=\frac{\sin\alpha}{\sin^2\alpha+\cos^2\alpha\sin^2 x} \leq \frac{1}{\sin\alpha}$$ 
for $\alpha, x \in (0; \pi/2]$.
This proves~(\ref{slte1}) for angles of the first type. We can also see that $f(x)$ is monotonously increasing and concave for $0\leq x \leq \pi/2$. 

If we have an angle of value $\lambda$ that is the union of two angles of the first type, then~(\ref{slte1}) holds by additivity. We call it an angle of the second type. 

Let an angle of value $\lambda$ be the difference of two angles of the first type of values $\lambda_1$ and $\lambda_2$. We call it an angle of the third type. Then by concativity of $f$ for the projections of these angles we have $$\omega = \omega_1-\omega_2 \leq (\lambda_1-\lambda_2) f'(\lambda_2) \leq \frac{\lambda}{ \sin\alpha}.$$ 

Now consider an arbitrary angle at $A_1$. Let $s$ be a line in $M^{\uparrow}$ through $A_1$ orthogonal to $AA_1$. Take the part of our angle that lies on the other side from $A$ with respect to $s$ and replace this part with its centrally symmetric image in $M^{\uparrow}$ with respect to $A_1$.  This does not change neither the total value of the angle nor the total value of the projection. Now we have at most two angles and each of them belongs to one of the three types described above. Then they satisfy~(\ref{slte1}) and so does the initial angle.

Consider a convex Fuchsian cone-manifold $P=P(d, V, h)$. Subdivide its face decomposition $\mathcal R(P)$ into triangles. For a vertex $v \in V$ incident to a triangle $T$ let $\lambda_{v, T}$, $\omega_{v, T}$ be the angles of the prism containing $T$ in the upper and lower faces corresponding to $v$. We showed that 
$$\omega_{v, T}\leq \frac{\lambda_{v,T}}{\sin(\sl(v, P, T))}.$$

Thus,
$$\omega_v(P) = \sum_{T{\rm~is~incident~}v} \omega_{v, T}\leq  \left( \sum_{T{\rm~is~incident~}v} \frac{\lambda_{v,T}}{\sin(\sl(v, P, T))} \right) \leq$$ $$\leq \frac{2\pi}{ \sin \alpha(d)}.$$
\end{proof}

\begin{lm}
\label{section}
Let $P=P(d, V, h)$ be a convex Fuchsian cone-manifold. Then for every $v \in V$ we have $$\sum\limits_{e{\rm~is~incident}~v} (\pi - \phi_e(P)) < \omega_v(P).$$
\end{lm}

\begin{proof}
Consider the spherical link of $P$ at $v$. It is a spherical polygon with a conical point of total angle $\omega_v(P)$ and boundary angles $\phi_e(P)$. Then its area is $$\omega_v(P) +\sum\limits_{e{\rm~is~incident}~v} (\phi_e(P)-\pi)>0.$$
\end{proof}

\begin{figure}
\begin{center}
\includegraphics[scale=0.55]{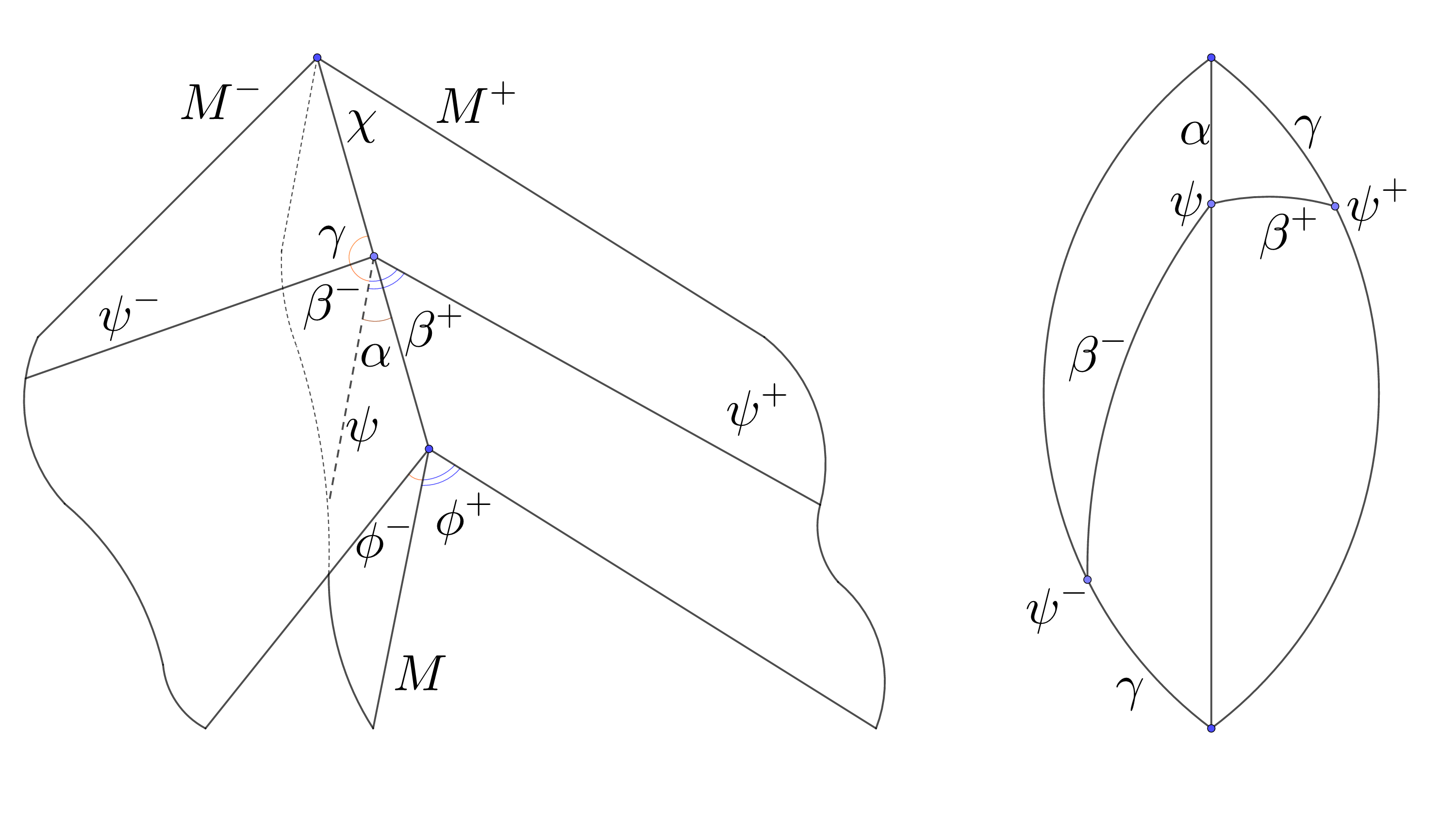}
\caption{To the statement of Lemma~\ref{impact}. \\ Left: In space. Right: in the spherical link of $p$.}
\label{Pic13}
\end{center}
\end{figure}

\begin{lm}
\label{impact}
Consider two half-planes $M^-,$ $M^+$ in $\H^3$ sharing a line $\chi$ (see Figure~\ref{Pic13}). Let $p \in \chi$ and $\psi^- \subset M^-$, $\psi^+ \subset M^+$ be geodesic rays from $p$ making an intrinsic geodesic in the union of $M^-$ and $M^+$. By $\gamma>0$ denote the (smallest) angle between this geodesic and~$\chi$.

Let $\psi$ be a geodesic ray from $p$ in the (convex) dihedral angle spanned by $M^-$ and $M^+$ and $M$ be the half plane spanned by $\chi$ and $\psi$. By $\phi^-$, $\phi^+$ denote the dihedral angles between $M$ and $M^-$, $M^+$ respectively. By $\beta^-$, $\beta^+$ denote the angles between $\psi$ and $\psi^-$, $\psi^+$ respectively. Define  $\phi=\phi^-+\phi^+,$ $\beta:=\beta^-+\beta^+$. By $\alpha$ denote the (smallest) angle between $\psi$ and $\chi$.

Assume that for some $\alpha_0$ with $0<\alpha_0<\pi/2$ we have
$$\alpha_0 \leq \alpha,\phi^-, \phi^+,\beta^-,\beta^+ \leq \pi-\alpha_0.$$ Then $$\pi -\phi \leq \frac{\pi-\beta}{\sin^2\alpha_0 \sin\gamma}.$$
\end{lm}

\begin{proof}
We consider the function $\frac{\pi-\beta}{(\pi-\phi)\sin\gamma}$ as the function of $\phi^-$, $\phi^+$, $\gamma$ and $\alpha$ and find its minimal value under the given restrictions. 

Note that $\beta \leq \pi$. As $0< \gamma, \alpha < \pi$, then $\phi=\pi$ if and only if $\beta=\pi$. So we assume that $\phi, \beta < \pi$.

Consider a spherical section at $p$, so geodesic rays become points and half-planes become half-circles. Assume also that $\beta^+$ belongs to the spherical triangle with the other sides $\gamma$ and $\alpha$ (otherwise, switch $M^-$ and $M^+$). From the cosine law we get
$$\cos\beta^+=\cos\gamma\cos\alpha+\sin\gamma\sin\alpha\cos\phi^+,$$
$$\cos\beta^-=-\cos\gamma\cos\alpha+\sin\gamma\sin\alpha\cos\phi^-,$$
$$\cos\beta^++\cos\beta^-=\sin\gamma\sin\alpha(\cos\phi^++\cos\phi^-).$$

As $\beta \leq \pi$, then from elementary trigonometry
$$\pi-\beta \geq \cos \beta^-+\cos\beta^+ \geq 0$$
with equality if and only if $\beta=\pi$. Thus,
$$\pi-\beta \geq \sin\gamma\sin\alpha(\cos\phi^++\cos\phi^-)$$
and
$$\frac{\pi-\beta}{(\pi-\phi)\sin\gamma}\geq \frac{\sin\alpha(\cos\phi^-+\cos\phi^+)}{\pi-\phi}.$$

Assume that $\phi^-\leq \phi^+$. Fix $\phi^-\leq \pi/2$ and consider the function $$f(\phi^+):=\frac{\cos\phi^-+\cos\phi^+}{\pi-\phi}$$ 
as a function of $\phi^+ \in [\phi^-; \pi-\phi^-)$. It suffices to prove $f(\phi^+)\geq \sin\alpha_0$. Using L'Hopital's rule we compute the limit of $f$ as $\phi^+ \rightarrow \pi-\phi^-$:
$$\lim_{\phi^+\rightarrow \pi-\phi^-}f(\phi^+)=\sin(\phi^-)\geq \sin\alpha_0.$$

Now we compute its derivative
$$\frac{\partial f}{\partial \phi^+}=\frac{1}{(\pi-\phi)^2}\left(-\sin\phi^+(\pi-\phi)+\cos\phi^-+\cos\phi^+\right).$$
We see that if $\phi^+$ is a critical point, then $f$ is equal to $\sin\phi^+\geq \sin\alpha_0$.

It remains to check what happens as $\phi^+=\phi^-$. From elementary trigonometry we have
$$f(\phi^-)=\frac{2\cos\phi^-}{\pi-2\phi^-}\geq \sin\phi^-\geq \sin\alpha_0.$$
In total we get 
$$\frac{\pi-\beta}{(\pi-\phi)\sin\gamma}\geq \sin^2\alpha_0,$$
which is equivalent to the desired inequality.
\end{proof}

\begin{lm}
\label{arlele}
Let $A_1A_2B_2B_1$ be a trapezoid in $\H^2$. Then
\begin{equation}
\label{arle}
\area(A_1A_2B_2B_1) \leq A_1A_2.
\end{equation}
\end{lm}

\begin{proof}
We claim that is enough to prove (\ref{arle}) for a right-angled trapezoid. Indeed, if both $\angle A_1$, $\angle A_2$ are less than $\pi/2$, then the trapezoid is ultraparallel and, moreover, can be cut into two right-angled ones. If we know (\ref{arle}) for both of them, then it holds also for the initial one. Suppose that $\angle A_1 > \pi/2$. Then we can rotate $A_1A_2$ around $A_2$ until $\angle A_1$ becomes $\pi/2$ so that the area only grows. 

Thereby, it is enough to prove (\ref{arle}) for $\angle A_1 = \pi/2$.  We use the notation from Subsection~\ref{tp}. The area of the trapezoid is $\pi/2 - \alpha_{21}$. From Corollary~\ref{rightangl} we have
$$\tan(\pi/2-\alpha_{21})=\cot\alpha_{21}=\tanh h_1 \sinh l_2 < \sinh l_2.$$
Therefore, the area of the trapezoid is less than $\arctan\sinh A_1A_2$. Consider the function $f(x)=\arctan(\sinh x)$. We have $$f'(x) = \frac{1}{\cosh x} \leq 1.$$ 

As $f(0)=0$, it follows that $f(x)\leq x$ for $x>0$. This finishes the proof.
\end{proof}

\begin{crl}
\label{bend}
Let $B_0A_0\ldots A_{m+1}B_{m+1}$ be a convex polygon in $\H^2$ with $\angle B_0 =\angle B_{m+1}=\pi/2$, $\alpha_k$ be the angle $A_{k-1}A_kA_{k+1}$ and $L$ be the sum of lengths $A_0A_1+\ldots +A_mA_{m+1}$. Then $$\sum_{k=1}^m(\pi - \alpha_k) \leq \pi + L.$$
\end{crl}

\begin{proof}
Decompose the polygon into trapezoids $A_{k-1}A_kB_{k-1}B_k$. Denote the angle $B_kA_kA_{k+1}$ by $\alpha^+_k$, the angle $B_kA_kA_{k-1}$ by $\alpha^-_k$ and the length $A_{k}A_{k+1}$ by $l_{k}$. Lemma~\ref{arlele} shows $\pi-\alpha^+_k-\alpha^-_{k+1} \leq l_{k}$. 

Then $$\sum_{k=1}^m(\pi - \alpha_k) = \pi/2 - \alpha^+_{1} + \sum_{k=1}^{m-1}(\pi - \alpha^-_k-\alpha^+_{k+1})+ \pi/2 - \alpha^-_m \leq$$ $$\leq \pi+\sum_{k=1}^{m-1}l_k \leq \pi+L.$$
\end{proof}

\subsubsection{Changing heights and boundary simultaneously}
\label{simult}


Lemma~\ref{metricch} gives us a Lipschitz bound on the change of $S$ when we change the boundary metric with fixed heights as we glue a bigon. Below we will obtain some other bounds of this type. Lemma~\ref{collect} combines such bounds with the results of Subsection~\ref{height} to get the full deformation. 

Recall that
$$H_S(d, V, m, K)=\{h \in H(d, V): \min_{v \in V}h_v \geq m, S(d, V, h) \geq K\}.$$

For $m \in \R_{>0}$ and $K \in \R$ define the sets
$$\mathfrak H_S(\mathfrak U, m, K) = \{ (d, h): d \in \mathfrak U, h \in H_S(d,V, m, K)\},$$ 
$$H_S(\mathfrak U, m, K) = \bigcup_{d \in \mathfrak U} H_S(d, V,m,K).$$
We claim

\begin{lm}
\label{H_U}
If $\mathfrak U$ is compact, then so are $\mathfrak H_S(\mathfrak U, m, K)$, $H_S(\mathfrak U, m, K)$. 
\end{lm}

\begin{proof}
By Corollary~\ref{H_S comp}, $H_S(d,V, m, K)$ is compact. The projection map $(d, h) \rightarrow d$ mapping $\mathfrak H_S\big(\mathfrak D_c(V), m, K\big)$ onto $\mathfrak D_c(V)$ is continuous, closed and have compact fibers. Therefore, it is proper. Then if $\mathfrak U$ is compact, $\mathfrak H_S(\mathfrak U, m, K)$ also is compact. The set $H_S(\mathfrak U, m, K)$ is compact as the image of $\mathfrak H_S(\mathfrak U, m, K)$ under the projection map $(d, h) \rightarrow h$.
\end{proof}

Note that except the proof of Lemma~\ref{collect} just below, Lemma~\ref{H_U} will be also prominent in our proof of Main Lemma IIIC. Now we are ready to show

\begin{lm}
\label{collect}
Let $d_t \subset \mathfrak D_{c}(V)$, $t \in [0; \tau]$, be a continuous path of metrics such that $d_t \in \mathfrak D_{sc}(V)$ for $t \in (0; \tau)$, and $P_0=P (d_0, V, h_0 )$ be a convex Fuchsian cone-manifold. 

Assume that there exists a constant $C\in\R$ with the following property. Let $[t'; t''] \subseteq [0; \tau]$ and let there be $h \in H(d_t, V)$ for all $t\in [t',t'']$ and a triangulation $\mathcal T_h$ compatible with Fuchsian cone-manifolds $P_t=P (d_t, V, h )$ for all $t\in [t',t'']$. Then $$S (d_{t''}, V, h )\geq S (d_{t'}, V, h ) - C(t''-t').$$

Here $C$ is independent on $h$, $\mathcal T_h$, $t', t''$.

Let $I \subseteq [0; \tau]$ be the set of $t$ such that there exists $h_t \in H(d_t, V)$ with the properties:

(1) $ h_{t, v} \geq  h_{0, v}$ for every $v \in V$;

(2) $S (d_t, V, h_t ) \geq S(d_0, V, h_0)-Ct$.

We claim that $I=[0;\tau].$
\end{lm}

\begin{proof}
Note that $0 \in I$. First suppose that $d_0 \in \mathfrak D_{sc}(V)$. Take $t \in I$ such that $t \neq \tau$. We prove that there is $\xi > 0$ such that $[t; t+\xi) \subset I$. Apply Lemma~\ref{open} to $P_t =P (d_t, V, h_t)$. We obtain $\xi >0$ and $h'$ such that:

(1) $  h'_v \geq  h_{t, v}\geq h_{0, v}$ for every $v \in V$;

(2) $S(d_t, V, h') \geq S(d_t, V, h_t)$;

(3) $h' \in H(d_{t'}, V)$ for every $t'\in [t; t+\xi)$.

We may also assume that Fuchsian cone-manifolds $P_{t'}=P(d_{t'}, V, h')$ for every $t'\in [t; t+\xi)$ have a common compatible triangulation provided $\xi$ is sufficiently small. From the assumption of the lemma we get $$S(d_{t'}, V, h') \geq S(d_t, V, h')-C(t'-t) \geq $$ $$\geq S(d_t, V, h_t)-C(t'-t) \geq S(d_0, V, h_0)-Ct'.$$
Therefore, $t'\in I$ and we can let $h_{t'}=h'$.

Suppose that $I\neq [0; \tau]$. Let $t \in (0; \tau]$ be the smallest value not in $I$. Consider a monotonous sequence $\{t_n\}$ converging to $t$ from the left. Since $t_n \in I$, there exists a sequence $\{h_n \in H(d_{t_n}, V)\}$ such that 

(1) $ h_{n, v} \geq   h_{0, v}$ for every $v \in V$;

(2) $S (d_{t_n}, V, h_n ) \geq S (d_0, V, h_0 )-Ct_n$.

Define $\mathfrak U = \{d_t: t\in [0; \tau]\} \subset \mathfrak D_c(V)$. Due to property (1) there also exists a lower bound $m>0$ for all $h_{n,v}$. Also define $K:=S(d_0, V, h_0)-C\tau$. We have all $h_n \in H_S(\mathfrak U, m, K)$. By Lemma~\ref{H_U} the set $H_S(\mathfrak U, m, K)$ is compact. Then there is a limit point $h_t \in H(d_t, V)$ for $\{h_{n}\}$, and for every $v \in V$ one has $ h_{t, v} \geq  h_{0,v}$. The discrete curvature $S$ is continuous over $\mathfrak H_S(\mathfrak U, m, K)$, therefore $S (d_t, V, h_t ) \geq S (d_0, V, h_0 )-Ct$.

Hence, $t \in I$ and we obtained a contradiction.

Now consider the case $d_0 \in \mathfrak D_c(V)$, but not in $\mathfrak D_{sc}(V)$. Apply Lemma~\ref{openV} to $d_0$ (with arbitrary $\mu>0$). For $t \in (0; \tau)$ we can still apply Lemma~\ref{open} to $d_t$. In this way for any $\mu>0$ and $t \in [0;\tau]$ we obtain $h_{t, \mu} \in H(d_t, V)$ such that 

(1) $ h_{t, \mu, v} \geq  h_{0, v}$ for every $v \in V$;

(2) $S (d_t, V, h_{t, \mu} ) \geq S(d_0, V, h_0)-Ct-\mu$.

Fix $t$ and consider a sequence $\mu_n \rightarrow 0$. The heights $h_{t, \mu_n}$ are uniformly bounded from below. The discrete curvature $S (d_t, V, h_{t, \mu_n} )$ is also uniformly bounded from below. Hence, by Corollary~\ref{H_S comp} they belong to a compact set and there exists the limit point $h_t \in H(d_t, V)$ such that

(1) $ h_{t, v} \geq  h_{0, v}$ for every $v \in V$;

(2) $S (d_t, V, h_t ) \geq S(d_0, V, h_0)-Ct$.
\end{proof}

Note that Lemma~\ref{collect} does not hold for paths in $\mathfrak D_c$ due to remarks in the proof of Lemma~\ref{openV}.

Now we return to the path of metrics $d_t$ described in Subsection~\ref{gluing}. Take a convex Fuchsian cone-manifold $P_0=P(d_0, V, h_0)$. Applying Lemma~\ref{collect} together with the lower bound of Lemma~\ref{metricch} we get (here $C=0$):

\begin{crl}
\label{glued1}
There exists $h_{\tau} \in H(d_{\tau}, V')$ such that

(1) $ h_{\tau, v} \geq  h_{0, v}$ for every $v \in V'$;

(2) $S(d_{\tau}, V', h_{\tau}) \geq S(d_0, V, h_0)$.
\end{crl}

Consider this path of metrics in the inverse direction, i.e., define $d'_t := d_{\tau-t}$ for $t \in [0;\tau]$. Let $P'_0=P(d'_0, V', h_0)$ be a convex Fuchsian cone-manifold. Applying the upper bound of Lemma~\ref{metricch} we have the following:

\begin{crl}
\label{glued2}
There exist $h_{\tau} \in H(d'_{\tau}, V')$ such that

(1) $ h_{\tau, v} \geq  h_{0, v}$ for every $v \in V'$;

(2) $$S(d'_{\tau}, V', h_{\tau}) \geq S(d'_0, V, h_0)-\left(\frac{(\pi+\Delta)\sinh\Delta}{\Delta}+1\right)\frac{2\pi\nu_{v}(d_0)\sinh\Delta}{\sin^{2}\alpha(d_{\tau})}.$$
\end{crl}

For the proof of Main Lemma IIIC we will also need the following variation of Lemma~\ref{collect}:

\begin{lm}
\label{collectC}
Let $d_t \subset \mathfrak D_{c}(V)$, $t \in [0; \tau]$, be a continuous path of metrics such that $d_t \in \mathfrak D_{sc}(V)$ for $t \in (0; \tau)$, and $P_0=P (d_0, V, h_0 )$ be a convex Fuchsian cone-manifold. 

Assume that there exists a constant $C\in\R$ with the following property. Let $[t'; t''] \subseteq [0; \tau]$ and let there be $h \in H(d_t, V)$ for all $t\in [t',t'']$ such that

(1*) $ h_{v} \geq  h_{0, v}$ for every $v \in V$;

(2*) $S (d_{t'}, V, h ) \geq S(d_0, V, h_0)-Ct'$.

Let there be a triangulation $\mathcal T_h$ compatible with Fuchsian cone-manifolds $P_t=P (d_t, V, h )$  for all $t\in [t',t'']$. Then $$S (d_{t''}, V, h )\geq S (d_{t'}, V, h ) - C(t''-t').$$

Here $C$ is independent on $h$ satisfying (1*) and (2*), $\mathcal T_h$, $t', t''$.

Let $I \subseteq [0; \tau]$ be the set of $t$ such that there exists $h_t \in H(d_t, V)$ with the properties:

(1) $ h_{t, v} \geq  h_{0, v}$ for every $v \in V$;

(2) $S (d_t, V, h_t ) \geq S(d_0, V, h_0)-Ct$.

We claim that $I=[0;\tau].$
\end{lm}

This Lemma with its monstrous statement is just a more restrictive version of Lemma~\ref{collect} and its proof is exactly the same: in the proof of Lemma~\ref{collect} we actually work only with heights $h$ satisfying (1*) and (2*). We need this version because in the proof of Main Lemma IIIC we can prove a Lipschitz bound that holds only on heights satisfying (1*) and (2*), but not globally.

\subsubsection{Proof of Main Lemmas IIIA-IIIB}
\label{proof}

\begin{lm}
\label{2pi}
Let $P=P(d, V, h)$ be a convex Fuchsian cone-manifold and $v \in V$ be a vertex with $\nu_v(d)=0$. Then $\kappa_v(P) \leq 0$.
\end{lm}

\begin{proof}
Assume that $\kappa_v(P) > 0$. Consider the spherical link of $P$ at $v$ and the solid cone $X$ (with particle singularity) defined by it. Let $x \in \partial X$ be a point such that the ray $vx$ determines the maximal angle with the axis (axis is the perpendicular from $v$ to $\partial_{\downarrow} P$). Consider the 2-dimensional cone  spanned by the rays $vx$ and the axis. It splits the dihedral angle at the ray $vx$ into two angles at most $\pi/2$ as $x$ formed the maximal angle with the axis.

Cut $X$ along this cone and glue there an orthogonal prism with dihedral angle $\kappa_v(P)$ at the axis. The resulting cone is convex. Its particle curvature will be zero, hence it can be embedded in $\H^3$. But its surface angle is greater than $2\pi$, which is a contradiction.
\end{proof}

\begin{proof}[Proof of Main Lemma IIIA]

Choose any $\Delta>0$ and $\delta_{\nu}, \delta_D$ with the help of (\ref{choice1}) and (\ref{choice2}). Put $\delta := \min\{\delta_{\nu}, \delta_D\}.$ Then we can merge all the curvature inside each triangle of $\mathcal T$ as described in the proof of Lemma~\ref{merge}. We are going to use the previous results of this section to deform the cone-manifold $P$. 

Let $v_0$ and $v_1$ be the first two vertices that we merge together, $w_1$ be the new vertex, $d_1$ be the obtained metric, $V'_1=V\cup \{w_1\}$, $V_1 = V(d_1) = V'_1 \backslash \{v_0, v_1\}$, $V'=V(\mathcal T)\cap V(d)$. Note that $V' \subset V_1 \subset V'_1$.
We transform the metric $d$ to $d_1$ continuously through $\mathfrak D_{sc}(V'_1)$ as described in Section~\ref{gluing} and transform respectively the cone-manifold $P=P(d, V, h)$. Using Corollary~\ref{glued1}  we obtain $h'_1 \in H(d_1, V'_1)$ such that 

(1) $ h'_{1,v} \geq  h_v$ for every $v \in V'$;

(2) $S(d_1, V'_1, h'_1) \geq S(d, V, h)$.

After this we have $\nu_{v_0}(d_1) = 0$. By Lemma~\ref{2pi} it has curvature ${\kappa_{v_0}(P'_1)\leq 0}$ in $P'_1=P(d_1, V'_1, h'_1)$. If it is less than zero, then by Corollary~\ref{decrease} we can decrease its height. The functional $S$ is increased under this deformation. Consider the minimal value of $h_{v_0}$ such that it can not be decreased (with all other heights fixed). It is greater than zero. Indeed, the number of triangulations compatible with cone-manifolds defined by $H(d, V)$ is finite. In every such triangulation there exists a prism containing $v_0$ and other vertex of $v$ (possibly, they are connected with flat edges). But for every prism, if we decrease sufficiently the height of one upper vertex with the fixed height of another vertex, then the prism will not be ultraparallel anymore. This contradicts to Lemma~\ref{ultrapar2}. Then for the minimal value of $h_{v_0}$ we get a convex cone-manifold with the particle curvature of $v_0$ equal to 0. 

Do the same with the vertex $v_1$. Let $h_1$ be the modified height function and $P_1=P(d_1, V'_1, h_1)$. It has no particle curvatures at points of $V'_1\backslash V_1$, therefore it can be rewritten as $P_1=P(d_1, V_1, h_1)$. We have 
$$S(P_1) \geq S(P'_1).$$

By repeating this procedure we obtain $\hat h \in H(\hat d, V(\hat d))$ such that

(1) $\hat h_v \geq h_v$ for every $v \in V'$;

(2) $S(\hat d, V(\hat d), \hat h) \geq S(d, V, h)$.

The property (1) holds because $V(\mathcal T)\cap V(d)=V' \subset V(d_i)$ for all intermediate metrics $d_i$. This finishes the proof.
\end{proof}

\begin{proof}[Proof of Main Lemma IIIB]

The proof of Main Lemma IIIB is just slightly more elaborate as now we have to control the decrement of $S$. We start again by choosing any $\Delta>0$, $\delta_{\nu}, \delta_D$ from (\ref{choice1}) and (\ref{choice2}) and putting $\delta(\Delta) := \min\{\delta_{\nu}, \delta_D\}$. The proof of Lemma~\ref{merge} shows that $\hat d$ can be obtained from $d$ via gluings bigons. Let $m$ be the number of steps to obtain $\hat d$. Now we do this process in the reverse order, i.e., we cut bigons from $\hat d$ to obtain $d$.

Let the metric $\hat d$ be obtained from the metric $d_{m}$ by merging $w_m$ with $v_m$ into $w_{m+1}$. Note that $\nu_{v_m}(d_{m})=\nu_{v_m}(d)$. Define $V'_{m}=V(\hat d)\backslash \{w_{m+1}\}$, $V_{m} = V(d_{m}) = V'_{m} \cup \{w_m, v_m\}$, $V'=V(\mathcal T)\cap V(d)$. Then $V'\subset V_m \subset V'_m$.

We have to make a remark on the diameters. Define
$$\hat D=\hat D(\Delta, D, A):=D+(A+2\pi(2-2g))\sinh\Delta.$$
Recall that Lemma~\ref{diamglobal} implies that $\hat d$ has the largest diameter among all metrics on our way and it satisfies
$$\diam(S_g, \hat d) \leq \diam(S_g, d) + \nu(S_g, d)\sinh\Delta \leq \hat D.$$

Set $$M=M(\Delta, D, A):=\frac{1}{\sin^{2}(\arccot(2\cosh\hat D))}.$$ 
Recall that for a metric $d$ on $S_g$ we defined $\alpha(d)=\arccot(2\cosh(\diam(S_g, d)).$ Hence, $M$ serves as an upper bound for $\frac{1}{\sin^{2}\alpha(d_i)}$ for all intermediate $d_i$. Also $V' \subset V(d_i)$ for all $i$.

We apply Corollary~\ref{glued2} and get $h'_{m} \in H(d_{m}, V'_{m})$ such that

(1) $ h'_{m,v} \geq \hat h_{m,v}$ for each $v \in V'$;

(2) $S(d_m, V'_m, P'_{m}) \geq S(\hat P)-\left(\frac{(\pi+\Delta)\sinh\Delta}{\Delta}+1\right)2\pi M\nu_{v_m}(d)\sinh\Delta.$

Now we have $\nu_{w_{m+1}}(d_{m})=0$. Thus, by Lemma~\ref{2pi} we get $\kappa_{w_{m+1}}(P'_{m})\leq 0$. Similarly to the proof of Main Lemma IIIA, we decrease the height of $w_{m+1}$ until we get $P_{m}=P(d_{m}, V'_{m}, h_{m})$ with $\kappa_{w_{m+1}}(P_{m})=0$, so $P_{m}=P(d_{m}, V_{m}, h_{m})$ and $S(P_{m}) \geq S(P'_{m})$.

After repeating these operations we obtain $h \in H(d, V)$ with 

(1) $h_v \geq \hat h_v$ for each $v \in V'$;

(2) $S(d, V, h) \geq S(\hat P) -\left(\frac{(\pi+\Delta)\sinh\Delta}{\Delta}+1\right)2\pi M\nu(S_g, d)\sinh\Delta$.

We have $$\nu(S_g, d)=\area(S_g, d)+2\pi(2-2g)\leq A+2\pi(2-2g).$$ We see that the decrement is bounded by $o(\Delta)$ as $\Delta \rightarrow 0$ (observe that $M$ is bounded for fixed $D,A$ and $\Delta \rightarrow 0$). So we can choose $\Delta$ such that $S(d,V,h) \geq S(\hat P)-\e$. This finishes the proof.
\end{proof}

\subsubsection{Dissolving curvature of swept triangles}
\label{flat}

In order to prove Main Lemma IIIC we should describe how to transform strict swept triangles with small curvature to hyperbolic ones. We apply the results from previous sections to obtain a simultaneous transformation of a Fuchsian cone-manifold and to bound the change of the discrete curvature.

\begin{dfn}
A swept triangle is called \emph{tetrahedral} if it is strict and is isometric to the surface of a tetrahedron (possibly degenerated) with a face excluded. A swept triangle is called \emph{short} if its metric is short.
\end{dfn}

A tetrahedral swept triangle might be non-short if it has some angles at least $\pi$. Denote a swept triangle by $OA_1A_2A_3$ as in Subsection~\ref{conesec}. Recall that $x_i$ is the length $OA_i$ and $\beta_i$ is $\angle OA_{i-1}A_{i+1}$. A strict swept triangle is tetrahedral if and only if $\beta_i$ satisfy the non-strict triangle inequalities. Indeed, if they do, then there exists a solid trihedral cone in $\H^3$ (possibly degenerated) with angles $\beta_i$ and we take the points at distances $x_i$ from the apex at the corresponding rays. The converse direction is obvious. It is clear that the tetrahedron is unique up to isometry.

\begin{lm}
A short strict swept triangle $OA_1A_2A_3$ is tetrahedral.
\end{lm}

\begin{proof}
Assume that $\beta_3> \beta_1 +\beta_2$. Develop the triangles $OA_2A_3$ and $OA_1A_3$ to $\H^2$ as the convex polygon $OA_2A_3A'_1$. Due to convexity, the segment $A'_1A_2 $ belongs to the union of the triangles. Then develop the triangle $OA_1A_2$ as the triangle $O A_2 A''_1$ such that $A''_1$ lies on the same side from $O A_2 $ as $A'_1$ (recall that $\beta_3<\pi$). As $\beta_3> \beta_1 +\beta_2$, it is easy to see that $A''_1A_2 >A'_1A_2 $, which contradicts that $A_1A_2$ is a shortest path.
\end{proof}

Let $l_1$, $l_2$, $l_3$ be three numbers satisfying the strict triangle inequalities and $\lambda_1^0$, $\lambda_2^0$, $\lambda_3^0$ be the angles of the corresponding hyperbolic triangle. Define $TCT(l_1,l_2,l_3)\subset \R^3$ as the set of triples $Z= (x_1, x_2, x_3)$ such that there exists a tetrahedral swept triangle $OA_1A_2A_3$, where the side lengths of $A_1A_2A_3$ are equal to $l_1$, $l_2$ and $l_3$ and the lengths $OA_1$, $OA_2$ and $OA_3$ are equal to $x_1$, $x_2$ and $x_3$. The boundary of $TCT(l_1,l_2,l_3)$ corresponds to the degenerations of the triangle inequalities for the angles $\beta_i$. The closure of $TCT(l_1,l_2,l_3)$ is formed by adding points with $\beta_1+\beta_2+\beta_3$, i.e., that correspond to hyperbolic triangles with $O$ belonging to $A_1A_2A_3$. We consider $\lambda_j$ as continuous functions over $TCT(l_1,l_2,l_3)$. Note that $\lambda_j(Z)\geq \lambda_j^0$ for each $Z$.

Consider $OA_1A_2A_3$ as a tetrahedron. Let $\tilde\tau$ be the distance from point $O$ to the triangle $A_1A_2A_3$. It defines another continuous function over $TCT(l_1,l_2,l_3)$. Moreover, it can be continuously extended to the closure of $TCT(l_1,l_2,l_3)$ in $\R^3$. It is equal to zero at all points from the closure that are not in $TCT(l_1,l_2,l_3)$.

\begin{lm}
\label{smallflat}
For every $\xi > 0$ there exists $\delta >0$ such that if for all $j=1,2,3$ we have $\lambda_j(Z)<\lambda_j^0+\delta$, then $\tilde\tau(Z)<\xi$.
\end{lm}

\begin{proof}
For any $Z \in TCT(l_1,l_2,l_3)$ we have $\lambda_j(Z) > \lambda_j^0$ for some $j$. Moreover, one can show that for every $\xi>0$ there exists $C>0$ such that if $x_i>C$ for some $i$, then $\lambda_{j}(Z)-\lambda_{j}^0>\xi$ for some $j$.

Now suppose the converse to the statement of the lemma. Then there exists $\xi>0$ and a sequence $Z_n$ such that $\tau(Z_n)\geq \xi$ for all $n$, but $\lambda_j(Z_n) \rightarrow \lambda_j^0$ for all $j$. Then $x_{n,i} \leq C$ for all $i$, all sufficiently large $n$ and a constant $C$ defined above. Then $Z_n$ belongs to a compact in the closure of $TCT(l_1,l_2,l_3)$ and has a limit point $Z$. We have $\lambda_j(Z)=\lambda_j^0$ for all $j$. This means that $Z$ does not belong to $TCT(l_1,l_2,l_3)$. However, this implies $\tilde\tau(Z)=0$.
\end{proof}

Consider $Z_0 \in TCT(l_1,l_2,l_3)$, the corresponding tetrahedron $OA_1A_2A_3$ and a point $O'$ in the interior of $A_1A_2A_3$. Let $\tau=OO'$ and $O_t \in OO'$ be the point at distance $t$ from $O$. This gives us a path of swept triangles parametrized by $t\in [0;\tau]$. Note that $\tau$ can be chosen arbitrarily close to $\tilde\tau(Z_0)$.

\begin{lm}
\label{angledecrease}
We have $\angle O_tA_1A_2+ \angle O_tA_1A_3 < \angle OA_1A_2+ \angle OA_1A_3$.
\end{lm}

\begin{proof}
Consider the spherical link of the tetrahedron $OA_1A_2A_3$ at point $A$. The inequality easily follows from the fact that the resulting spherical triangle for $O_tA_1A_2A_3$ is strictly contained in the spherical triangle for $OA_1A_2A_3$ and the fact that the perimeter of convex figures decreases under inclusion.
\end{proof}

In other words, the angles of swept triangles determined by $Z_t$ are strictly decreasing.

Consider now a Fuchsian cone-manifold $P=P(d_0, V, h)$, $d_0 \in \mathfrak D_{sc}(V)$. Assume that $T$ is a convex short swept triangle in $(S_g, d_0)$ with the conical point $w$ and vertices from $V$. We perform the deformation of $(T, d_0)$ transforming it to a hyperbolic triangle as described above and obtain a path of metrics $d_t \in \mathfrak D(V)$, $t\in [0;\tau]$. Until the end of this section we denote the realization of $T$ in $d_t$ by $X_t$. Because the angles of $X_t$ are strictly decreasing and all other angles remain the same, we have that $d_t \in \mathfrak D_{sc}(V)$ except the last point $d_{\tau} \in \mathfrak D_c(V)$ because $\nu_w(d_{\tau})=0$. We would like to transform simultaneously the Fuchsian cone-manifold with the help of Lemma~\ref{collect} and to control the change of $S$. To this purpose we need an analogue of Lemma~\ref{metricch}.

It is easy to see that the diameter of metrics $d_t$ does not increase. Let $\Delta$ be an upper bound for the diameter of $X_0$, so it is an upper bound for the diameters of $X_t$. We denote the three vertices of $T$ by $v_1$, $v_2$, $v_3$, the shortest paths from $w$ to them by $\psi_1$, $\psi_2$ and $\psi_3$, the edges of $T$ by $e_1$, $e_2$ and $e_3$ and the angles of $X_t$ by $\lambda_1(t)$, $\lambda_2(t)$, $\lambda_3(t)$ respectively. Define $\Lambda(t):=\lambda_1(t)+\lambda_2(t)+\lambda_3(t)$.

\begin{lm}
\label{flatSbound}
Let $0\leq t' < t'' \leq \tau$. Assume that there is 
$$h \in \bigcap_{t \in [t'; t'']} H(d_t, V)$$
and a triangulation $\mathcal T_h$ with $V(\mathcal T_h)=V$ such that every convex Fuchsian cone-manifold $P_t=P(d_t, V, h)$ is compatible with $\mathcal T_h$. Then
$$\frac{-1}{\sin^{2}\alpha(d_0)}\left((\pi+\Delta)\Delta\left(\Lambda(t')-\Lambda(t'')\right)+2\pi(t''-t')\right)\leq$$ $$  \leq S(d_{t''}, V, h) - S(d_{t'}, V, h)\leq \frac{ 2\pi(t''-t')}{\sin\alpha(d_0)}.$$
\end{lm}

\begin{proof}
The proof is similar to the proof of Lemma~\ref{metricch}. By $S(t)$ denote $S(d_t, V, h)$. We have
$$S(t'')-S(t')= \int_{t'}^{t''} \sum_{e \in E(\mathcal T_h)} (\pi  - \theta_e)\dot l_e dt.$$

First, we say that we consider \emph{the simple case} and assume that each edge of $\mathcal T_h$ that changes its the length during the deformation belongs to only one of the following four types:

(1) edges that intersect $\psi_1$ once;

(2) edges that intersect $\psi_2$ once;

(3) edges that intersect $\psi_3$ once;

(4) edges that have $w$ as one (and only one) of the endpoints.

Denote these types by $E_1$, $E_2$, $E_3$ and $E_4$. Triangulate $X_0$ naturally into 3 triangles and extend this to a triangulation $\mathcal T_0$ of $(S_g, d_0)$ with vertex set $V$. Then $d_t \in \mathcal D_{sc}(V, \mathcal T_0)$ for all $t \in [0; \tau)$, and only three edges of $\mathcal T_0$ change their lengths: $\psi_1, \psi_2$ and $\psi_3$.

Similarly to Lemma~\ref{metricch} for $e \in E_1$ we consider the path of triangles of $\mathcal T_0$ along $e$, develop it to $\H^2$ and denote by $y_e(t)$ the hyperbolic sine of the length of the perpendicular from the image of $v_1$ under the developing map to the line containing the image of $e$. Then we obtain 
$$\frac{\partial l_e}{\partial \lambda_1}  = y_e ,~~~ \frac{\partial l_e}{\partial \lambda_2}  =\frac{\partial l_e}{\partial \lambda_3}  =0.$$

By Lemma~\ref{angledecrease}, $\lambda_1$ is strictly decreasing, therefore $\dot l_e(t)\leq 0$. Then as in the proof of Lemma~\ref{metricch} using the bound $\len(\psi_1, d_t)<\Delta$ we get the estimate
$$-\frac{(\pi+\Delta)\Delta}{\sin^2\alpha(d_0)}(\lambda_1(t')-\lambda_1(t'')) \leq \int_{t'}^{t''} \sum_{e \in E_1} (\pi-\theta_e) \dot l dt \leq 0.$$

Similar estimates hold for edges of the second and the third type. 

We proceed with the fourth type. By Lemma~\ref{section} and Lemma~\ref{vertexangle} we have $$\sum_{e \in E_4} (\pi - \theta_e)\leq \omega_w \leq \frac{2\pi}{\sin\alpha(d_0)}.$$

We claim that for $e \in E_4$ we have $$-1\leq \dot l_e(t) \leq 1.$$

Assume that $e$ leaves $X_t$ through $e_3$. Let $v_4$ be the second vertex of $e$. Develop to $\H^2$ the path of triangles of $\mathcal T_0$ intersecting $e$ in the moment $t$. Let $A_1$ and $A_2$ be the images of $v_1$ and $v_2$, $B_1$ and $C$ be the images of $w$ and $v_4$ respectively. For some $\Delta t$ develop the path of triangles intersecting $e$ in the moment $t+\Delta t$ in such a way that the image of each triangle coincides with its image under the previous developments, except the first one, which is developed to a triangle $A_1A_2B_2$. We have $l_e(t)=B_1C$, $l_e(t+\Delta t)=B_2C$, so $$-B_1B_2 \leq l_e(t+\Delta t)-l_e(t) \leq B_1B_2.$$

One can see that $B_1B_2 \leq \Delta t$. Indeed, there is a hyperbolic tetrahedron $A_1A_2O_1O_2$ with $A_1A_2O_1=A_1A_2B_1$, $A_1A_2O_2=A_1A_2B_2$ and $O_1O_2=\Delta t$ (this tetrahedron is a subset of the tetrahedron $OA_1A_2A_3$ defining our deformation). Clearly if we rotate the triangle $A_1A_2O_2$ around the edge $A_1A_2$, the distance $O_1O_2$ is minimized, when $O_2$ is in the plane $A_1A_2O_1$. But in this moment $O_1O_2=B_1B_2$. Therefore, $B_1B_2 \leq \Delta t$ and we obtain the estimate of $\dot l_e(t)$. Thus, we have
$$\frac{-2\pi}{\sin\alpha(d_0)} \leq \int_{t'}^{t''} \sum_{e \in E_4} (\pi - \theta_e) \dot l_e dt \leq \frac{2\pi}{\sin\alpha(d_0)}.$$

The proof of the general case is done via elementary deformations of the four types as in the proof of Lemma~\ref{metricch}.
Consider the general case and an edge $e \in E(\mathcal T_h)$  with an orientation. As before, develop the path of triangles intersecting $e$ at moment $t$ as a triangulated polygon $R_0$ in $\H^2$. Enumerate triangles of $R_0$ with respect to the orientation. By $k$ denote the number of edges in $R_0$ and by $z_{1}(t), \ldots, z_{k}(t)$ denote their lengths. So the length $l_e(t)$ can be considered as a function $l_e(t)=l_e(z_{1}(t), \ldots, z_{k}(t)).$ 

All $z_i$ are constant except those that are images of some $\psi_j$ under the developing map. The length of each $\psi_j$ is completely determined by $l_j$ and $\lambda_j$ due to Lemma~\ref{conerigid}. Assume that $l_j$ are fixed, but $\lambda_j$ vary. Then $z_i$ can be considered as functions of $\lambda_j$. In this way we view each $z_i(t)$ as a function $z_i(\lambda_1(t), \lambda_2(t), \lambda_3(t)).$ We consider triple coordinate derivatives 
$$\frac{\partial l_e}{\partial z_{i}}\frac{\partial z_{i}}{\partial \lambda_j}\frac{\partial \lambda_j}{\partial t}(t).$$ 

Let $e$ enter $X_t$ from outside and then leaves $X_t$ again. This corresponds to some development of $X_t$ in $R_0$. Let $I$ be the set of indices of the edges of $R_0$ corresponding to the interior edges of $X_t$ during such event. Define $$L(e, I, j)(t):=\sum_{i\in I}\frac{\partial l_e}{\partial z_{i}}\frac{\partial z_{i}}{\partial \lambda_j}\frac{\partial \lambda_j}{\partial t}(t).$$ 
We say that $L(e, I, j)$ is an elementary deformation of the type $(j)$, $j=1,2,3$. Similarly to the proof of Lemma~\ref{metricch} one can see that an elementary deformation of type $(j)$ can be treated similarly to the edge of type $(j)$ in the simple case. We note that if, e.g., $e$ enters $X_t$, crosses $\psi_1$ and then leaves $X_t$, then
$$L(e, I, 2)=L(e, I, 3)=0$$
similarly to edges of type (1) in the simple case. We also remark that when $e$ enters $X_t$, it can cross each $\psi_j$ at most once, as otherwise it self-intersects inside $X_t$.

Finally, if $e$ starts or ends in $w$ and $i_1$, $i_2$ are the indices of the respective edges of $R_0$, then for $I=\{i_1, i_2\}$ we define
$$L(e, I)(t):=\sum_{i \in I}\frac{\partial l_e}{\partial z_{i}}\frac{\partial z_{i}}{\partial t}(t)$$
and say that it is an elementary deformation of type (4). 

We have
$$\dot l_e(t) = \sum_{i=1}^k\sum_{j=1}^3 \frac{\partial l_e}{\partial z_{i}}\frac{\partial z_{i}}{\partial \lambda_j}\frac{\partial \lambda_j}{\partial t}(t).$$
Some of these triple derivatives form elementary deformations, all others are equal zero.
Consider elementary deformations of every edge and decompose $$\sum_{e \in E(\mathcal T_h)} (\pi  - \theta_e(t))\dot l_e(t)$$ into sums of elementary deformations of each type. Inside each type their sum is estimated in the same way as edges of this type in the simple case. This gives the same estimate for $S$ as in the simple case.
\end{proof}

Now we can transform a convex Fuchsian cone-manifold $P_0=P(d_0, V, h_0)$ along the path of metrics $d_t$, $t \in [0;\tau]$, described above. With the help of Lemma~\ref{collect} and the lower bound of Lemma~\ref{flatSbound} we get

\begin{crl}
\label{flatned1}
Assume that for some $\xi>0$ and for all $j=1,2,3$ we have $\lambda_j(0)-\lambda_j(\tau)<\xi$ and $\tau<\xi$.

Then there exists $h_{\tau} \in H(d_{\tau}, V)$ such that

(1) $ h_{\tau}(v) \geq  h_0(v)$ for every $v \in V$;

(2) $S(d_{\tau}, V, h_{\tau}) \geq S(d_0, V, h_0)-\frac{\left(3(\pi+\Delta)\Delta+2\pi\right)\xi}{\sin^{2}\alpha(d_0)}$.
\end{crl}

Now consider the same path of metrics in the inverse direction, i.e., $d'(t) := d(\tau-t)$ for $t \in [0;\tau]$. Let $P'_0=P(d'_0, V, h_0)$ be a convex Fuchsian cone-manifold. Then from Lemma~\ref{collect} and the upper bound Lemma~\ref{flatSbound} we get 

\begin{crl}
\label{flatned2}
There exists $h_{\tau} \in H(d'_{\tau}, V)$ such that

(1) $ h_{\tau}(v) \geq  h_0(v)$ for every $v \in V$;

(2) $S(d'_{\tau}, V, h_{\tau}) \geq S(d'_0, V, h_0)-\frac{2\pi\tau}{\sin\alpha(d'_{\tau})}.$
\end{crl}

\subsubsection{Proof of Main Lemma IIIC}

Let $\mathfrak S(\mathcal T)$ be the set of convex cone-metrics swept with respect to $\mathcal T$ (defined up to isometry isotopic to identity with respect to $V(\mathcal T)$). Hence, $\hat d \in \mathfrak S(\mathcal T)$. By $\mathfrak S(\mathcal T, \hat d, \delta) \subset \mathfrak S(\mathcal T)$ denote the set of metrics $d \in \mathfrak S(\mathcal T)$ such that for each triangle $T$ of $\mathcal T$ we have $||T(d)-T(\hat d)||_{\infty}<\delta$. By the assumptions of Main Lemma IIIC we have $d^1$, $d^2 \in \mathfrak S(\mathcal T, \hat d, \delta).$

Define $$\hat V:=V(\hat d)\cup V(\mathcal T),~~~D:=2\diam(S_g, \hat d),$$ $$M:=\sin^{-2}(\arccot D),~~~\Delta:=2\max_{T \in \mathcal T}\diam(T, \hat d).$$ By $\hat T$ we denote the triangulation obtained by refining each triangle of $\mathcal T$ that has a conical point in $\hat d$ into three. Naturally $\hat d \in \mathfrak D_c(\hat V, \hat{\mathcal T})$. Recall from Subsection~\ref{polmetr} that with the help of the edge-length map, $\mathfrak D_c(\hat V, \hat{\mathcal T})$ is considered as a subset of $\R^N$, where $N=|E(\hat{\mathcal T})|$, endowed with $l_{\infty}$-metric.

Our strategy is as follows. If a triangle $T$ of $\mathcal T$ is non-strict in $\hat d$, then its curvature in $d^1$, $d^2$ is small provided that $\delta$ is small. Then we use the results of Subsection~\ref{flat} to dissolve this curvature. Doing this for each triangle of this kind we obtain two metrics in $\mathfrak D_{c}(\hat V, \hat{\mathcal T})$ that are very close with respect to the $l_{\infty}$-metric. By an indirect argument we connect them with a short path and transform the Fuchsian cone-manifold along the path.

There exists $\delta_0$ such that for each $d \in \mathfrak S(\mathcal T, \hat d, \delta_0)$ the following hold

(1) $V(\hat d) \subseteq V(d)$;

(2) $\diam (S_g, d) < D$;

(3) $\max_{T \in \mathcal T}\diam(T, d)< \Delta$.

Let $k$ be the number of triangles of $\mathcal T$ that are non-strict in $\hat d$. Choose $\xi>0$ small enough so that $$M\left(3(\pi+\Delta)\Delta+2\pi\right)\xi<\frac{\e}{3k}.$$
Then also $2\pi M \xi < \e/3k$.

Take a triangle $T$ of $\mathcal T$ that is non-strict in $\hat d$. By Lemma~\ref{smallflat} there exists $\delta_1(T)$ such that if $d \in \mathfrak S(\mathcal T, \hat d, \delta_1(T))$, then $\tau < \xi$, where $\tau$ is the distance from the conical point $O$ to an interior point of the triangle $A_1A_2A_3$ for a tetrahedral realization of the swept triangle $(T, d)$. We also assume that $\delta_1(T) < \xi$. Hence, if $d \in \mathfrak S(\mathcal T, \hat d, \delta_1(T))$, $\lambda(d)$ is an angle of $(T, d)$ and $\lambda(\hat d)$ is the same angle of $(T, \hat d)$, then $\lambda(d) < \lambda(\hat d)+\xi$.

Take $\delta_1$ as the minimum of $\delta_1(T)$ over all such triangles $T$.

For $\sigma_0>0$ define $\mathfrak B_{c}(\hat d, \sigma_0)$, $\mathfrak B_{sc}(\hat d, \sigma_0)$ as in Section~\ref{polmetr} to be the intersections of the open ball in $(\R^N, l_{\infty})$ of radius $\sigma_0$ with center at $\hat d$ with $\mathfrak D_{c}(\hat V, \hat{\mathcal T})$, $\mathfrak D_{sc}(\hat V, \hat{\mathcal T})$ respectively. By Corollary~\ref{connect} the set $\mathfrak B_{sc}(\hat d, \sigma_0)$ is connected for sufficiently small $\sigma_0$. 
Its boundary is locally piecewise analytic, thus $\mathfrak B_{sc}(\hat d, \sigma_0)$ is an open, connected and bounded subset of $\R^N$ with Lipschitz boundary. From~\cite[Chapters 2.5.1-2.5.2]{BrBrBook} it is quasiconvex: there exists a constant $C_1'=C_1'(\hat d, \mathcal T, \sigma_0)$ such that every two points of $\mathfrak B_{sc}(\hat d, \sigma_0)$ at the distance $\sigma$ can be connected through $\mathfrak B_{sc}(\hat d, \sigma_0)$ by a path of length at most $C_1'\sigma$. By Lemma~\ref{decreaseangles} the set $\mathfrak B_{c}(\hat d, \sigma_0)$ belongs to the closure of $\mathfrak B_{sc}(\hat d, \sigma_0)$ Then for any $C_1>C_1'$ the set $\mathfrak B_{c}(\hat d, \sigma_0)$ is strictly $C_1$-quasiconvex, i.e., every two points of $\mathfrak B_{c}(\hat d, \sigma_0)$ at the distance $\sigma$ can be connected through $\mathfrak B_{sc}(\hat d, \sigma_0)$ by a path of length at most $C_1\sigma$. Take any such $C_1$. We also assume that the closure $\overline{\mathfrak B}_{c}(\hat d, \sigma_0)$ is in $\mathfrak D_c(\hat V, \hat{\mathcal T})$ (so no triangles become degenerate in the closure). 

Recall that 
$$H_S(d, V, m, K)=\{h \in H(d, V): \min_{v \in V}h_v \geq m, S(d, V, h) \geq K\},$$
$$\mathfrak H_S(\overline{\mathfrak B}_{c}(\hat d, \sigma_0), m, K) = \{ (d, h): d \in \overline{\mathfrak B}_{c}(\hat d, \sigma_0), h \in H_S(d,V, m, K)\}\subset \mathfrak H(\mathfrak D_c(\hat V, \hat{\mathcal T})).$$ 
By Lemma~\ref{S1der} $S$ is continuously differentiable over $\mathfrak H(\mathfrak D_c(\hat V, \hat{\mathcal T}))$. Put
$$K:=S(d^1, V(d^1), h^1)-\e.$$ As $\overline{\mathfrak B}_{c}(\hat d, \sigma_0)$ is compact, it follows from Lemma~\ref{H_U} that so is $\mathfrak H_S(\overline{\mathfrak B}_{c}(\hat d, \sigma_0), m, K)$. Hence, there exists a constant $C_2=C_2(\hat d, \mathcal T)$ such that $S$ is $C_2$-Lipschitz over $\mathfrak H_S(\overline{\mathfrak B}_{c}(\hat d, \sigma_0), m, K)$. 

Choose $\sigma>0$ such that $C_1C_2\sigma < \e/3$ and $\sigma<\sigma_0$. By Corollary~\ref{homeo} we can choose $\delta_2$ sufficiently small such that if $d \in \mathfrak S(\mathcal T, \hat d, \delta_2)$, then $d \in {\mathfrak B}_{c}(\hat d, \sigma/2)$.

Now put $\delta:=\min\{\delta_0, \delta_1, \delta_2\}$ and take $d^1, d^2 \in \mathfrak S(\mathcal T, \hat d, \delta)$. We have $$V(\hat d)\subseteq (V(d^1)\cap V(d^2)).$$

Let $T$ be a triangle of $\mathcal T$ such that $T$ is strict in $d^1$, but non-strict in $\hat d$. Denote the conical point of $T$ in $d_1$ by $w$.
Define $d_0:=d^1$, $V_0:=V(d^1)$. Consider a deformation $d_t$ through $\mathfrak D_{sc}(V_0)$, $t \in [0; \tau)$, described in Section~\ref{flat} that dissolves the curvature of $w$. It follows from the discussion above that $\tau<\xi$ and the variations of angles of $T$ are smaller than $\xi$. Denote the resulting metric by $d_1$. By Corollary~\ref{flatned1} we obtain a convex Fuchsian cone-manifold $P'_1=P(d_1, V_0, h'_1)$ such that

(1) $ h'_{1,v} \geq  h^1_v$ for every $v \in V(\hat d)$;

(2) $S(d_1, V_0,  h'_1) \geq S(d^1, V^1, h^1) - \frac{\left(3(\pi+\Delta)\Delta+2\pi\right)\xi}{\sin^{2}\alpha(d_0)}\geq S(d^1, V^1, h^1) - \frac{\e}{3k}$.

Then we reduce the height of $w$ until its particle curvature disappears as in the proofs of Main Lemmas IIIA-B. This increases $S$ and produces $h_1 \in H(d_1, V_1)$, where $V_1=V_0 \backslash \{w\}$, such that
$$S(d_1, V_1,  h_1) \geq S(d_1, V_0,  h'_1).$$

We do this for every such triangle and obtain a metric $\hat d^1 \in \mathfrak D_{c}(\hat V, \hat{\mathcal T})$. We note that the diameters of all intermediate metrics are smaller than $D$. We get $\hat h^1 \in H(\hat d^1, \hat V)$ such that

(1) $ \hat h^1_v \geq  h^1_v$ for every $v \in V(\hat d)$;

(2) $S(\hat d^1, \hat V, \hat h^1) \geq S(d^1, V^1, h^1) - \e/3$.

Now we take $d^2$ and replace each strict swept triangle of $\mathcal T$ that is non-strict in $\hat d$ by the hyperbolic triangle with the same side lengths (of $d^2$). Denote the obtained metric by $\hat d^2$. Note that $\hat d^1, \hat d^2 \in \mathfrak S(\mathcal T, \hat d, \delta)$. Thus, $\hat d^1, \hat d^2 \in {\mathfrak B}_{c}(\hat d, \sigma/2)$. Hence, they can be connected through $\mathfrak B_{sc}(\hat d, \sigma_0)$ by a curve $d_t$ of length at most $C_1\sigma$.
 
We take the Fuchsian cone-manifold $\hat P^1 = P(\hat d^1, \hat V, \hat h^1)$, which we obtained before. We want to show that there exist a height function $\hat h^{21} \in H(\hat d^2, \hat V)$ such that 

(1) $ \hat h^{21}_v\geq   \hat h^1_v$ for every $v \in V(\hat d)$;

(2) $S(\hat d^{2}, \hat V, \hat h^{21}) \geq S(\hat d^{1}, \hat V, \hat h^{1}) - C_1C_2\sigma \geq S(\hat d^{1}, \hat V, \hat h^{1}) - \e/3$.

This follows from Lemma~\ref{collectC} together with the fact that $S$ is $C_2$-Lipschitzian over $\mathfrak H_S(\overline{\mathfrak B}_{c}(\hat d, \sigma_0), m, K)$.

It remains to transform $\hat d^2$ to $d^2$ as in Section~\ref{flat} and to transform the Fuchsian cone-manifold $P^{21}=P(\hat d^{2}, \hat V, \hat h^{21})$ with the help of Corollary~\ref{flatned2}. We obtain $h^{21} \in H(d^2, V^2)$ such that

(1) $  h^{21}_v\geq  \hat h^{21}_v$ for every $v \in V(\hat d)$;

(2) $S(d^{2}, V^2, h^{21}) \geq S(\hat d^{2}, \hat V, \hat h^{21})-k2\pi M \xi\geq S(\hat d^{2}, \hat V, \hat h^{21}) - \e/3$.

In total we get

(1) $h^{21}_v \geq h^1_v$ for each $v \in V(\hat d)$;

(2) $S(d^2, V^2, h^{21}) \geq S(d^1, V^1, h^1) - \e$.

This finishes the proof.

\appendix
\section{Appendix: A discussion of Volkov's paper~\cite{Vol}}
\label{appendix}

Due to several requests, here we attach a discussion of gaps in Volkov's paper~\cite{Vol}. We refer to its English translation in Subsection 12.1 of~\cite{Ale2}.

Volkov's paper is divided in two sections: about convex caps and about convex polyhedra. The latter is claimed to be reduced to the former. The proof in the first section is based on two main lemmas: Lemma 1 and Lemma 2 (corresponding to our Main Lemma I and Main Lemmas IIIA-B respectively). Our main concerns are about Lemma 2. There Volkov has two convex caps with close boundary metrics (``developments'' in the translation). As in our proofs of Main Lemmas III, he transforms the first boundary metric to the second by some elementary operations. Along this he transforms the first convex cap in the class of convex caps with cone-singularities (``generalized caps'' in the translation) in a controlled way, and he bounds the change of the discrete curvature (``the total curvature'').

The main trouble is located on p. 471 in paragraph 7 (``Now, the required deformation proceeds...''). Volkov writes ``using the standard trick of cutting-off digons, we ``drive'' the curvature from all the triangles splitting $S^0$ and $S^1$.'' He gives no further explanations of this procedure and we do not have a clear understanding what he means here. Volkov adds further ``The deformation of the strips of one of these developments to the strips of the other is a nondifficult, though laborious, essentially planar problem. We will not dwell upon its solution.''

We were unable to decipher, what was Volkov's process of ``driving'' the curvature and how he planned to deform the obtained strips. Instead we propose to merge curvature in each triangle to a point and then to transform the obtained cone-triangles of one metric to those of the second with the help of an indirect argument. We realized this in Main Lemma II and Main Lemma IIIC.

The next mysterious place in Volkov's work is p. 482 starting from paragraph 2 (bounding the change of discrete curvature). We do not understand Volkov's procedure of ``leveling'' and we are unable to confirm his observations starting from the words ``It is easy to show...''. We spent Subsections~\ref{boundchange}--\ref{slopes} to obtain bounds that play a similar role in our paper.

Except this, in Volkov's proof of Lemma 5 he skips the most non-trivial case when some $OB_k > \pi/2$. We remark that his Lemma 5 follows from our Lemma~\ref{spherarea} and the isoperimetric inequality for cone-surfaces of non-positive curvature~\cite{Izm2}.

Finally, Volkov's reduction of the case of convex polyhedra to the case of convex caps in Subsection 2.4 is very sketchy and requires additional explanations.

\bibliographystyle{abbrv}
\bibliography{biblio}

\begin{thebibliography}{10}

\bibitem{AlBi}
S.~{Alexander} and R.~L. {Bishop}.
\newblock {\(\mathcal FK\)-convex functions on metric spaces.}
\newblock {\em {Manuscr. Math.}}, 110(1):115--133, 2003.

\bibitem{AKP}
S.~Alexander, V.~Kapovitch, and A.~Petrunin.
\newblock Alexandrov geometry: preliminary version no. 1, 2019.
\newblock ArXiv e-print 1903.08539.

\bibitem{Ale1}
A.~D. Alexandrov.
\newblock Existence of a convex polyhedron and of a convex surface with a given
  metric.
\newblock {\em Rec. Math. [Mat. Sbornik] N.S.}, 11(53):15--65, 1942.

\bibitem{Ale2}
A.~D. Alexandrov.
\newblock {\em Convex polyhedra}.
\newblock Springer Monographs in Mathematics. 2005.

\bibitem{Ale3}
A.~D. Alexandrov.
\newblock {\em Intrinsic geometry of convex surfaces}.
\newblock Chapman \& Hall/CRC, Boca Raton, FL, 2006.
\newblock A. {D}. {A}lexandrov selected works. {P}art {II}.

\bibitem{AlZa}
A.~D. Alexandrov and V.~A. Zalgaller.
\newblock {\em Intrinsic geometry of surfaces}.
\newblock American Mathematical Society, Providence, R.I., 1967.

\bibitem{BeCa}
J.~Bertrand and P.~Castillon.
\newblock Prescribing the gauss curvature of convex bodies in hyperbolic space,
  2019.
\newblock ArXiv e-print 1903.06502.

\bibitem{BoIz}
A.~I. Bobenko and I.~Izmestiev.
\newblock Alexandrov's theorem, weighted {D}elaunay triangulations, and mixed
  volumes.
\newblock {\em Ann. Inst. Fourier (Grenoble)}, 58(2):447--505, 2008.

\bibitem{BPS}
A.~I. Bobenko, U.~Pinkall, and B.~A. Springborn.
\newblock Discrete conformal maps and ideal hyperbolic polyhedra.
\newblock {\em Geom. Topol.}, 19(4):2155--2215, 2015.

\bibitem{BLP}
M.~Boileau, B.~Leeb, and J.~Porti.
\newblock Geometrization of 3-dimensional orbifolds.
\newblock {\em Ann. of Math. (2)}, 162(1):195--290, 2005.

\bibitem{BCM}
J.~F. Brock, R.~D. Canary, and Y.~N. Minsky.
\newblock The classification of {K}leinian surface groups, {II}: {T}he ending
  lamination conjecture.
\newblock {\em Ann. of Math. (2)}, 176(1):1--149, 2012.

\bibitem{Bro3}
K.~Bromberg.
\newblock Projective structures with degenerate holonomy and the {B}ers density
  conjecture.
\newblock {\em Ann. of Math. (2)}, 166(1):77--93, 2007.

\bibitem{BrBrBook}
A.~{Brudnyi} and Y.~{Brudnyi}.
\newblock {\em {Methods of geometric analysis in extension and trace problems.
  Vol. 1.}}, volume 102.
\newblock Basel: Birkh\"auser, 2012.

\bibitem{Bru}
L.~Brunswic.
\newblock Alexandrov theorem for 2+1 flat radiant spacetimes, 2020.
\newblock ArXiv e-print 2012.01275.

\bibitem{BBI}
D.~Burago, Y.~Burago, and S.~Ivanov.
\newblock {\em A course in metric geometry}, volume~33 of {\em Graduate Studies
  in Mathematics}.
\newblock American Mathematical Society, Providence, RI, 2001.

\bibitem{CoV}
S.~{Cohn-Vossen}.
\newblock {Zwei S\"atze \"uber die Starrheit der Eifl\"achen.}
\newblock {\em {Nachr. Ges. Wiss. G\"ottingen, Math.-Phys. Kl.}},
  1927:125--134, 1927.

\bibitem{CoV2}
S.~{Cohn-Vossen}.
\newblock {Bending of surfaces in the large}.
\newblock {\em {Usp. Mat. Nauk}}, 1:33--76, 1936.

\bibitem{Fil1}
F.~Fillastre.
\newblock Polyhedral realisation of hyperbolic metrics with conical
  singularities on compact surfaces.
\newblock {\em Ann. Inst. Fourier (Grenoble)}, 57(1):163--195, 2007.

\bibitem{Fil2}
F.~Fillastre.
\newblock Polyhedral hyperbolic metrics on surfaces.
\newblock {\em Geom. Dedicata}, 134:177--196, 2008.

\bibitem{Fil3}
F.~Fillastre.
\newblock Fuchsian polyhedra in {L}orentzian space-forms.
\newblock {\em Math. Ann.}, 350(2):417--453, 2011.

\bibitem{FiIz1}
F.~Fillastre and I.~Izmestiev.
\newblock Hyperbolic cusps with convex polyhedral boundary.
\newblock {\em Geom. Topol.}, 13(1):457--492, 2009.

\bibitem{FiIz2}
F.~Fillastre and I.~Izmestiev.
\newblock Gauss images of hyperbolic cusps with convex polyhedral boundary.
\newblock {\em Trans. Amer. Math. Soc.}, 363(10):5481--5536, 2011.

\bibitem{FIV}
F.~Fillastre, I.~Izmestiev, and G.~Veronelli.
\newblock Hyperbolization of cusps with convex boundary.
\newblock {\em Manuscripta Math.}, 150(3-4):475--492, 2016.

\bibitem{FiSe}
F.~Fillastre and A.~Seppi.
\newblock Spherical, hyperbolic, and other projective geometries: convexity,
  duality, transitions.
\newblock In {\em Eighteen essays in non-{E}uclidean geometry}, volume~29 of
  {\em IRMA Lect. Math. Theor. Phys.}, pages 321--409. Eur. Math. Soc.,
  Z\"{u}rich, 2019.

\bibitem{FiSl}
F.~Fillastre and D.~Slutskiy.
\newblock Embeddings of non-positively curved compact surfaces in flat
  {L}orentzian manifolds.
\newblock {\em Math. Z.}, 291(1-2):149--178, 2019.

\bibitem{Gro}
M.~Gromov.
\newblock {\em Partial differential relations}, volume~9.
\newblock Springer-Verlag, Berlin, 1986.

\bibitem{GGLSW}
X.~Gu, R.~Guo, F.~Luo, J.~Sun, and T.~Wu.
\newblock A discrete uniformization theorem for polyhedral surfaces {II}.
\newblock {\em J. Differential Geom.}, 109(3):431--466, 2018.

\bibitem{GuLi}
P.~{Guan} and Y.~Y. {Li}.
\newblock {The Weyl problem with nonnegative Gauss curvature.}
\newblock {\em {J. Differ. Geom.}}, 39(2):331--342, 1994.

\bibitem{Hei}
E.~Heinz.
\newblock On {W}eyl's embedding problem.
\newblock {\em J. Math. Mech.}, 11:421--454, 1962.

\bibitem{Her}
G.~{Herglotz}.
\newblock {\"Uber die Starrheit der Eifl\"achen.}
\newblock {\em {Abh. Math. Semin. Univ. Hamb.}}, 15:127--129, 1943.

\bibitem{HoKe2}
C.~D. Hodgson and S.~P. Kerckhoff.
\newblock Universal bounds for hyperbolic {D}ehn surgery.
\newblock {\em Ann. of Math. (2)}, 162(1):367--421, 2005.

\bibitem{HoZu}
J.~{Hong} and C.~{Zuily}.
\newblock {Isometric embedding of the 2-sphere with non negative curvature in
  \(\mathbb{R}^ 3\).}
\newblock {\em {Math. Z.}}, 219(3):323--334, 1995.

\bibitem{ILTC}
C.~Indermitte, T.~Liebling, M.~Troyanov, and H.~Clemencon.
\newblock Voronoi diagrams on piecewise flat surfaces and an application to
  biological growth.
\newblock {\em Theoretical Computer Science}, 263(1):263 -- 274, 2001.

\bibitem{Izm1}
I.~Izmestiev.
\newblock A variational proof of {A}lexandrov's convex cap theorem.
\newblock {\em Discrete Comput. Geom.}, 40(4):561--585, 2008.

\bibitem{Izm3}
I.~Izmestiev.
\newblock Variational properties of the discrete {H}ilbert-{E}instein
  functional.
\newblock {\em Actes des rencontres du CIRM}, 3(1):151--157, 11 2013.

\bibitem{Izm2}
I.~{Izmestiev}.
\newblock {A simple proof of an isoperimetric inequality for Euclidean and
  hyperbolic cone-surfaces.}
\newblock {\em {Differ. Geom. Appl.}}, 43:95--101, 2015.

\bibitem{Lab}
H.~Labeni.
\newblock Realizing metrics of curvature $\geq -1$ on closed surfaces in
  fuchsian anti-de sitter manifolds.
\newblock {\em J. of the Australian Math. Soc., to appear}, 2021.

\bibitem{Lab2}
F.~Labourie.
\newblock M\'{e}triques prescrites sur le bord des vari\'{e}t\'{e}s
  hyperboliques de dimension {$3$}.
\newblock {\em J. Differential Geom.}, 35(3):609--626, 1992.

\bibitem{LaSc}
F.~Labourie and J.-M. Schlenker.
\newblock Surfaces convexes fuchsiennes dans les espaces lorentziens \`a
  courbure constante.
\newblock {\em Math. Ann.}, 316(3):465--483, 2000.

\bibitem{Lew}
H.~{Lewy}.
\newblock {On the existence of a closed convex surface realizing a given
  Riemannian metric.}
\newblock {\em {Proc. Natl. Acad. Sci. USA}}, 24:104--106, 1938.

\bibitem{Lib}
J.~Liberman.
\newblock Geodesic lines on convex surfaces.
\newblock {\em C. R. (Doklady) Acad. Sci. URSS (N.S.)}, 32:310--313, 1941.

\bibitem{Mac}
Y.~Machigashira.
\newblock The {G}aussian curvature of {A}lexandrov surfaces.
\newblock {\em J. Math. Soc. Japan}, 50(4):859--878, 1998.

\bibitem{Mard}
A.~{Marden}.
\newblock {The geometry of finitely generated Kleinian groups.}
\newblock {\em {Ann. Math. (2)}}, 99:383--462, 1974.

\bibitem{Mil2}
A.~D. Milka.
\newblock The lemma of {B}usemann and {F}eller in spherical and hyperbolic
  spaces.
\newblock {\em Ukrain. Geometr. Sb.}, (10):40--49, 1971.

\bibitem{Mil}
A.~D. Milka.
\newblock Unique determinacy of general closed convex surfaces in
  {L}oba\v{c}evski\u{\i} space.
\newblock {\em Ukrain. Geom. Sb.}, (23):99--107, iii, 1980.

\bibitem{Mil3}
A.~D. Milka.
\newblock Space-like convex surfaces in pseudo-{E}uclidean spaces.
\newblock In {\em Some questions of differential geometry in the large}, volume
  176 of {\em Amer. Math. Soc. Transl. Ser. 2}, pages 97--150. Amer. Math.
  Soc., Providence, RI, 1996.

\bibitem{Min}
Y.~Minsky.
\newblock The classification of {K}leinian surface groups. {I}. {M}odels and
  bounds.
\newblock {\em Ann. of Math. (2)}, 171(1):1--107, 2010.

\bibitem{Mos}
G.~D. Mostow.
\newblock Quasi-conformal mappings in {$n$}-space and the rigidity of
  hyperbolic space forms.
\newblock {\em Inst. Hautes \'{E}tudes Sci. Publ. Math.}, (34):53--104, 1968.

\bibitem{Nir}
L.~{Nirenberg}.
\newblock {The Weyl and Minkowski problems in differential geometry in the
  large.}
\newblock {\em {Commun. Pure Appl. Math.}}, 6:337--394, 1953.

\bibitem{Olo}
S.~Olovianishnikoff.
\newblock G\'{e}n\'{e}ralisation du th\'{e}or\`eme de {C}auchy sur les
  poly\`edres convexes.
\newblock {\em Rec. Math. [Mat. Sbornik] N.S.}, 18(60):441--446, 1946.

\bibitem{Pet}
A.~{Petrunin}.
\newblock {Applications of quasigeodesics and gradient curves.}
\newblock In {\em {Comparison geometry}}, pages 203--219. Cambridge: Cambridge
  University, 1997.

\bibitem{Pog1}
A.~V. Pogorelov.
\newblock The rigidity of convex surfaces.
\newblock {\em Trudy Mat. Inst. Steklov}, 29:99, 1949.

\bibitem{Pog2}
A.~V. Pogorelov.
\newblock {\em Extrinsic geometry of convex surfaces}.
\newblock American Mathematical Society, Providence, R.I., 1973.

\bibitem{Pra}
G.~Prasad.
\newblock Strong rigidity of {${\bf Q}$}-rank {$1$} lattices.
\newblock {\em Invent. Math.}, 21:255--286, 1973.

\bibitem{Pro1}
R.~Prosanov.
\newblock Ideal polyhedral surfaces in {F}uchsian manifolds.
\newblock {\em Geometria Dedicata}, 206(1):151--179, 2020.

\bibitem{Sac}
R.~{Sacksteder}.
\newblock {The rigidity of hypersurfaces.}
\newblock {\em {J. Math. Mech.}}, 11:929--939, 1962.

\bibitem{Sch2}
J.-M. Schlenker.
\newblock Surfaces convexes dans des espaces lorentziens \`a courbure
  constante.
\newblock {\em Comm. Anal. Geom.}, 4(1-2):285--331, 1996.

\bibitem{Sch4}
J.-M. Schlenker.
\newblock Hyperbolic manifolds with convex boundary.
\newblock {\em Invent. Math.}, 163(1):109--169, 2006.

\bibitem{Sch}
R.~Schneider.
\newblock {\em Convex bodies: the {B}runn-{M}inkowski theory}, volume 151 of
  {\em Encyclopedia of Mathematics and its Applications}.
\newblock Cambridge University Press, Cambridge, expanded edition, 2014.

\bibitem{Sec}
S.~Sechelmann.
\newblock {\em Alexandrov's polyhedron}.
\newblock
  https://webstart.discretization.de/Sechelmann-AlexandrovPolyhedra-2016-12-19/.

\bibitem{Sen}
E.~P. Sen'kin.
\newblock Bending of convex surfaces.
\newblock In {\em Problems in geometry, {V}ol. 10 ({R}ussian)}, pages 193--222,
  224. VINITI, Moscow, 1978.

\bibitem{Slu}
D.~Slutskiy.
\newblock Compact domains with prescribed convex boundary metrics in
  quasi-{F}uchsian manifolds.
\newblock {\em Bull. Soc. Math. France}, 146(2):309--353, 2018.

\bibitem{Spr}
B.~Springborn.
\newblock Hyperbolic polyhedra and discrete uniformization.
\newblock {\em Discrete and Computational Geometry}, 64:63--108, 2020.

\bibitem{Str}
S.~{Straszewicz}.
\newblock {\"Uber exponierte Punkte abgeschlossener Punktmengen.}
\newblock {\em {Fundam. Math.}}, 24:139--143, 1935.

\bibitem{Tam}
A.~Tamburelli.
\newblock Prescribing metrics on the boundary of anti--de {S}itter 3-manifolds.
\newblock {\em Int. Math. Res. Not. IMRN}, (5):1281--1313, 2018.

\bibitem{Thu1}
W.~P. Thurston.
\newblock {\em The geometry and topology of 3-manifold}.
\newblock Princeton University Press, Princeton, NJ, 1978.

\bibitem{Thu2}
W.~P. Thurston.
\newblock Three-dimensional manifolds, {K}leinian groups and hyperbolic
  geometry.
\newblock {\em Bull. Amer. Math. Soc. (N.S.)}, 6(3):357--381, 1982.

\bibitem{Vol2}
Y.~A. Volkov.
\newblock Existence of convex polyhedra with prescribed development {I}.
\newblock {\em Vestn. Leningr. Univ}, 15:75--86, 1960.

\bibitem{Vol}
Y.~A. Volkov.
\newblock An estimate of the deformation of a convex surface as a function of
  the change in its intrinsic metric.
\newblock {\em Ukrain. Geometr. Sb. Vyp.}, 5--6:44--69, 1968.

\bibitem{VolThe}
Y.~A. Volkov.
\newblock Existence of a polyhedron with prescribed development.
\newblock {\em Zap. Nauchn. Sem. S.-Peterburg. Otdel. Mat. Inst. Steklov.
  (POMI)}, 476(13):50--78, 2018.

\bibitem{Wey}
H.~{Weyl}.
\newblock {\"Uber die Bestimmung einer geschlossenen konvexen Fl\"ache durch
  ihr Linienelement.}
\newblock {Z\"urich. Naturf. Ges. 61, 40-72 (1916).}, 1916.

\bibitem{Zhi}
O.~K. {Zhitomirsky}.
\newblock {Sur la non-flexibilite des ovaloides.}
\newblock {\em {C. R. (Dokl.) Acad. Sci. URSS, n. Ser.}}, 25:347--349, 1939.

\end{thebibliography}

\end{document}